\documentclass[review]{elsarticle}

\usepackage{lineno,hyperref}

\usepackage{amsmath} 
\usepackage[subnum]{cases}

\usepackage{mathtools,cases}

\usepackage{graphicx}
\usepackage{bm}
\usepackage{makeidx}
\usepackage{multirow}
\usepackage{multicol}
\usepackage[dvipsnames,svgnames,table]{xcolor}

\usepackage[colorinlistoftodos,prependcaption]{todonotes}

\usepackage{epstopdf}
\usepackage{ulem}
\usepackage{hyperref}
\usepackage{amsmath,mathtools}
\usepackage{amssymb}
\usepackage{array}
\usepackage[caption=false]{subfig}
\usepackage{float}
\usepackage{amsfonts}
\usepackage{calrsfs}
\usepackage[title]{appendix}
\usepackage{comment}
\usepackage{soul}

\newcolumntype{C}[1]{>{\centering\arraybackslash}p{#1}}
\usepackage[]{algorithm2e}
\usepackage{setspace}

\newdefinition{rmk}{Remark}

\modulolinenumbers[5]

\journal{Journal of Computer Methods in Applied Mechanics and Engineering}

\usepackage{empheq}
\usepackage[many]{tcolorbox}

\tcbset{
	myhlight/.style={
		colback=yellow,
		arc=0pt,
		outer arc=0pt,
		boxrule=0pt,
		top=2pt,
		bottom=2pt,
		left=2pt,
		right=2pt,
	},
	highlight math style={myhlight}
}

\newtcolorbox{myhl}{
	breakable,
	myhlight
}

\newcommand{\splitatcommas}[1]{%
	\begingroup
	\begingroup\lccode`~=`, \lowercase{\endgroup
		\edef~{\mathchar\the\mathcode`, \penalty0 \noexpand\hspace{0pt plus 1em}}%
	}\mathcode`,="8000 #1%
	\endgroup
}

\begin{document}

\begin{frontmatter}

\title{Phase-field material point method for dynamic brittle fracture with isotropic and anisotropic surface energy}

\author[mymainaddress]{Emmanouil G. Kakouris}
\author[mymainaddress]{Savvas P. Triantafyllou \corref{mycorrespondingauthor}}
\cortext[mycorrespondingauthor]{Corresponding author}
\ead{savvas.triantafyllou@nottingham.ac.uk}

\address[mymainaddress]{Centre for Structural Engineering and Informatics, The University of Nottingham, Nottingham, NG7 2RD, UK}

\begin{abstract}

A novel phase field material point method is introduced for robust simulation of dynamic fracture in elastic media considering the most general case of anisotropic surface energy. Anisotropy is explicitly introduced through a properly defined crack density functional. The particular case of impact driven fracture is treated by employing a discrete field approach within the material point method setting. In this, the equations of motion and phase field governing equations are solved independently for each discrete field using a predictor-corrector algorithm. Contact at the interface is resolved through frictional contact conditions. The proposed method is verified using analytical predictions. The influence of surface energy anisotropy and loading conditions on the resulting crack paths is assessed through a set of benchmark problems. Comparisons are made with the standard Phase Field Finite Element Method and experimental observations.

\end{abstract}


\end{frontmatter}

\linenumbers
\section{Introduction} \label{sec:Introduction}

Failure of materials subjected to dynamic loading is commonly associated with complex yet intriguing phenomena, i.e., crack merging, branching and arrest \cite{Ravi_1984,Ravi_1998}. These phenomena become even more pronounced in the case of anisotropy and high-rate loading conditions, e.g., impact. Anisotropy governs the fracture response of both natural and manufactured materials as in the case of granitic rocks \cite{Chandler_2016}, biological tissues \cite{Holzapfel_2000}, single crystals \cite{Ast_2014}, and composite sheets \cite{Takei_2013}. Furthermore, the response of such materials under impact loading is being receiving considerable attention as it pertains to numerous industrial applications particularly within the automotive and aerospace sector, see, e.g., \cite{Turner_2018}. Numerical simulation of fracture propagation under such conditions can provide valuable insight into the underlying mechanical processes while also providing a framework for optimum design of materials considering their post-fracture response under impact loading. However, robust and accurate simulation of impact driven dynamic fracture is a challenging task as it requires the fusion of robust fracture propagation modelling with contact induced non-linearities and large displacement kinematics.

Within the framework of Computational Fracture Mechanics, a variety of methods has been introduced to address the problem of crack propagation. Among the most commonly used mesh based methods are the element deletion method \cite{Song_2008}, the Cohesive Zone Method \cite{Xu_1994,Park_2011} , the eXtended Finite Element Method (XFEM) \cite{moes_1999, Agathos_2017} and crack-driving configurational force approaches \cite{Miehe_2007, Kaczmarczyk_2017}. In these methods, algorithmic tracking of individual cracks is required as these evolve, merge, or branch. This results in considerable increase of the underlying computational complexity especially in the three dimensional case. Furthermore, an ad-hoc crack growth criterion is required for crack evolution. 

Francfort and Marigo \cite{francfort_1998} introduced a framework for avoiding these issues by establishing brittle fracture as an energy minimization problem within a robust variational structure. More recently, Bourdin et al. \cite{bourdin_2008} provided a regularization of the variational formulation which is more suitable for numerical solution schemes using as point of departure the phase field approximation of the Mumford-Shah potential presented in \cite{ambrosio_1990}. Within this variational setting, brittle fracture is formulated as a coupled, i.e., displacement and phase field problem, and the crack path naturally emerges from the solution of corresponding field equations. 

As a result, standard re-meshing or enrichment strategies near the crack tip and the requirement for algorithmic tracking of the crack front are avoided. Complex crack topologies, e.g., crack merging and/or branching as well as applications to three dimensional domains are efficiently resolved in the same manner, see, e.g., \cite{borden_2012}. Finite element based phase field formulations have been introduced to treat brittle \cite{aldakheel2018phase,moutsanidis2018hyperbolic}, ductile fracture \cite{ambati_2016,borden_2016} and hydraulic fracture \cite{Wilson_2016,miehe_2016a}. Phase field models for anisotropic fracture have been presented \cite{Li_2015,Teichtmeister_2017, Gultekin_2018,Nguyen_2017b} although not within a dynamic setting. Very recently, Hesch et al. \cite{Hesch_2016} have developed a method to resolve contact problems involving isotropic phase field fracture. In this formulation, a finite element based mortar contact algorithm in conjunction with a hierarchical refinement scheme is employed that reduces computational costs although relying on the predefinition of contact areas. Therefore, an adaptive hierarchical refinement is required for arbitrary impact fracture problems to resolve the local contact features.

The Material Point Method (MPM) \cite{sulsky_1994} has been introduced as a promising alternative to computationally expensive particle based methods that can efficiently deal with contact and large displacement problems. MPM is an extension of Particle-In-Cell (PIC) methods that efficiently treats history-dependent variables. In MPM, the continuum is represented by a set of Lagrangian particles, i.e., the material points, that are mapped onto a non-deforming Eulerian mesh (computational grid) where the governing equations are solved. This combined Eulerian-Lagrangian approach has been proven particularly advantageous in problems pertinent to high material and geometric nonlinearities since the distortion error is minimized \cite{zhang_2016,Charlton_2017,sofianos_2018}. Within this context, MPM has already been used to simulate very challenging engineering problems e.g. penetration \cite{Huang_2011}, cutting process simulations \cite{Ambati_2012} and solid-fluid interaction problems \cite{bandara_2015, Hamad_2017}.

To this point, few research has been conducted in damage simulation utilizing MPM using either discrete \cite{nairn_2003,Liang_2017}, cohesive \cite{daphalapurkar_2007, bardenhagen_2011}, or continuum damage models \cite{nairn_2017,Homel_2017}. Taking advantage of the good qualities of phase field modelling in naturally resolving complex crack paths, a Phase Field Material Point Method (PF-MPM) has been successfully introduced by the authors in \cite{kakouris_2017} for quasi-static brittle fracture problems while a variant accounting for anisotropy in the quasi-static regime has been developed in \cite{kakouris_2017b}.

Moving beyond the state-of-the-art, we present a phase field MPM method for the solution of dynamic fracture considering materials with anisotropic fracture energy; isotropy emerges as a special case of the proposed formulation. Following, the method is extended to also account for frictional contact fracture problems. We use as point of departure the MPM contact algorithm introduced in Bardenhagen et al. \cite{bardenhagen_2000} where multiple fields, termed discrete fields, are introduced in the non-deforming Eulerian mesh so that each contact body corresponds to a different field. We define the variational structure of our phase field implementation of impact driven fracture at each discrete field from which the coupled weak form of the contact problem naturally emerges. Finally, we develop a predictor corrector solution algorithm for the solution of the governing equations over time.

This paper is organized as follows. In section \ref{sec:Preliminaries} phase-field modelling is briefly described in both isotropic and anisotropic brittle fracture. The discrete field formulation for phase field fracture due to impact is presented in section \ref{sec:DFFPFIF}. The Material Point Method implementation for frictional contact fracture is presented in section \ref{sec:NumericalImplementation}. Finally, in section \ref{sec:NumericalExamples}, a set of benchmark problems are examined to demonstrate the accuracy and robustness of the proposed method.
\section{Preliminaries} \label{sec:Preliminaries}

\subsection{Phase field modelling} \label{sec:PFmodelling}

\begin{figure}
	\centering
	\begin{tabular}{cc}
		\subfloat[\label{fig:Theory_ContGamma}]{
			\includegraphics[width=0.50\columnwidth]{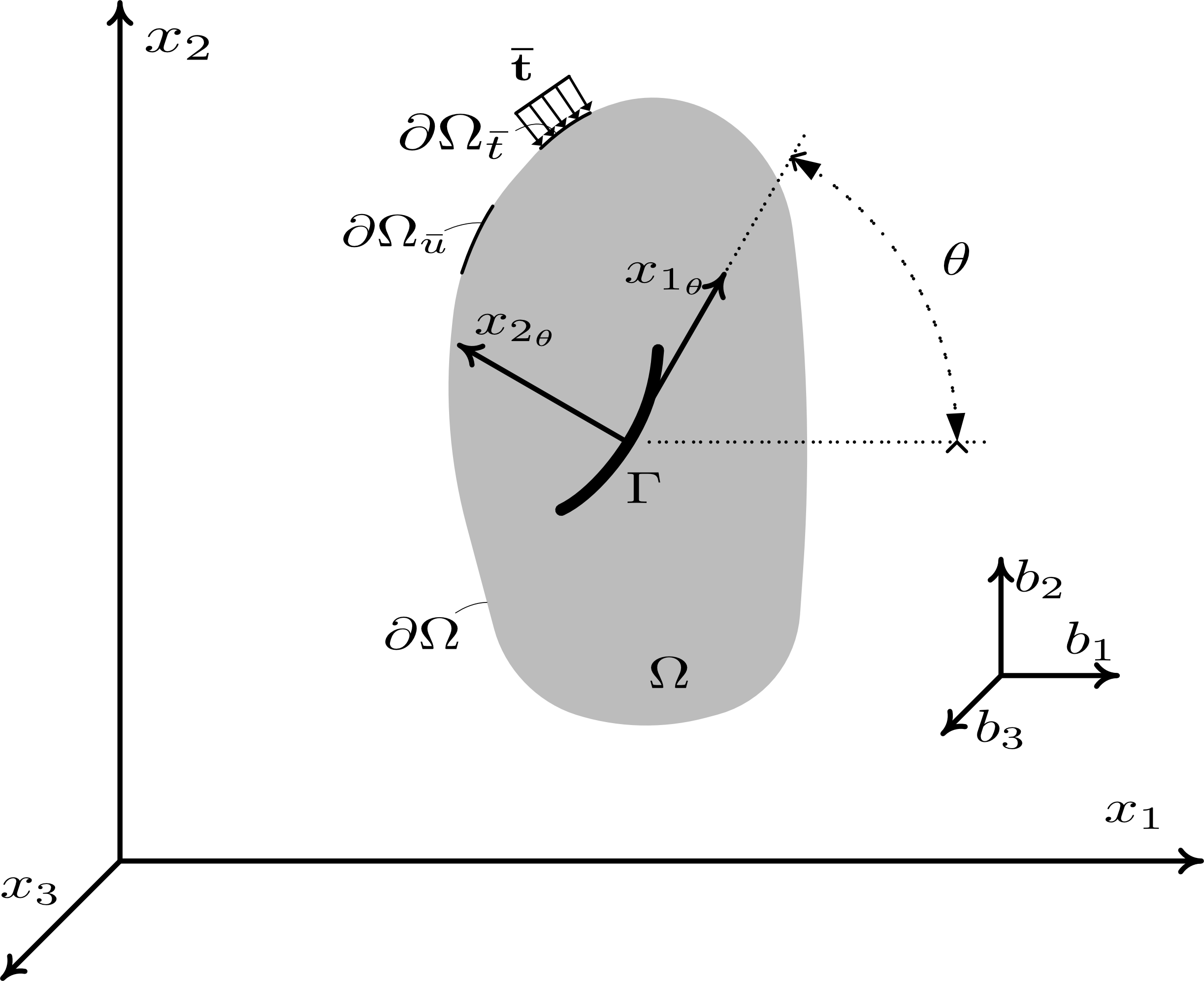}} &
		{\subfloat[\label{fig:Theory_ContGPF}]{
				\includegraphics[width=0.50\columnwidth]{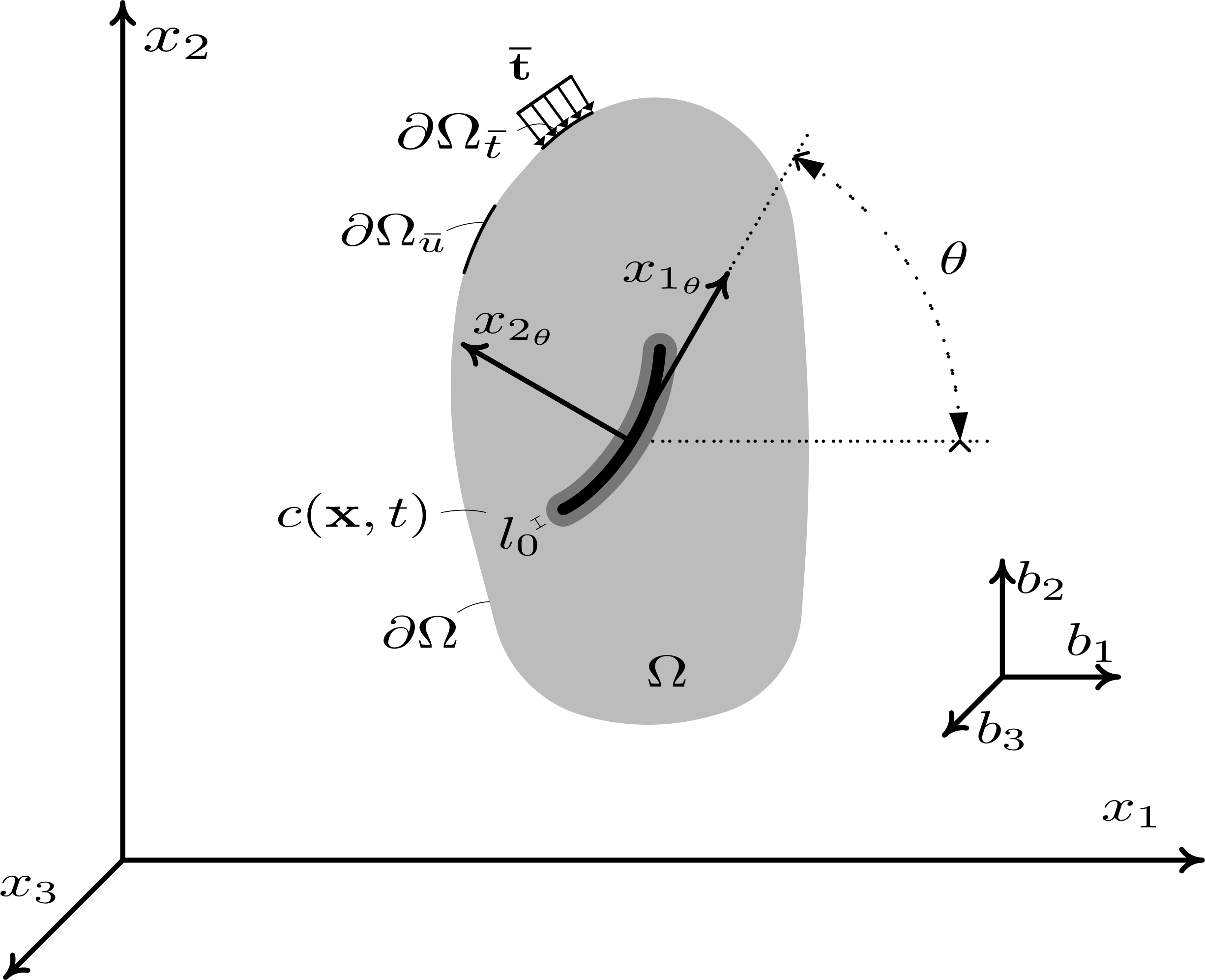}}}
	\end{tabular}
	\caption[]{\subref{fig:Theory_ContGamma} Solid body $\Omega$ with a crack path $\Gamma$ and \subref{fig:Theory_ContGPF} Phase field approximation of the crack path $\Gamma$.}
	\label{fig:Theory_ContGamma_ContGPF}
\end{figure}

In the following, the case of an arbitrary deformable domain $\Omega$ is considered, with an external boundary $\partial \Omega$ and a crack path $\Gamma$ as shown in Fig. \eqref{fig:Theory_ContGamma}. The deformable domain $\Omega$ with domain volume $V$, is subjected to body forces $\bm{b}=\left[\begin{matrix}b_1 & b_2 & b_3\end{matrix} \right]^T$. Furthermore, a set of traction/pressure loads $\bm{\bar{t}}$ is applied on the boundary $\partial\Omega_{\bar{t}} \subseteq \partial \Omega$. A prescribed displacement field, denoted as $\mathbf{\bar{u}}$, is imposed on the boundary $\partial\Omega_{\bar{u}} \subseteq \partial \Omega$. 

According to Griffith's  theory \cite{Griffith_163} the stored energy ${\Psi _{s}}$ of the body $\Omega$ can be expressed as 
\begin{equation}
\label{eqn:PotentialEnergyDecom}
{\Psi _{s}} = {\Psi _{el}} + {\Psi _{f}} = \int\limits_\Omega  {\psi _{el} \left( \boldsymbol{\varepsilon} \right) dV } + \int\limits_\Gamma  \mathcal{G}_c \left( \theta \right) d \Gamma
\end{equation}
where ${\Psi _{el}}$ and ${\Psi _{f}}$ are the elastic strain energy and the fracture energy (surface energy), respectively. Moreover, $\psi _{el} \left( \boldsymbol{\varepsilon} \right)$ corresponds to the elastic energy density and $\boldsymbol{\varepsilon}$ is the symmetric strain tensor which under the small strain assumption is defined as
\begin{equation}
\boldsymbol \varepsilon = \frac{1}{2} \Big(\nabla \mathbf{u} + \nabla \mathbf{u}^{T} \Big)
\label{eqn:SmallStrainsApprox}
\end{equation}
The $\left(\nabla\right)$ symbol in Eq. \eqref{eqn:SmallStrainsApprox} stands for the gradient operator and $\mathbf{u}(\mathbf{x},t)$ for the displacement field of a point $\mathbf{x}=\left[\begin{matrix}x_1 & x_2 & x_3\end{matrix} \right]^T$ at time $t$.

Due to material anisotropy, the critical fracture energy density $\mathcal{G}_c \left( \theta \right)$ in equation \eqref{eqn:PotentialEnergyDecom} explicitly depends on the orientation angle of the crack $\theta\left(s\right)$, $s \in \Gamma$. In the 2D case, the orientation angle is defined as the angle between the tangent vector at any point to the crack path $\Gamma$ and the horizontal. In the 3D case, the orientation can be defined by considering the direction cosines of the normal to the tangent plane of the fracture surface with respect to the global coordinate system. 

In principle, the stored energy $\Psi _{s}$ is known provided that both $\Gamma$ and $\theta$ at the current configuration are known. Hence, the computational treatment of elastic fracture mechanics gives rise to a nonlinear problem whereby standard procedures revert to path tracking and optimization algorithms to predict and resolve the crack path as this evolves. In the phase field approximation, the path dependent fracture energy surface integral is transformed into a volume integral defined over the entire domain $\Omega$ (Bourdin et al. \cite{bourdin_2008}) - see also, Fig. \ref{fig:Theory_ContGPF}. Hence, the phase field approximation gives rise to equation \eqref{eqn:FracEnergPFApprox}

\begin{equation}
{\Psi _{f}} = \int\limits_\Gamma  {\mathcal{G}_c \left( \theta \right) d \Gamma}  \approx \int\limits_\Omega  {\bar{\mathcal{G}}_c \mathcal{Z}_{c,Anis} dV}
\label{eqn:FracEnergPFApprox}
\end{equation}
where the functional $\mathcal{Z}_{c,Anis}=\mathcal{Z}_{c,Anis}\left(c, \theta \right)$ and $c$ is the phase field. Parameter $\bar{\mathcal{G}}_c$ in Eq. \eqref{eqn:FracEnergPFApprox} corresponds to the energy required to create a unit area of fracture surface $\mathcal{A}_c \left( \theta \right)$. This is assumed to be constant for all directions even though the actual surface energy of an anisotropic material is direction dependent. To resolve this inconsistency, directionality of the fracture toughness is accounted for in the definition of the functional $\mathcal{Z}_{c,Anis}$ by introducing a fourth-order anisotropic sensor $\gamma$ as discussed in section \ref{subsec:ACDFunc}. This formulation enables for a direct comparison to be drawn between isotropic and anisotropic models by controlling the value of a single parameter as further discussed in section \ref{sec:NumericalExamples}. Although strain rate dependence of the fracture toughness has been reported in the literature \cite{Grady_1979}, such concepts are beyond the scope of this work. 

\subsection{Anisotropic crack density functional} \label{subsec:ACDFunc}

To account for the general case of anisotropic material behaviour, $\mathcal{Z}_{c,Anis}$ is defined as the fourth-order functional utilized in \cite{Li_2015,kakouris_2017b} according to Eq.  \eqref{eqn:CrackDensityFunctional}
\begin{equation}
\mathcal{Z}_{c,Anis} = \left[ {\frac{{{{\left( {c - 1} \right)}^2}}}{{4{l_0}}} + {l_0}  | \nabla c |^2 } + l_0^{3} \sum_{\substack{ijkl}} \gamma_{ijkl} \frac{\partial^2 c}{\partial x_{i} \partial x_{j}} \frac{\partial^2 c}{\partial x_{k} \partial x_{l}} \right]
\label{eqn:CrackDensityFunctional}
\end{equation}
where $c(\mathbf{x},t) \in \left[0,1\right]$ is the phase field defined over the domain $\Omega$, $l_0 \in \mathbb{R}^{+}$ is a length scale parameter and $\gamma_{ijkl}$, $i,j,k,l=1\dots3$ are the components of the fourth-order anisotropic tensor corresponding to the anisotropic constitutive behaviour of the material. Phase field values of $c=1$ correspond to uncracked regions of the domain $\Omega$. Conversely, values of $c=0$ correspond to cracked regions. The length scale parameter $l_0$ controls the width of the regularized crack topology.

\begin{rmk} 
	A second-order functional can be employed to model anisotropy on the fracture properties. However, on the modelling side, fourth-order phase field functionals have been shown to successfully and robustly account for strong anisotropies \cite{li_2016,Teichtmeister_2017} while avoiding ill-posedness associated with second-order anisotropic models \cite{Li_2011}. On the simulation side, they improve the convergence rate of the underlying Newton procedure \cite{borden_2014}.		
	\label{rmk:fourthorder}
\end{rmk}

\begin{rmk} 
	A mathematical proof on the $\Gamma-$~convergence of the fourth-order anisotropic theory presented in this work has not been yet established. The extensive numerical studies performed in this work, within the bounds of the anisotropic tensors consider, hint that solutions provided by the fourth-order functional used in this work indeed converge. A relevant discussion on the aspect of $\Gamma-$~convergence for the fourth-order isotropic functional can be found in \cite{borden_2014}.		
	\label{rmk:GammaConvergence}
\end{rmk}

The anisotropic tensor $\bm{\gamma}$ is conveniently defined in the three dimensional space utilizing Voigt notation as
\setcounter{MaxMatrixCols}{20}
\begin{equation}
\bm{\gamma} = \begin{bmatrix}
\gamma_{1111} & \gamma_{1122} & \gamma_{1133} & \gamma_{1112} & \gamma_{1123} & \gamma_{1113}  \\
\gamma_{2211} & \gamma_{2222} & \gamma_{2233} & \gamma_{2212} & \gamma_{2223} & \gamma_{2213}  \\
\gamma_{3311} & \gamma_{3322} & \gamma_{3333} & \gamma_{3312} & \gamma_{3323} & \gamma_{3313}  \\
\gamma_{1211} & \gamma_{1222} & \gamma_{1233} & \gamma_{1212} & \gamma_{1223} & \gamma_{1213}  \\
\gamma_{2311} & \gamma_{2322} & \gamma_{2333} & \gamma_{2312} & \gamma_{2323} & \gamma_{2313}  \\
\gamma_{1311} & \gamma_{1322} & \gamma_{1333} & \gamma_{1312} & \gamma_{1323} & \gamma_{1313}  \\
\end{bmatrix}. 
\label{eqn:GammaTensor_3D}
\end{equation}
To demonstrate the versatility of the anisotropic functional in describing different material symmetries we focus on the 2D case for brevity in which case, $\bm{\gamma}$  reduces to
\setcounter{MaxMatrixCols}{20}
\begin{equation}
\bm{\gamma} = \begin{bmatrix}
\gamma_{1111} & \gamma_{1122} & \gamma_{1112} \\
\gamma_{2211} & \gamma_{2222} & \gamma_{2212} \\
\gamma_{1211} & \gamma_{1222} & \gamma_{1212} \\
\end{bmatrix}. 
\label{eqn:GammaTensor_2D}
\end{equation}
The direction angle of the crack path $\theta$ can be explicitly introduced in the expression of $\mathcal{Z}_{c,Anis}$ through a coordinate transformation,i.e., by transforming the Cartesian coordinate system $\mathbf{x}=\left[\begin{matrix} x_1 & x_2\end{matrix} \right]^T$  to $\mathbf{x}_{\theta} =\left[\begin{matrix} x_{\theta_1} & x_{\theta_2} \end{matrix} \right]^T$ where axis ${x}_{1_{\theta}}$ is defined along the crack path $\Gamma$ and axis ${x}_{2_{\theta}}$ is normal to the crack interface (see Fig. \ref{fig:Theory_ContGamma_ContGPF}). Thus, the transformation relation \eqref{eqn:CoordRotateCrackPath} holds
\begin{equation}
\mathbf{x}_{\theta} = \mathbf{R}_{\theta} \mathbf{x}
\label{eqn:CoordRotateCrackPath}
\end{equation}
where $\theta$ is the counter-clockwise angle between the ${x}_{1}$-axis and ${x}_{1_{\theta}}$ and $\mathbf{R}_{\theta}$ is the standard 2D rotation matrix. Detailed derivations on the transformation are provided by the authors in Appendix B of \cite{kakouris_2017b}. Eventually, the surface energy density $\mathcal{G}_{c} \left( \theta \right)$ for each angle $\theta$ is cast in the following form
\begin{equation}
\begin{aligned}
\mathcal{G}_{c} \left( \theta \right) = \int\limits_\Gamma  {\mathcal{G}_c \left( \theta \right) d \Gamma}  \approx \int_{-\infty}^{+\infty} \bar{\mathcal{G}}_c \mathcal{Z}_{c,Anis} dx_{2_{\theta}} \approx \int_{-x_{lb}}^{+x_{lb}} \bar{\mathcal{G}}_c \mathcal{Z}_{c,Anis} dx_{2_{\theta}} 
\end{aligned}
\label{TotalFracEnergPC}
\end{equation}
where $x_{lb}$ is the distance from the crack. Eq. \eqref{TotalFracEnergPC} enables numerical evaluation and visualization of the anisotropic surface energy density in polar form. Integration is performed along the normal to the crack path where the phase field variations are significant; assuming a value $x_{lb}=50 l_0$ yields a reasonable approximation.

\begin{figure}
	\centering
	\begin{tabular}{cc}
		\subfloat[\label{fig:1DExp_SurfEng}]{
			\includegraphics[width=0.45\columnwidth]{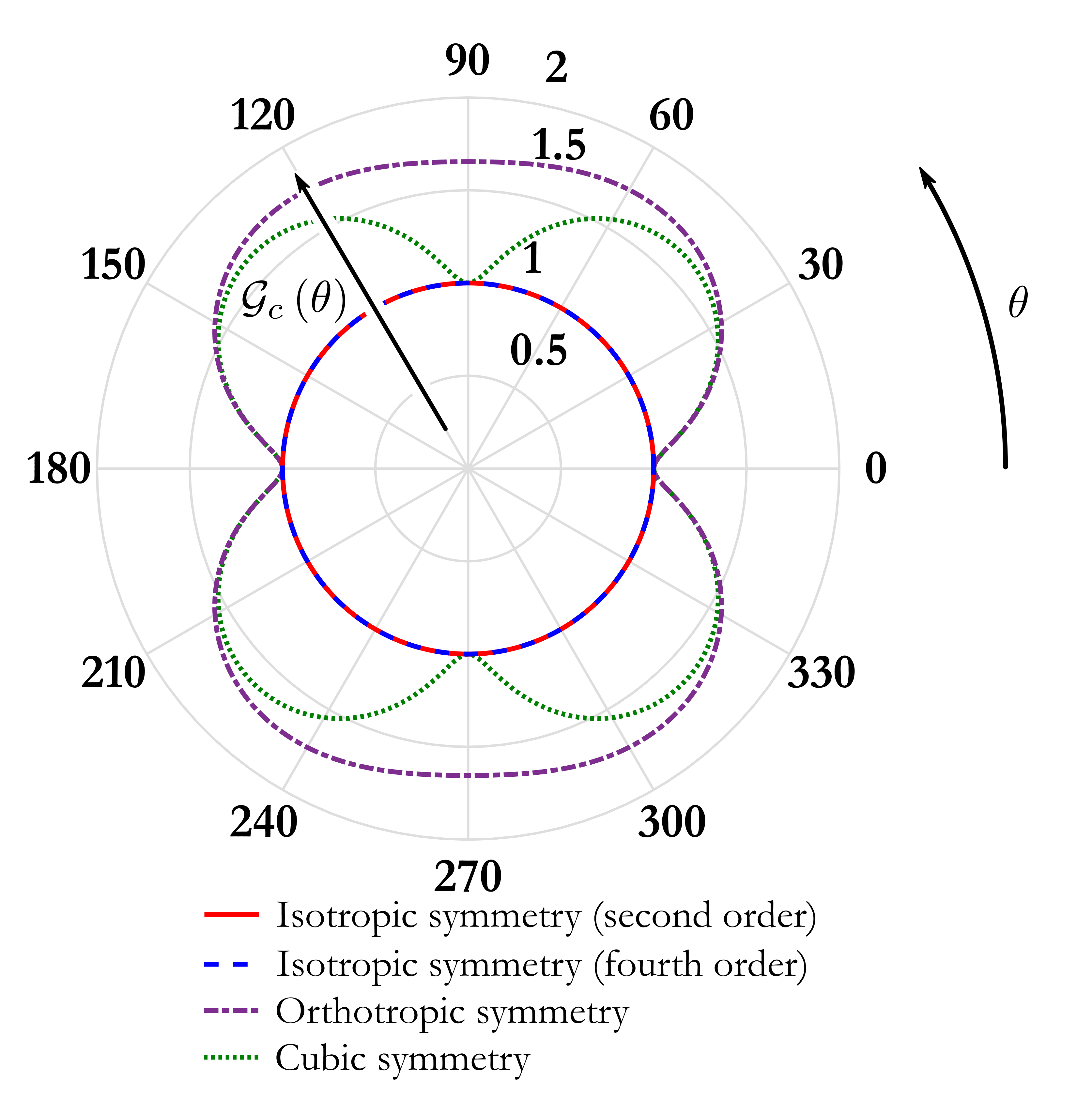}} &
		{\subfloat[\label{fig:1DExp_InvSurfEng}]{
				\includegraphics[width=0.45\columnwidth]{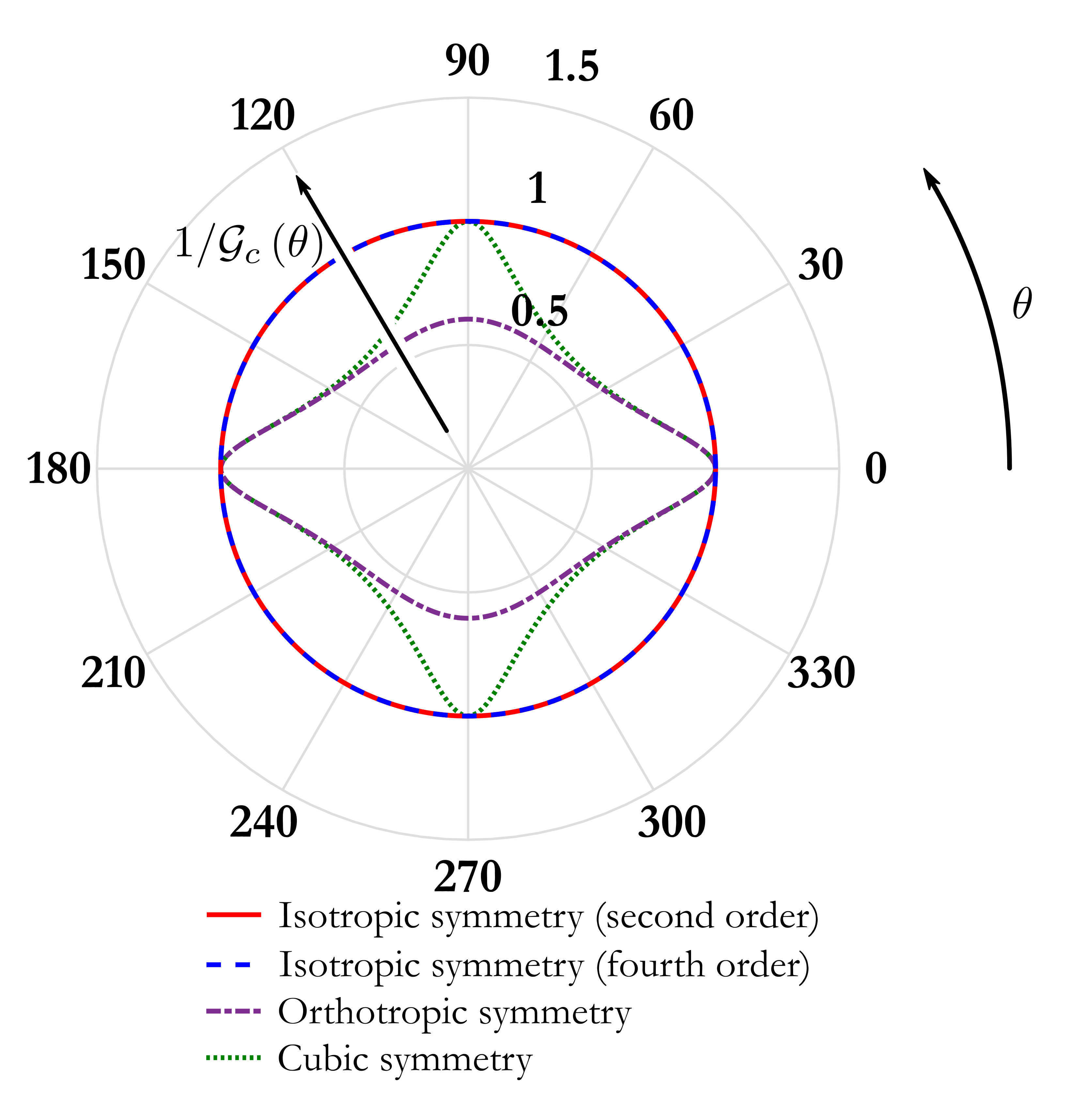}}}
	\end{tabular}
	\caption[]{Polar plots: \subref{fig:1DExp_SurfEng} Surface energy densities $\mathcal{G}_{c} \left( \theta \right)$ and \subref{fig:1DExp_InvSurfEng} their reciprocals $1/\mathcal{G}_{c} \left( \theta \right)$ in polar coordinates.}
	\label{fig:fig:1DExp_PolarPlots}
\end{figure}

Figs. \ref{fig:1DExp_SurfEng} and \ref{fig:1DExp_InvSurfEng} illustrate the surface energy densities $\mathcal{G}_{c} \left( \theta \right)$ and their reciprocals $1/ \mathcal{G}_{c} \left( \theta \right)$, respectively in polar coordinates for isotropic symmetry with the second and fourth order phase field models. The cases of cubic and orthotropic symmetry are also shown in Fig. \ref{fig:fig:1DExp_PolarPlots}. 
To derive these polar plots, the parameter ${{{\bar{\mathcal{G}}_c}}}$ is chosen to be ${{{\bar{\mathcal{G}}_c}}}=0.70710$ kN/m for the fourth order isotropic, cubic and orthotropic symmetry whereas  ${{{\bar{\mathcal{G}}_c}}}=1$ kN/m for second order isotropic symmetry. The parameter ${{{\bar{\mathcal{G}}_c}}}$ is chosen so that all previously mentioned models have the same minimum value of surface energy density $\mathcal{G}_{c_{min}}=1$ kN/m.

\begin{rmk}
Similar polar plots can be derived for the 3D case by considering the transformation of the global coordinate system to the coordinate system defined by the tangent plane at the fracture surface and its normal. It is useful to note that the 3D equivalent of Eq. \ref{TotalFracEnergPC} is still a line integral as integration is performed along the normal to the fracture surface. Such aspects are beyond the scope of this work; an intuitive approach on the rotation of anisotropic tensors is provided in \cite{nordmann2018visualising}.
\end{rmk}

\section{Governing equations for phase field fracture due to impact} \label{sec:DFFPFIF}

\subsection{Derivation of the coupled strong form for impact-fracture problems} \label{subsec:DCSF}

\begin{figure}
	\centering
	\begin{tabular}{lr}
		\multicolumn{2}{c}{\subfloat[\label{fig:Theory_Gamma}]{
				\includegraphics[width=0.40\columnwidth]{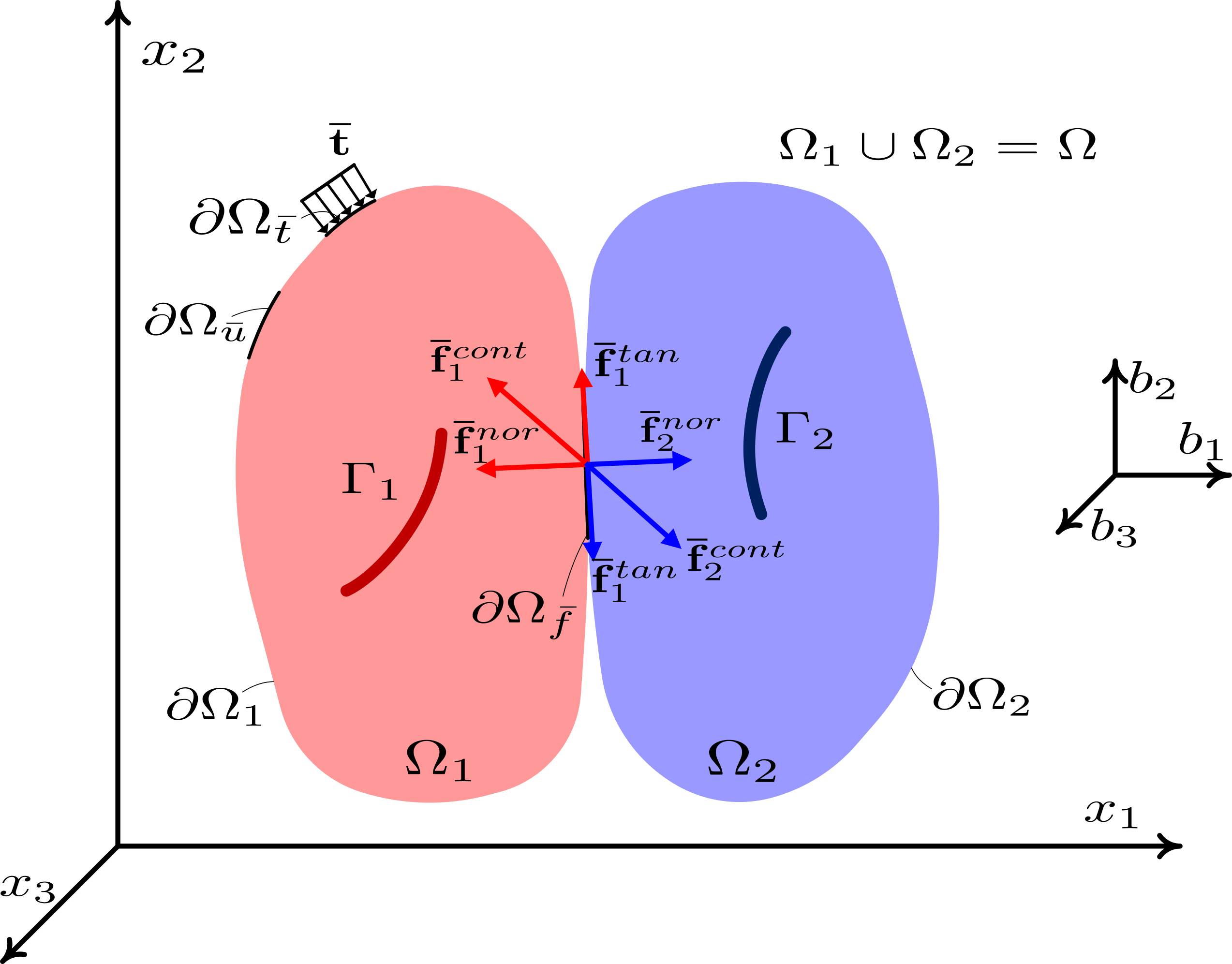}}} \\
		\subfloat[\label{fig:Theory_Cont}]{
			\includegraphics[width=0.40\columnwidth]{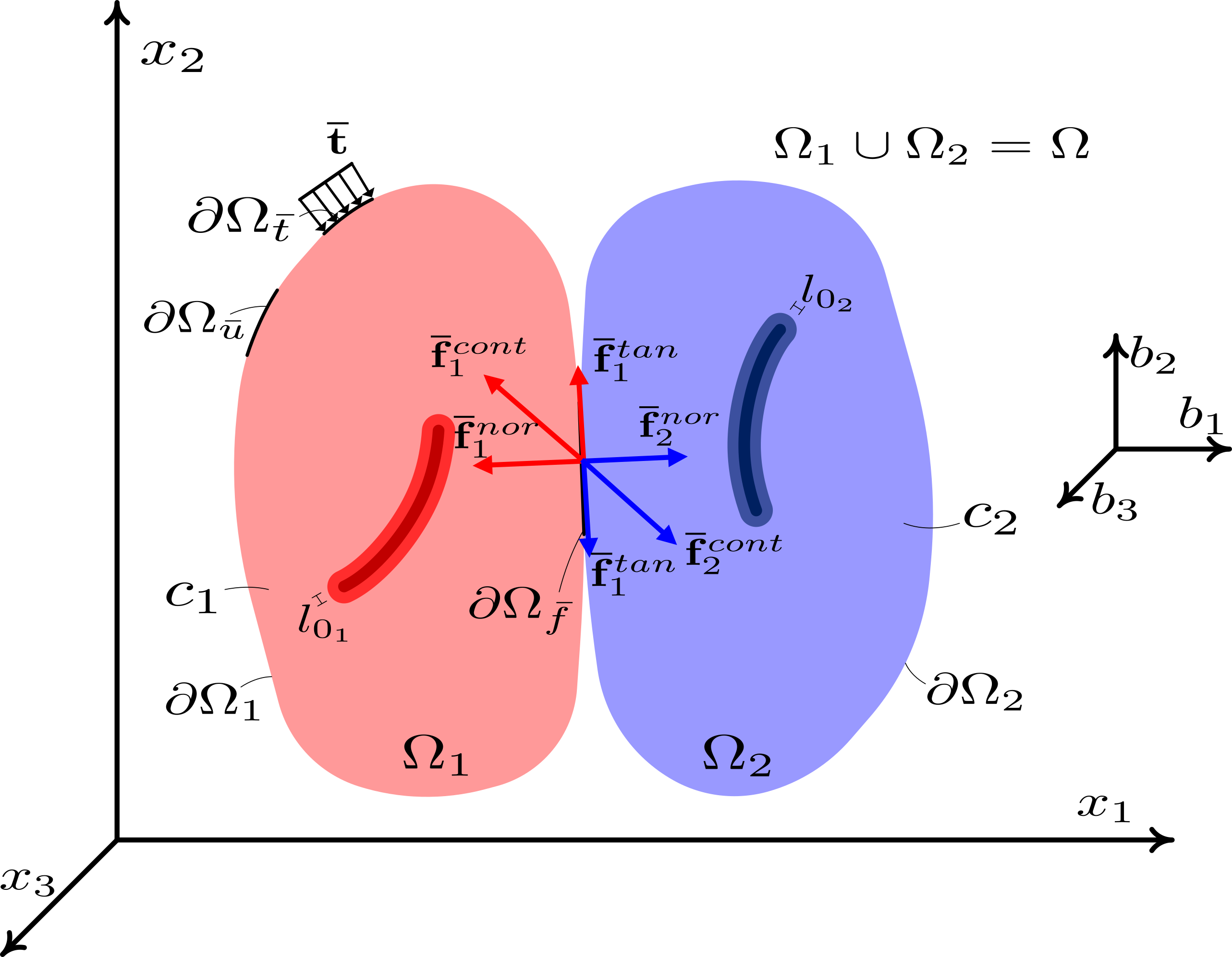}} &
		{\subfloat[\label{fig:Theory_Disc}]{
				\includegraphics[width=0.40\columnwidth]{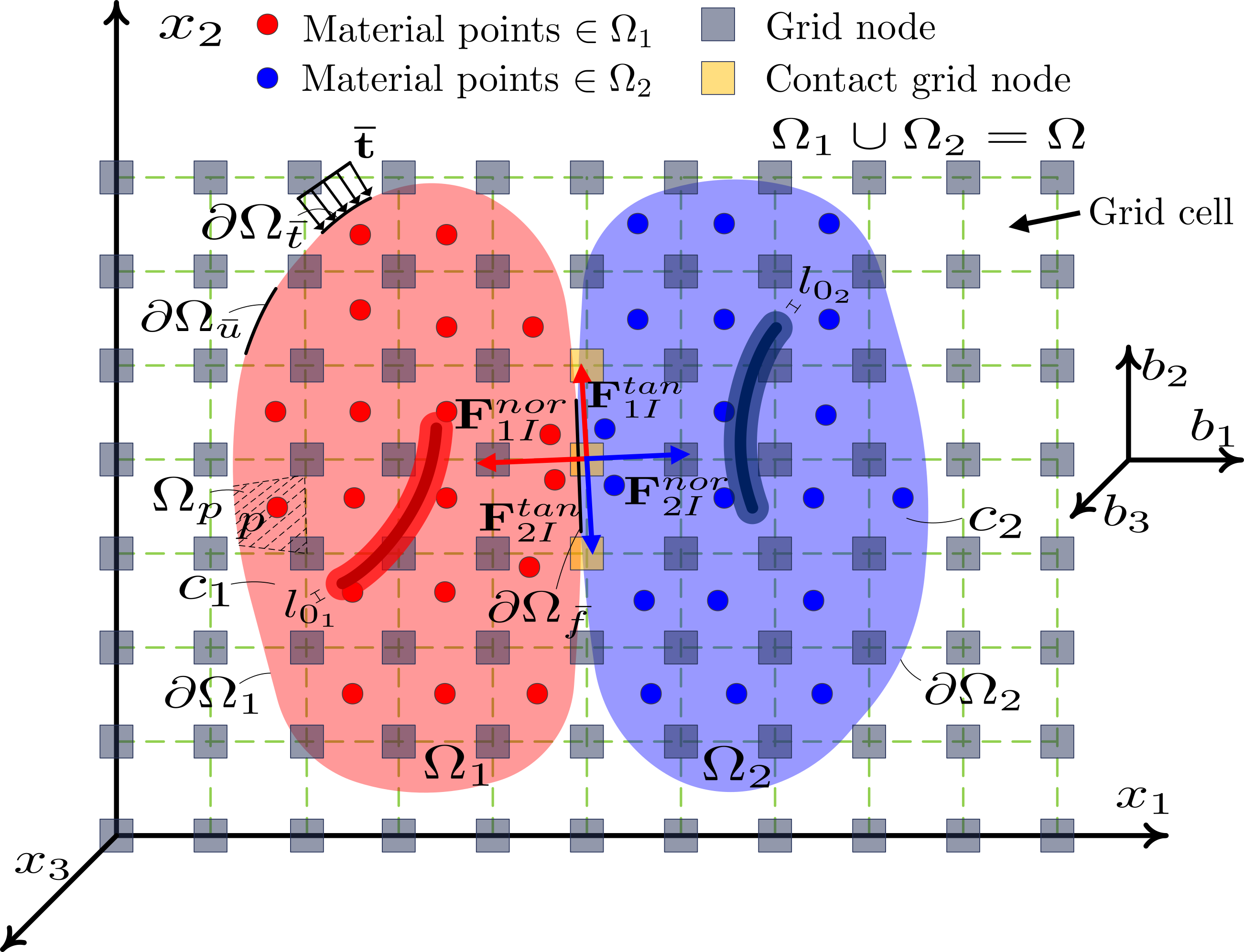}}}
	\end{tabular}
	\caption[]{\subref{fig:Theory_Gamma} Two bodies ($\Omega_{1} \cup \Omega_{2} = \Omega$)  into contact with two crack paths $\Gamma_{1}$ and $\Gamma_{2}$ \subref{fig:Theory_Cont} Phase field approximation of the crack paths and \subref{fig:Theory_Disc} Phase field material point method approximation.}
	\label{fig:Theory_Cont_Disc}
\end{figure}

In this section, the governing equations for contact induced brittle fracture are introduced. For brevity, the case of two bodies is presented herein. In Fig. \ref{fig:Theory_Gamma}, two deformable domains, namely $\Omega_1$ and $\Omega_2$ are considered such that $\Omega_{1} \cup \Omega_{2} = \Omega$. Their external boundaries are defined as $\partial \Omega_1$ and $\partial \Omega_2$, respectively.

At time $t$ the two bodies are in contact along the surface $\partial\Omega_{1\bar{f}} = \partial\Omega_{2\bar{f}} = \partial\Omega_{\bar{f}} \subseteq \partial\Omega$. Furthermore, two crack paths are defined, i.e., $\Gamma_{1}$ and $\Gamma_{2}$ at $\Omega_1$ and $\Omega_2$, respectively, under the action of a set of tractions $\bm{\bar{t}}$ and body forces $\bm{b}=\left[\begin{matrix}b_1 & b_2 & b_3\end{matrix} \right]^T$. 

When the two bodies are in contact, a contact force $\bm{\bar{f}}_{1}^{cont}$ is applied to body $\Omega_{1}$ from body $\Omega_{2}$. This is defined in component form according to Eq. \eqref{eqn:ContForce} 
\begin{equation}
\bm{\bar{f}}_{1}^{cont} = \bm{\bar{f}}_{1}^{nor} + \bm{\bar{f}}_{1}^{tan} = \bar{f}_{1}^{nor} \cdot \mathbf{n}_{1}^{cont} + \bar{f}_{1}^{tan} \cdot \mathbf{s}_{1}^{cont}
\label{eqn:ContForce}
\end{equation}
where $\bm{\bar{f}}_{1}^{nor}$ and $\bm{\bar{f}}_{1}^{tan}$ stand for the normal and tangential contact force vectors whereas $\bar{f}_{1}^{nor}$ and $\bar{f}_{1}^{tan}$ are their corresponding components. The normal and tangential surface unit vectors on contact surface $\partial\Omega_{1\bar{f}}$ are denoted as $\mathbf{n}_{1}^{cont}$ and $\mathbf{s}_{1}^{cont}$, respectively. Similarly, a contact force $\bm{\bar{f}}_{2}^{cont}$ is applied from $\Omega_{1}$ to $\Omega_{2}$ with components $\bm{\bar{f}}_{2}^{nor}$, $\bm{\bar{f}}_{2}^{tan}$ and normal and tangential surface unit vectors $\mathbf{n}_{2}^{cont}$ and $\mathbf{s}_{2}^{cont}$ being defined accordingly. 

Employing a phase field representation of fracture allows for a robust derivation of the impact-fracture strong form by considering the energy balance equation \eqref{eqn:PrincipleOfBalanceEnergy}
\begin{equation}
\dot{\mathcal{K}} \left( \dot{\mathbf{u}} \right) + \dot{\mathcal{W}}^{int} \left( \dot{\mathbf{u}},\dot{c},\nabla \dot{c} \right) - \dot{\mathcal{W}}^{ext} \left( \dot{\mathbf{u}} \right) - \dot{\mathcal{W}}^{cont} \left( \dot{\mathbf{u}} \right) = 0
\label{eqn:PrincipleOfBalanceEnergy}
\end{equation}
where $\dot{\mathcal{K}} \left( \dot{\mathbf{u}} \right)$ is the rate of the kinetic energy, $\dot{\mathcal{W}}^{int} \left( \dot{\mathbf{u}},\dot{c},\nabla \dot{c} \right)$ is the rate of internal work, $\dot{\mathcal{W}}^{ext} \left( \dot{\mathbf{u}} \right)$ is the rate of the work done by external forces, and $\dot{\mathcal{W}}^{cont} \left( \dot{\mathbf{u}} \right)$ is the rate of work done by contact forces. Furthermore, $\dot{\mathbf{u}}={d \mathbf{u}}/{d t}$ corresponds to the velocity field, $\dot{c}={d c}/{d t}$ is the phase field time derivative, and $\nabla \dot{c}$ corresponds to the rate of the phase field spatial derivative, i.e.,
\begin{equation}
\nabla \dot{c}=\frac{d}{d t} \left( \frac{\partial c}{\partial x_i} \right)
\end{equation}
for $i=1,...,3$. 

The kinetic energy rate functional $\dot{\mathcal{K}} \left( \dot{\mathbf{u}} \right)$ is expressed as 
\begin{equation}
\dot{\mathcal{K}} \left( \dot{\mathbf{u}} \right) = \frac{d}{dt} \int\limits_\Omega \frac{1}{2} \rho |\dot{\mathbf{u}}|^{2} dV 
\label{eqn:RateKineticEnergyFunc}
\end{equation}
where $\rho$ corresponds to the mass density.

The rate of internal work is expressed as 
\begin{equation}
\dot{\mathcal{W}}^{int} \left( \dot{\mathbf{u}},\dot{c},\nabla \dot{c} \right) = \frac{d \Psi _{s}}{dt} = \frac{d}{dt} \int\limits_\Omega \left( \psi_{el} +  \bar{\mathcal{G}}_c \mathcal{Z}_{c,Anis} \right) dV
\label{eqn:RateInternalWorkFunc}
\end{equation}
where the phase field approximation introduced in Eq. \eqref{eqn:FracEnergPFApprox} is employed to define the fracture energy corresponding to the crack paths $\Gamma_{1}$ and $\Gamma_{2}$ - see also Fig. \ref{fig:Theory_Cont}).

In this work, $\psi_{el}$ in Eq. \eqref{eqn:RateInternalWorkFunc} is decomposed into a purely tensile $\psi_{el}^{+}$ and a purely compressive $\psi_{el}^{-}$ parts according to the strain energy density decomposition introduced by Miehe et al. \cite{miehe_2010b} that is based on a spectral decomposition of the strain tensor. In this, the strain energy density is defined as
\begin{equation}
\begin{aligned}
\psi_{el} = g(c) \psi_{el}^{+} + \psi_{el}^{-}
\end{aligned}
\label{eqn:ElasticStrainEnergyDesn}
\end{equation}
where $g(c) \in \left[0,1\right]$ is a degradation function expressed as
\begin{equation}
\begin{aligned}
g = (1-k_{f}) c^2 + k_{f}
\end{aligned}
\label{eqn:DegradFunc}
\end{equation}
and $0 \le k_{f} \ll 1$ is a model parameter to treat potential ill-conditioning. In this work the model parameter is considered to be $k_{f}=0$ with no impact on the results as also highlighted by Braides \cite{braides_1998}. The stress field is derived from Eq. \eqref{eqn:Stress} as
\begin{equation}
\boldsymbol{\sigma} = \frac{\partial \psi_{el}}{\partial \boldsymbol {\varepsilon}}.
\label{eqn:Stress}
\end{equation}

Eq. \eqref{eqn:ElasticStrainEnergyDesn} is adopted herein for the purpose of verification however different schemes, also with significant computational advantages, can be found in the literature (see, e.g., \cite{ambati_2015} and \cite{li_2016} for a review of existing models). 

The rate of the external work functional $\dot{\mathcal{W}}^{ext} \left( \dot{\mathbf{u}} \right)$ is defined as
\begin{equation}
\dot{\mathcal{W}}^{ext} \left( \dot{\mathbf{u}} \right) = \int_{\partial \Omega_{\bar{t}}} ( \bm{\bar{t}} \cdot \dot{\mathbf{u}}) \,d \partial \Omega_{\bar{t}} + \int_\Omega (\bm{b} \cdot \dot{\mathbf{u}}) \,dV.
\label{eqn:RateExternalWorkFunc}
\end{equation}
Finally, the rate of work done by contact forces $\dot{\mathcal{W}}^{cont} \left( \dot{\mathbf{u}} \right)$ is expressed as
\begin{equation}
\begin{aligned}
\dot{\mathcal{W}}^{cont} \left( \dot{\mathbf{u}} \right) = \int_{\partial \Omega_{\bar{f}}} ( \bm{\bar{f}}^{cont} \cdot \dot{\mathbf{u}}) \,d \partial \Omega_{\bar{f}} = \int_{\partial \Omega_{\bar{f}}} \left( ( \bm{\bar{f}}^{nor} + \bm{\bar{f}}^{tan} ) \cdot \dot{\mathbf{u}}\right) \,d \partial \Omega_{\bar{f}}.
\end{aligned}
\label{eqn:RateContactWorkFunc}
\end{equation}
Clearly, $\dot{\mathcal{W}}^{cont} \left( \dot{\mathbf{u}} \right)$ must vanish as the contributing forces are always opposite. However, $\dot{\mathcal{W}}^{cont} \left( \dot{\mathbf{u}} \right)$ is retained in the energy balance equation and is further decomposed into discrete field components; this greatly facilitates numerical approximation as will be highlighted in section \ref{sec:NumericalImplementation}. Therefore,  Eq. \eqref{eqn:RateContactWorkFunc} is expressed as
\begin{equation}
\begin{aligned}
\dot{\mathcal{W}}^{cont} \left( \dot{\mathbf{u}} \right) = \dot{\mathcal{W}}_{1}^{cont} \left( \dot{\mathbf{u}}_{1} \right) + \dot{\mathcal{W}}_{2}^{cont} \left( \dot{\mathbf{u}}_{2} \right) = \\ \int_{\partial \Omega_{{1}\bar{f}}} ( \bm{\bar{f}}_{1}^{cont} \cdot \dot{\mathbf{u}}_{1}) \,d \partial \Omega_{{1}\bar{f}} + \int_{\partial \Omega_{{2}\bar{f}}} ( \bm{\bar{f}}_{2}^{cont} \cdot \dot{\mathbf{u}}_{2}) \,d \partial \Omega_{{2}\bar{f}}= 0
\end{aligned}
\label{eqn:RateContactWorkFunc_Fields}
\end{equation}
where $\dot{\mathbf{u}}_{1}$ and $\dot{\mathbf{u}}_{2}$ are the velocity fields at body $\Omega_{1}$ and $\Omega_{2}$, respectively. 
  
Applying the divergence theorem in Eq. \eqref{eqn:PrincipleOfBalanceEnergy}, performing the necessary algebraic manipulations, and finally considering that the resulting expression must hold for arbitrary values of $\dot{\mathbf{u}}$ and $\dot{c}$, the strong form of the problem is derived as (see \cite{kakouris_2017b} for details)
\begin{equation}
\label{AnisStrongForm}
\begin{cases}
\nabla \cdot \boldsymbol {\sigma}  + \bm{b} = \rho \ddot{\mathbf{u}} & \text{on } \Omega\\\\
\begin{aligned}
\left( {\frac{{4{l_0}\left( {1 - k_{f}} \right) \mathcal{H} }}{{{\bar{\mathcal{G}}_c}}} + 1} \right)c &- 4l_0^{2} \Delta c \\&+ 4 l_0^{4} \sum_{\substack{ijkl}} \gamma_{ijkl} \frac{\partial^4 c}{\partial x_{i} \partial x_{j} \partial x_{k} \partial x_{l}}= 1
\end{aligned} & \text{on } \Omega
\end{cases}
\end{equation}
where $\mathcal{H}$ is a history field defined as the maximum value of the tensile part of the elastic energy density $\psi _{el}^ +$ obtained in time domain $[0,t]$ and $\ddot{\mathbf{u}}={d \dot{\mathbf{u}}}/{d t}$ corresponds to the acceleration field. 

In our implementation, we use the history field (see, e.g., Miehe et al. \cite{miehe_2010}) to enforce the irreversibility condition pertinent to the crack propagation problem, i.e., $\prescript{(t)}{}{  \Gamma   }^{}_{} \subseteq \prescript{(t+\Delta t)}{}{  \Gamma   }^{}_{} $ according to the following Kuhn-Tucker conditions for loading and unloading, i.e.,
\begin{equation*}
\begin{array}{*{20}{c}}
{\psi _e^ +  - \mathcal{H} \le 0}&{\dot{\mathcal{H}} \ge 0}&{\mathcal{ \dot{H}}\left( {\psi _e^ +  - \mathcal{H}} \right) = 0}.
\end{array}
\end{equation*}

The coupled field equations \eqref{AnisStrongForm} are subject to the set of boundary and initial conditions defined in Eq. \eqref{AnisStrongFormBC}
\begin{equation}
\label{AnisStrongFormBC}
\begin{cases}
\boldsymbol{\sigma} \cdot \mathbf{n} = \bm{\bar{t}}, & \text{on } \partial\Omega_{\bar{t}} \\
\mathbf{u}=\mathbf{\bar{u}}, & \text{on } \partial\Omega_{\bar{u}} \\
\mathbf{u}=\prescript{(0)}{}{  \mathbf{u}   }^{}_{}, & \text{on } \prescript{(0)}{}{  \Omega   }^{}_{} \\
\dot{\mathbf{u}}=\prescript{(0)}{}{  \dot{\mathbf{u}}   }^{}_{}, & \text{on } \prescript{(0)}{}{  \Omega   }^{}_{} \\
\ddot{\mathbf{u}}=\prescript{(0)}{}{  \ddot{\mathbf{u}}   }^{}_{}, & \text{on } \prescript{(0)}{}{  \Omega   }^{}_{} \\
\left[ 4 l_0^{2} \nabla c - 2 l_0^{4} \sum_{\substack{ijkl}} \gamma_{ijkl} \left( \frac{\partial^3 c}{\partial x_{j} \partial x_{k} \partial x_{l}} \right) - 2 l_0^{4} \sum_{\substack{ijkl}} \gamma_{ijkl} \left( \frac{\partial^3 c}{\partial x_{i} \partial x_{j} \partial x_{k}} \right) \right] \cdot \mathbf{n} = 0, & \text{on }  \partial \Omega \\
2 l_0^{4} \sum_{\substack{ijkl}} \gamma_{ijkl} \left(  \frac{\partial^2 c}{\partial x_{k} \partial x_{l}} \right) + 2 l_0^{4} \sum_{\substack{ijkl}} \gamma_{ijkl} \left( \frac{\partial^2 c}{\partial x_{i} \partial x_{j}} \right) = 0, & \text{on } \partial \Omega \\
c=\prescript{(0)}{}{  c   }^{}_{}, & \text{on } \prescript{(0)}{}{  \Omega   }^{}_{} \\
\end{cases}
\end{equation}
where $\mathbf{n}$ stands for the outward unit normal vector on the boundary.

Furthermore, the coupled field equations \eqref{AnisStrongForm} are subjected to the kinematic constraints presented in Eqs. \eqref{CollinearityNormalCond} to \eqref{ComplementaryNormalCond} and \eqref{CollinearityTangCond} to \eqref{ComplementaryTangCond} at contact surface $\partial\Omega_{\bar{f}}$ \cite{Yastrebov_2013}. The kinematic constraints of Eqs. \eqref{CollinearityNormalCond} to \eqref{ComplementaryNormalCond} correspond to the normal contact laws
\begin{numcases}{}
\mathbf{n}^{cont}_{1} = - \mathbf{n}^{cont}_{2}, & collinearity, \space \space \space \space \space \space \space \space on $\partial\Omega_{\bar{f}}$  \label{CollinearityNormalCond} \\
\bm{\bar{f}}_{1}^{nor} = - \bm{\bar{f}}_{2}^{nor}, & collinearity, \space \space \space \space \space \space \space \space on $\partial\Omega_{\bar{f}}$  \label{ContactForceNormalCond} \\
\bar{f}^{nor} \leq 0, & non-tension, \qquad \space on $\partial\Omega_{\bar{f}}$ \label{NonTensionNormalCond} \\
\gamma_{n} \leq 0, & impenetrability, \space \space \space  on $\partial\Omega_{\bar{f}}$ \label{ImpenetratibilityCond} \\
\gamma_{n} \bar{f}^{nor} = 0, & complementarity, \space on $\partial\Omega_{\bar{f}}$ \label{ComplementaryNormalCond}
\end{numcases}
whereas Eq. \eqref{CollinearityTangCond} to \eqref{ComplementaryTangCond} correspond to the tangential contact and friction laws, where the Coulomb friction model is adopted.
\begin{numcases}{}
\mathbf{s}^{cont}_{1} = - \mathbf{s}^{cont}_{2}, & collinearity, \qquad \space \space on $\partial\Omega_{\bar{f}}$ \label{CollinearityTangCond} \\
\bm{\bar{f}}_{1}^{tan} = - \bm{\bar{f}}_{2}^{tan}, & collinearity, \qquad \space \space on $\partial\Omega_{\bar{f}}$ \label{ContactForceTangCond} \\
| \bar{f}^{tan} | \leq \mu_{f} | \bar{f}^{nor} |, & coulomb friction, \space on $\partial\Omega_{\bar{f}}$ \label{CoulombCond} \\
| \gamma_{s} | \geq 0, & slip/non-slip, \space \space \space \space \space \space  on $\partial\Omega_{\bar{f}}$ \label{gsCond} \\
| \gamma_{s} | \left(  | \bar{f}^{tan} | - \mu_{f} | \bar{f}^{nor} |  \right) = 0 , & complementarity, \space on $\partial\Omega_{\bar{f}}$. \label{ComplementaryTangCond}
\end{numcases}

Kinematic constraints \eqref{CollinearityNormalCond}, \eqref{CollinearityTangCond} and \eqref{ContactForceNormalCond} and \eqref{ContactForceTangCond} are imposed to satisfy Newton's third law at the contact surface $\partial\Omega_{\bar{f}}$. Condition \eqref{NonTensionNormalCond} is imposed on the normal component of the contact force that is defined according to Eq. \eqref{NonStick}
\begin{equation}
\bar{f}^{nor} = \bm{\bar{f}}^{cont}_{1} \cdot \mathbf{n}^{cont}_{1} = \bm{\bar{f}}^{cont}_{2} \cdot \mathbf{n}^{cont}_{2}
\label{NonStick}
\end{equation}
and implies a non-tension, i.e., non-stick, condition at the contact surface $\partial \Omega_{{}\bar{f}}$.

The tangential component $\bar{f}^{tan}$ is defined accordingly as
\begin{equation}
\bar{f}^{tan} = \bm{\bar{f}}^{cont}_{1} \cdot \mathbf{s}^{cont}_{1} = \bm{\bar{f}}^{cont}_{2} \cdot \mathbf{s}^{cont}_{2}.
\label{TangComp}
\end{equation}

The impenetrability condition \eqref{ImpenetratibilityCond} is imposed to ensure no penetration between the contact surfaces $\partial \Omega_{{1}\bar{f}}$ and $\partial \Omega_{{2}\bar{f}}$ when the two bodies are in contact.

\subsubsection{Discrete field formulation for the coupled governing equations} \label{subsec:DFFCGE}

In this work, a discrete field approach is adopted for the robust and efficient numerical treatment of contact dynamics between deformable bodies whereby each body is treated independently as discrete field. In the general case, it is assumed that the entire domain consists of a set of independent discrete fields $\{\mathcal{D}\mid \mathcal{D}=1,2,\dots,N_{\mathcal{D}}\}$, where $N_{\mathcal{D}} \in \mathbb{Z}^{+}$ stands for the total number of discrete fields and $\mathcal{D}$ indexes the $\mathcal{D}^{th}$ discrete field. Furthermore, all corresponding quantities that belong to discrete field $\mathcal{D}$, i.e., body $\Omega_{\mathcal{D}}$, are denoted with the subscript $\mathcal{D}$. Hence, in the two body case considered in this section $\{\mathcal{D}\mid \mathcal{D}=1,2\}$.

Within the discrete field setting, the contact forces arising from the interaction of the discrete fields are treated as additional external forces. Hence, the energy balance equation \eqref{eqn:PrincipleOfBalanceEnergy}, is re-defined for each discrete field $\mathcal{D}$ as
\begin{equation}
\dot{\mathcal{K}}_{\mathcal{D}} \left( \dot{\mathbf{u}}_{\mathcal{D}} \right) + \dot{\mathcal{W}}^{int}_{\mathcal{D}} \left( \dot{\mathbf{u}}_{\mathcal{D}},\dot{c}_{\mathcal{D}},\nabla \dot{c}_{\mathcal{D}} \right) - \dot{\mathcal{W}}^{ext}_{\mathcal{D}} \left( \dot{\mathbf{u}}_{\mathcal{D}} \right) - \dot{\mathcal{W}}^{cont}_{\mathcal{D}} \left( \dot{\mathbf{u}}_{\mathcal{D}} \right) = 0
\label{eqn:PrincipleOfBalanceEnergy_Disc}
\end{equation}
where  $\dot{\mathcal{W}}^{cont}_{\mathcal{D}} \left( \dot{\mathbf{u}}_{\mathcal{D}} \right)$ is the rate of work done by contact forces and is expressed as
\begin{equation}
\begin{aligned}
\dot{\mathcal{W}}^{cont}_{\mathcal{D}} \left( \dot{\mathbf{u}}_{\mathcal{D}} \right) = \int_{\partial \Omega_{{{\mathcal{D}}}\bar{f}}} ( \bm{\bar{f}}_{{\mathcal{D}}}^{cont} \cdot \dot{\mathbf{u}}_{{\mathcal{D}}}) \,d \partial \Omega_{{{\mathcal{D}}}\bar{f}}.
\end{aligned}
\label{eqn:RateContactWorkFunc_Fields_Disc}
\end{equation}
Thus, the coupled strong form introduced in Eqs. \eqref{AnisStrongForm} is now defined for the discrete field $\mathcal{D}$ as
\begin{equation}
\label{AnisStrongForm_Disc}
\begin{cases}
\nabla \cdot \boldsymbol \sigma_{\mathcal{D}}  + \bm{b}_{\mathcal{D}} = \rho_{\mathcal{D}} \ddot{\mathbf{u}}_{\mathcal{D}} & \text{on } \Omega_{\mathcal{D}}\\\\
\begin{aligned}
\left( {\frac{{4{l_{0_{\mathcal{D}}}}\left( {1 - k_{f_{\mathcal{D}}}} \right) \mathcal{H}_{\mathcal{D}} }}{{{\bar{\mathcal{G}}_{c_{\mathcal{D}}}}}} + 1} \right) c_{\mathcal{D}} & - 4l_{0_{\mathcal{D}}}^{2} \Delta c_{\mathcal{D}} \\&+ 4 l_{0_{\mathcal{D}}}^{4} \sum_{\substack{ijkl}} \gamma_{{ijkl}_{\mathcal{D}}} \frac{\partial^4 c_{\mathcal{D}}}{\partial x_{i} \partial x_{j} \partial x_{k} \partial x_{l}}= 1
\end{aligned} & \text{on } \Omega_{\mathcal{D}}.
\end{cases}
\end{equation}
The set of boundary and initial conditions introduced in Eqs. \eqref{AnisStrongFormBC} are modified for each discrete field $\mathcal{D}$ accordingly as
\begin{equation}
\label{AnisStrongFormBC_Disc}
\begin{cases}
\boldsymbol{\sigma}_{\mathcal{D}} \cdot \mathbf{n}_{\mathcal{D}} = \bm{\bar{t}}_{\mathcal{D}}, & \text{on } \partial\Omega_{\bar{t}_{\mathcal{D}}} \\
\mathbf{u}_{\mathcal{D}}=\mathbf{\bar{u}}_{\mathcal{D}}, & \text{on } \partial\Omega_{\bar{u}_{\mathcal{D}}} \\
\mathbf{u}_{\mathcal{D}}=\prescript{(0)}{}{  \mathbf{u}   }^{}_{{\mathcal{D}}}, & \text{on } \prescript{(0)}{}{  \Omega   }^{}_{{\mathcal{D}}} \\
\dot{\mathbf{u}}_{\mathcal{D}}=\prescript{(0)}{}{  \dot{\mathbf{u}}   }^{}_{{\mathcal{D}}}, & \text{on } \prescript{(0)}{}{  \Omega   }^{}_{{\mathcal{D}}} \\
\ddot{\mathbf{u}}_{\mathcal{D}}=\prescript{(0)}{}{  \ddot{\mathbf{u}}   }^{}_{{\mathcal{D}}}, & \text{on } \prescript{(0)}{}{  \Omega   }^{}_{{\mathcal{D}}} \\
\left[ 4 l_{0_{\mathcal{D}}}^{2} \nabla c_{\mathcal{D}} - 2 l_{0_{\mathcal{D}}}^{4} \sum_{\substack{ijkl}} \gamma_{{ijkl}_{{\mathcal{D}}}} \left( \frac{\partial^3 c}{\partial x_{j} \partial x_{k} \partial x_{l}} \right) - 2 l_{0_{\mathcal{D}}}^{4} \sum_{\substack{ijkl}} \gamma_{{ijkl}_{{\mathcal{D}}}} \left( \frac{\partial^3 c_{\mathcal{D}}}{\partial x_{i} \partial x_{j} \partial x_{k}} \right) \right] \cdot \mathbf{n}_{\mathcal{D}} = 0, & \text{on }  \partial \Omega_{\mathcal{D}} \\
2 l_{0_{\mathcal{D}}}^{4} \sum_{\substack{ijkl}} \gamma_{{ijkl}_{{\mathcal{D}}}} \left(  \frac{\partial^2 c_{\mathcal{D}}}{\partial x_{k} \partial x_{l}} \right) + 2 l_{0_{\mathcal{D}}}^{4} \sum_{\substack{ijkl}} \gamma_{{ijkl}_{{\mathcal{D}}}} \left( \frac{\partial^2 c_{\mathcal{D}}}{\partial x_{i} \partial x_{j}} \right) = 0, & \text{on } \partial \Omega_{\mathcal{D}} \\
c_{\mathcal{D}}=\prescript{(0)}{}{  c   }^{}_{{\mathcal{D}}}, & \text{on } \prescript{(0)}{}{  \Omega   }^{}_{{\mathcal{D}}} \\
\boldsymbol{\sigma}_{{\mathcal{D}}} \cdot \mathbf{n}_{{\mathcal{D}}}^{cont} = \bm{\bar{f}}_{{\mathcal{D}}}^{cont}, & \text{on } \partial\Omega_{{\mathcal{D}}\bar{f}} \\
\end{cases}
\end{equation}
where the last boundary condition is due to the contact forces that in this implementation are considered as forces applied externally to the discrete field $\mathcal{D}$.
\section{Material Point Method for dynamic anisotropic fracture} \label{sec:NumericalImplementation}

Dynamic fracture under impact naturally involves large displacement kinematics especially in the pre- and post-fracture regime, e.g., in the case of high
velocity projectile impact problems. To accurately resolve the pre and post fracture kinematics, the Material Point Method \cite{sulsky_1994} is used in this work to solve the system of coupled governing Eqs. \eqref{AnisStrongForm_Disc}.

In the Material Point Method framework employed herein, the entire domain $\Omega=\Omega_{1} \cup \Omega_{2}$ is discretized into a set of material points $\mathcal{P}=\{p\mid p=1,2,\dots,N_p\}$, where $N_p \in \mathbb{Z}^{+}$ is the total number of material points whereas $p$ indexes the $p^{th}$ material point. It is assumed herein that $N_{\mathcal{D}p}$ material points belong to discrete field $\mathcal{D}$, i.e. body $\Omega_{\mathcal{D}}$ (see Fig. \ref{fig:Theory_Disc}). 

According to the MPM approximation, the mass density $\rho_{\mathcal{D}}$ and domain volume $V_{\mathcal{D}}$ corresponding to the discrete field $\mathcal{D}$ are additively decomposed into the corresponding material point contributions according to Eqs. \eqref{MassDensity} and \eqref{Volume}, respectively, i.e., 
\begin{equation}
\label{MassDensity}
\rho_{\mathcal{D}} \left(\mathbf{x}_{\mathcal{D}},t\right)=\sum\limits_{p = 1}^{{N_p}} {{\rho_{\mathcal{D}p}}V_{\mathcal{D}p}\delta \left( {\mathbf{x}_{\mathcal{D}} - \mathbf{x}_{\mathcal{D}p}} \right)} 
\end{equation}
and
\begin{equation}
\label{Volume}
V_{\mathcal{D}} \left(\mathbf{x}_{\mathcal{D}},t\right)=\sum\limits_{p = 1}^{{N_p}} {{V_{\mathcal{D}p}}\delta \left( {\mathbf{x}_{\mathcal{D}} - \mathbf{x}_{\mathcal{D}p}} \right)} 
\end{equation}
where $\mathbf{x}_{\mathcal{D}}$ is the position vector of discrete field $\mathcal{D}$  and $\delta$ is the Dirac delta function. The material point mass density is defined as $\rho_{\mathcal{D}p}=M_{\mathcal{D}p}/V_{\mathcal{D}p}$ where $M_{\mathcal{D}p}$ and $V_{\mathcal{D}p}$ are the material point mass and volume, respectively. Furthermore,  $\mathbf{x}_{\mathcal{D}p}$ corresponds to the position vector of material point $p$ at discrete field $\mathcal{D}$.

\begin{figure}
	\centering
	\includegraphics[width=0.80\columnwidth]{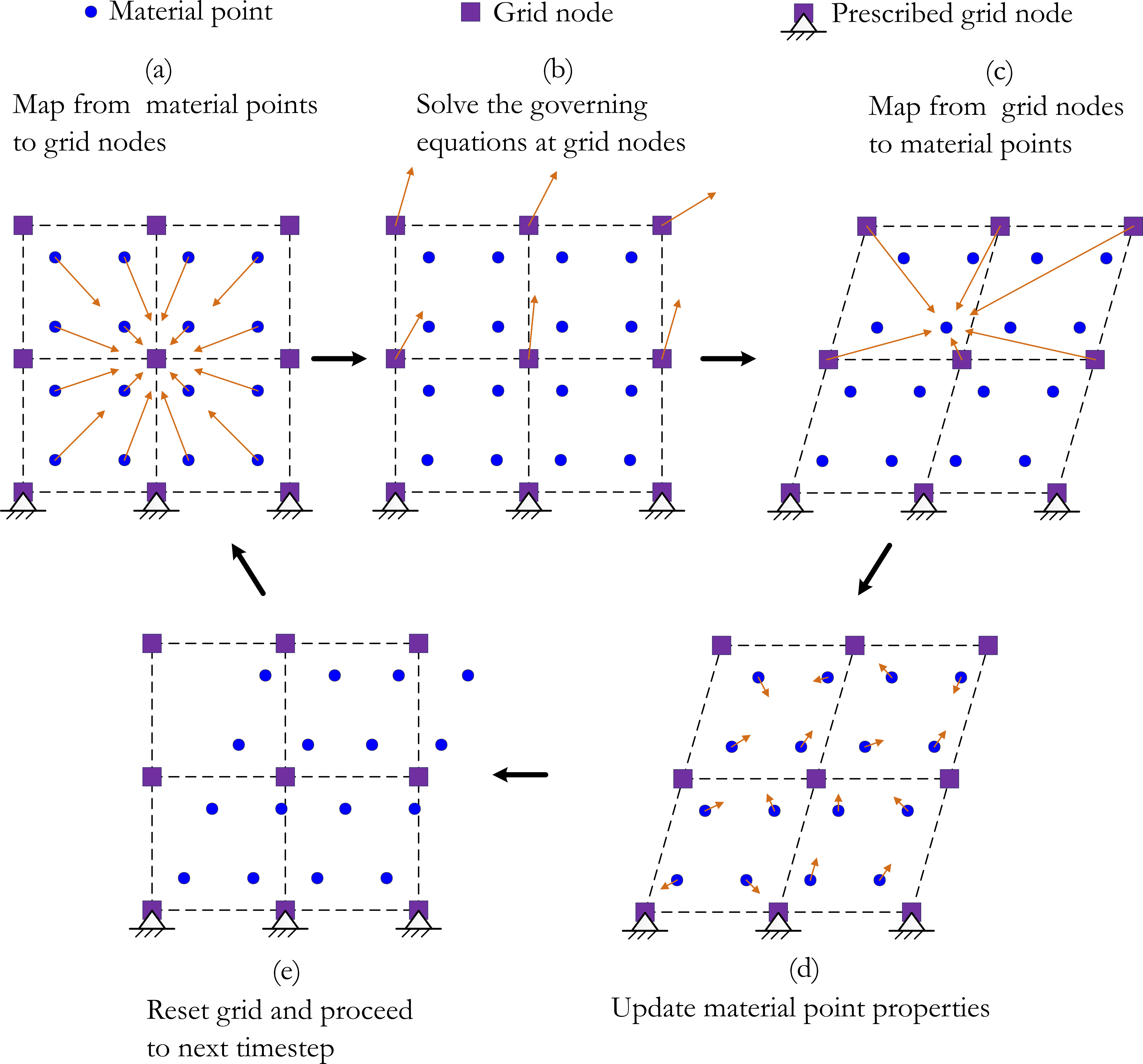}
	\caption[]{Material point method computational cycle}
	\label{fig:ComCy}
\end{figure}

These material points are moving within a fixed computational grid, i.e. an Eulerian grid. The Eulerian grid is a non-deforming mesh that consists a set of $N_{n} \in \mathbb{Z}^{+}$ grid nodes and $N_{cells} \in \mathbb{Z}^{+}$ grid cells (see Fig. \ref{fig:Theory_Disc}). The material points are mapped onto the Eulerian grid where the governing equations are solved. The updated solution is mapped back from grid nodes to the material points. Finally, the background grid is reset and the computational cycle proceeds. The steps of the MPM are shown in Fig. \ref{fig:ComCy}. In this work, mapping from material points to grid nodes and vice versa is implemented by utilizing higher order B-splines interpolation functions. 

\subsection{Discrete equilibrium equations for contact dynamics} \label{subsec:EDE}

Defining appropriate trial solution and weighting function spaces for the displacement field, i.e., 
\begin{equation*}
\mathcal{V}=\{\mathbf{u}  \in \left( H^{1}\left(\Omega\right) \right)^{d} \text{ } \mid \mathbf{u}=\mathbf{\bar{u}} \text{ on } \partial\Omega_{\bar{u}}\}
\end{equation*}
and 
\begin{equation*}
\mathcal{U}=\{\mathbf{w}  \in \left( H^{1}\left(\Omega\right) \right)^{d} \text{ } \mid \mathbf{w}=0 \text{ on } \partial\Omega_{\bar{u}}\},
\end{equation*}
respectively, the discrete form of the equations of motion introduced in the first of Eqs. \eqref{AnisStrongForm_Disc} is expressed for each discrete field $\mathcal{D}$ as
\begin{equation}
\label{eqn:SemiDiscreGalerkinDF}
\begin{aligned}
\int_{\Omega_{\mathcal{D}}}  ( \rho_{\mathcal{D}} \ddot{\mathbf{u}}_{\mathcal{D}} \cdot \mathbf{w}_{\mathcal{D}} ) \,d V_{\mathcal{D}} + \int_{\Omega_{\mathcal{D}}}  ( \boldsymbol{\sigma}_{\mathcal{D}} : \mathbf{\nabla w}_{\mathcal{D}}) \,d V_{\mathcal{D}} = \int_{\partial \Omega_{\mathcal{D} \bar{t}}} ( \bm{\bar{t}}_{\mathcal{D}} \cdot \mathbf{w}_{\mathcal{D}}) \,d \partial \Omega_{\mathcal{D}\bar{t}} + \\ \int_{\Omega_{\mathcal{D}}} (\bm{b}_{\mathcal{D}} \cdot \mathbf{w}_{\mathcal{D}}) \,d V_{\mathcal{D}} + \int_{\partial \Omega_{\mathcal{D} \bar{f}}} ( \bm{\bar{f}}_{\mathcal{D}}^{cont} \cdot \mathbf{w}_{\mathcal{D}}) \,d \partial \Omega_{\mathcal{D}\bar{f}}
\end{aligned}
\end{equation}
where $\mathbf{w}$ are weighting functions that satisfy the homogeneous essential boundary conditions of the problem \cite{hughes_2000}. 

Substituting the material point approximation introduced in Eqs. \eqref{MassDensity} and \eqref{Volume} into Eq. \eqref{eqn:SemiDiscreGalerkinDF}, Eq. \eqref{eqn:DiscreGalerkinDF} is established
\begin{equation}
\label{eqn:DiscreGalerkinDF}
\begin{aligned}
\sum_{p=1}^{N_p} ( \rho_{\mathcal{D}p} \ddot{\mathbf{u}}_{\mathcal{D}p} \cdot \mathbf{w}_{\mathcal{D}p}  ) V_{\mathcal{D}p} + \sum_{p=1}^{N_p} (\boldsymbol{\sigma}_{\mathcal{D}p} : \mathbf{\nabla w }_{\mathcal{D}p}) V_{\mathcal{D}p} = \int_{\partial \Omega_{\mathcal{D}\bar{t}}} ( \bm{\bar{t}}_{\mathcal{D}} \cdot \mathbf{w_{\mathcal{D}}}) \,d \partial \Omega_{\mathcal{D}\bar{t}} + \\ \sum_{p=1}^{N_p} ( \mathbf{b}_{\mathcal{D}p} \cdot \mathbf{ w }_{\mathcal{D}p}) V_{\mathcal{D}p} + \int_{\partial \Omega_{\mathcal{D}\bar{f}}} ( \bm{\bar{f}}_{\mathcal{D}}^{cont} \cdot \mathbf{w_{\mathcal{D}}}) \,d \partial \Omega_{\mathcal{D}\bar{f}}.
\end{aligned}
\end{equation}

Next, the weighting functions $\mathbf{w}_{\mathcal{D}p}$ and their spatial derivatives $\mathbf{\nabla w }_{\mathcal{D}p}$ are interpolated in the Galerkin sense according to relations \eqref{DispTestFunc}
\begin{equation}
\label{DispTestFunc}
\mathbf{w}_{\mathcal{D}p} = \sum_{I=1}^{N_n} N_{I}(\mathbf{x}_{\mathcal{D}p}) \mathbf{w}_{\mathcal{D}I} 
\end{equation}
and \eqref{DerDispTestFunc}, respectively, 
\begin{equation}
\label{DerDispTestFunc}
\mathbf{\nabla w }_{\mathcal{D}p} = \sum_{I=1}^{N_n} \nabla N_{I}(\mathbf{x}_{\mathcal{D}p}) \mathbf{w}_{\mathcal{D}I} 
\end{equation}
where $N_{I}(\mathbf{x}_{p})$ are the higher-order B-spline interpolation functions evaluated at the material point positions $\mathbf{x}_{\mathcal{D}p}$ \cite{kakouris_2017b,hughes_2005}. Furthermore, $\mathbf{w}_{\mathcal{D}I}$ are the weighting function, values evaluated at the background grid nodes; $I$ refers to the $I^{th}$ grid node. 

Similar expressions are established for the displacement, velocity, and acceleration field, i.e.,
\begin{equation}
\label{DispFunc}
\mathbf{u}_{\mathcal{D}p} = \sum_{I=1}^{N_n} N_{I}(\mathbf{x}_{p}) \mathbf{u}_{\mathcal{D}I},
\end{equation}
\begin{equation}
\label{VelFunc}
\dot{\mathbf{u}}_{\mathcal{D}p} = \sum_{I=1}^{N_n} N_{I}(\mathbf{x}_{p}) \dot{\mathbf{u}}_{\mathcal{D}I},
\end{equation}
and
\begin{equation}
\label{AccFunc}
\ddot{\mathbf{u}}_{\mathcal{D}p} = \sum_{I=1}^{N_n} N_{I}(\mathbf{x}_{p}) \ddot{\mathbf{u}}_{\mathcal{D}I},
\end{equation}
respectively, where $\mathbf{u}_{\mathcal{D}I}$, $\dot{\mathbf{u}}_{\mathcal{D}I}$, and $\ddot{\mathbf{u}}_{\mathcal{D}I}$ are the components of the nodal displacement, velocity and acceleration vectors, respectively, evaluated at node $I$.

Substituting Eqs. \eqref{DispTestFunc} and \eqref{DerDispTestFunc} in relation \eqref{eqn:DiscreGalerkinDF} and performing the necessary algebraic manipulations, the following expression is established 
\begin{equation}
\label{eqn:DiscreGalerkinDFGridArbi}
\sum_{I=1}^{N_n} \mathbf{w}_{\mathcal{D}I} \cdot \left[\bm{F}^{irt}_{\mathcal{D}I} + \bm{F}^{int}_{\mathcal{D}I} - \bm{F}^{ext}_{\mathcal{D}I} - \bm{F}^{cont}_{\mathcal{D}I}  \right]=0
\end{equation}
where $\bm{F}^{irt}_{\mathcal{D}I}$ are the nodal components of the inertia forces evaluated as
\begin{equation}
\label{eqn:InertiaForcesNodal}
\bm{F}^{irt}_{\mathcal{D}I} = \sum_{p=1}^{N_p} (\rho_{\mathcal{D}p} \ddot{\mathbf{u}}_{\mathcal{D}p} \cdot N_I(\mathbf{x}_{p})) V_{\mathcal{D}p}
\end{equation}
whereas $\bm{F}^{int}_{\mathcal{D}I}$ are the nodal components of the internal forces
\begin{equation}
\label{eqn:InternalForcesNodal}
\bm{F}^{int}_{\mathcal{D}I} = \sum_{p=1}^{N_p} (\boldsymbol{\sigma}_{\mathcal{D}p} \cdot \nabla N_I(\mathbf{x}_{p})) V_{\mathcal{D}p}.
\end{equation}
Similarly, the nodal components of the external force vector $\bm{F}^{ext}_{\mathcal{D}I}$ assume the following form
\begin{equation}
\label{eqn:ExternalForces}
\bm{F}^{ext}_{\mathcal{D}I} = \int_{\partial \Omega_{\mathcal{D}\bar{t}}} ( \bm{\bar{t}_{\mathcal{D}}} N_I(\mathbf{x})) \,d \partial \Omega_{\mathcal{D}\bar{t}} + \sum_{p=1}^{N_p} \bm{b}_{\mathcal{D}p} N_I(\mathbf{x}_{p}) V_{\mathcal{D}p}.
\end{equation}
Finally, $\bm{F}^{cont}_{\mathcal{D}I}$ corresponds to the contact force nodal vector defined as
\begin{equation}
\label{eqn:ContactForces}
\bm{F}^{cont}_{\mathcal{D}I} = \int_{\partial \Omega_{\mathcal{D}\bar{f}}} ( \bm{\bar{f}_{\mathcal{D}}}^{cont} N_I(\mathbf{x})) \,d \partial \Omega_{\mathcal{D}\bar{f}}.
\end{equation}
As the weighting functions in Eq. \eqref{eqn:DiscreGalerkinDF} are chosen arbitrarily, Eq. \eqref{eqn:DiscreGalerkinDFGridArbi} should hold for every set of nodal values $\mathbf{w}_{\mathcal{D}I}$. Hence, the following equilibrium equation is finally established
\begin{equation}
\bm{R}^{u}_{\mathcal{D}I} (\mathbf{u}_{\mathcal{D}}) = \bm{F}^{irt}_{\mathcal{D}I} + \bm{F}^{int}_{\mathcal{D}I} - \bm{F}^{ext}_{\mathcal{D}I} - \bm{F}^{cont}_{\mathcal{D}I} = 0, \quad I=1\dots,{N_n}
\label{eqn:ResidualVec_DF}
\end{equation}
where $\bm{R}^{u}_{\mathcal{D}I}$ is the nodal residual force vector at grid node $I$. 

Finally, substituting Eq. \eqref{AccFunc} in Eq. \eqref{eqn:InertiaForcesNodal}, Eq. \eqref{eqn:ResidualVec_DF} is rewritten in the following form
\begin{equation}
\bm{M}^{u}_{\mathcal{D}} \ddot{\mathbf{u}}_{\mathcal{D}} + \bm{F}^{int}_{\mathcal{D}} = \bm{F}^{ext}_{\mathcal{D}} + \bm{F}^{cont}_{\mathcal{D}}
\label{eqn:EqOfMotion}
\end{equation}
where $\bm{M}^{u}_{\mathcal{D}}$ is the global lumped mass matrix of the structure whose $M_{\mathcal{D}I}^{u}$ component is expressed as
\begin{equation}
\begin{aligned}
M_{\mathcal{D}I}^{u} = \sum_{p=1}^{N_p} \Big( \rho_{\mathcal{D}p} N_{I}(\mathbf{x}_{p}) \Big) V_{\mathcal{D}p}.
\end{aligned}
\label{eqn:MassLumpedMatrix}
\end{equation}
Eq. \eqref{eqn:EqOfMotion} lends itself conveniently into an explicit predictor-corrector time integration scheme as will be further discussed in Section \ref{subsec:SolProc}.

\subsection{Discrete phase field equations } \label{subsec:PFDE}

The discrete form of the anisotropic phase field governing equations introduced in the second of Eqs. \eqref{AnisStrongForm_Disc} can be also derived on the basis of the Material Point setting and Galerkin approximation. Similarly to the case of the displacement field, the phase field $c$ and the corresponding weighting functions $q$ are defined with respect to the following spaces, i.e.,
\begin{equation*}
\mathcal{Y}=\{c  \in H^{1}\left(\Omega\right)\}
\end{equation*}
and
\begin{equation*}
\mathcal{Q}=\{q  \in H^{1}\left(\Omega\right) \},
\end{equation*}
respectively.

Hence, the weak form of the phase field governing equations is expressed for each discrete field $\mathcal{D}$ as 
\begin{equation}
\begin{aligned}
\int_{\Omega_{\mathcal{D}}} \Big( \frac{4l_{0_{\mathcal{D}}}(1-k_{f_{\mathcal{D}}})\mathcal{H}_{\mathcal{D}}}{\bar{\mathcal{G}}_{c_{\mathcal{D}}}} + 1 \Big) c_{\mathcal{D}} q_{\mathcal{D}} \,d V_{\mathcal{D}} &+ \int_{\Omega_{\mathcal{D}}} 4l_{0_{\mathcal{D}}}^2 (\nabla c_{\mathcal{D}} : \nabla q_{\mathcal{D}}) \,d V_{\mathcal{D}} \\&+ \int_{\Omega_{\mathcal{D}}} 4l_{0_{\mathcal{D}}}^4 \sum_{\substack{ijkl}} \gamma_{ijkl_{\mathcal{D}}} \left( \frac{\partial^2 c_{\mathcal{D}}}{\partial x_{i} \partial x_{j}}  \frac{\partial^2 q_{\mathcal{D}}}{\partial x_{k} \partial x_{l}} \right) \,d V_{\mathcal{D}} \\ &= \int_{\Omega_{\mathcal{D}}} q_{\mathcal{D}} \,d V_{\mathcal{D}}.
\label{eqn:SemiDiscreGalerkinPF}
\end{aligned}
\end{equation}
Introducing the MPM approximation (Eq. \eqref{Volume}) into Eq. \eqref{eqn:SemiDiscreGalerkinPF}, the following expression is obtained
\begin{equation}
\begin{aligned}
\sum_{p=1}^{N_p} \mathcal{F}_{\mathcal{D}p} c_{\mathcal{D}p} q_{\mathcal{D}p} V_{\mathcal{D}p} &+ \sum_{p=1}^{N_p} 4l_{0_{\mathcal{D}p}}^2 (\nabla c_{\mathcal{D}p} : \nabla q_{\mathcal{D}p}) V_{\mathcal{D}p} \\&+ \sum_{p=1}^{N_p} 4l_{0_{\mathcal{D}p}}^4 \sum_{\substack{ijkl}} \gamma_{ijkl_{\mathcal{D}p}} \left( \frac{\partial^2 c_{\mathcal{D}p}}{\partial x_{i} \partial x_{j}}  \frac{\partial^2 q_{\mathcal{D}p}}{\partial x_{k} \partial x_{l}} \right) \Omega_{\mathcal{D}p}= \sum_{p=1}^{N_p} q_{\mathcal{D}p} V_{\mathcal{D}p} 
\end{aligned}
\label{eqn:DiscreGalerkinPF}
\end{equation}
where $c_{\mathcal{D}p}$, $q_{\mathcal{D}p}$ and $\gamma_{ijkl_{\mathcal{D}p}}$ are the phase field, weighting functions and anisotropic tensor components evaluated at the material point $p$. Parameter $\mathcal{F}_{\mathcal{D}p}$ in Eq. \eqref{eqn:DiscreGalerkinPF} is expressed as
\begin{equation}
\mathcal{F}_{\mathcal{D}p}=\frac{4l_{0_{\mathcal{D}p}}(1-k_{f_{\mathcal{D}p}})\mathcal{H}_{\mathcal{D}p}}{\bar{\mathcal{G}}_{c_{\mathcal{D}p}}} + 1
\label{eqn:ZeroTermCoefDiscreGalerkinPF}
\end{equation}
where $l_{0_{\mathcal{D}p}}$, $k_{\mathcal{D}p}$, $\mathcal{H}_{\mathcal{D}p}$ and $\bar{\mathcal{G}}_{c_{\mathcal{D}p}}$ are the length scale parameter, model parameter, history field and critical fracture energy density of material point $\mathbf{x}_p$.
Both $c_{\mathcal{D}p}$ and $q_{\mathcal{D}p}$ are interpolated at the nodal points of the background mesh, similarly to the case of the displacement field. Following the same procedure as in Section \ref{subsec:EDE} the nodal residual vector for the phase field is established as    
\begin{equation}
R^{c}_{\mathcal{D}I} (c_{\mathcal{D}}) = S^{c}_{\mathcal{D}I} - F^{c}_{\mathcal{D}I} = 0, \quad I=1\dots,{N_n}
\label{eqn:DiscreGalerkinPFGridArbi_5}
\end{equation}
where
\begin{equation}
\begin{aligned}
S^{c}_{\mathcal{D}I} = \sum_{p=1}^{N_p} \mathcal{F}_{\mathcal{D}p} c_{\mathcal{D}p} N_I(\mathbf{x}_{p}) V_{\mathcal{D}p} &+ \sum_{p=1}^{N_p} 4l_{0_{\mathcal{D}p}}^2 (\nabla c_{\mathcal{D}p} \cdot \nabla N_I(\mathbf{x}_{p})) V_{\mathcal{D}p} \\&+ \sum_{p=1}^{N_p} 4l_{0_{\mathcal{D}p}}^4 \sum_{\substack{ijkl}} \gamma_{ijkl_{\mathcal{D}p}} \left( \frac{\partial^2 c_{\mathcal{D}p}}{\partial x_{i} \partial x_{j}}  \frac{\partial^2 N_I(\mathbf{x}_{p})}{\partial x_{k} \partial x_{l}} \right) V_{\mathcal{D}p}
\end{aligned}
\label{eqn:Test}
\end{equation}
and
\begin{equation}
\label{eqn:PhaseFieldVolume}
F^{c}_{\mathcal{D}I}=\sum_{p=1}^{N_p} N_I(\mathbf{x}_{p}) V_{\mathcal{D}p}.
\end{equation}

Finally, applying the phase field interpolation and its spatial derivatives as
\begin{equation}
c_{\mathcal{D}p} = \sum_{I=1}^{N_n} N_{I}(\mathbf{x}_{p}) c_{\mathcal{D}I} 
\label{eqn:PhaseField}
\end{equation}
\begin{equation}
\nabla c_{\mathcal{D}p} = \sum_{I=1}^{N_n} \nabla N_{I}(\mathbf{x}_{p}) c_{\mathcal{D}I} 
\label{eqn:DerPFApprox}
\end{equation}
and
\begin{equation}
\Delta c_{\mathcal{D}p} = \sum_{I=1}^{N_n} \Delta N_{I}(\mathbf{x}_{p}) c_{\mathcal{D}I}. 
\label{eqn:Der2PFApprox}
\end{equation}
Eq. \eqref{eqn:DiscreGalerkinPFGridArbi_5} can be rewritten in the following convenient form as
\begin{equation}
\bm{K}^{c}_{\mathcal{D}} \bm{c}_{\mathcal{D}}= \bm{F}^c_{\mathcal{D}}
\label{eqn:PhaseFieldSolve}
\end{equation}
where $\bm{K}^{c}_{\mathcal{D}}$ is an $(N_n \times N_n )$ coefficient matrix whose ${K}^{c}_{I,J}$ component is expressed as
\begin{equation}
\begin{aligned}
K_{\mathcal{D} I,J}^{c} = \sum_{p=1}^{N_p} \Bigg( \mathcal{F}_{\mathcal{D}p} N_{J}(\mathbf{x}_{p}) N_{I}(\mathbf{x}_{p}) & + 4l_{0_{\mathcal{D}p}}^2  \Big( \nabla N_{J}(\mathbf{x}_{p}) \cdot \nabla N_I(\mathbf{x}_{p}) \Big) \\ &+ 4l_{0_{\mathcal{D}p}}^4\sum_{\substack{ijkl}} \gamma_{ijkl_{\mathcal{D}p}} \left( \frac{\partial^2 N_J(\mathbf{x}_{p})}{\partial x_{i} \partial x_{j}} \frac{\partial^2 N_I(\mathbf{x}_{p})}{\partial x_{k} \partial x_{l}} \right) \Bigg)  V_{\mathcal{D}p}.
\end{aligned}
\label{eqn:PhaseFieldMatrix}
\end{equation}
The $(N_n \times 1 )$ vector $\bm{c}_{\mathcal{D}}$ holds the nodal values of the phase field, defined at the background grid nodes,  and $\bm{F}^c_{\mathcal{D}}$ is the $(N_n \times 1 )$  vector whose $F^{c}_{\mathcal{D}I}$ component is defined from Eq. \eqref{eqn:PhaseFieldVolume}.

\subsection{Solution procedure}
\label{subsec:SolProc}
A staggered solution procedure \cite{miehe_2010} is employed to numerically solve the coupled Eqs. \eqref{eqn:EqOfMotion} and \eqref{eqn:PhaseFieldSolve}. In this, the two sets of equations are treated independently, by allowing the equation of motion to be solved either implicitly or explicitly \cite{borden_2012}. Although an explicit time integration scheme is utilized herein to integrate Eqs. \eqref{eqn:EqOfMotion} in the time domain \cite{Huang_2011}, an implicit time integration scheme can also be employed in a straightforward manner \cite{Chen_2017}. 
\subsubsection{Explicit time integration scheme} \label{subsubsec:ETI}

To numerically solve the equation of motion \eqref{eqn:EqOfMotion}, we employ the momentum formulation of the Material Point Method algorithm \cite{sulsky_1995}. Eq. \eqref{eqn:EqOfMotion}, is rewritten at the grid node $I$ at time $t$ as 
\begin{equation}
\prescript{(t)}{}{M}^{u}_{\mathcal{D}I}  \prescript{(t)}{}{  \ddot{\mathbf{u}}  }^{ }_{\mathcal{D}I} + \prescript{(t)}{}{  \bm{F}  }^{ int }_{\mathcal{D}I} = \prescript{(t)}{}{  \bm{F}  }^{ ext }_{\mathcal{D}I} + \prescript{(t)}{}{  \bm{F}  }^{ cont }_{\mathcal{D}I}
\label{eqn:ExpEqMotion_1}
\end{equation}
and considering, a forward Euler integration scheme, the acceleration field is expressed as 
\begin{equation}
\prescript{(t)}{}{  \ddot{\mathbf{u}}  }^{ }_{\mathcal{D}I} = ( \prescript{(t+\Delta t)}{}{  \dot{\mathbf{u}}  }^{ }_{\mathcal{D}I} - \prescript{(t)}{}{  \dot{\mathbf{u}}  }^{ }_{\mathcal{D}I} ) / \Delta t
\label{forwardEulerAcc}
\end{equation}
where $\Delta t$ stands for the corresponding time step. 

In view of Eq. \eqref{forwardEulerAcc}, Eq. \eqref{eqn:ExpEqMotion_1} is rewritten as
\begin{equation}
\begin{aligned}
\prescript{(t)}{}{M}^{u}_{\mathcal{D}I}  \prescript{(t+\Delta t)}{}{  \dot{\mathbf{u}}  }^{ }_{\mathcal{D}I} = \prescript{(t)}{}{M}^{u}_{\mathcal{D}I}  \prescript{(t)}{}{  \dot{\mathbf{u}}  }^{ }_{\mathcal{D}I} + \Delta t \left(  \prescript{(t)}{}{  \bm{F}  }^{ ext }_{\mathcal{D}I} + \prescript{(t)}{}{  \bm{F}  }^{ cont }_{\mathcal{D}I} - \prescript{(t)}{}{  \bm{F}  }^{ int }_{\mathcal{D}I} \right) \Leftrightarrow \\ \prescript{(t+\Delta t)}{}{\mathbf{p}}^{}_{\mathcal{D}I} = \prescript{(t)}{}{\mathbf{p}}^{}_{\mathcal{D}I} + \Delta t \left(  \prescript{(t)}{}{  \bm{F}  }^{ ext }_{\mathcal{D}I} + \prescript{(t)}{}{  \bm{F}  }^{ cont }_{\mathcal{D}I} - \prescript{(t)}{}{  \bm{F}  }^{ int }_{\mathcal{D}I} \right)
\end{aligned}
\label{eqn:ExpEqMotion_2}
\end{equation}
where $\prescript{(t+\Delta t)}{}{\mathbf{p}}^{}_{\mathcal{D}I}$ and $\prescript{(t)}{}{\mathbf{p}}^{}_{\mathcal{D}I}$ are the  nodal momentum at time $t+\Delta t$ and $t$, respectively. 

At time $t$ the nodal momentums $\prescript{(t)}{}{\mathbf{p}}_{\mathcal{D}I}$ are unknown; hence, these are mapped from material points to grid node $I$ using Eq. \eqref{eqn:MomentumProj}  
\begin{equation}
\prescript{(t)}{}{\mathbf{p}}^{}_{\mathcal{D}I} = \prescript{(t)}{}{M}^{u}_{\mathcal{D}I}  \prescript{(t)}{}{  \dot{\mathbf{u}}  }^{ }_{\mathcal{D}I} = \sum_{p=1}^{N_p} N_I(\prescript{(t)}{}{\mathbf{x}}^{}_{p}) M_{\mathcal{D}p} \prescript{(t)}{}{  \dot{\mathbf{u}}  }^{ }_{\mathcal{D}p}.
\label{eqn:MomentumProj}
\end{equation}

Similarly, the nodal internal forces $\prescript{(t)}{}{  \bm{F}  }^{ int }_{\mathcal{D}I}$ are evaluated as
\begin{equation}
\label{eqn:InternalForcesNodal_Exp}
\prescript{(t)}{}{  \bm{F}  }^{ int }_{\mathcal{D}I} = \sum_{p=1}^{N_p} ( \prescript{(t)}{}{   \boldsymbol{\sigma}   }^{}_{   \mathcal{D}p  } \cdot   \nabla N_I(\prescript{(t)}{}{\mathbf{x}}^{}_{p})) \prescript{(t)}{}{   V }^{}_{   \mathcal{D}p  }.
\end{equation}

Eq. \eqref{eqn:ExpEqMotion_2} is numerically solved by extending a predictor-corrector algorithm introduced by \cite{bardenhagen_2000} for granular media and further improved by \cite{Huang_2011} for the case of impact induced plasticity. In this, the trial momentums are initially evaluated for each discrete field $\mathcal{D}$, neglecting the contact forces $\prescript{(t)}{}{  \bm{F}  }^{ cont }_{\mathcal{D}I}$, as  
\begin{equation}
\begin{aligned}
\prescript{(t+\Delta t)}{}{\mathbf{p}}^{trl}_{\mathcal{D}I} = \prescript{(t)}{}{\mathbf{p}}^{}_{\mathcal{D}I} + \Delta t \left(  \prescript{(t)}{}{  \bm{F}  }^{ ext }_{\mathcal{D}I} - \prescript{(t)}{}{  \bm{F}  }^{ int }_{\mathcal{D}I} \right).
\end{aligned}
\label{eqn:ExpEqMotionTrl_1}
\end{equation}
The corresponding trial nodal velocities $\prescript{(t+\Delta t)}{}{\dot{\mathbf{u}}}^{trl}_{\mathcal{D}I}$ are then computed accordingly as 
\begin{equation}
\begin{aligned}
\prescript{(t+\Delta t)}{}{\dot{\mathbf{u}}}^{trl}_{\mathcal{D}I} = \frac{ \prescript{(t+\Delta t)}{}{\mathbf{p}}^{trl}_{\mathcal{D}I}}{ \prescript{(t)}{}{M}^{u}_{\mathcal{D}I} }.
\end{aligned}
\label{eqn:VelTrl}
\end{equation}
The trial velocities correspond to the velocities of each discrete field $\mathcal{D}$ when no contact force is exerted between them. 

The predicted trial velocities $\prescript{(t+\Delta t)}{}{\dot{\mathbf{u}}}^{trl}_{\mathcal{D}I}$ (evaluated from Eq. \eqref{eqn:VelTrl}) are then corrected according to Eq. \eqref{eqn:CorVel}
\begin{equation}
\label{eqn:CorVel}
\begin{aligned}
\prescript{(t+\Delta t)}{}{\dot{\mathbf{u}}}^{}_{\mathcal{D}I} = \prescript{(t+\Delta t)}{}{\dot{\mathbf{u}}}^{trl}_{\mathcal{D}I} + \Delta t \frac{ \prescript{(t)}{}{  \bm{F}  }^{ cont }_{\mathcal{D}I} }{\prescript{(t)}{}{M}^{u}_{\mathcal{D}I}}
\end{aligned}
\end{equation}
where $\prescript{(t+\Delta t)}{}{\dot{\mathbf{u}}}^{}_{\mathcal{D}I}$ is the vector of corrected nodal velocities at time $t+\Delta t$. 
To evaluate the corrected nodal velocities using Eq. \eqref{eqn:CorVel}, the contact forces $\prescript{(t)}{}{  \bm{F}  }^{ cont }_{\mathcal{D}I}$ must be evaluated first. The procedure for evaluating the contact forces between two discrete fields is presented in section  \ref{ContForceEval}.

\subsubsection{Contact force evaluation} \label{ContForceEval}

The contact force vector $\prescript{(t)}{}{  \bm{F}  }^{ cont }_{\mathcal{D}I}$ is the sum of a normal  $\prescript{(t)}{}{  \bm{F}  }^{ nor }_{\mathcal{D}I}$ and a tangential $\prescript{(t)}{}{  \bm{F}  }^{ tan }_{\mathcal{D}I}$ force vector. Hence, the corresponding components of these vectors, i.e. $\prescript{(t)}{}{  F  }^{ nor }_{\mathcal{D}I}$ and $\prescript{(t)}{}{  F  }^{ tan }_{\mathcal{D}I}$, should be initially computed taking into account the kinematic contact constraints presented in Eqs. \eqref{CollinearityNormalCond} to \eqref{ComplementaryNormalCond} and \eqref{CollinearityTangCond} to \eqref{ComplementaryTangCond}. Their evaluation is performed through the following procedure.

The nodal centre of mass velocities are calculated using Eq. \eqref{eqn:CentreMassVel} below
\begin{equation}
\label{eqn:CentreMassVel}
\prescript{(t+\Delta t)}{}{\dot{\mathbf{u}}}^{cm}_{I} = \frac{\sum_{\mathcal{D}=1}^{N_{\mathcal{D}}} \prescript{(t+\Delta t)}{}{\mathbf{p}}^{trl}_{\mathcal{D}I} }{\sum_{\mathcal{D}=1}^{N_{\mathcal{D}}} \prescript{(t)}{}{M}^{u}_{\mathcal{D}I}} = \frac{\sum_{\mathcal{D}=1}^{N_{\mathcal{D}}} \prescript{(t)}{}{M}^{u}_{\mathcal{D}I} \prescript{(t+\Delta t)}{}{\dot{\mathbf{u}}}^{trl}_{\mathcal{D}I} }{\sum_{\mathcal{D}=1}^{N_{\mathcal{D}}} \prescript{(t)}{}{M}^{u}_{\mathcal{D}I}}.
\end{equation}

These correspond to the velocities that each discrete field $\mathcal{D}$ would have if these were to move as a single field (non-slip contact). The normal component of the contact force $\prescript{(t)}{}{  F  }^{ nor,s }_{\mathcal{D}I}$ is evaluated considering the impenetrability condition defined in Eq. \eqref{ImpenetratibilityCond},  at contact grid node $I$ as

\begin{equation}
\label{eqn:ImpCon0}
\begin{aligned}
\prescript{(t+\Delta t)}{}{\gamma}^{}_{nI}= \left( \prescript{(t+\Delta t)}{}{\dot{\mathbf{u}}}^{}_{1I} - \prescript{(t+\Delta t)}{}{\dot{\mathbf{u}}}^{}_{2I} \right) \cdot \prescript{(t)}{}{\mathbf{n}}^{cont}_{1I} = 0. 
\end{aligned}
\end{equation}


As aforementioned in section \ref{subsec:DCSF}, when two bodies come into contact at contact grid node $I$ it holds that $\gamma_{nI}=0$. Substituting relation \eqref{eqn:CorVel} into \eqref{eqn:ImpCon0}, considering the equilibrium of contact forces on the contact surface, i.e.,
\begin{equation}
\label{eqn:Fcont_1}
\prescript{(t)}{}{  \bm{F}  }^{ cont }_{1I} = -\prescript{(t)}{}{  \bm{F}  }^{ cont }_{2I}
\end{equation}
and also Eq. \eqref{eqn:CentreMassVel}, the normal component of contact force is expressed as
\begin{equation}
\label{eqn:Fnors}
\prescript{(t)}{}{  F  }^{ nor,s }_{\mathcal{D}I} = \frac{\prescript{(t)}{}{M}^{u}_{\mathcal{D}I}}{\Delta t} \left( \prescript{(t+\Delta t)}{}{\dot{\mathbf{u}}}^{cm}_{I} - \prescript{(t+\Delta t)}{}{\dot{\mathbf{u}}}^{trl}_{\mathcal{D}I} \right) \cdot \prescript{(t)}{}{\mathbf{n}}^{cont}_{\mathcal{D}I}.
\end{equation}
The surface unit normal vector is computed by using the mass gradients \cite[see, e.g.,]{Homel_2017} as
\begin{equation}
\label{eqn:SurfNornalVec_Un}
\prescript{(t)}{}{\mathbf{\hat{n}}}^{cont}_{\mathcal{D}I} =\frac{\sum_{p=1}^{N_p} \nabla N_I(\prescript{(t)}{}{\mathbf{x}}^{}_{p}) M_{\mathcal{D}p} }{  \| \sum_{p=1}^{N_p} \nabla N_I(\prescript{(t)}{}{\mathbf{x}}^{}_{p}) M_{\mathcal{D}p} \|}.
\end{equation}
However, as also mentioned in \cite{Huang_2011}, Eq. \eqref{eqn:SurfNornalVec_Un} should be modified to satisfy the collinearity conditions \eqref{CollinearityNormalCond} and \eqref{CollinearityTangCond} at the contact surface $\partial\Omega_{\bar{f}}$ as 
\begin{equation}
\label{eqn:SurfNornalVec}
\prescript{(t)}{}{\mathbf{n}}^{cont}_{1I} = -\prescript{(t)}{}{\mathbf{n}}^{cont}_{2I} = \frac{ \prescript{(t)}{}{\mathbf{\hat{n}}}^{cont}_{1I} - \prescript{(t)}{}{\mathbf{\hat{n}}}^{cont}_{2I} }{ \| \prescript{(t)}{}{\mathbf{\hat{n}}}^{cont}_{1I} - \prescript{(t)}{}{\mathbf{\hat{n}}}^{cont}_{2I}  \| }
\end{equation}
to insure conservation of momentum.

To satisfy the non-tensional constraint (Eq. \eqref{NonTensionNormalCond}) during contact, the normal component should be modified as 
\begin{equation}
\label{eqn:Fnor}
\prescript{(t)}{}{  F  }^{ nor }_{\mathcal{D}I}=\min(0,\prescript{(t)}{}{  F  }^{ nor,s }_{\mathcal{D}I}). 
\end{equation}
Similarly, the tangential component of the contact force is evaluated considering the non-slip condition introduced in Eq. \eqref{gsCond} as
\begin{equation}
\label{eqn:NonSlip}
\begin{aligned}
\prescript{(t+\Delta t)}{}{\gamma}^{}_{sI}= \left( \prescript{(t+\Delta t)}{}{\dot{\mathbf{u}}}^{}_{1I} - \prescript{(t+\Delta t)}{}{\dot{\mathbf{u}}}^{}_{2I} \right) \cdot \prescript{(t)}{}{\mathbf{s}}^{cont}_{1I} = 0. 
\end{aligned}
\end{equation}
Substituting relation \eqref{eqn:CorVel} into \eqref{eqn:NonSlip} and then making use of Eqs. \eqref{eqn:Fcont_1} and \eqref{eqn:CentreMassVel}, the tangential component of contact force is expressed as
\begin{equation}
\label{eqn:Ftans}
\prescript{(t)}{}{  F  }^{ tan,s }_{\mathcal{D}I}=\frac{\prescript{(t)}{}{M}^{u}_{\mathcal{D}I}}{\Delta t} \left( \prescript{(t+\Delta t)}{}{\dot{\mathbf{u}}}^{cm}_{I} - \prescript{(t+\Delta t)}{}{\dot{\mathbf{u}}}^{trl}_{\mathcal{D}I} \right) \cdot \prescript{(t)}{}{\mathbf{s}}^{cont}_{\mathcal{D}I}   
\end{equation}
where the surface unit tangential vector $\mathbf{s}^{cont}_{\mathcal{D}I}$  can be derived as the unit vector that forms an orthogonal basis with $\mathbf{n}^{cont}_{\mathcal{D}I}$. The tangential component can be further modified to account for sliding at contact grid node $I$, considering the Coulomb friction model, as
\begin{equation}
\label{eqn:Ftan}
\prescript{(t)}{}{  F  }^{ tan }_{\mathcal{D}I} = \min \left( \mu_{f} | \prescript{(t)}{}{  F  }^{ nor }_{\mathcal{D}I} | , | \prescript{(t)}{}{  F  }^{ tan,s }_{\mathcal{D}I} | \right) \text{sign} \left( \prescript{(t)}{}{  F  }^{ tan,s }_{\mathcal{D}I} \right).
\end{equation}
Therefore, the contact force is eventually evaluated as
\begin{equation}
\label{eqn:Fcont}
\begin{aligned}
\prescript{(t)}{}{  \bm{F}  }^{ cont }_{\mathcal{D}I} = \prescript{(t)}{}{  F  }^{ nor }_{\mathcal{D}I} \cdot \prescript{(t)}{}{\mathbf{n}}^{cont}_{\mathcal{D}I} + \prescript{(t)}{}{  F  }^{ tan }_{\mathcal{D}I} \cdot \prescript{(t)}{}{\mathbf{s}}^{cont}_{\mathcal{D}I}   
\end{aligned}
\end{equation}
when the impenetrability condition
\begin{equation}
\label{eqn:ImpCond}
\begin{aligned}
\left( \prescript{(t+\Delta t)}{}{\dot{\mathbf{u}}}^{trl}_{\mathcal{D}I} - \prescript{(t+\Delta t)}{}{\dot{\mathbf{u}}}^{cm}_{\mathcal{D}I}  \right) \cdot \prescript{(t)}{}{\mathbf{n}}^{cont}_{\mathcal{D}I} > 0 
\end{aligned}
\end{equation}
is satisfied at contact grid node $I$. 

Finally, once the contact force vector is computed from Eq. \eqref{eqn:Fcont}, the initially predicted nodal velocities $\prescript{(t+\Delta t)}{}{\dot{\mathbf{u}}}^{trl}_{\mathcal{D}I}$ should be corrected according to Eq. \eqref{eqn:CorVel}.

\subsubsection{Material point properties update} \label{MPPU}

The corrected nodal velocities $\prescript{(t+\Delta t)}{}{\dot{\mathbf{u}}}^{}_{\mathcal{D}I}$ are utilized to update the material point properties. Hence, the total strains at the $p^{th}$ material point are evaluated as
\begin{equation}
\label{eqn:StrainsVelForm}
\begin{aligned}
\prescript{(t+\Delta t)}{}{\boldsymbol {\varepsilon}}^{}_{\mathcal{D}p} =  \prescript{(t)}{}{\boldsymbol {\varepsilon}}^{}_{\mathcal{D}p} + \frac{1}{2} \Delta t \sum_{I=1}^{N_n} \left(  \nabla N_{I}(\prescript{(t)}{}{\mathbf{x}}^{}_{p}) \prescript{(t+\Delta t)}{}{\dot{\mathbf{u}}}^{}_{\mathcal{D}I} + \left( \nabla N_{I}(\prescript{(t)}{}{\mathbf{x}}^{}_{p}) \prescript{(t+\Delta t)}{}{\dot{\mathbf{u}}}^{}_{\mathcal{D}I} \right)^{T} \right).
\end{aligned}
\end{equation}
The total stresses are evaluated from Eq. \eqref{eqn:Stress}. Finally, the displacement, velocity and acceleration of all material points are updated as
\begin{equation}
\label{eqn:MpDisp}
\begin{aligned}
\prescript{(t+\Delta t)}{}{  \mathbf{u}   }^{}_{\mathcal{D}p} = \prescript{(t)}{}{  \mathbf{u}   }^{}_{\mathcal{D}p} + \Delta t \sum_{I=1}^{N_n} \left( N_{I}(\prescript{(t)}{}{\mathbf{x}}^{}_{p}) \prescript{(t+\Delta t)}{}{\dot{\mathbf{u}}}^{}_{\mathcal{D}I}  \right) 
\end{aligned}
\end{equation}
\begin{equation}
\label{eqn:MpVel}
\begin{aligned}
\prescript{(t+\Delta t)}{}{  \dot{\mathbf{u}}   }^{}_{\mathcal{D}p} = \prescript{(t)}{}{  \dot{\mathbf{u}}   }^{}_{\mathcal{D}p} + \Delta t \sum_{I=1}^{N_n} \left( N_{I}(\prescript{(t)}{}{\mathbf{x}}^{}_{p}) \frac{\prescript{(t)}{}{  \bm{F}  }^{ ext }_{\mathcal{D}I} + \prescript{(t)}{}{  \bm{F}  }^{ cont }_{\mathcal{D}I} - \prescript{(t)}{}{  \bm{F}  }^{ int }_{\mathcal{D}I}}{\prescript{(t)}{}{M}^{u}_{\mathcal{D}I}}    \right)
\end{aligned}
\end{equation}
and 
\begin{equation}
\label{eqn:MpAcc}
\begin{aligned}
\prescript{(t)}{}{  \ddot{\mathbf{u}}   }^{}_{\mathcal{D}p} = \sum_{I=1}^{N_n} \left( N_{I}(\prescript{(t)}{}{\mathbf{x}}^{}_{p}) \frac{\prescript{(t)}{}{  \bm{F}  }^{ ext }_{\mathcal{D}I} + \prescript{(t)}{}{  \bm{F}  }^{ cont }_{\mathcal{D}I} - \prescript{(t)}{}{  \bm{F}  }^{ int }_{\mathcal{D}I}}{\prescript{(t)}{}{M}^{u}_{\mathcal{D}I}}    \right),
\end{aligned}
\end{equation}
respectively. The material point positions are also updated as 
\begin{equation}
\label{eqn:MpPosition}
\begin{aligned}
\prescript{(t+\Delta t)}{}{  \mathbf{x}   }^{}_{\mathcal{D}p} = \prescript{(t)}{}{  \mathbf{x}   }^{}_{\mathcal{D}p} + \Delta t \sum_{I=1}^{N_n} \left( N_{I}(\prescript{(t)}{}{\mathbf{x}}^{}_{p}) \prescript{(t+\Delta t)}{}{\dot{\mathbf{u}}}^{}_{\mathcal{D}I}  \right).
\end{aligned}
\end{equation}

\subsubsection{Staggered solution algorithm} \label{subsubsec:SSA}

The solution procedure is summarized in Algorithm \ref{alg:PFMPMStaggExplicit} where $E_{\mathcal{D}p}$ and $\nu_{\mathcal{D}p}$ are the Young's modulus and Poisson ratio at the material points. At each time increment $\left( m=0,...,N_{steps}-1 \right)$, the active part of the Eulerian grid is determined and the total number of active grid nodes $N_{n}$, unconstrained degrees of freedom $N_{dofs}$ and active cells $N_{cells}$ are evaluated (see \cite{kakouris_2017} for details). 

The background grid basis functions and their derivatives, i.e., $ \bm{N}_{} \left( \prescript{(m)}{}{ \mathbf{x} }^{}_{p} \right)$, $\nabla \bm{N}_{} \left( \prescript{(m)}{}{ \mathbf{x} }^{}_{p} \right)$ and $\Delta \bm{N}_{} \left( \prescript{(m)}{}{ \mathbf{x} }^{}_{p} \right)$ are evaluated at the material points with respect to the global coordinate system. To account for the arbitrary material orientation, the first and second spatial derivatives of the basis functions, i.e., $\nabla \bm{N}_{\phi_p} \left( \prescript{(m)}{}{ \mathbf{x} }^{}_{p} \right)$,  and $\Delta \bm{N}_{\phi_p} \left( \prescript{(m)}{}{ \mathbf{x} }^{}_{p} \right)$, respectively are also evaluated in the local material coordinate system. In the 2D cases examined herein, the principal material orientation at the $p^{th}$ material point is defined with regards to the angle $\phi_{\mathcal{D}p}$ between the global axis $x_1$ and the principal material axis.

Following, the contact grid nodes are detected among the discrete fields according to Remark \ref{rmk:DetectContactNodes}.
\begin{rmk} 
	Two discrete fields are in contact at grid node $I$ when at least one material point from both discrete fields is projected into grid node $I$. In this case, the grid node $I$ is a contact grid node for this pair of discrete fields.	
	\label{rmk:DetectContactNodes}
\end{rmk}

Next, the outward normal $\prescript{(m)}{}{\mathbf{n}}^{cont}_{\mathcal{D}I}$ and tangential $\prescript{(m)}{}{\mathbf{s}}^{cont}_{\mathcal{D}I}$ unit vectors are computed at the contact grid nodes. Mass, momentum and internal forces are projected from material points to grid nodes; thus, the quantities $\prescript{ (m)  }{   }{   M     }^{  u }_{ \mathcal{D}I  }$, $\prescript{ (m)  }{   }{  \mathbf{p}      }^{   }_{ \mathcal{D}I  }$ and $\prescript{ (m)  }{   }{  \bm{F}      }^{ int  }_{ \mathcal{D}I  }$ are obtained. Finally, the solution of the coupled Eqs. \eqref{eqn:ExpEqMotion_2} and \eqref{eqn:PhaseFieldSolve} is obtained within a set of $N_{staggs}$ staggered iterations $\left( k=1,...,N_{staggs} \right)$.

In the employed staggered scheme, the phase field Eq. \eqref{eqn:PhaseFieldSolve} is initially solved for a specific value of the history field $\prescript{ (m)  }{   }{   \mathcal{H}     }^{ (k)  }_{ \mathcal{D}p  }$. The basis functions $ \bm{N}_{} \left( \prescript{(m)}{}{ \mathbf{x} }^{}_{p} \right)$ and their spatial derivatives evaluated with respect to the material principal axes are utilized to compute the phase field coefficient matrix $\bm{K}^{c}_{\mathcal{D}}$ from relation \eqref{eqn:PhaseFieldMatrix}. Thus, the phase field nodal values $\prescript{ (m)  }{   }{   c     }^{ (k)  }_{  \mathcal{D}I }$ are obtained for each discrete field $\mathcal{D}$. Next, the phase field nodal values are mapped back onto the material points and the degradation function $\prescript{ (m)  }{   }{   g     }^{ (k)  }_{ \mathcal{D}p  }$ is computed at each material point. Next, the equation of motion \eqref{eqn:ExpEqMotion_2} is integrated in time employing the predictor-corrector algorithm described in Section \ref{subsubsec:ETI} and updated values for the history field $\prescript{ (m)  }{   }{   \mathcal{H}     }^{ (k)  }_{ \mathcal{D}p  }$ are obtained. 

Finally, the phase field nodal residual vector $\prescript{ (m)  }{   }{    R    }^{  c(k) }_{ I  }$ is evaluated according to the updated value of the history field $\prescript{ (m)  }{   }{   \mathcal{H}     }^{ (k)  }_{ \mathcal{D}p  }$ and convergence is checked as $ \| \prescript{ (m)  }{   }{    \bm{R}    }^{ c(k)  }_{   } \| \leq tol_{c} $ or $k \geq N_{staggs}$ where $tol_{c}$ and $\| \cdot \|$ stand for the phase field tolerance value and the Euclidean norm, respectively. After convergence, the material point properties are updated and the algorithm proceeds to the next increment $m$. 

Four conditions, namely \textbf{C.1} to \textbf{C.4}, are also included in Algorithm \ref{alg:PFMPMStaggExplicit}. These are employed to verify that the kinematic constraints introduced in Eqs. \eqref{CollinearityNormalCond}-\eqref{ComplementaryNormalCond} and \eqref{CollinearityTangCond}-\eqref{ComplementaryTangCond} are satisfied at contact grid nodes.

\begin{algorithm}
	\begin{small}
		\setstretch{0.60}
		\KwData{Define dynamic parameters, computational grid, material point properties ($\Delta t$, $\prescript{ (0)  }{   }{   \mathbf{x}     }^{   }_{  \mathcal{D}p }$, $\prescript{ (0)  }{   }{   V     }^{   }_{  \mathcal{D}p }$, $E_{\mathcal{D}p}$, $\nu_{\mathcal{D}p}$, $\prescript{ (0)  }{   }{   \rho     }^{   }_{  \mathcal{D}p }$, $l_{0_{\mathcal{D}p}}$, $k_{f_{\mathcal{D}p}}$, $\bar{\mathcal{G}}_{c_{\mathcal{D}p}}$, $\gamma_{ijkl_{\mathcal{D}p}}$, $\phi_{\mathcal{D}p}$, $\prescript{ (0)  }{   }{   \mathcal{H}     }^{   }_{  \mathcal{D}p }$, $\prescript{ (0)  }{   }{   \boldsymbol {\varepsilon}     }^{   }_{  \mathcal{D}p }$, $\prescript{ (0)  }{   }{   \boldsymbol{\sigma}     }^{   }_{  \mathcal{D}p }$, $\prescript{ (0)  }{   }{   \mathbf{u}     }^{   }_{  \mathcal{D}p }$, $\prescript{ (0)  }{   }{   \dot{\mathbf{u}}     }^{   }_{  \mathcal{D}p }$, $\prescript{ (0)  }{   }{   \ddot{\mathbf{u}}     }^{   }_{  \mathcal{D}p }$, $\prescript{ (0)  }{   }{   \mathbf{x}     }^{   }_{  \mathcal{D}p }$)}
		\For{each time step $m=0,..,N_{steps}-1$ }{
			Reset the computational grid: Find active part of Eulerian Grid, $N_n$, $N_{dofs}$, $N_{cells}$\;
			Compute: $\bm{N}_{} (\prescript{(m)}{}{ \mathbf{x} }^{}_{p} )$, $\nabla \bm{N}_{} ( \prescript{(m)}{}{ \mathbf{x} }^{}_{p} )$ and $\Delta \bm{N}_{} ( \prescript{(m)}{}{ \mathbf{x} }^{}_{p} )$ , for all material points. \;
			Compute: $\nabla \bm{N}_{\phi_p} ( \prescript{(m)}{}{ \mathbf{x} }^{}_{p} )$ and $\Delta \bm{N}_{\phi_p} ( \prescript{(m)}{}{ \mathbf{x} }^{}_{p} )$, for all material points. \;
			Detect contact grid nodes (see Remark \ref{rmk:DetectContactNodes})\;
			Compute: $\prescript{(m)}{}{\mathbf{n}}^{cont}_{\mathcal{D}I}$ and $\prescript{(m)}{}{\mathbf{s}}^{cont}_{\mathcal{D}I}$ vectors (see Eqs. \eqref{eqn:SurfNornalVec})\;
			\textbf{C.1:} Check collinearity conditions for all contact nodes (see Eqs. \eqref{CollinearityNormalCond} and \eqref{CollinearityTangCond})\;
			Map mass, momentum and internal forces from material points to grid nodes: $\prescript{ (m)  }{   }{   M     }^{  u }_{ \mathcal{D}I  }$, $\prescript{ (m)  }{   }{  \mathbf{p}      }^{   }_{ \mathcal{D}I  }$ and $\prescript{ (m)  }{   }{  \bm{F}      }^{ int  }_{ \mathcal{D}I  }$ (see Eqs. \eqref{eqn:MassLumpedMatrix}, \eqref{eqn:MomentumProj} and \eqref{eqn:InternalForcesNodal_Exp}) \;
			\For{each staggered iteration $k=1,2,..,N_{staggs}$ }{
				Compute: $\prescript{ (m)  }{   }{   \bm{F}     }^{  c(k) }_{ \mathcal{D}  }$ (see Eq. \eqref{eqn:PhaseFieldVolume} according to $\bm{N}_{} (\prescript{(m)}{}{ \mathbf{x} }^{}_{p} )$. \;
				Compute: $\prescript{ (m)  }{   }{   \bm{K}     }^{  c(k) }_{ \mathcal{D}  }$ (see Eq. \eqref{eqn:PhaseFieldMatrix}) according to $\bm{N}_{} (\prescript{(m)}{}{ \mathbf{x} }^{}_{p} )$, $\nabla \bm{N}_{\phi_p} ( \prescript{(m)}{}{ \mathbf{x} }^{}_{p} )$, $\Delta \bm{N}_{\phi_p} ( \prescript{(m)}{}{ \mathbf{x} }^{}_{p} )$ and $\prescript{ (m)  }{   }{   \mathcal{H}     }^{ (k)  }_{ \mathcal{D}p  }$. \;
				Solve: $\prescript{ (m)  }{   }{   \bm{K}     }^{  c(k) }_{ \mathcal{D}  } \prescript{ (m)  }{   }{   \bm{c}     }^{  (k) }_{ \mathcal{D}  } = \prescript{ (m)  }{   }{   \bm{F}     }^{  c(k) }_{ \mathcal{D}  }$ \;
				Map phase field ($\prescript{ (m)  }{   }{   \bm{c}     }^{  (k) }_{ \mathcal{D}  }$) from grid nodes to material points. Evaluate: $\prescript{ (m)  }{   }{   c     }^{  (k) }_{ \mathcal{D}p  }$, $\prescript{ (m)  }{   }{   \nabla c     }^{  (k) }_{ \mathcal{D}p  }$, $\prescript{ (m)  }{   }{   \Delta c     }^{  (k) }_{ \mathcal{D}p  }$, $\prescript{ (m)  }{   }{   g     }^{  (k) }_{ \mathcal{D}p  }$, for all material points (see Eqs. \eqref{eqn:PhaseField}, \eqref{eqn:DerPFApprox},\eqref{eqn:Der2PFApprox} and  
				\eqref{eqn:DegradFunc}). \;
				Update trial momentum: $\prescript{ (m+1)  }{   }{   \mathbf{p}     }^{  trl(k) }_{ \mathcal{D}I  }$ (see Eq. \eqref{eqn:ExpEqMotionTrl_1})\;
				Compute: $\prescript{ (m+1)  }{   }{   \dot{\mathbf{u}}    }^{  trl(k) }_{ \mathcal{D}I  }$ and $\prescript{ (m+1)  }{   }{   \dot{\mathbf{u}}    }^{  cm(k) }_{ \mathcal{D}I  }$ (see Eqs. \eqref{eqn:VelTrl} and \eqref{eqn:CentreMassVel} ). \;
				\If{ Eq. \eqref{eqn:ImpCond} is satisfied at contact grid node $I$ }{
					Compute: $\prescript{ (m)  }{   }{   F    }^{  nor(k) }_{ \mathcal{D}I  }$, $\prescript{ (m)  }{   }{   F    }^{  tan(k) }_{ \mathcal{D}I  }$ and $\prescript{ (m)  }{   }{   \bm{F}    }^{  cont(k) }_{ \mathcal{D}I  }$ (see Eqs. \eqref{eqn:Fnor}, \eqref{eqn:Ftan} and \eqref{eqn:Fcont}) \;
				}
				\textbf{C.2:} Check collinearity conditions for all contact nodes (see Eq. \eqref{ContactForceNormalCond} and Eq. \ref{ContactForceTangCond})\;
				Correct velocities: $\prescript{ (m+1)  }{   }{   \dot{\mathbf{u}}    }^{  (k) }_{ \mathcal{D}I  }$ (see Eq. \eqref{eqn:CorVel})\;  
				\textbf{C.3:} Check impenetrability and complementarity (normal) conditions for all contact nodes (see Eqs. \eqref{ImpenetratibilityCond} and \eqref{ComplementaryNormalCond}) \;
				\textbf{C.4:} Check slip/non-slip and complementarity (tangential) conditions for all contact nodes (see Eqs. \eqref{gsCond} and \eqref{ComplementaryTangCond}) \;
				Compute: $\prescript{ (m+1)  }{   }{   \boldsymbol{\varepsilon}    }^{  (k) }_{ \mathcal{D}p  }$ and $\prescript{ (m+1)  }{   }{   \boldsymbol{\sigma}    }^{  (k) }_{ \mathcal{D}p  }$, for all material points (see Eq. \eqref{eqn:StrainsVelForm} and \eqref{eqn:Stress} ) \;	
				Compute: $\prescript{ (m+1)  }{   }{   \psi    }^{  +(k) }_{ el_{\mathcal{D}p}  }$, for all material points $\rightarrow \prescript{ (m)  }{   }{   \mathcal{H}     }^{ (k)  }_{ \mathcal{D}p  } = \begin{cases}  \prescript{ (m+1)  }{   }{   \psi    }^{  +(k) }_{ el_{\mathcal{D}p}  } , & \text{for } \prescript{ (m+1)  }{   }{   \psi    }^{  +(k) }_{ el_{\mathcal{D}p}  }>\prescript{ (m)  }{   }{   \mathcal{H}     }^{ (k)  }_{ \mathcal{D}p  } \\ \prescript{ (m)  }{   }{   \mathcal{H}     }^{ (k)  }_{ \mathcal{D}p  } , & \text{otherwise }\\ \end{cases}$\;
				Compute Residual (Phase-Field): $\prescript{ (m)  }{   }{    \bm{R}    }^{ c(k)  }_{   }$ (see Eq. \eqref{eqn:DiscreGalerkinPFGridArbi_5}) according to $\prescript{ (m)  }{   }{   c     }^{  (k) }_{ \mathcal{D}p  }$, $\prescript{ (m)  }{   }{   \nabla c     }^{  (k) }_{ \mathcal{D}p  }$, $\prescript{ (m)  }{   }{   \Delta c     }^{  (k) }_{ \mathcal{D}p  }$, $\prescript{ (m)  }{   }{   g     }^{  (k) }_{ \mathcal{D}p  }$ \;
				Convergence Check (Phase Field): If $ \| \prescript{ (m)  }{   }{    \bm{R}    }^{ c(k)  }_{   } \| \leq tol_{c} $ or $k \geq N_{staggs}$ then "exit" from loop else $k =k + 1$ go to next stagger iteration. \;
			}
			Update material point properties: $\prescript{ (m+1)  }{   }{   \mathbf{u}     }^{  c(k) }_{ \mathcal{D}p  }$, $\prescript{ (m+1)  }{   }{   \dot{\mathbf{u}}     }^{  c(k) }_{ \mathcal{D}p  }$ and $\prescript{ (m)  }{   }{   \ddot{\mathbf{u}}     }^{  c(k) }_{ \mathcal{D}p  }$ (see Eqs. \eqref{eqn:MpDisp}, \eqref{eqn:MpVel}, \eqref{eqn:MpAcc} and \eqref{eqn:MpPosition}). \;
			Update material point history field: $\prescript{ (m+1)  }{   }{   \mathcal{H}     }^{   }_{ \mathcal{D}p  }=\prescript{ (m)  }{   }{   \mathcal{H}     }^{   }_{ \mathcal{D}p  }$ \;
		}
		\caption{Anisotropic Phase-Field Material Point Method pseudo-code for impact-fracture problems (Staggered Solution Algorithm with Explicit time integration).}
		\label{alg:PFMPMStaggExplicit}
	\end{small}
\end{algorithm}

\section{Numerical examples} \label{sec:NumericalExamples}

In this section, a set of two-dimensional numerical examples is presented. The numerical examples demonstrate the accuracy of the proposed PF-MPM against the standard PF-FEM as well as its computational efficiency in impact-fracture problems. Both the isotropic and anisotropic phase field models are examined within both single and multi discrete-field examples. Quadratic B-splines ( $C^{1}$ ) are utilized for the background grid. The initial cell density is chosen to be at least $3x3=9$ material points per cell element. Extensive numerical experiments performed in this work have demonstrated that this cell density results in accurate estimates when quadratic basis functions are utilized. Higher order B-splines are employed not only to compute the anisotropic phase field matrix in Eq. \eqref{eqn:PhaseFieldMatrix}, but also to treat the ``cell-crossing error'' of the Material Point Method \cite{Gan_2017}. In all cases examined in this section, stability of the explicit integration scheme is established on the basis of the following upper bound for the time increment $\Delta t$ 
\begin{equation}
\Delta t \le \tilde{\Delta t_{cr}}
\label{eqn:TimeCrFact}
\end{equation}
where 
\begin{equation}\label{eqn:TimeCrFact2}
\tilde{\Delta t_{cr}}=\alpha_{c} \cdot \Delta t_{cr}
\end{equation}
and $\Delta t_{cr}$ corresponds to the critical time step prescribed by the Courant – Friedrichs – Lewy (CFL) condition. Parameter $\alpha_{c} \in \left[0.8, 0.98\right]$  in Eq. \eqref{eqn:TimeCrFact} depends on the nonlinearities of the system \cite{zhang_2016}. In all the numerical experiments, we consider $\alpha_{c}=0.80$. 

In all numerical experiments presented the phase field residual tolerance was set to $tol_c = 10^{-6}$ and a single stagger iteration was required for solution convergence. This is due to the small time step, imposed by the stability requirements (Eq. \eqref{eqn:TimeCrFact}).
\subsection{Plate under impact loading} \label{lbl:PuIL}

\begin{figure}
	\centering
	\begin{tabular}{lcr}
		{\subfloat[\label{fig:PuIL_GeometryBCs}]{
				\includegraphics[width=0.30\columnwidth]{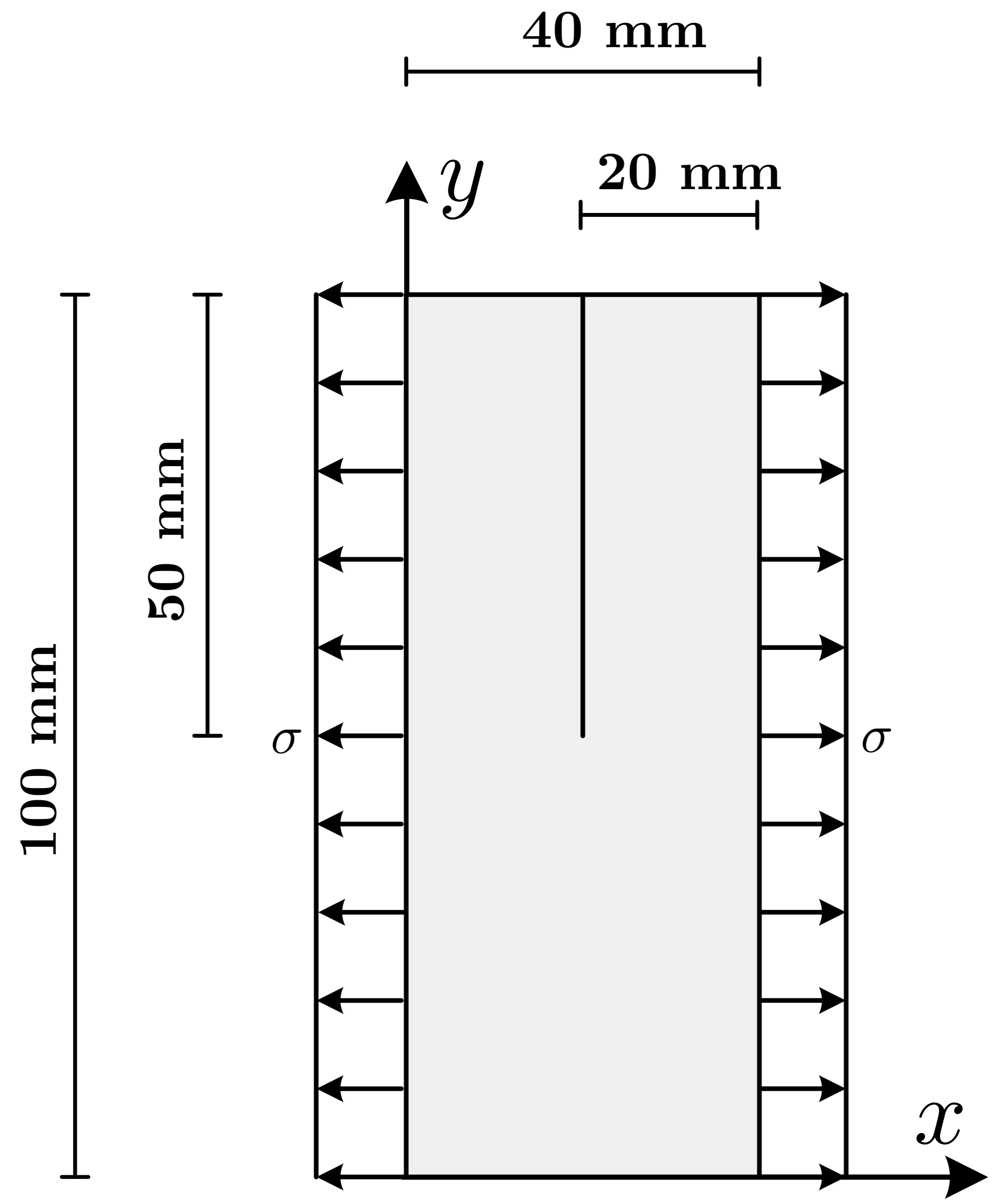}}} &
		\subfloat[\label{fig:PuIL_SurfEng}]{
			\includegraphics[width=0.30\columnwidth]{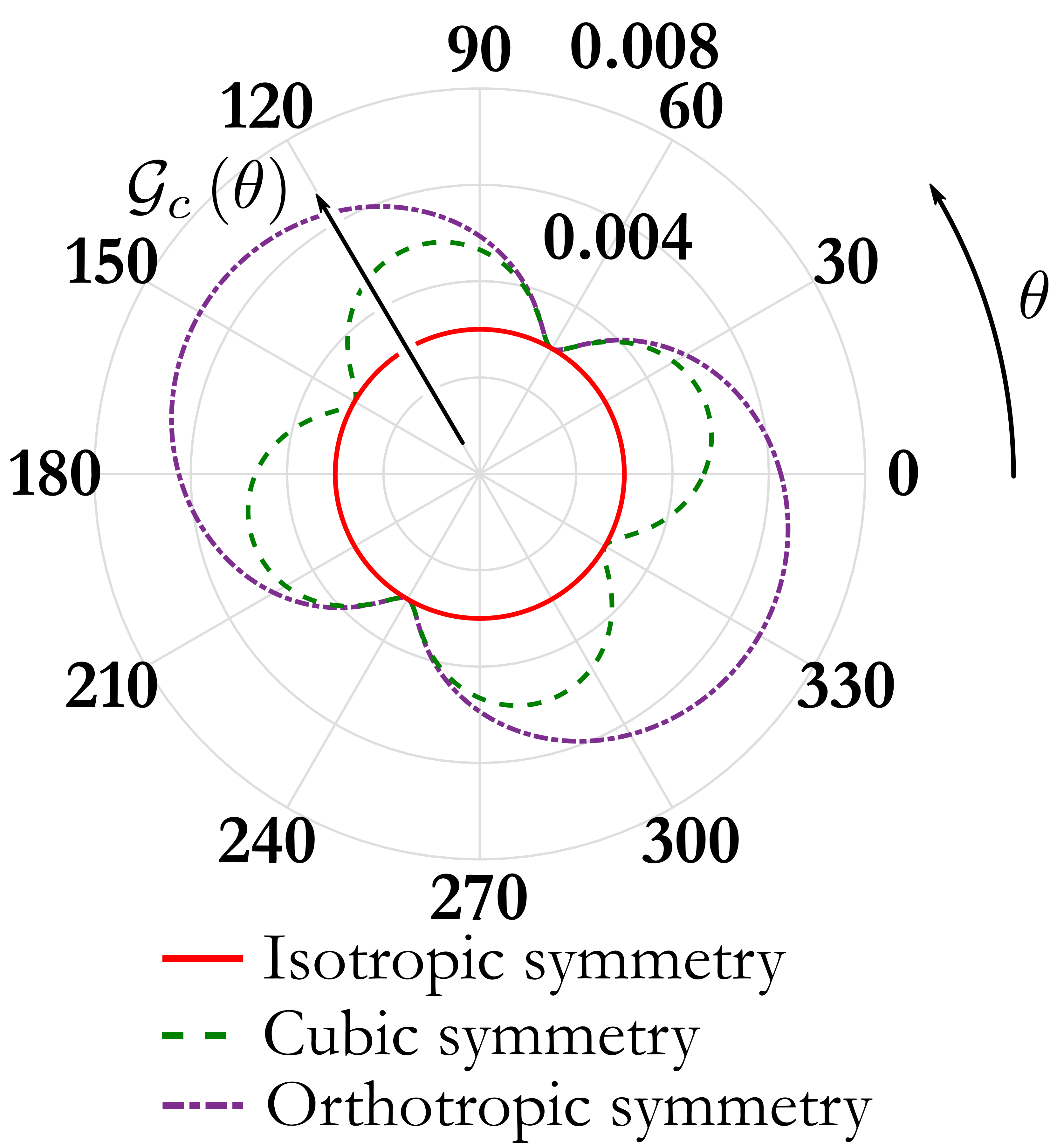}} &
		{\subfloat[\label{fig:PuIL_InvSurfEng}]{
				\includegraphics[width=0.30\columnwidth]{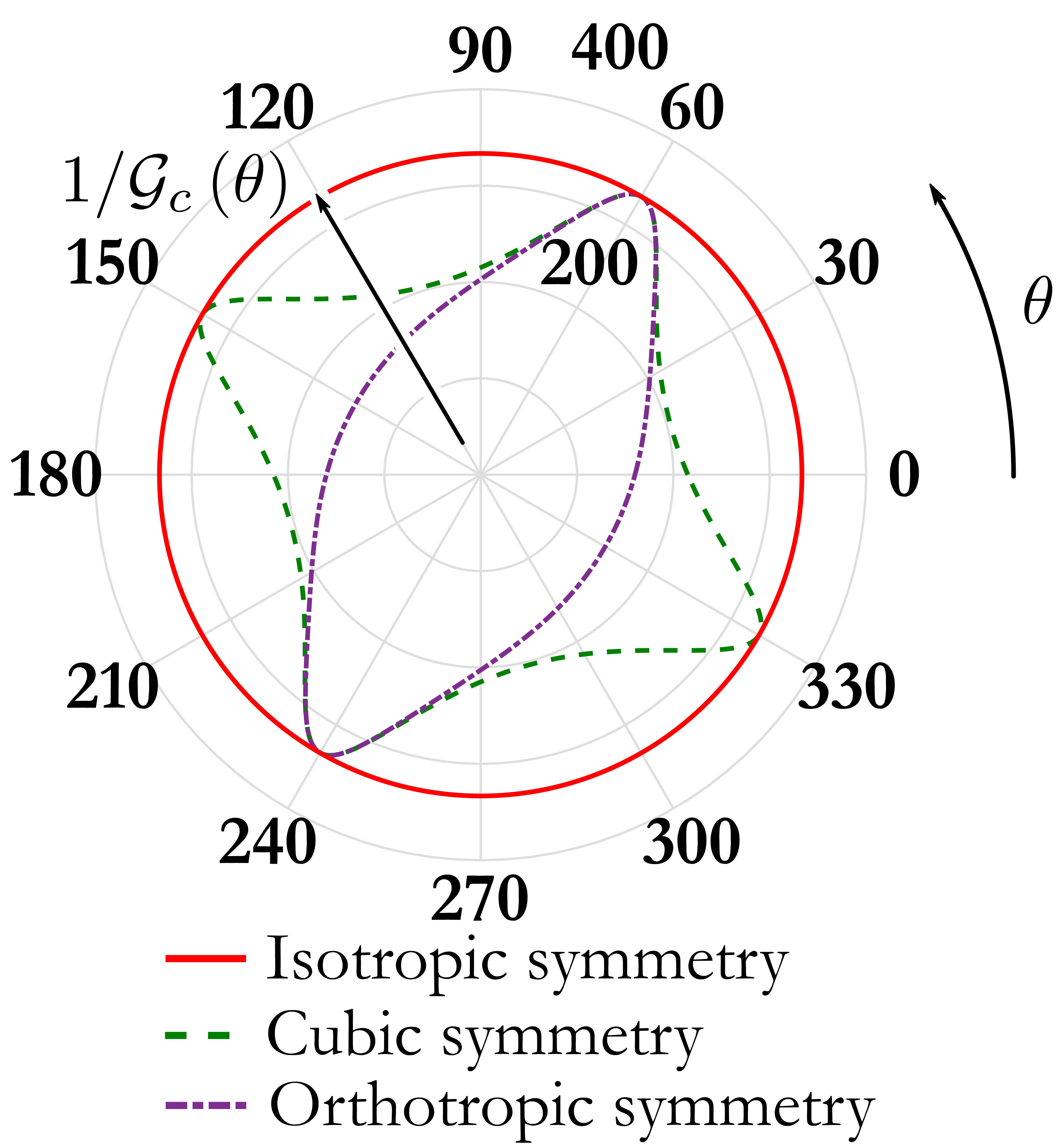}}}
	\end{tabular}
	\caption[]{Plate under impact loading: \subref{fig:PuIL_GeometryBCs} Geometry and boundary conditions. \subref{fig:PuIL_SurfEng} Surface energy densities $\mathcal{G}_c \left( \theta \right)$  and \subref{fig:PuIL_InvSurfEng} their reciprocals $1/\mathcal{G}_c \left( \theta \right)$ for material orientation $\phi=+30^{o}$ (with respect to $x$ axis (clockwise)) in polar coordinates.}
	\label{fig:PuIL_1}
\end{figure}

A plate under impact loading is examined. The same problem has been previously analysed by Borden et al. \cite{borden_2012} with a finite element phase field implementation, considering a second order isotropic phase field formulation. The geometry and boundary conditions are presented in Fig. \ref{fig:PuIL_GeometryBCs}. Herein, three cases are considered, i.e., (i) isotropic symmetry, (ii) cubic symmetry, and (iii) orthotropic symmetry. The material orientation is considered to be $\phi=+30^{o}$ with respect to $x$ axis (clockwise) as shown in Fig. \ref{fig:PuIL_ElastEng_Borden_C1SecondMPM}.

The cell (patch) spacing is $h=0.125$ mm and plane strain conditions are assumed. The grid is formed by two knot vectors $\splitatcommas{\Xi=\{0,0,0,0.001240,0.002481,...,0.997518,0.998759,1,1,1\}}$ and $\splitatcommas{H = \{0,0,0,0.003067,0.006134,...,0.993865,0.996932,1,1,1\}}$, $265024$ control points and $806x326=262756$ cells. The total number of material points is $2304000$ and the elastic material properties are $E=32000$ N/mm\textsuperscript{2}, $\nu=0.20$ and $\rho=2450$ kg/m\textsuperscript{3}.
	
The length scale parameter is chosen to be $l_0=0.25$ mm and $k_f=0.00$. In case (i), all the anisotropic material parameters are chosen such that $\gamma_{ijkl}=0$. Hence, the anisotropic phase field model reduces to the second order isotropic case. The maximum and minimum surface energy densities are equal to $\mathcal{G}_{c} \left( \theta \right)=\bar{\mathcal{G}}_c=\mathcal{G}_{c_{max}}=\mathcal{G}_{c_{min}}=0.003$ N/mm. In case (ii) cubic symmetry of the surface energy density is considered with $\bar{\mathcal{G}}_c=0.002121$ N/mm and anisotropic parameters $\gamma_{1111}=\gamma_{2222}=1.00$, $\gamma_{1122}=0.00$ and $\gamma_{1212}=74.00$. These parameters result into maximum and minimum surface energy densities $\mathcal{G}_{c_{max}}=0.0049$ N/mm and $\mathcal{G}_{c_{min}}=0.003$ N/mm, respectively. In case (iii) the anisotropic parameter $\gamma_{2222}$ is increased to $\gamma_{2222}=80.00$ giving rise to orthotropic symmetry with maximum and minimum surface energy densities $\mathcal{G}_{c_{max}}=0.0067$ N/mm and $\mathcal{G}_{c_{min}}=0.003$ N/mm, respectively. 

The surface energy densities and their reciprocals for material orientation $\phi=+30^{o}$ are shown in Fig. \ref{fig:PuIL_SurfEng} and Fig. \ref{fig:PuIL_InvSurfEng}, respectively. In cases (ii) and (iii), the parameter $\bar{\mathcal{G}}_c$ is chosen so that $\mathcal{G}_{c_{min}}=0.003$ N/mm and to facilitate comparisons between all cases (see also, Figs. \ref{fig:PuIL_SurfEng} and \ref{fig:PuIL_InvSurfEng}).

A single discrete field is considered in this example. The solution procedure is implemented with a time step $\Delta t = 0.025$ $\mu$s for $N_{steps}=3200$ steps. The critical time step is $\tilde{\Delta t_{cr}}=0.026$ $\mu$s. The traction is considered to be constant $\sigma = 1$ N/mm\textsuperscript{2} during the analysis. The initial crack is modelled by introducing an initial history field at the corresponding material points as in \cite{borden_2012}. The Rayleigh wave speed is $\dot{u}_{R}=2125$ m/s for the material parameters of that specimen \cite{Freund_1998}.

\subsubsection{Case (i): Isotropy}

Initially, the PF-MPM is compared against the Phase Field Finite Element Method (PF-FEM) with the results obtained in \cite{borden_2012} for the same cell (patch) spacing $h=0.125$ mm. 

\begin{figure}
	\centering
	\begin{tabular}{ccc}
		\subfloat[\label{fig:PuIL_ElastEng_Borden_C1SecondMPM}]{
			\includegraphics[width=0.37\columnwidth]{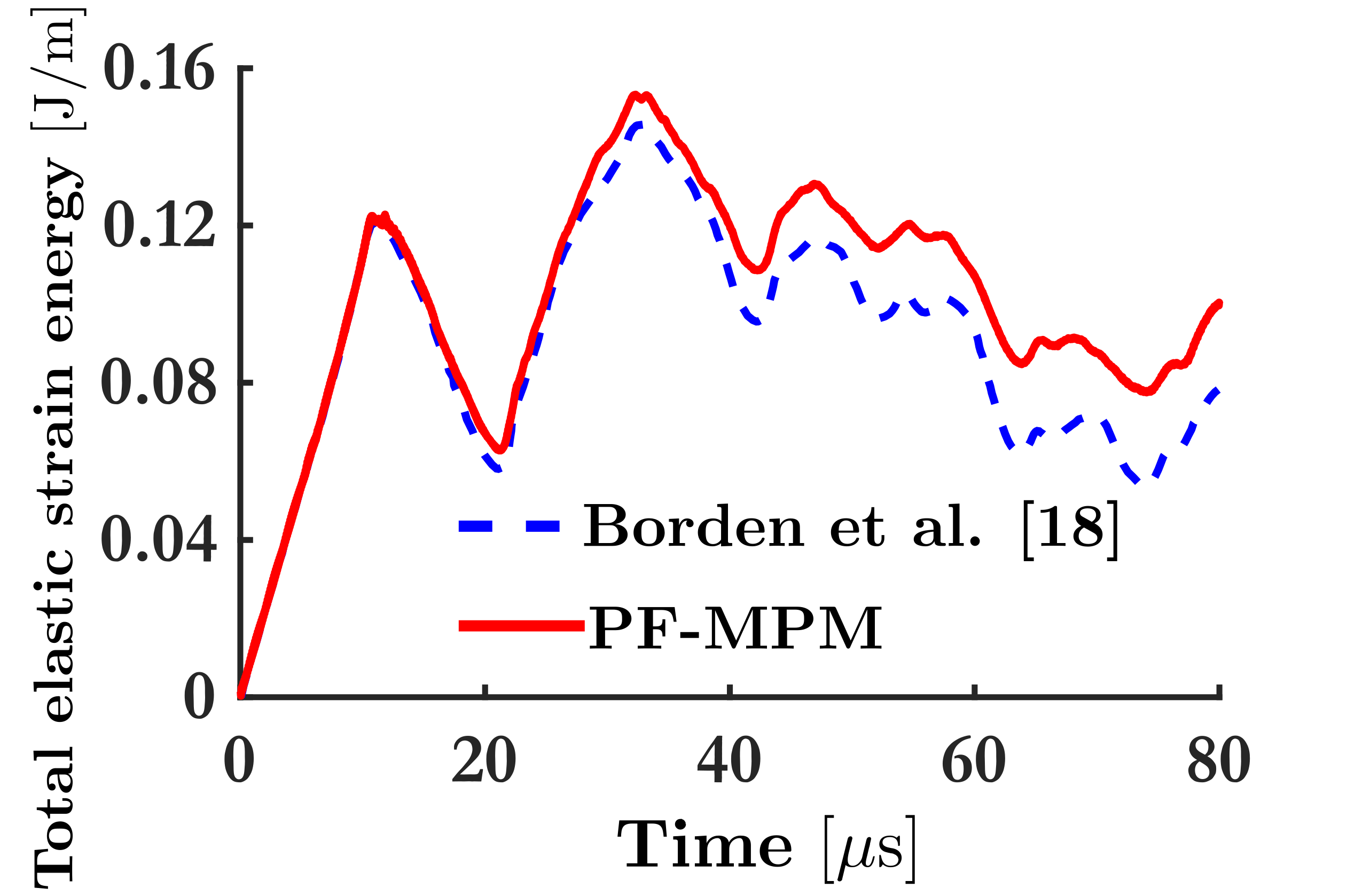}} &
		{\subfloat[\label{fig:PuIL_FractEng_Borden_C1SecondMPM}]{
				\includegraphics[width=0.37\columnwidth]{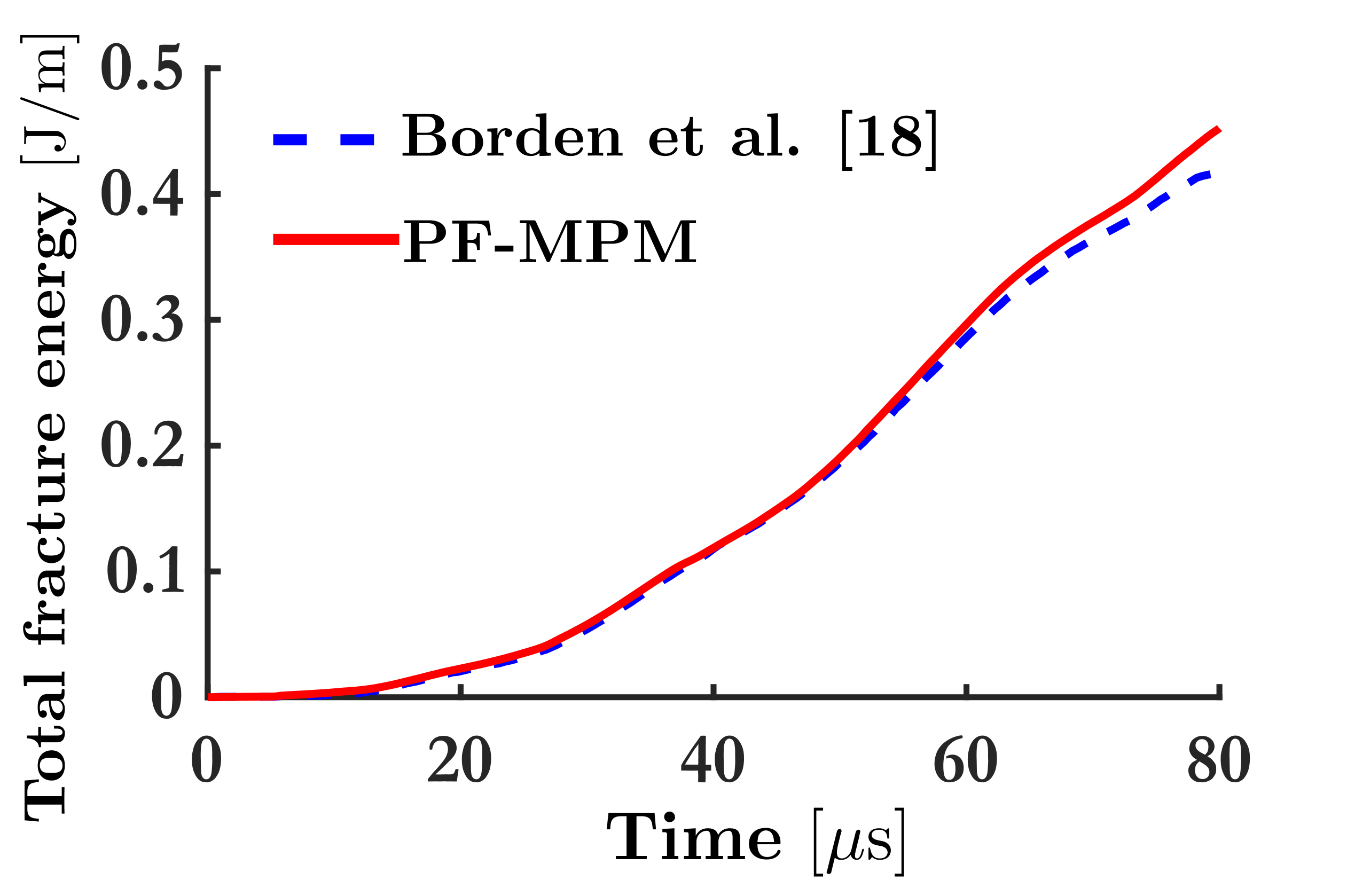}}} &
		{\subfloat[\label{fig:PuIL_CrackVel_Borden_C1SecondMPM}]{
				\includegraphics[width=0.37\columnwidth]{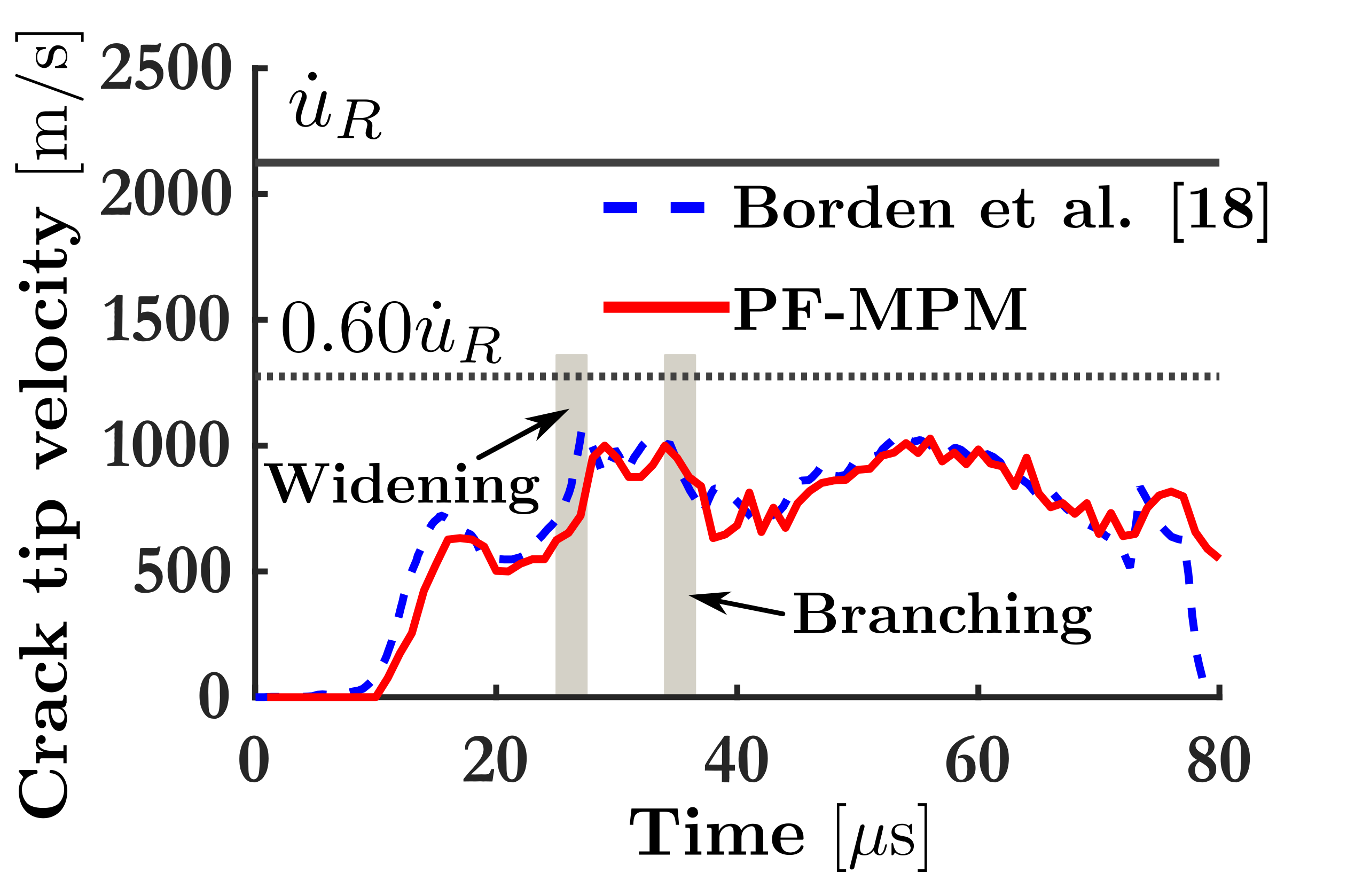}}} 
	\end{tabular}
	\caption[]{Plate under impact loading: \subref{fig:PuIL_ElastEng_Borden_C1SecondMPM} Total elastic strain energies, \subref{fig:PuIL_FractEng_Borden_C1SecondMPM} Total fracture energies and \subref{fig:PuIL_CrackVel_Borden_C1SecondMPM} Crack tip velocities time histories for Borden et al. \cite{borden_2012}  and PF-MPM 2nd order isotropic model (case (i)). The traction is considered to be $\sigma = 1$ N/mm\textsuperscript{2}.}
	\label{fig:Borden_C1SecondMPM}
\end{figure}

The total energy time-histories for the two solutions are shown in Figs. \ref{fig:PuIL_ElastEng_Borden_C1SecondMPM} and \ref{fig:PuIL_FractEng_Borden_C1SecondMPM} where the 2 methods demonstrate a very good agreement. The total fracture energy results (see Fig. \ref{fig:PuIL_FractEng_Borden_C1SecondMPM}) are in perfect agreement with the results reported in \cite{borden_2012}. Minor differences are observed, especially for time $t>50$ $\mu$s. The total elastic strain energies (see Fig. 10a) also demonstrate very good agreement with some differences after $t>30$ $\mu$s. 

\begin{figure}
	\centering
	\begin{tabular}{cccc}
		\subfloat[\label{fig:PuIL_C1SecondMPM_0sec}]{
			\includegraphics[width=0.20\columnwidth]{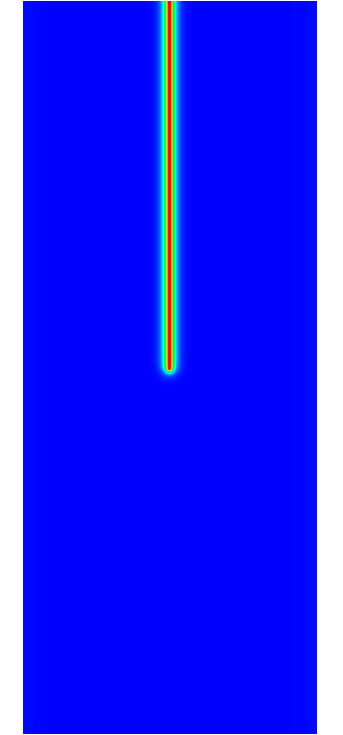}} &
		\subfloat[\label{fig:PuIL_C1SecondMPM_50sec}]{
			\includegraphics[width=0.20\columnwidth]{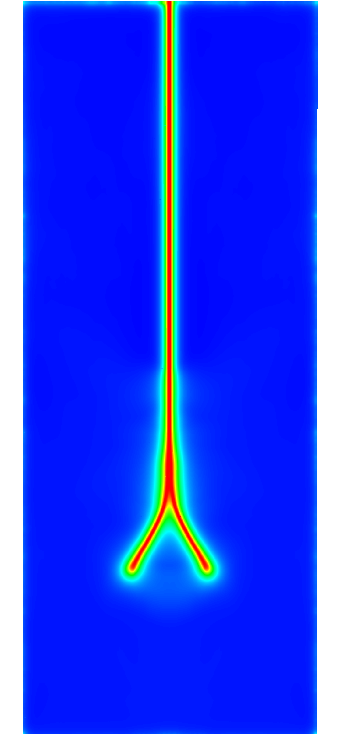}} &
		\subfloat[\label{fig:PuIL_C1SecondMPM_65sec}]{
			\includegraphics[width=0.20\columnwidth]{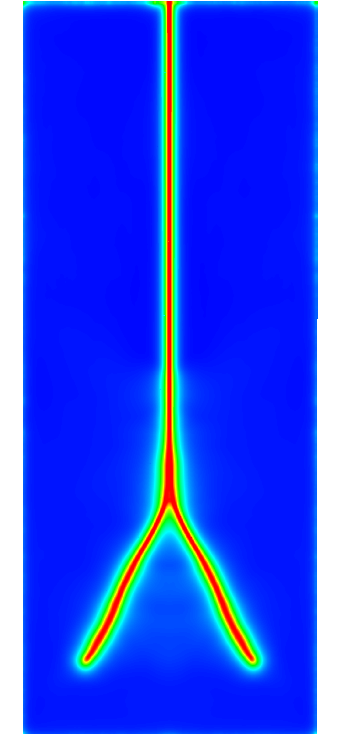}} &
		\subfloat[\label{fig:PuIL_C1SecondMPM_80sec}]{
			\includegraphics[width=0.20\columnwidth]{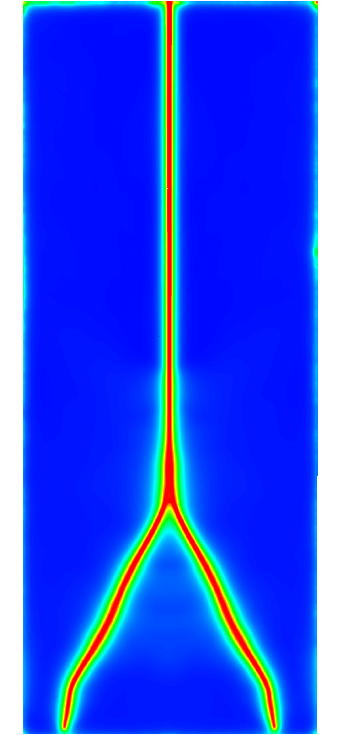}} \\
		\multicolumn{4}{c}{\subfloat{
				\includegraphics[width=0.75\columnwidth]{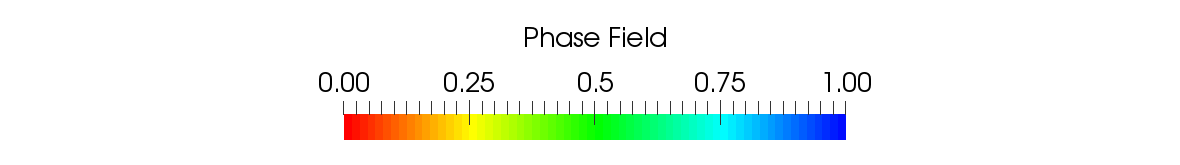}}}		
	\end{tabular}
	\caption[]{Plate under impact loading: Phase field for time steps \subref{fig:PuIL_C1SecondMPM_0sec} t=0 $\mu$s \subref{fig:PuIL_C1SecondMPM_50sec} t=50 $\mu$s \subref{fig:PuIL_C1SecondMPM_50sec} t=65 $\mu$s and \subref{fig:PuIL_C1SecondMPM_50sec} t=80 $\mu$s. Results for case (i): 2nd order isotropic phase field model and $\sigma = 1$ N/mm\textsuperscript{2}.}
	\label{fig:PuIL_PF_C1SecondMPM}
\end{figure}

\begin{figure}
	\centering
	\begin{tabular}{lrrr}
		\subfloat[\label{fig:PuIL_HS_C1SecondMPM_0sec}]{
			\includegraphics[width=0.20\columnwidth]{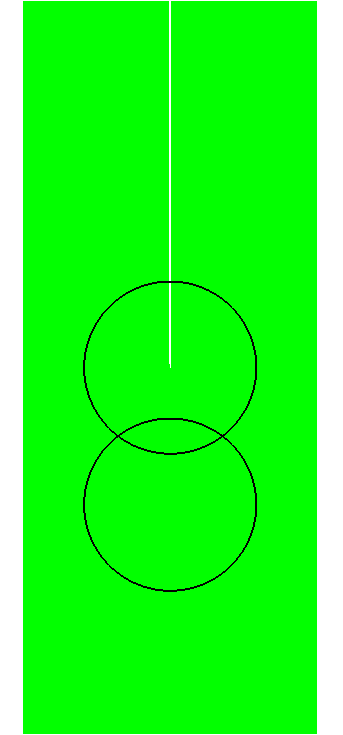}} &
		\subfloat[\label{fig:PuIL_HS_C1SecondMPM_50sec}]{
			\includegraphics[width=0.20\columnwidth]{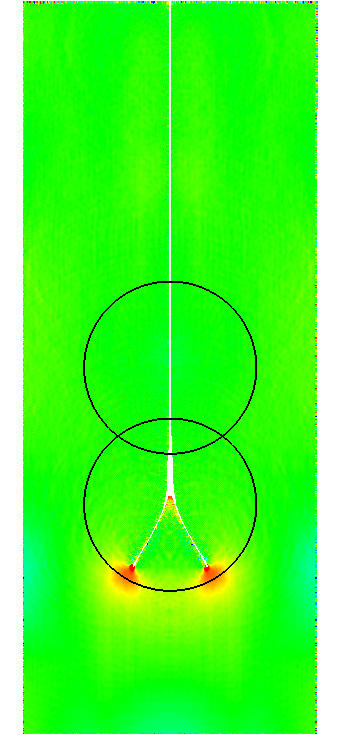}} &
		\subfloat[\label{fig:PuIL_HS_C1SecondMPM_65sec}]{
			\includegraphics[width=0.20\columnwidth]{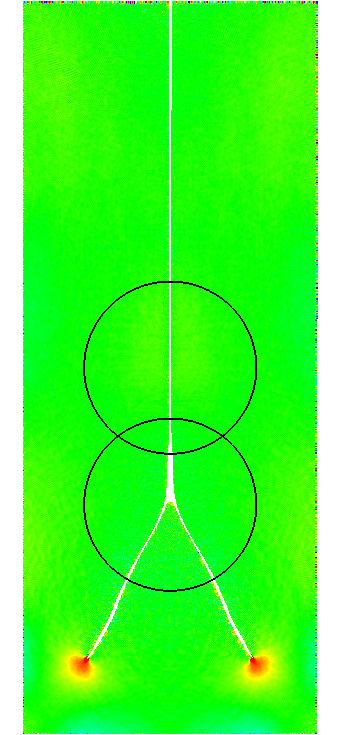}} &
		\subfloat[\label{fig:PuIL_HS_C1SecondMPM_80sec}]{
			\includegraphics[width=0.20\columnwidth]{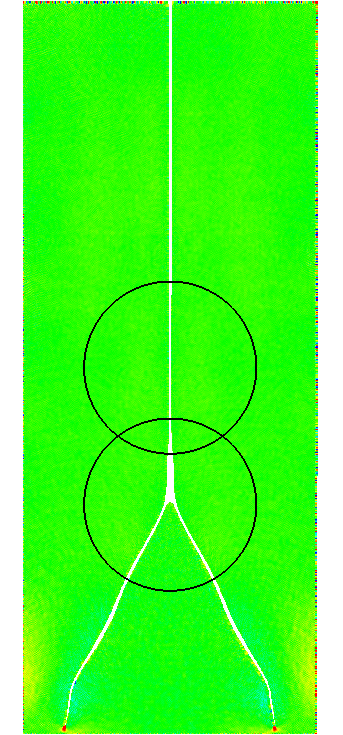}} \\
		\multicolumn{4}{c}{\subfloat{
				\includegraphics[width=0.75\columnwidth]{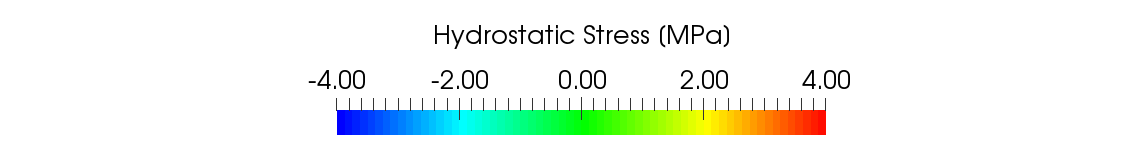}}}				
	\end{tabular}
	\caption[]{Plate under impact loading: Hydrostatic stress for time steps \subref{fig:PuIL_HS_C1SecondMPM_0sec} t=0 $\mu$s \subref{fig:PuIL_HS_C1SecondMPM_50sec} t=50 $\mu$s \subref{fig:PuIL_HS_C1SecondMPM_65sec} t=65 $\mu$s and \subref{fig:PuIL_HS_C1SecondMPM_80sec} t=80 $\mu$s. Results for case (i): 2nd order isotropic phase field model and $\sigma = 1$ N/mm\textsuperscript{2}. Material points with $c_p<0.10$ have been removed.}
	\label{fig:PuIL_HS_C1SecondMPM}
\end{figure}

The evolution of the phase field is presented in Fig. \ref{fig:PuIL_PF_C1SecondMPM} for specific time steps. In these, the occurrence of a branched crack is observed at approximately $t=35$ $\mu$s. The evolution of the hydrostatic stress for the same timesteps is shown in Fig. \ref{fig:PuIL_HS_C1SecondMPM}. To demonstrate the influence of the surface energy density into the resulting crack paths, the reciprocal of the surface energy density is also plotted (black circle); the hydrostatic stresses are also shown in the same figure. Since the surface energy density is isotropic, hence independent of the material orientation, the crack naturally propagates along the vertical axis (see Figs. \ref{fig:PuIL_C1SecondMPM_50sec} and \ref{fig:PuIL_HS_C1SecondMPM_50sec}) until branching occurs. Crack branching is perfectly symmetrical due to structure, load symmetry and the isotropic phase field model.

The crack tip velocities for the two methods are presented in Fig. \ref{fig:PuIL_CrackVel_Borden_C1SecondMPM}. As already mentioned in numerous works (see \cite{borden_2012}, \cite{hofacker_2013}, \cite{Schluter_2014}), the crack tip and the exact location of crack branching cannot be identified uniquely due to the smooth description of the crack. Therefore, the crack tip velocity is measured with the methodology employed in \cite{borden_2012} to facilitate verification. The results of both methods illustrate very good agreement. The crack widening and branching regions are almost the same for the two solutions and they are also shown in Fig. \ref{fig:PuIL_CrackVel_Borden_C1SecondMPM}. Crack widening here refers to the broadening of the damage zone prior to branching in accordance with the definition introduced in \cite{borden_2012}.

In Fig. \ref{fig:PuIL_CrackVel_Borden_C1SecondMPM}, the crack tip velocity is clearly below the Rayleigh wave speed which stands for the crack speed limit as elaborated by \cite{Freund_1998} \cite{Ravi_1998}. However, experimental studies have shown that cracks rarely propagate at speeds close to the Rayleigh wave speed. In fact, they propagate at a fraction of the Rayleigh wave speed, i.e. 60\% $\dot{u}_{R}$ \cite{Ravi_1984}. As shown in Fig. \ref{fig:PuIL_CrackVel_Borden_C1SecondMPM}, the resulting crack tip velocities are below this limit.

\subsubsection{Case (ii): Cubic Symmetry}

\begin{figure}
	\centering
	\begin{tabular}{cccc}
		\subfloat[\label{fig:PuIL_C1FourthCubicMPM_0sec}]{
			\includegraphics[width=0.20\columnwidth]{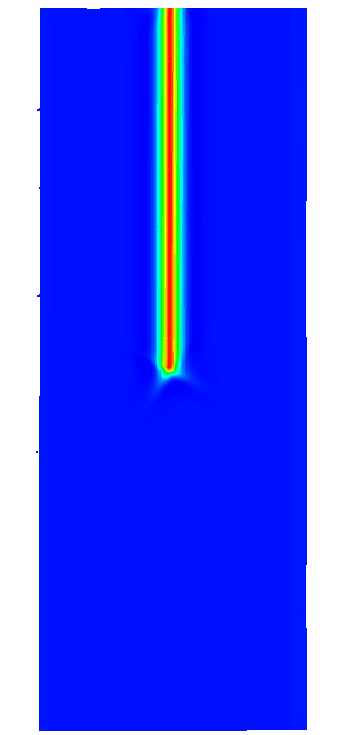}} &
		\subfloat[\label{fig:PuIL_C1FourthCubicMPM_50sec}]{
			\includegraphics[width=0.20\columnwidth]{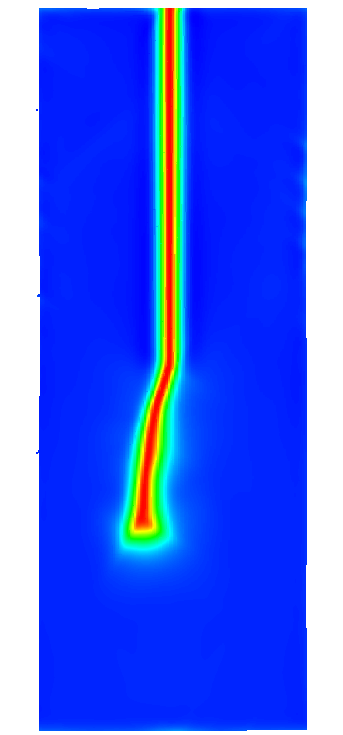}} &
		\subfloat[\label{fig:PuIL_C1FourthCubicMPM_65sec}]{
			\includegraphics[width=0.20\columnwidth]{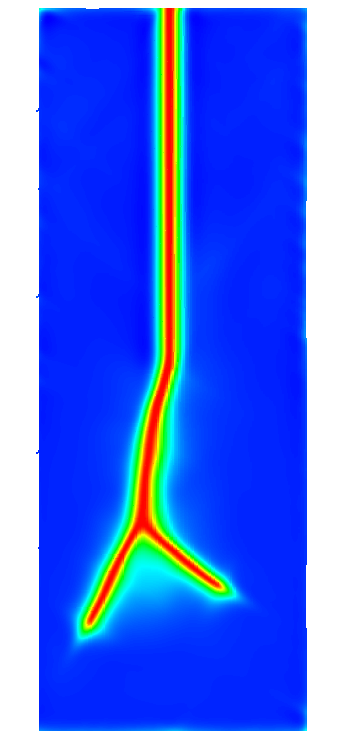}} &
		\subfloat[\label{fig:PuIL_C1FourthCubicMPM_80sec}]{
			\includegraphics[width=0.20\columnwidth]{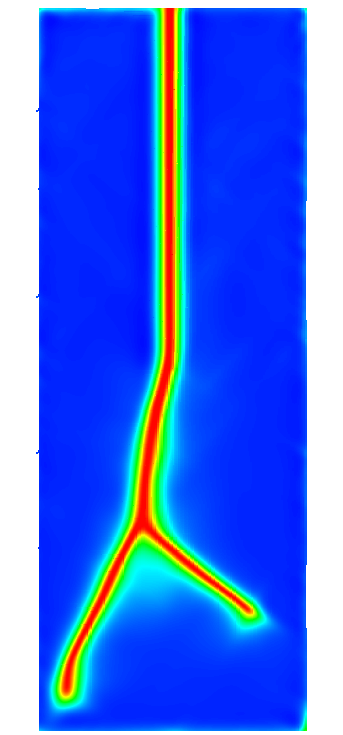}} \\
		\multicolumn{4}{c}{\subfloat{
						\includegraphics[width=0.75\columnwidth]{PF_ColorBar.png}}}			
	\end{tabular}
	\caption[]{Plate under impact loading: Phase field for time steps \subref{fig:PuIL_C1FourthCubicMPM_0sec} t=0 $\mu$s \subref{fig:PuIL_C1FourthCubicMPM_50sec} t=50 $\mu$s \subref{fig:PuIL_C1FourthCubicMPM_65sec} t=65 $\mu$s and \subref{fig:PuIL_C1FourthCubicMPM_80sec} t=80 $\mu$s. Results for case (ii): 4th order anisotropic cubic phase field model and $\sigma = 1$ N/mm\textsuperscript{2}.}
	\label{fig:PuIL_PF_C1FourthCubicMPM}
\end{figure}

\begin{figure}
	\centering
	\begin{tabular}{cccc}
		\subfloat[\label{fig:PuIL_HS_C1FourthCubicMPM_0sec}]{
			\includegraphics[width=0.20\columnwidth]{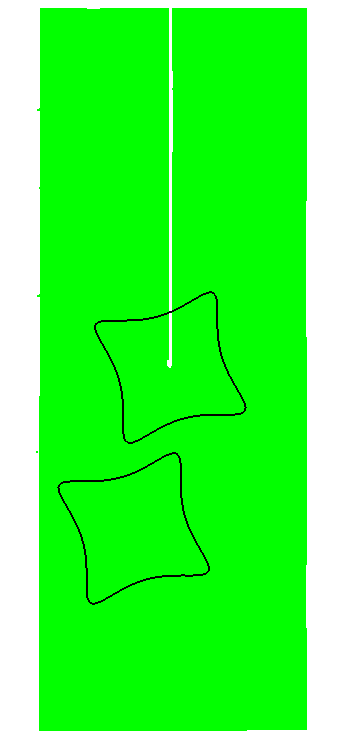}} &
		\subfloat[\label{fig:PuIL_HS_C1FourthCubicMPM_50sec}]{
			\includegraphics[width=0.20\columnwidth]{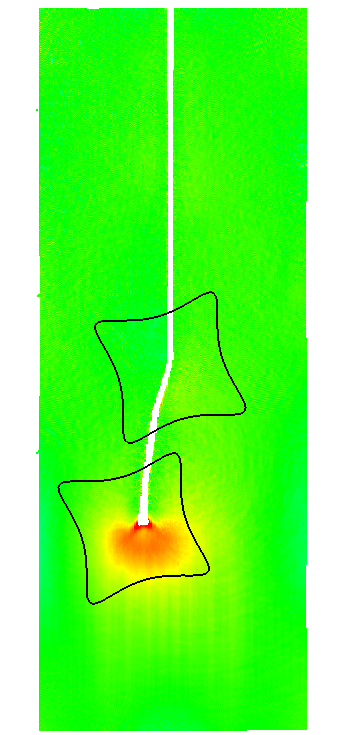}} &
		\subfloat[\label{fig:PuIL_HS_C1FourthCubicMPM_65sec}]{
			\includegraphics[width=0.20\columnwidth]{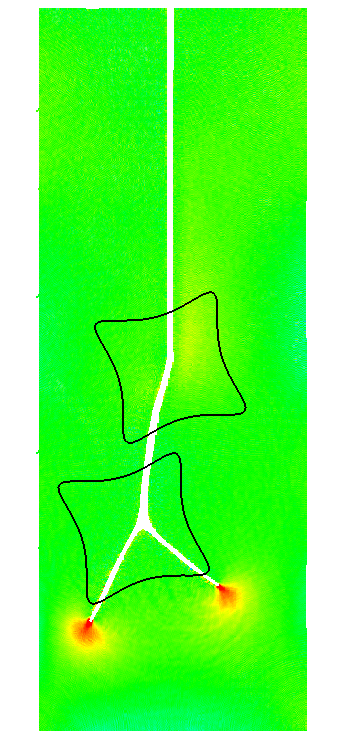}} &
		\subfloat[\label{fig:PuIL_HS_C1FourthCubicMPM_80sec}]{
			\includegraphics[width=0.20\columnwidth]{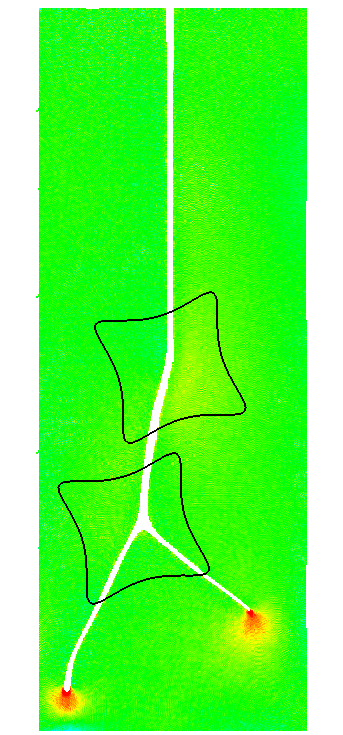}} \\
		\multicolumn{4}{c}{\subfloat{
				\includegraphics[width=0.75\columnwidth]{HS_4_ColorBar.png}}}						
	\end{tabular}
	\caption[]{Plate under impact loading: Hydrostatic stress for time steps \subref{fig:PuIL_HS_C1FourthCubicMPM_0sec} t=0 $\mu$s \subref{fig:PuIL_HS_C1FourthCubicMPM_50sec} t=50 $\mu$s \subref{fig:PuIL_HS_C1FourthCubicMPM_65sec} t=65 $\mu$s and \subref{fig:PuIL_HS_C1FourthCubicMPM_80sec} t=80 $\mu$s. Results for case (ii): 4th order anisotropic cubic phase field model and $\sigma = 1$ N/mm\textsuperscript{2}. Material points with $c_p<0.10$ have been removed.}
	\label{fig:PuIL_HS_C1FourthCubicMPM}
\end{figure}

The evolution of the phase field and the hydrostatic stress for specific time steps is presented in Fig. \ref{fig:PuIL_PF_C1FourthCubicMPM} and \ref{fig:PuIL_HS_C1FourthCubicMPM}, respectively. From Figs. \ref{fig:PuIL_C1FourthCubicMPM_50sec} and \ref{fig:PuIL_HS_C1FourthCubicMPM_50sec}, it is observed that due to the anisotropic cubic model and the material orientation $\phi=+30^{o}$ the crack does not initiate along the vertical axis. The crack propagates until the crack branches at approximately $t=50$ $\mu$s. In this case, the branched crack is not symmetrical and it branches along its two preferential weak directions (see Figs. \ref{fig:PuIL_HS_C1FourthCubicMPM_65sec} and \ref{fig:PuIL_HS_C1FourthCubicMPM_80sec}). To further illustrate this, the reciprocal of the surface energy density is also plotted in Fig. \ref{fig:PuIL_HS_C1FourthCubicMPM}.

\subsubsection{Case (iii): Orthotropic Symmetry}

\begin{figure}
	\centering
	\begin{tabular}{cccc}
		\subfloat[\label{fig:PuIL_C1FourthOrthoMPM_0sec}]{
			\includegraphics[width=0.20\columnwidth]{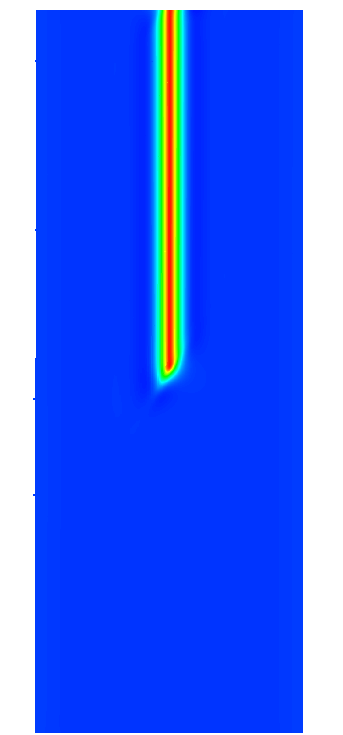}} &
		\subfloat[\label{fig:PuIL_C1FourthOrthoMPM_50sec}]{
			\includegraphics[width=0.20\columnwidth]{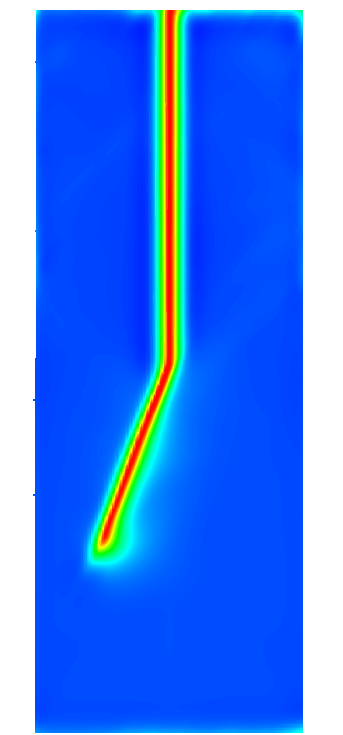}} &
		\subfloat[\label{fig:PuIL_C1FourthOrthoMPM_65sec}]{
			\includegraphics[width=0.20\columnwidth]{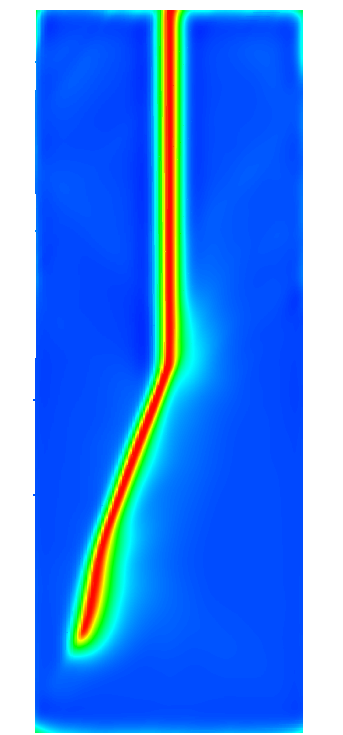}} &
		\subfloat[\label{fig:PuIL_C1FourthOrthoMPM_80sec}]{
			\includegraphics[width=0.20\columnwidth]{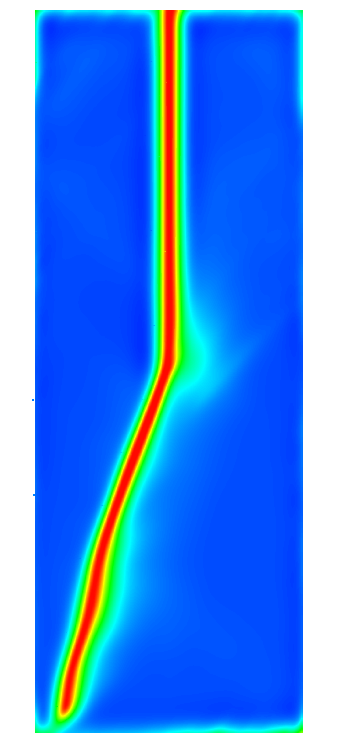}} \\
		\multicolumn{4}{c}{\subfloat{
				\includegraphics[width=0.75\columnwidth]{PF_ColorBar.png}}}		
	\end{tabular}
	\caption[]{Plate under impact loading: Phase field for time steps \subref{fig:PuIL_C1FourthOrthoMPM_0sec} t=0 $\mu$s \subref{fig:PuIL_C1FourthOrthoMPM_50sec} t=50 $\mu$s \subref{fig:PuIL_C1FourthOrthoMPM_65sec} t=65 $\mu$s and \subref{fig:PuIL_C1FourthOrthoMPM_80sec} t=80 $\mu$s. Results for case (iii): 4th order anisotropic orthotropic phase field model and $\sigma = 1$ N/mm\textsuperscript{2}.}
	\label{fig:PuIL_PF_C1FourthOrthoMPM}
\end{figure}

\begin{figure}
	\centering
	\begin{tabular}{lrrr}
		\subfloat[\label{fig:PuIL_HS_C1FourthOrthoMPM_0sec}]{
			\includegraphics[width=0.20\columnwidth]{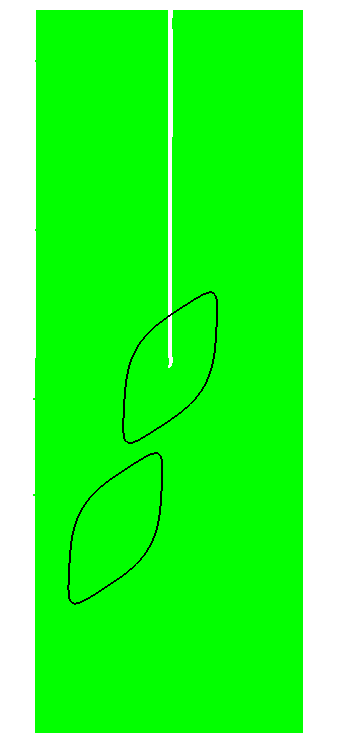}} &
		\subfloat[\label{fig:PuIL_HS_C1FourthOrthoMPM_50sec}]{
			\includegraphics[width=0.20\columnwidth]{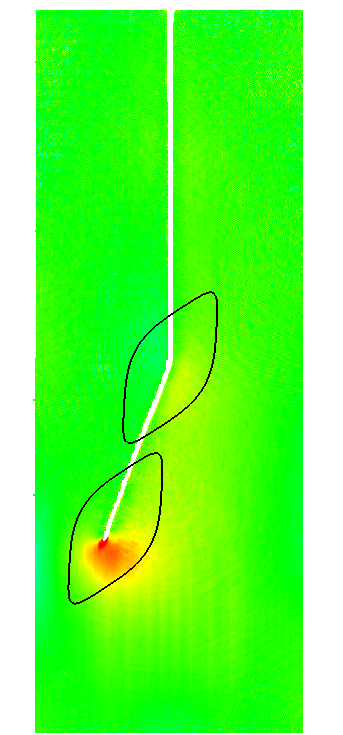}} &
		\subfloat[\label{fig:PuIL_HS_C1FourthOrthoMPM_65sec}]{
			\includegraphics[width=0.20\columnwidth]{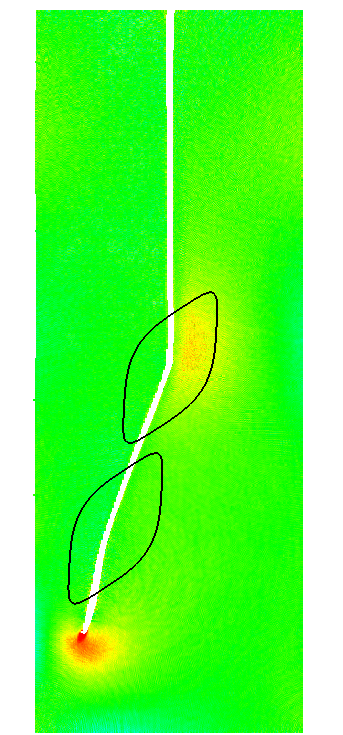}} &
		\subfloat[\label{fig:PuIL_HS_C1FourthOrthoMPM_80sec}]{
			\includegraphics[width=0.20\columnwidth]{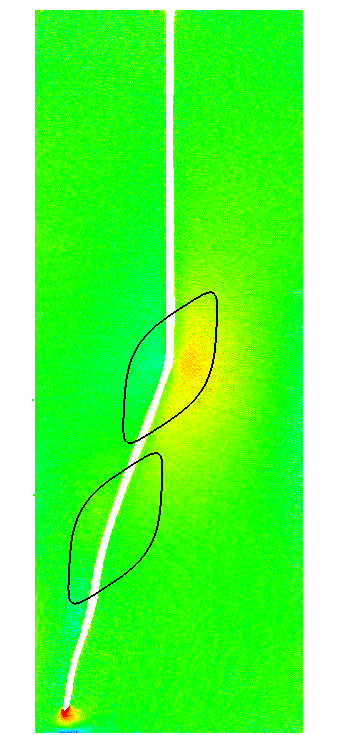}} \\
		\multicolumn{4}{c}{\subfloat{
				\includegraphics[width=0.75\columnwidth]{HS_4_ColorBar.png}}}					
	\end{tabular}
	\caption[]{Plate under impact loading: Hydrostatic stress for time steps \subref{fig:PuIL_HS_C1FourthOrthoMPM_0sec} t=0 $\mu$s \subref{fig:PuIL_HS_C1FourthOrthoMPM_50sec} t=50 $\mu$s \subref{fig:PuIL_HS_C1FourthOrthoMPM_65sec} t=65 $\mu$s and \subref{fig:PuIL_HS_C1FourthOrthoMPM_80sec} t=80 $\mu$s. Results for case (iii): 4th order anisotropic orthotropic phase field model and $\sigma = 1$ N/mm\textsuperscript{2}. Material points with $c_p<0.10$ have been removed.}
	\label{fig:PuIL_HS_C1FourthOrthoMPM}
\end{figure}

The evolution of the phase field and the hydrostatic stresses for several time steps in case (iii) are shown in Figs. \ref{fig:PuIL_PF_C1FourthOrthoMPM} and \ref{fig:PuIL_HS_C1FourthOrthoMPM}, respectively. Contrary to the isotropic case, crack initiation does not occur along the the vertical axis as shown in Figs. \ref{fig:PuIL_C1FourthOrthoMPM_50sec} and \ref{fig:PuIL_HS_C1FourthOrthoMPM_50sec}. Indeed, the anisotropic orthotropic model and the material orientation trigger the crack to propagate along a weak direction that is not aligned with vertical axis, similar to case (ii). However, contrary to case (ii) no crack branching is observed in case (iii). This can be justified by the fact that there is only one preferential weak direction. As a result, the crack path continues to propagate at one half of the plate. To further illustrate the effect of anisotropy on the resulting crack path, the reciprocal of the surface energy density is also plotted in Fig. \ref{fig:PuIL_HS_C1FourthOrthoMPM}. 

\begin{figure}
	\centering
	\begin{tabular}{ccc}
		\subfloat[\label{fig:PuIL_Elast_IsoCubicOrtho}]{
			\includegraphics[width=0.37\columnwidth]{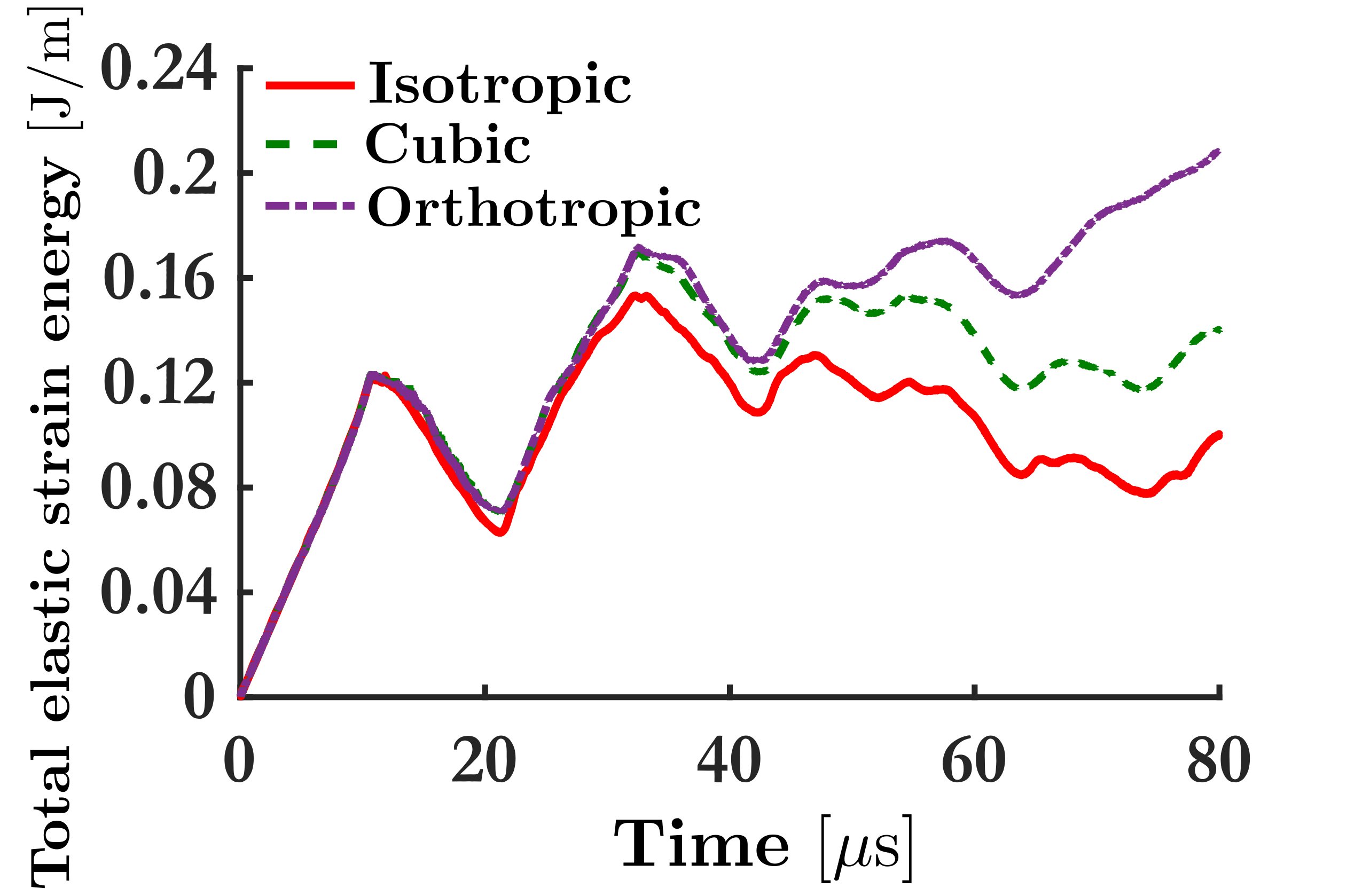}} &
		{\subfloat[\label{fig:PuIL_Fract_IsoCubicOrtho}]{
				\includegraphics[width=0.37\columnwidth]{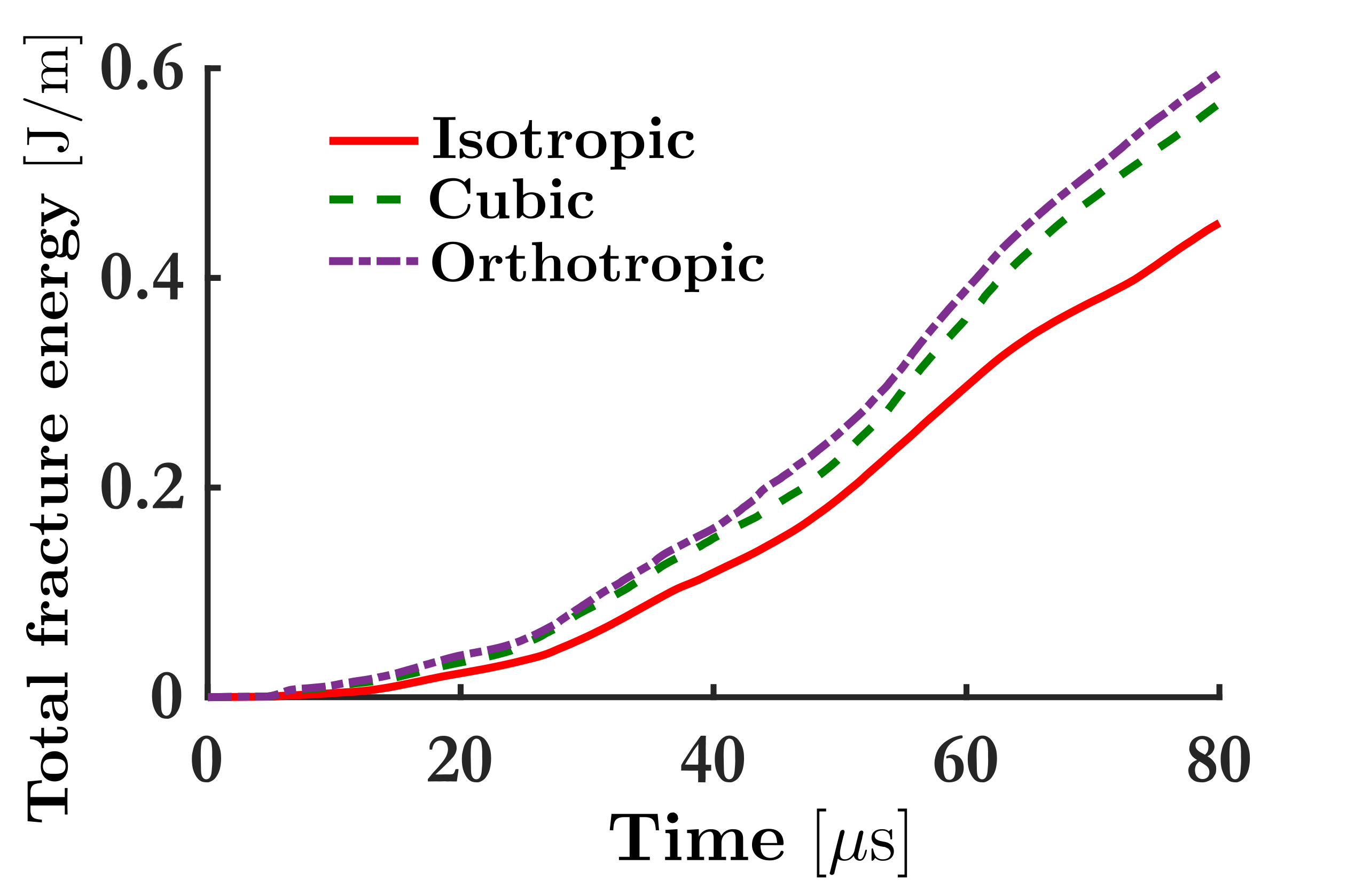}}} &
		{\subfloat[\label{fig:PuIL_CrackVel_IsoCubicOrtho}]{
				\includegraphics[width=0.37\columnwidth]{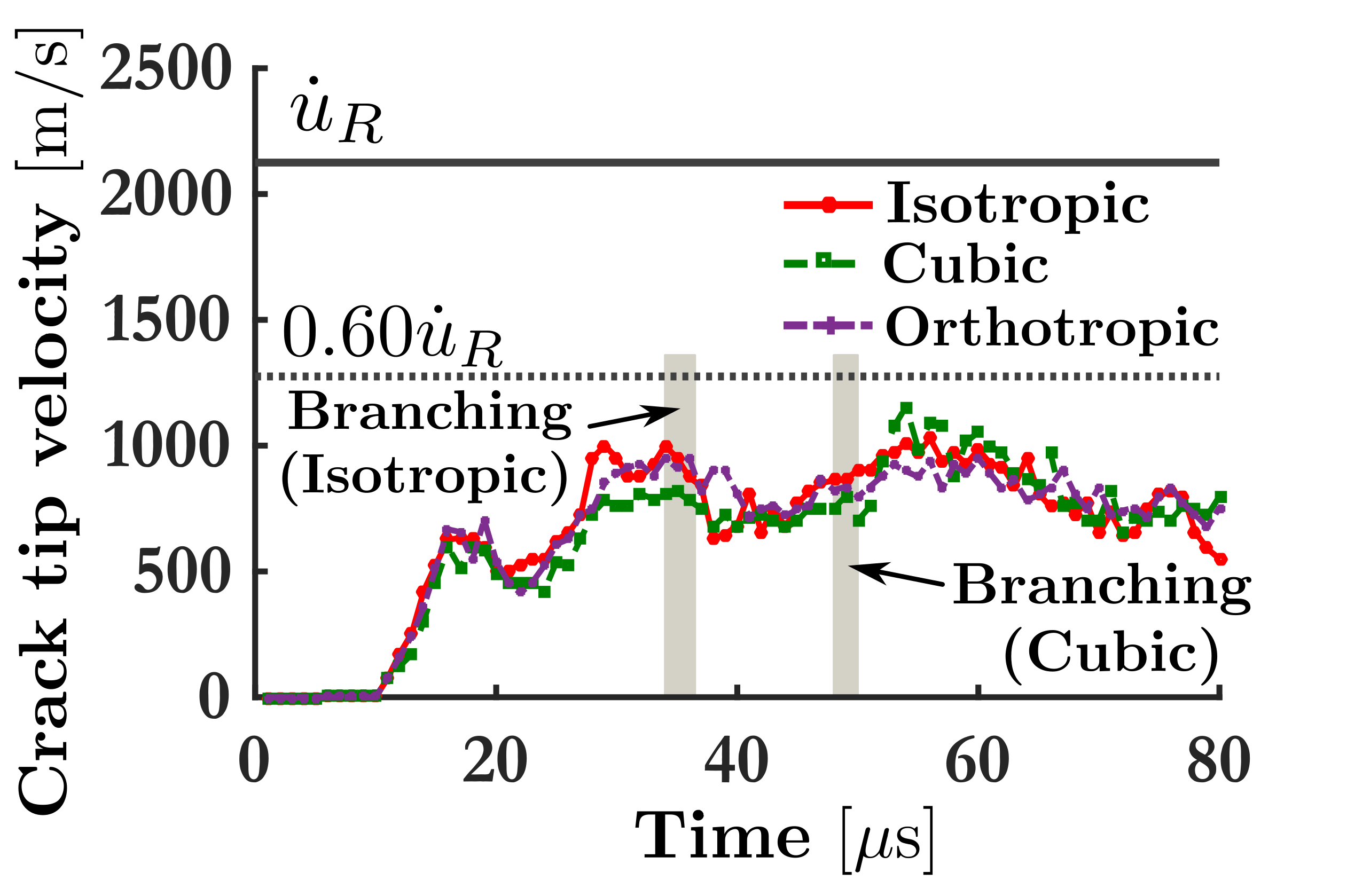}}} 
	\end{tabular}
	\caption[]{Plate under impact loading: \subref{fig:PuIL_Elast_IsoCubicOrtho} Total elastic strain energies, \subref{fig:PuIL_Fract_IsoCubicOrtho} Total fracture energies and \subref{fig:PuIL_CrackVel_IsoCubicOrtho} Crack tip velocities over time for PF-MPM 2nd order isotropic model (case (i)), PF-MPM 4th order cubic model (case (ii)) and PF-MPM 4th order orthotropic model (case (iii)). The traction is considered to be $\sigma = 1$ N/mm\textsuperscript{2}.}
	\label{fig:IsoCubicOrtho}
\end{figure}

The energy time-histories for all cases are shown in Figs. \ref{fig:PuIL_Elast_IsoCubicOrtho} and \ref{fig:PuIL_Fract_IsoCubicOrtho}. In all cases the crack initiates when the total elastic energy becomes approximately equal to $0.12$ J/m. It should be stressed that the total elastic strain energy evolves in an almost identical fashion in cases (ii) and (iii) until approximately $t=50$ $\mu$s. After that, the two models diverge as a result of the crack branching in the case of cubic symmetry only. 

The crack tip velocities for the three cases are shown in Fig. \ref{fig:PuIL_CrackVel_IsoCubicOrtho}. The branched regions of both case (i) (isotropic) and  (ii) (cubic) are also highlighted in Fig. \ref{fig:PuIL_CrackVel_IsoCubicOrtho}. The crack branching in isotropic symmetry is observed earlier than in cubic symmetry with a corresponding decrease in crack speed. In all cases, the crack initiates at approximately $t=10$ $\mu$s; the resulting crack tip velocities in all cases are comparable. This can be justified by the fact that the impact energy imposed as well as $\mathcal{G}_{c_{min}}$ are identical in all cases. In the orthotropic case, all the results are below the 60\% $\dot{u}_{R}$ limit.

\subsubsection{Crack branching and merging}

\begin{figure}
	\centering
	\begin{tabular}{cccc}
		\subfloat[\label{fig:PuIL_C1SecondMPMDem_0sec}]{
			\includegraphics[width=0.20\columnwidth]{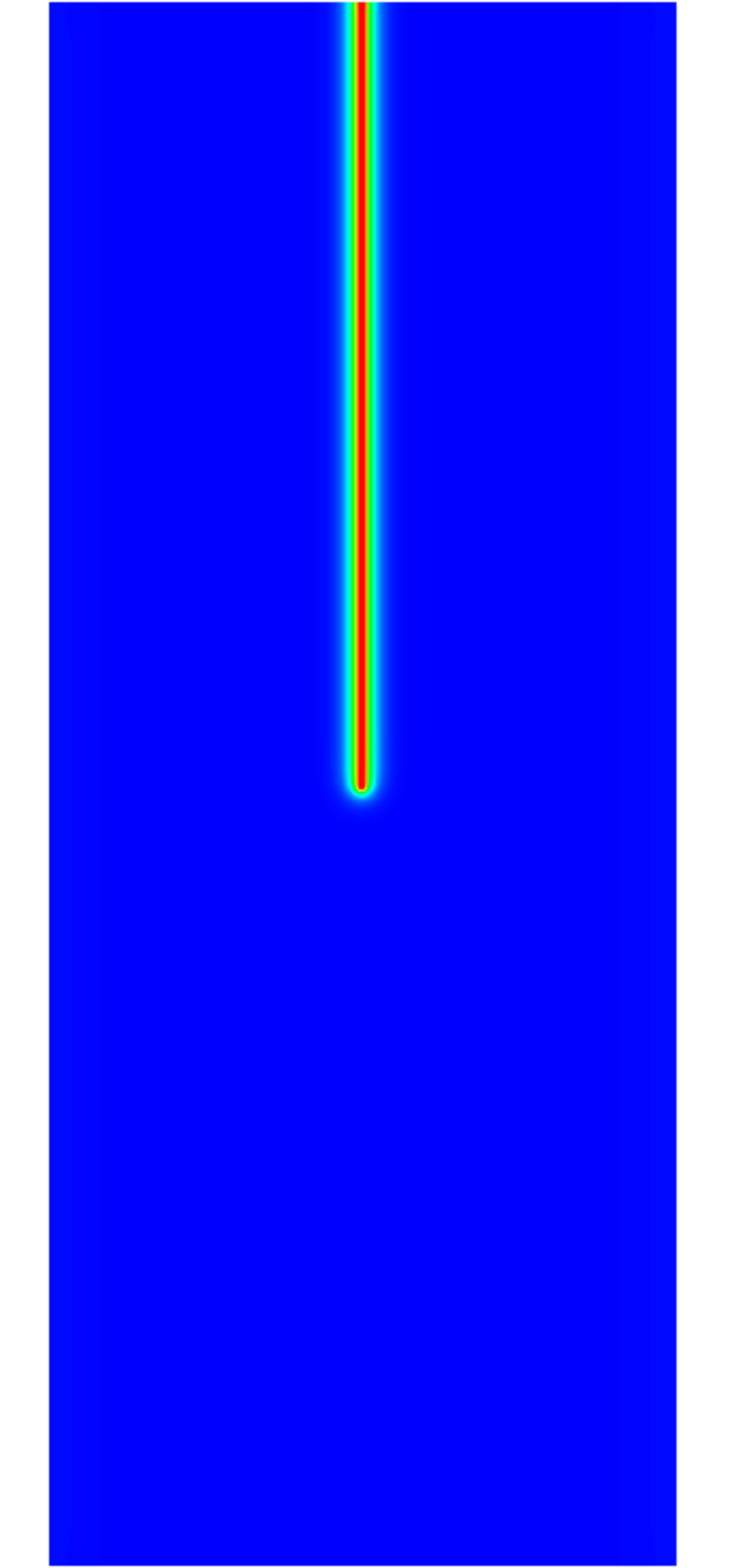}} &
		\subfloat[\label{fig:PuIL_C1SecondMPMDem_50sec}]{
			\includegraphics[width=0.20\columnwidth]{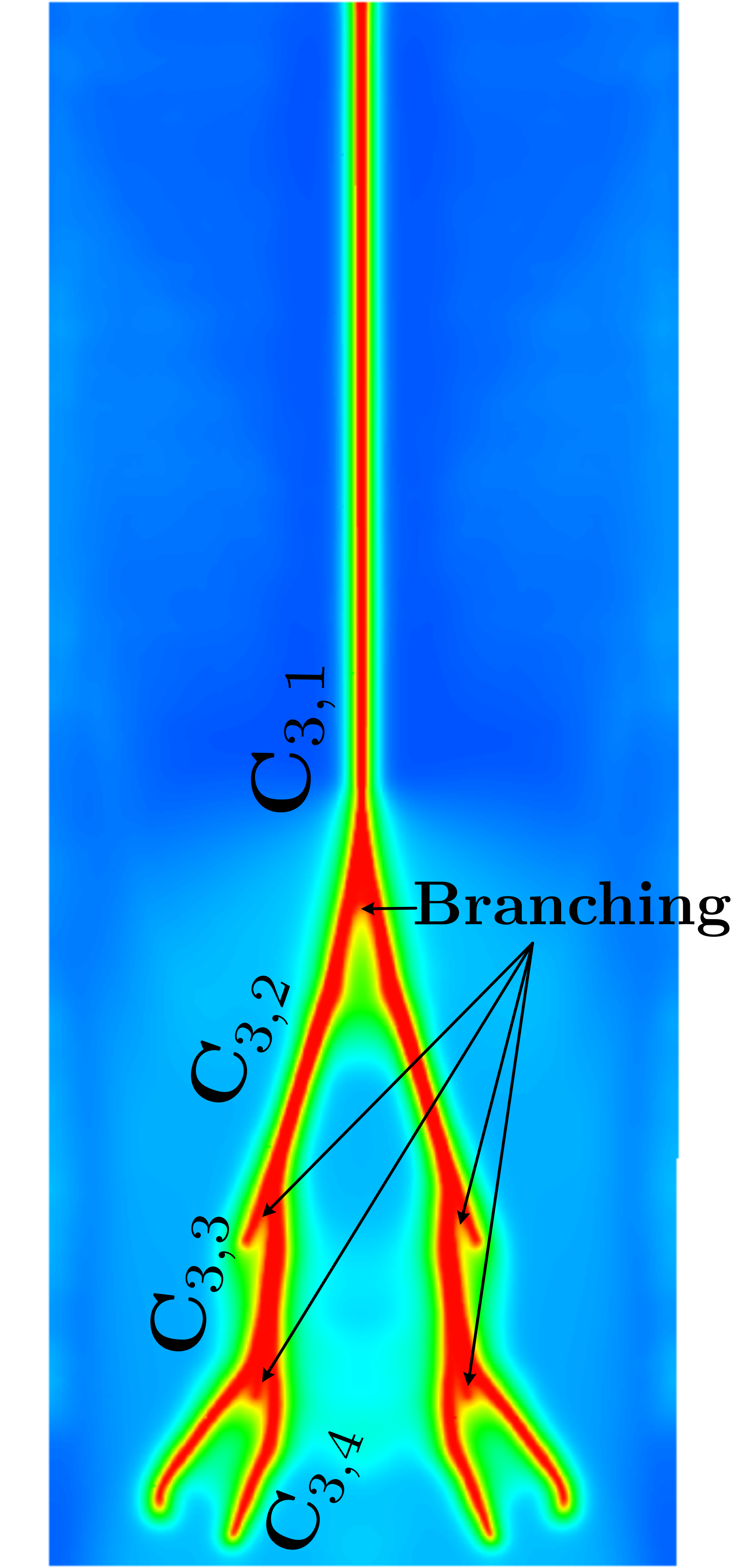}} &
		\subfloat[\label{fig:PuIL_C1SecondMPMDem_110sec}]{
			\includegraphics[width=0.20\columnwidth]{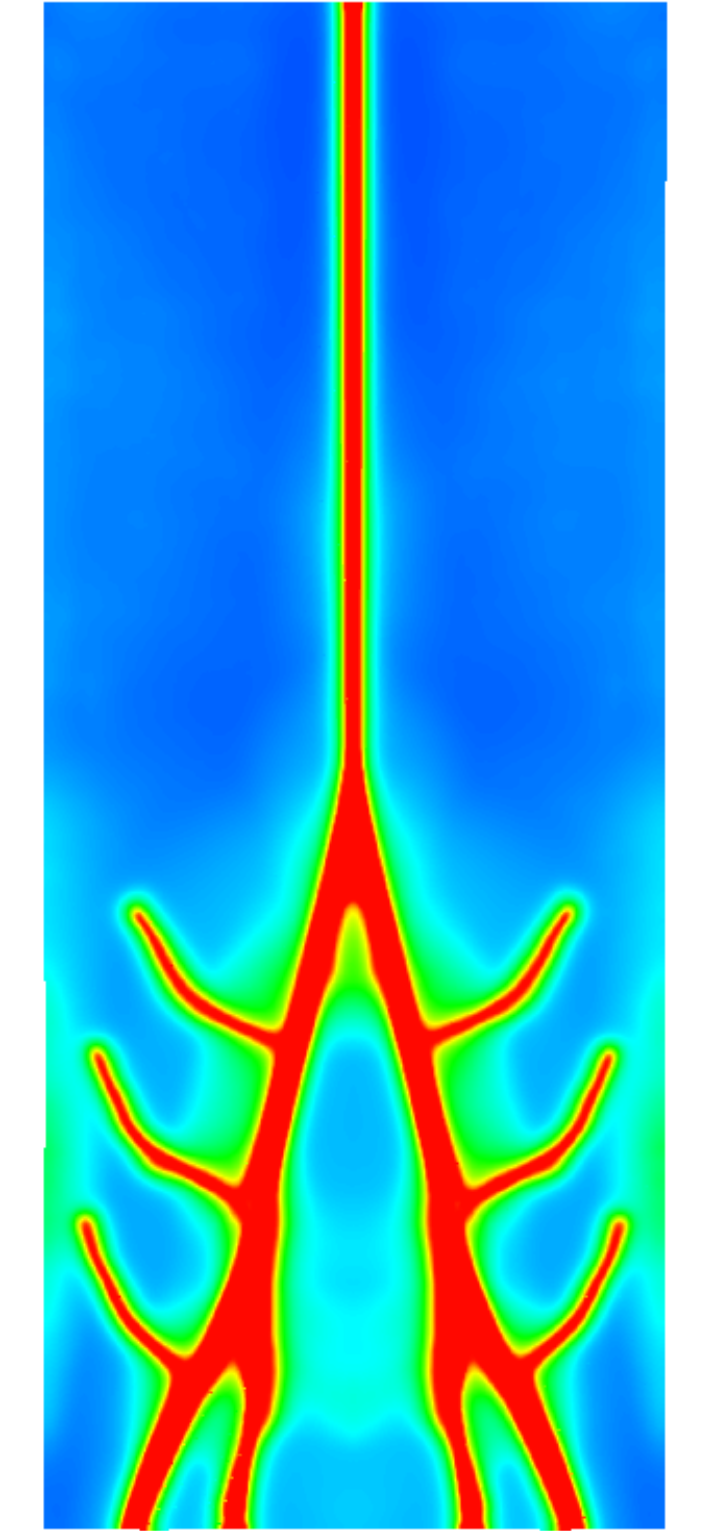}} &
		\subfloat[\label{fig:PuIL_C1SecondMPMDem_130sec}]{
			\includegraphics[width=0.20\columnwidth]{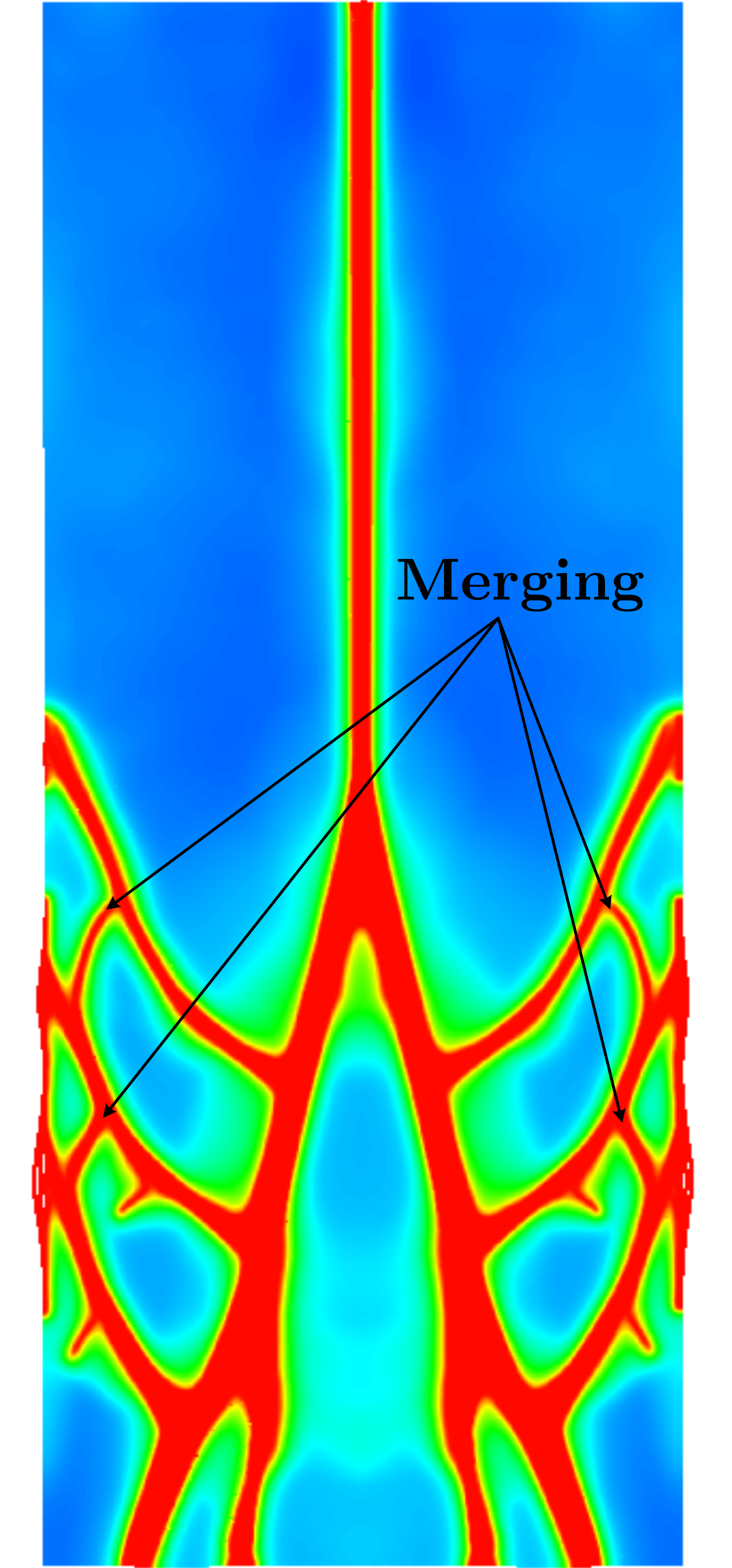}} \\
		\multicolumn{4}{c}{\subfloat{
				\includegraphics[width=0.75\columnwidth]{PF_ColorBar.png}}}		
	\end{tabular}
	\caption[]{Plate under impact loading: Phase field for time steps \subref{fig:PuIL_C1SecondMPMDem_0sec} t=0 $\mu$s \subref{fig:PuIL_C1SecondMPMDem_50sec} t=50 $\mu$s \subref{fig:PuIL_C1SecondMPMDem_110sec} t=110 $\mu$s and \subref{fig:PuIL_C1SecondMPMDem_130sec} t=130 $\mu$s. Results for 2nd order isotropic phase field model and $\sigma = 2.3$ N/mm\textsuperscript{2}.}
	\label{fig:PuIL_PF_C1SecondDemMPM}
\end{figure}

\begin{figure}
	\centering
	\begin{tabular}{cccc}
		\subfloat[\label{fig:PuIL_HS_C1SecondMPMDem_0sec}]{
			\includegraphics[width=0.20\columnwidth]{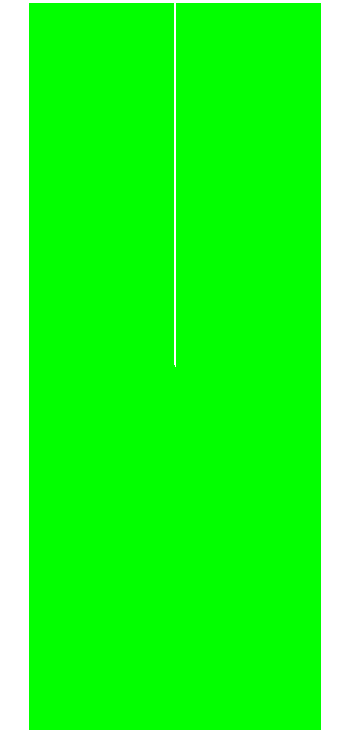}} &
		\subfloat[\label{fig:PuIL_HS_C1SecondMPMDem_50sec}]{
			\includegraphics[width=0.20\columnwidth]{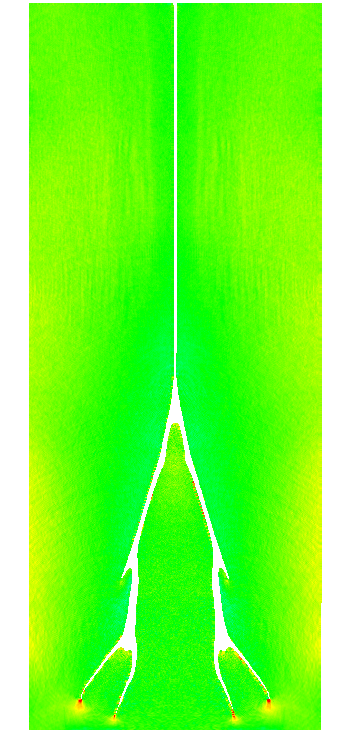}} &
		\subfloat[\label{fig:PuIL_HS_C1SecondMPMDem_110sec}]{
			\includegraphics[width=0.20\columnwidth]{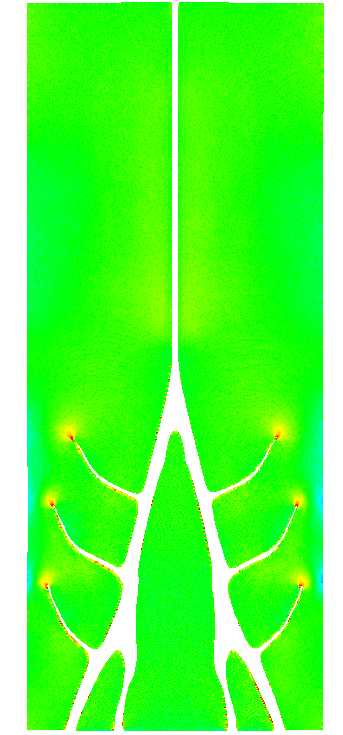}} &
		\subfloat[\label{fig:PuIL_HS_C1SecondMPMDem_130sec}]{
			\includegraphics[width=0.20\columnwidth]{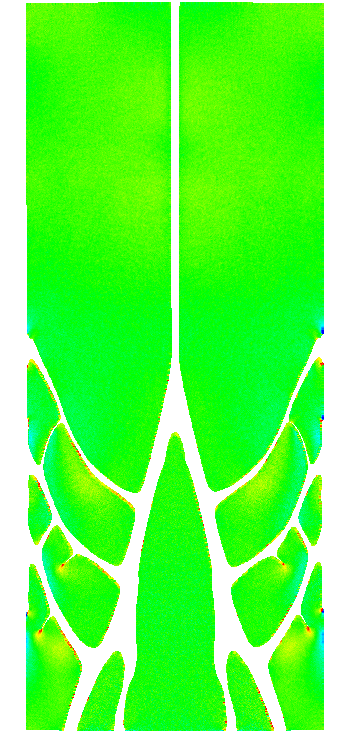}} \\
		\multicolumn{4}{c}{\subfloat{
				\includegraphics[width=0.75\columnwidth]{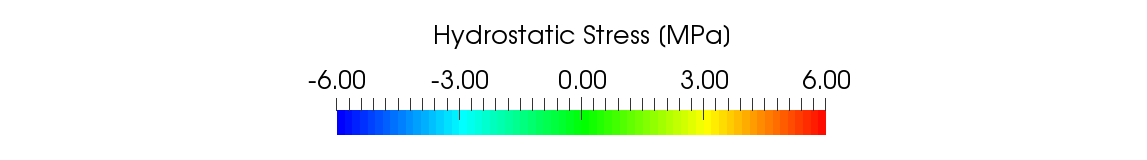}}}		
	\end{tabular}
	\caption[]{Plate under impact loading: Hydrostatic stress for time steps \subref{fig:PuIL_HS_C1SecondMPMDem_0sec} t=0 $\mu$s \subref{fig:PuIL_HS_C1SecondMPMDem_50sec} t=50 $\mu$s \subref{fig:PuIL_HS_C1SecondMPMDem_110sec} t=110 $\mu$s and \subref{fig:PuIL_HS_C1SecondMPMDem_130sec} t=130 $\mu$s. Results for 2nd order isotropic phase field model and $\sigma = 2.3$ N/mm\textsuperscript{2}. Material points with $c_p<0.10$ have been removed.}
	\label{fig:PuIL_HS_C1SecondDemMPM}
\end{figure}

\begin{figure}
	\centering
	\begin{tabular}{ccc}
		\subfloat[\label{fig:PuIL_ElastEng_C1SecondMPMDem}]{
			\includegraphics[width=0.37\columnwidth]{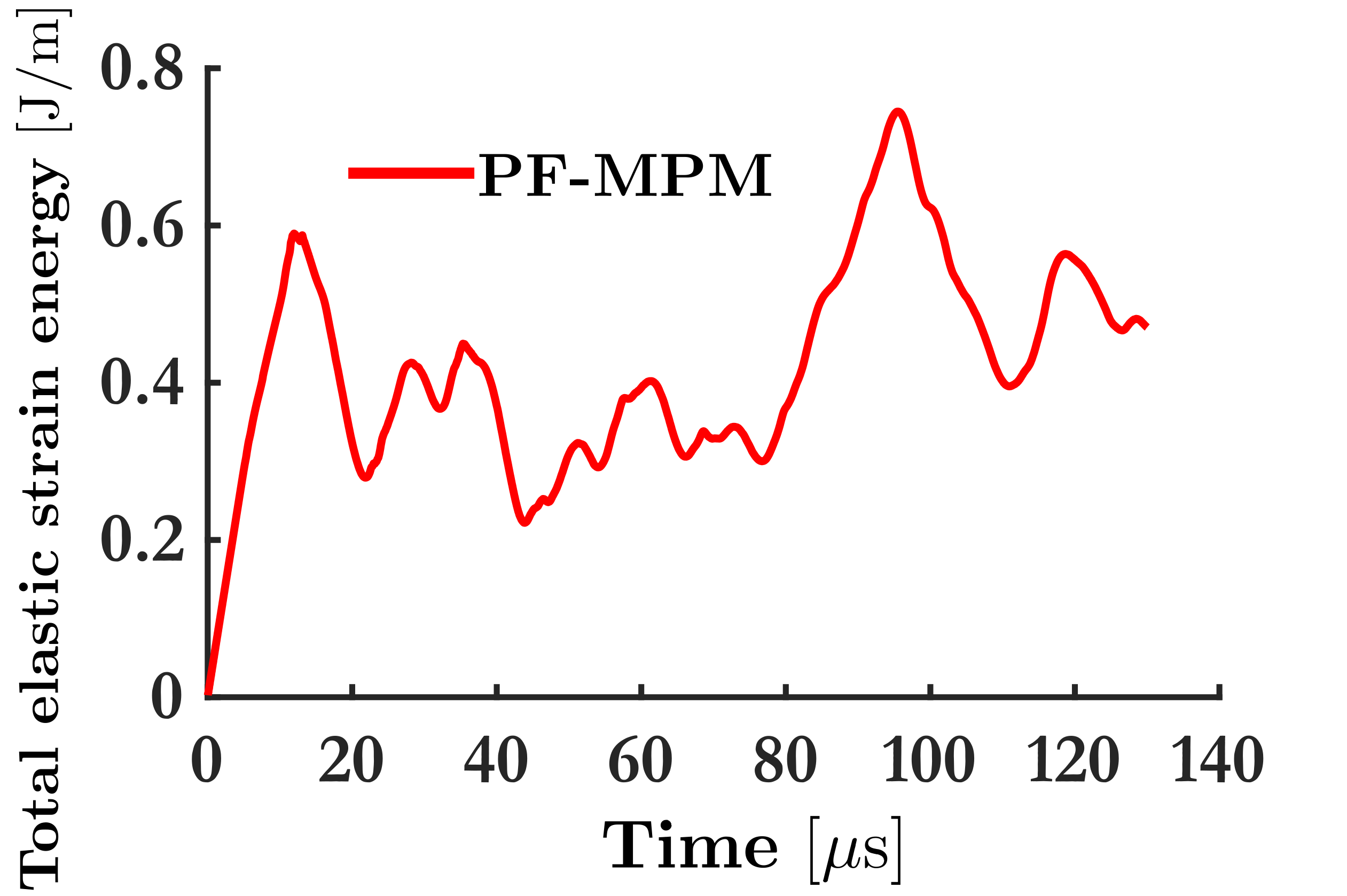}} &
		{\subfloat[\label{fig:PuIL_FractEng_C1SecondMPMDem}]{
				\includegraphics[width=0.37\columnwidth]{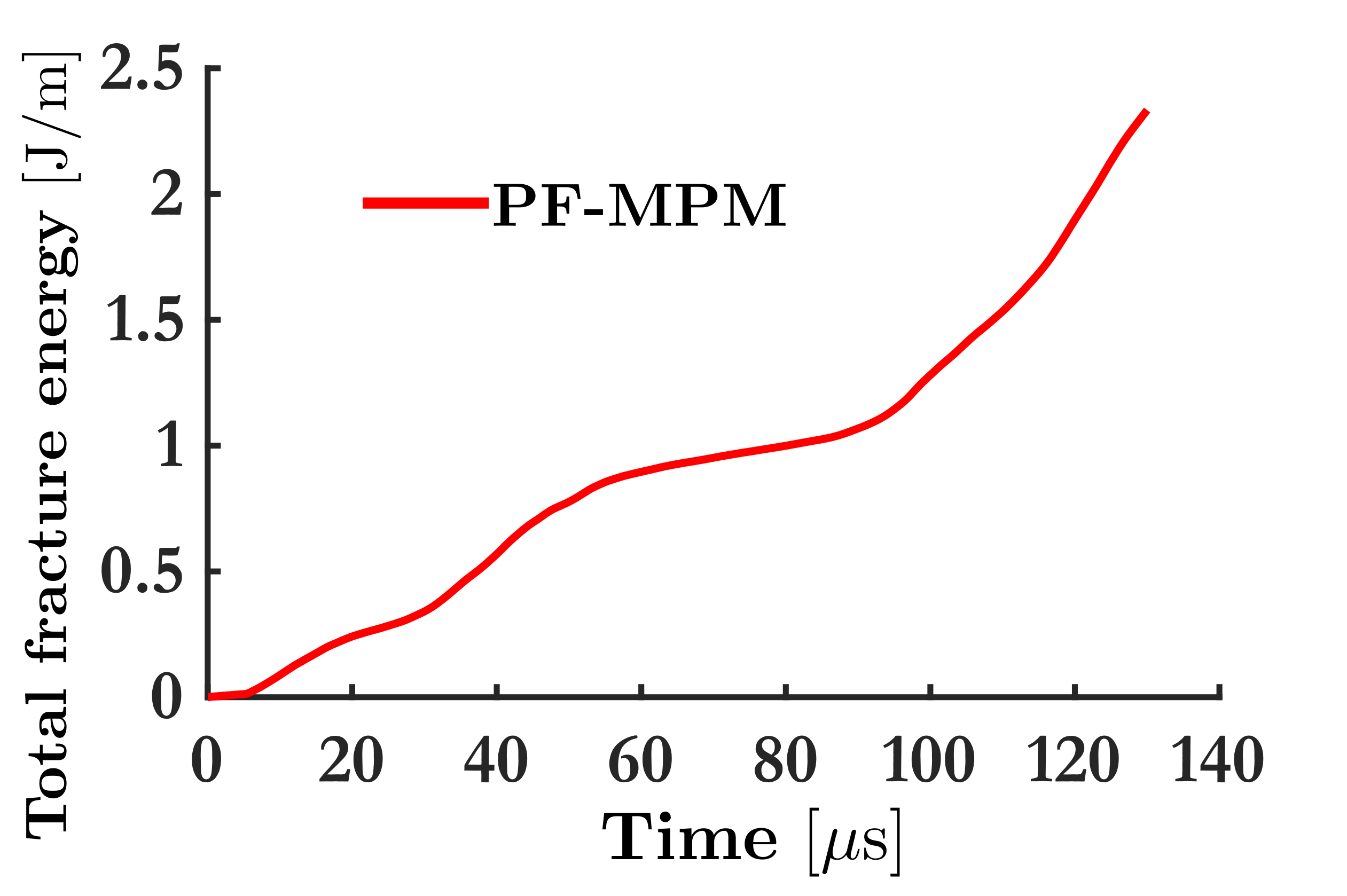}}} &
		{\subfloat[\label{fig:PuIL_CrackVel_C1SecondMPMDem}]{
				\includegraphics[width=0.37\columnwidth]{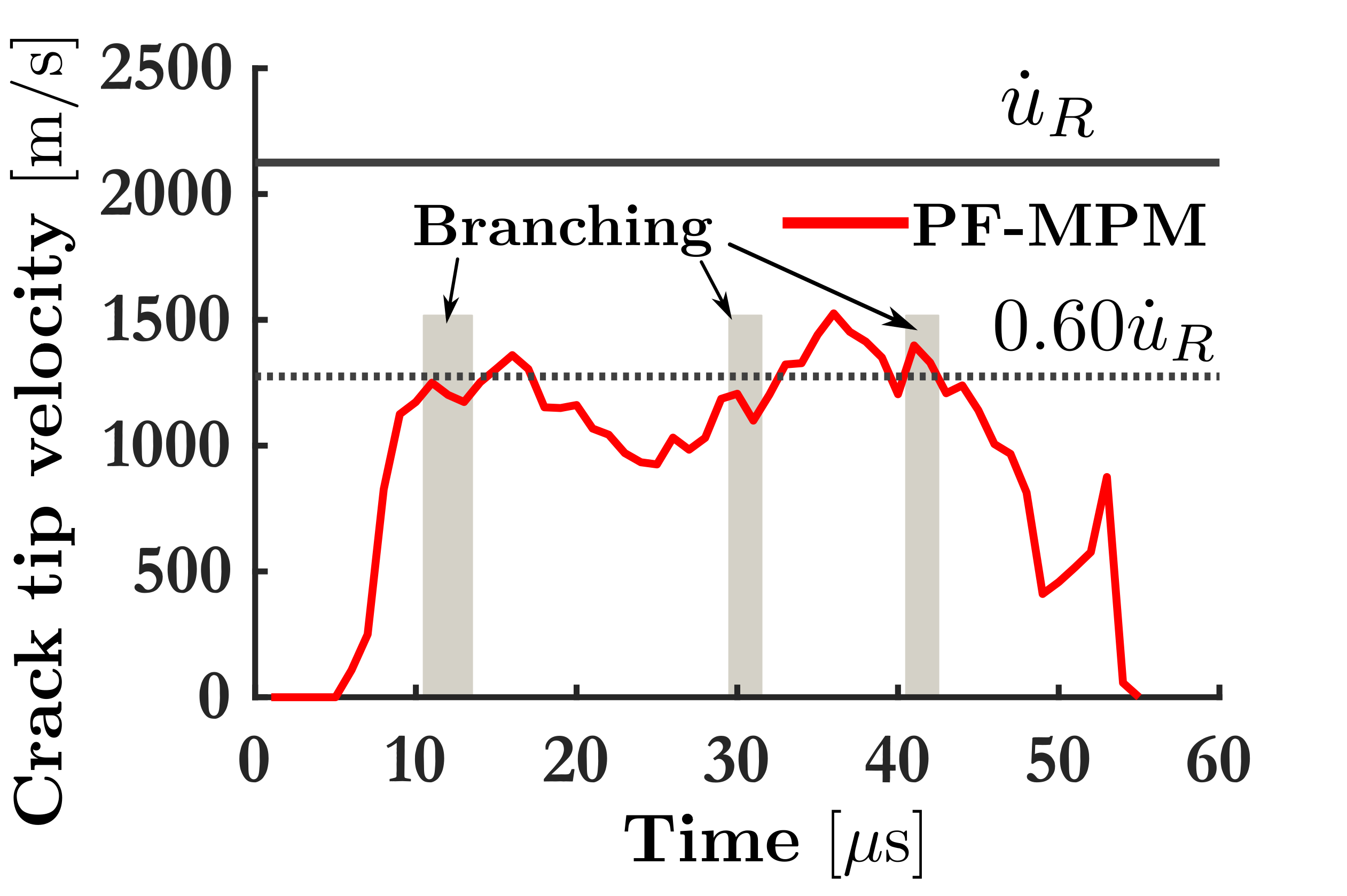}}} 
	\end{tabular}
	\caption[]{Plate under impact loading: \subref{fig:PuIL_ElastEng_C1SecondMPMDem} Total elastic strain energy, \subref{fig:PuIL_FractEng_C1SecondMPMDem} Total fracture energy and \subref{fig:PuIL_CrackVel_C1SecondMPMDem} Crack tip velocity over time for PF-MPM 2nd order isotropic model. The traction is considered to be $\sigma = 2.3$ N/mm\textsuperscript{2}.}
	\label{fig:C1SecondMPMDem}
\end{figure}

The efficiency of phase field models to deal with complex crack paths, i.e. including crack branching and crack merging is demonstrated herein. The traction is increased to $\sigma = 2.3$ N/mm\textsuperscript{2}. All other model parameters are kept constant. The second order isotropic phase field model is used. The total duration of the analysis is $t_{tot} = 130$ $\mu$s within $N_{steps}=5200$ steps. The evolution of the phase field for several time steps is presented in Fig. \ref{fig:PuIL_PF_C1SecondDemMPM}. In particular, Fig. \ref{fig:PuIL_C1SecondMPMDem_50sec} illustrates 5 branched cracks, i.e., 1 main, 2 secondary and 2 tertiary branches. More branched cracks are observed in Fig. \ref{fig:PuIL_C1SecondMPMDem_110sec}  while in Fig. \ref{fig:PuIL_C1SecondMPMDem_130sec} 4 merged cracks are presented. The total strain energy together with the total fracture energy is shown in Fig. \ref{fig:C1SecondMPMDem} whereas The evolution of the hydrostatic stress for that case is shown in Fig. \ref{fig:PuIL_HS_C1SecondDemMPM} for several time steps 

The crack tip velocity for that case is shown in Fig. \ref{fig:PuIL_CrackVel_C1SecondMPMDem}. The crack tip velocity is measured along the paths $\mathbf{C}_{3,1}$, $\mathbf{C}_{3,2}$, $\mathbf{C}_{3,3}$ and $\mathbf{C}_{3,4}$ that are marked in Fig. \ref{fig:PuIL_C1SecondMPMDem_50sec}. The increased impact loading, i.e. $\sigma = 2.3$ N/mm\textsuperscript{2}, leads to a crack initiation at approximately $t=5$ $\mu$s. This occurs earlier than in case (i), where the first crack initiates at approximately $t=10$ $\mu$s. The crack tip rapidly accelerates to the 60\% $\dot{u}_{R}$ limit. Although, some points exceed the 60\% $\dot{u}_{R}$ limit, the majority of measured points satisfy this condition while all points are clearly below the Rayleigh wave speed. In particular, the crack tip propagates with an average speed close to the 60\% $\dot{u}_{R}$ limit. The three branched regions are also illustrated in Fig. \ref{fig:PuIL_CrackVel_C1SecondMPMDem} where a decrease in crack speed is observed.

\subsection{Collision of two rings} \label{lbl:CTR}

\begin{figure}
	\centering
	\begin{tabular}{cc}
		\subfloat[\label{fig:CTR_GeometryBCs}]{
			\includegraphics[width=0.50\columnwidth]{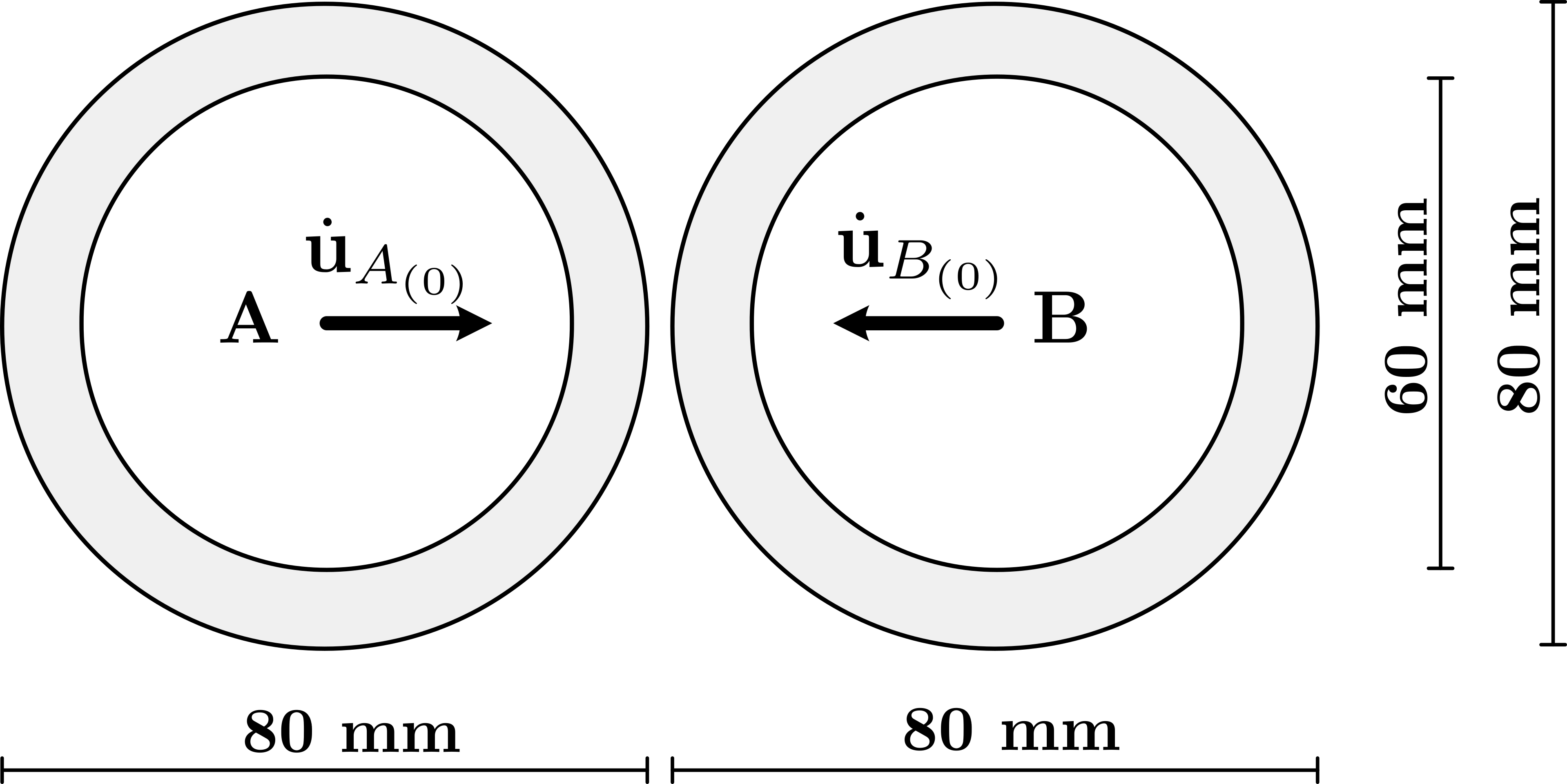}} &
		{\subfloat[\label{fig:CTR_FractEngs}]{
				\includegraphics[width=0.50\columnwidth]{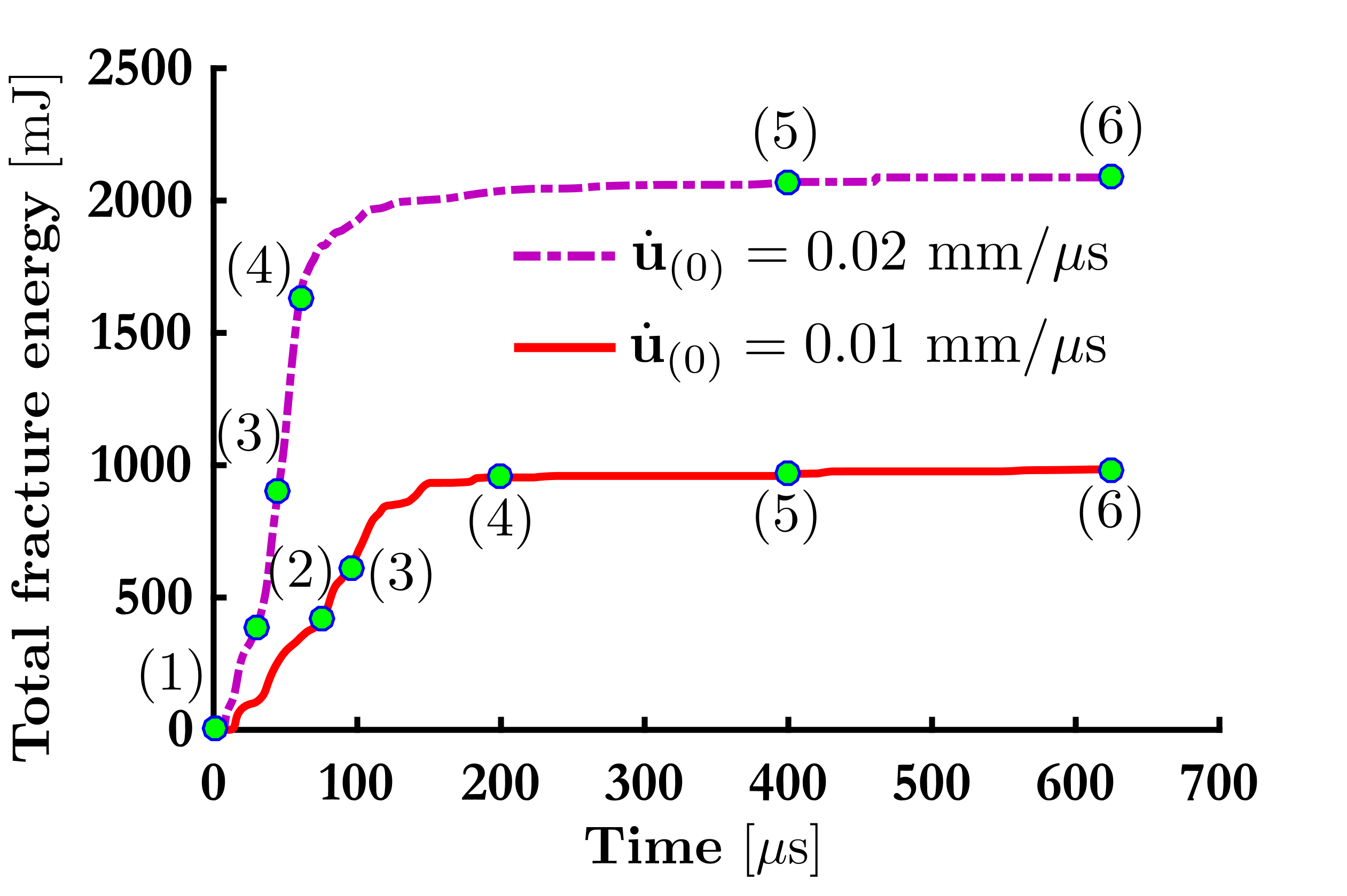}}}
	\end{tabular}
	\caption[]{Collision of two rings: \subref{fig:CTR_GeometryBCs} Geometry and boundary conditions \subref{fig:CTR_FractEngs} Total fracture energy over time for case (i) $\dot{\mathbf{u}}_{{(0)}}=$ $0.01$ mm/$\mu$s and case (ii) $\dot{\mathbf{u}}_{{(0)}}=$ $0.02$ mm/$\mu$s.}
	\label{fig:CTR_GeometryBCsFractEngs}
\end{figure}

Herein, the collision of two rings is analysed with the geometry and boundary conditions of the problem shown in Fig. \ref{fig:CTR_GeometryBCs}. The aim of this example is to demonstrate the robustness of the proposed  method into resolving fragmentation problems and the interactions occurring between fragments.

The cell (patch) spacing is chosen to be $h=0.50$ mm for the numerical implementation and plane stress conditions are assumed with thickness $2$ mm. The grid is formed by two knot vectors $\splitatcommas{\Xi=\{0,0,0,0.00263,0.0052,...,0.9947,0.9973,1,1,1\}}$ and $\splitatcommas{H=\{0,0,0,0.00357,0.0071,...,0.9928,0.9964,1,1,1\}}$, $107724$ control points and $380x280=106400$ cells. Two discrete fields are considered, i.e., field A (left ring) and B (right ring). The corresponding friction coefficient is $\mu_{f}=0.65$. The total number of material points is $325620$. 

The elastic material parameters are chosen to be $E=190000$ N/mm\textsuperscript{2}, $\nu=0.30$ and $\rho=8000$ kg/m\textsuperscript{3} for both bodies. A time step $\Delta t = 0.0125$ $\mu$s for $N_{steps} = 50000$ is considered. The critical time step is $\tilde{\Delta t_{cr}}=0.071$ $\mu$s. The initial distance between the two rings is assumed to be $2h=1.00$ mm. An initial velocity is applied to the material points of the two rings as $\dot{\mathbf{u}}_{Ap_{(0)}}=\dot{\mathbf{u}}_{{(0)}}$ and $\dot{\mathbf{u}}_{Bp_{(0)}}=-\dot{\mathbf{u}}_{{(0)}}$. To examine the influence of the initial velocity into the resulting crack paths, two cases are considered, namely (i) $\dot{\mathbf{u}}_{{(0)}}=0.01$ mm/$\mu$s and (ii) $\dot{\mathbf{u}}_{{(0)}}=0.02$ mm/$\mu$s, respectively. The second order isotropic model is utilized for that problem, therefore $\gamma_{ijkl} = 0$ with length scale parameter $l_0=1.00$ mm, $k_f=0.00$ and $\mathcal{G}_{c} \left( \theta \right)=\bar{\mathcal{G}}_c=\mathcal{G}_{c_{max}}=\mathcal{G}_{c_{min}}=6.00$ N/mm.

The total fracture energy time-history for both cases is shown in Fig. \ref{fig:CTR_FractEngs}. The evolution of the phase field and the hydrostatic stress for points ((1)-(6)) labelled in Fig. \ref{fig:CTR_FractEngs} is shown in Figs. \ref{fig:CTR_PF_V10} and \ref{fig:CTR_HS_V10} for case (i) and in Figs. \ref{fig:CTR_PF_V20} and \ref{fig:CTR_HS_V20} for case (ii), respectively. 

\subsubsection{Case (i): $\dot{\mathbf{u}}_{{(0)}}=0.01$ {\normalfont mm/$\mu$s}}

\begin{figure}
	\centering
	\begin{tabular}{ccc}
		\subfloat[\label{fig:CTR_PF_V10_0sec}]{
			\includegraphics[width=0.3\columnwidth]{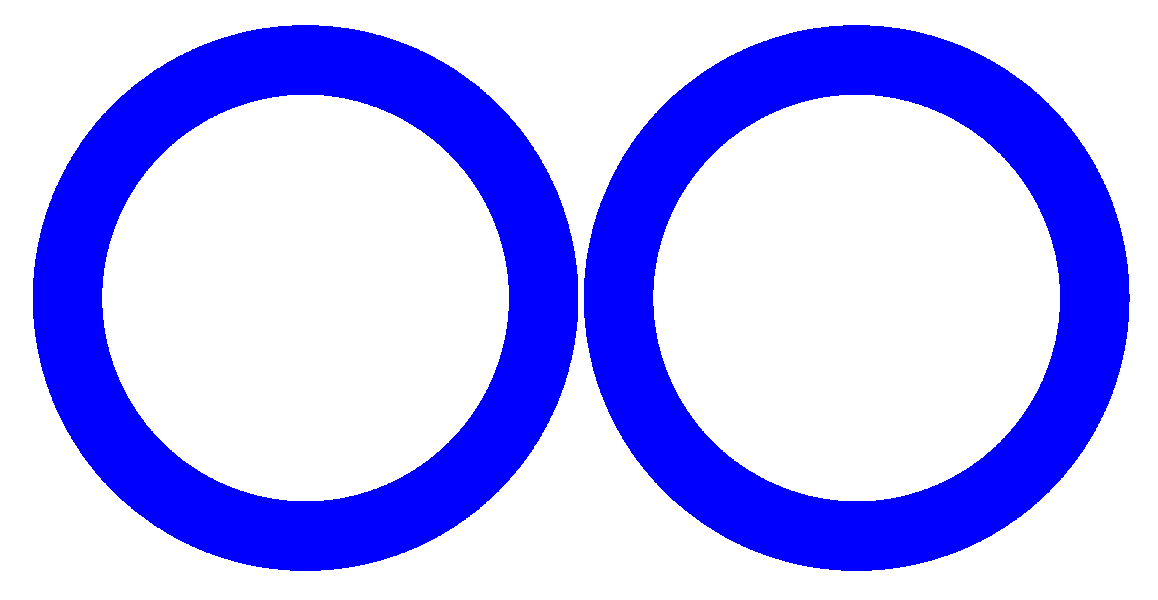}} &
		\subfloat[\label{fig:CTR_PF_V10_75sec}]{
			\includegraphics[width=0.3\columnwidth]{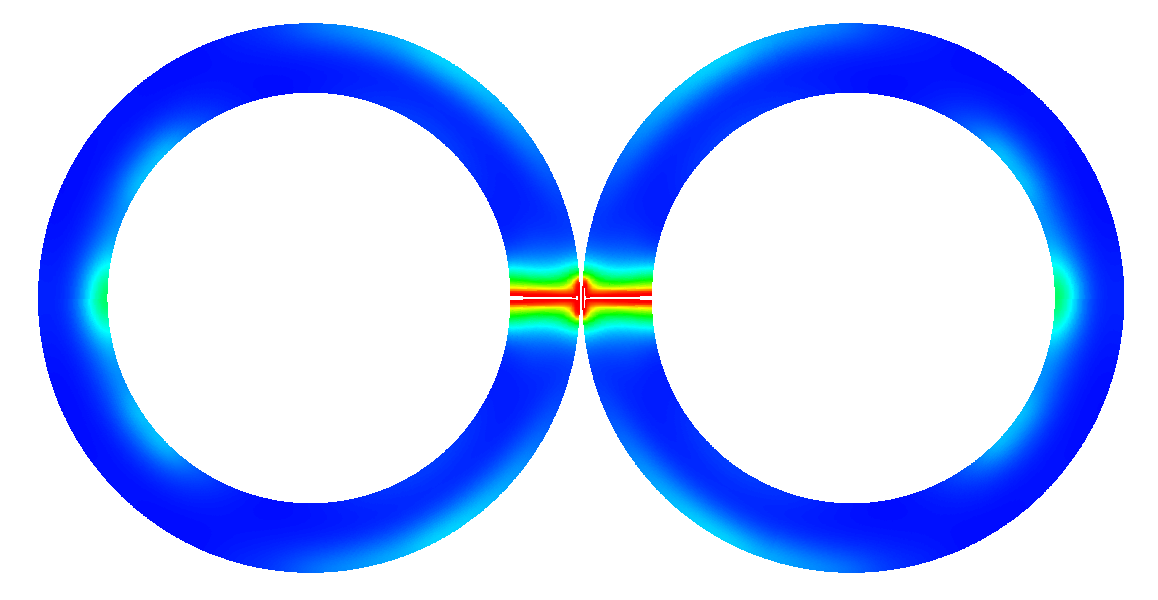}} &
		\subfloat[\label{fig:CTR_PF_V10_95sec}]{
			\includegraphics[width=0.3\columnwidth]{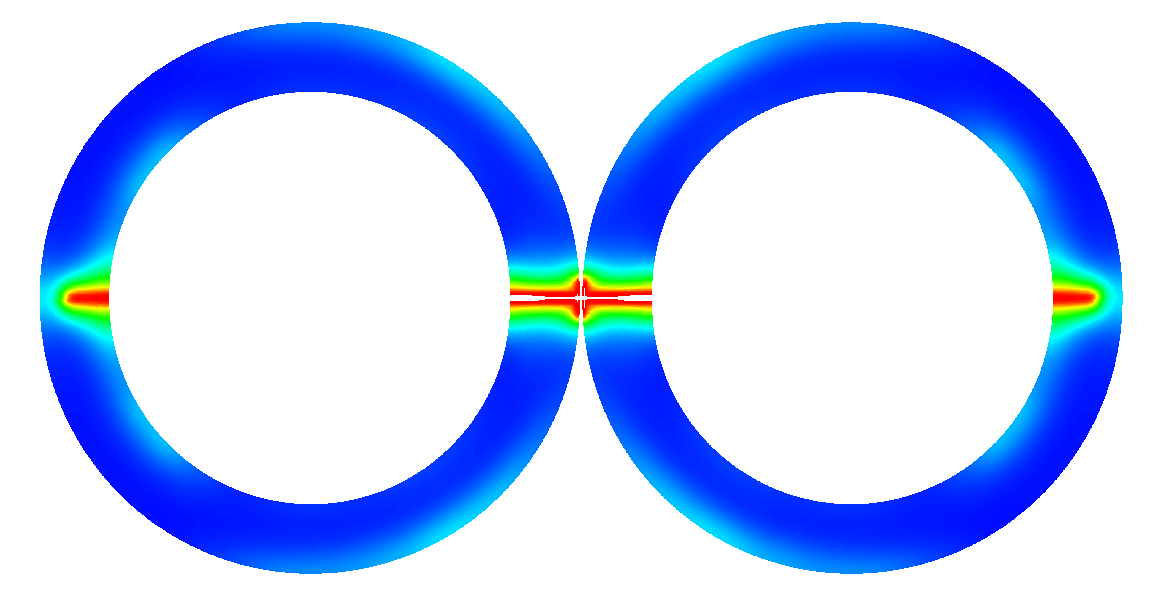}} \\
		\subfloat[\label{fig:CTR_PF_V10_200sec}]{
			\includegraphics[width=0.3\columnwidth]{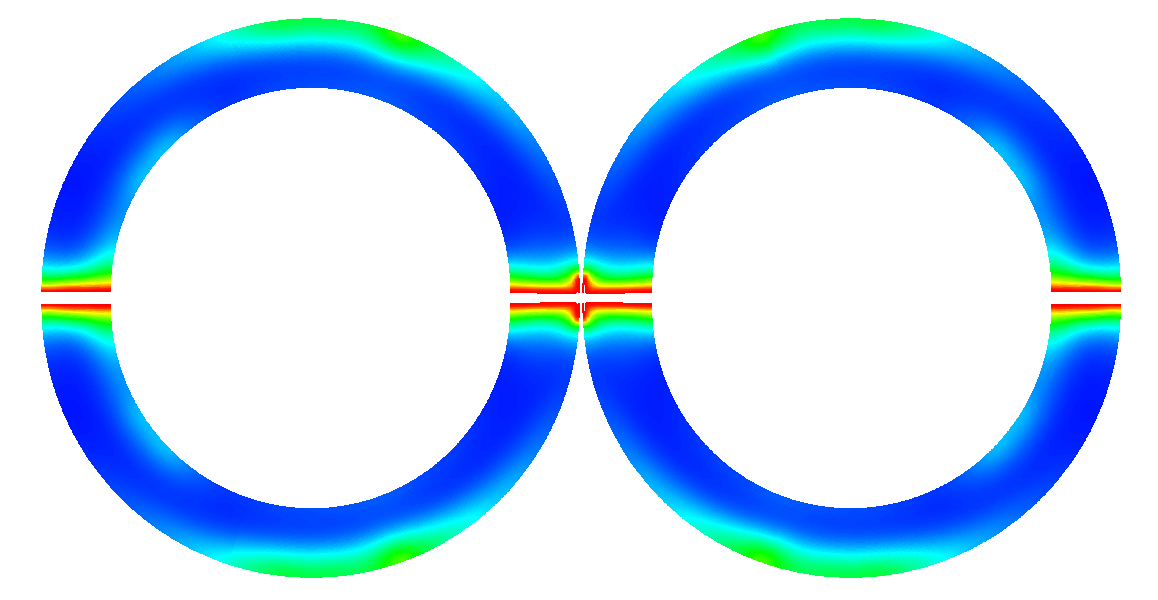}} &
		\subfloat[\label{fig:CTR_PF_V10_400sec}]{
			\includegraphics[width=0.3\columnwidth]{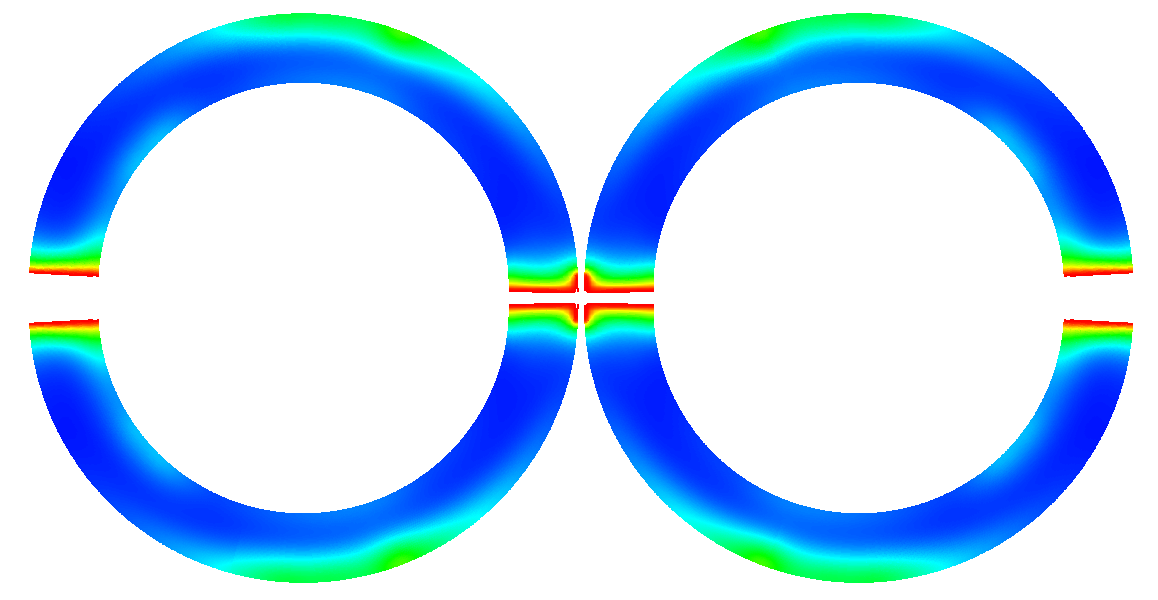}} &
		\subfloat[\label{fig:CTR_PF_V10_625sec}]{
			\includegraphics[width=0.3\columnwidth]{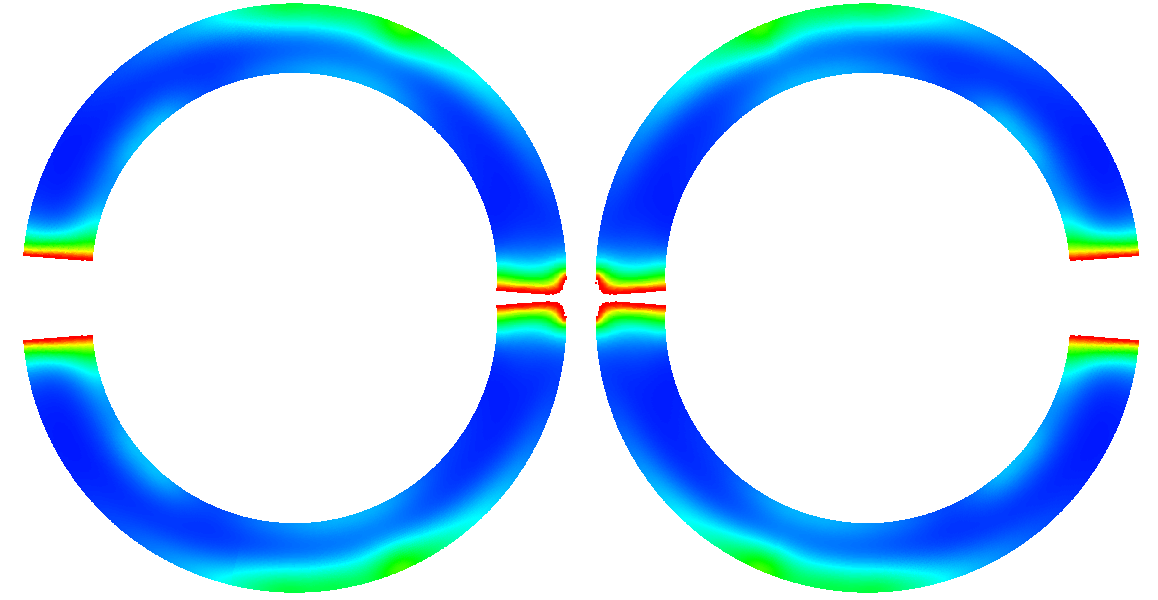} } \\
		\multicolumn{3}{c}{\subfloat{
				\includegraphics[width=1\columnwidth]{PF_ColorBar.png}}} 
	\end{tabular}
	\caption[]{Collision of two rings - Case (i): Phase field for time steps \subref{fig:CTR_PF_V10_0sec} t=0 $\mu$s \subref{fig:CTR_PF_V10_75sec} t=75 $\mu$s \subref{fig:CTR_PF_V10_95sec} t=95 $\mu$s \subref{fig:CTR_PF_V10_200sec} t=200 $\mu$s \subref{fig:CTR_PF_V10_400sec} t=400 $\mu$s and \subref{fig:CTR_PF_V10_625sec} t=625 $\mu$s.}
	\label{fig:CTR_PF_V10}
\end{figure}

\begin{figure}
	\centering
	\begin{tabular}{ccc}
		\subfloat[\label{fig:CTR_HS_V10_0sec}]{
			\includegraphics[width=0.3\columnwidth]{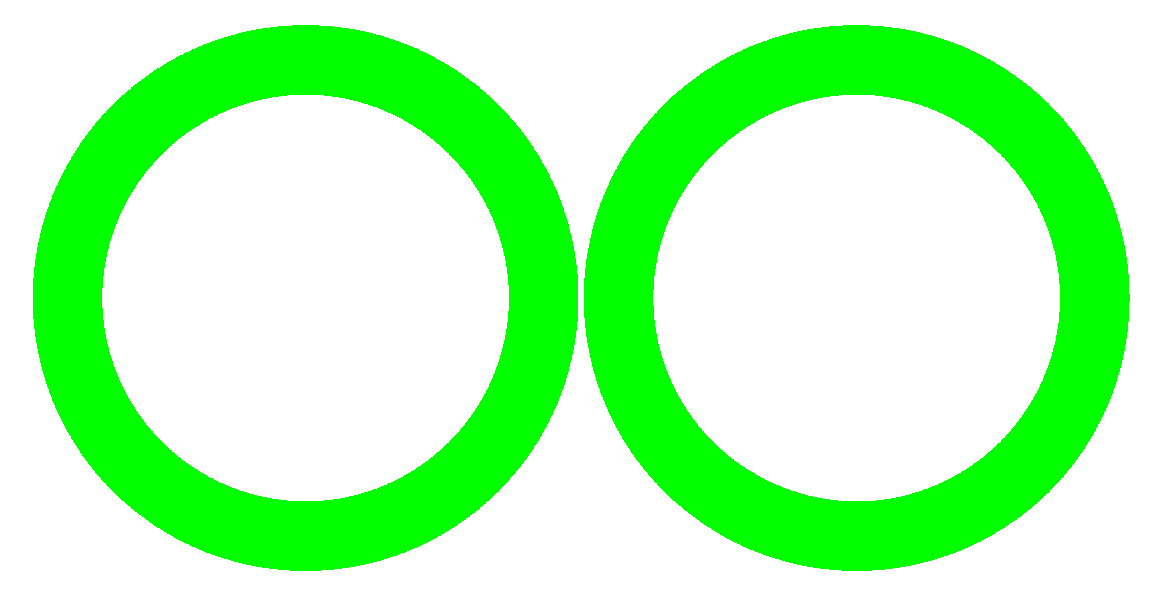}} &
		\subfloat[\label{fig:CTR_HS_V10_75sec}]{
			\includegraphics[width=0.3\columnwidth]{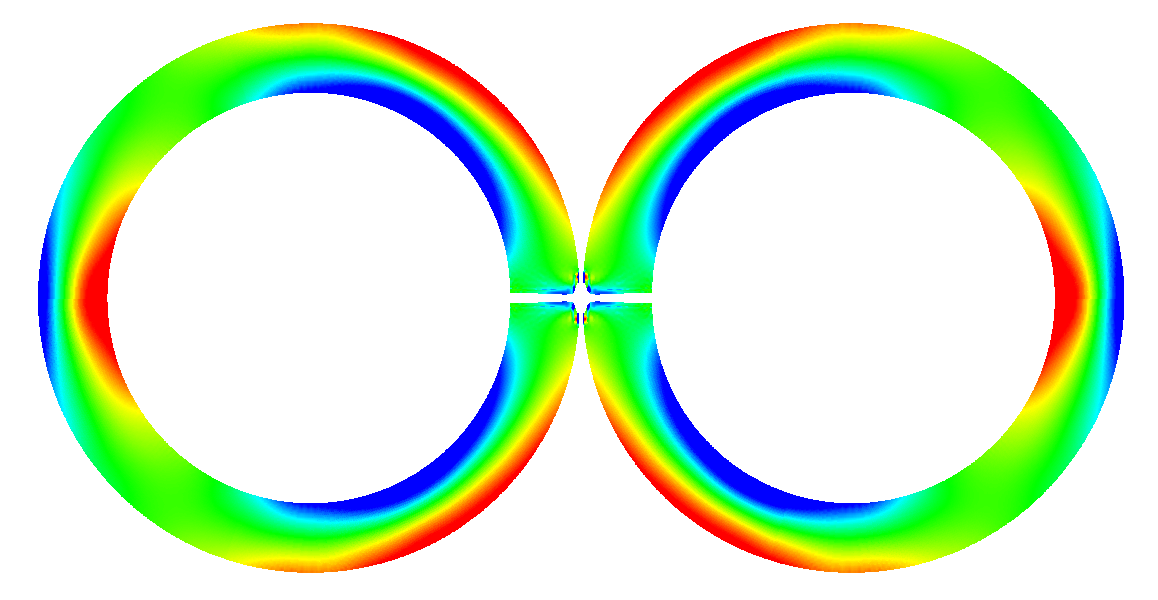}} &
		\subfloat[\label{fig:CTR_HS_V10_95sec}]{
			\includegraphics[width=0.3\columnwidth]{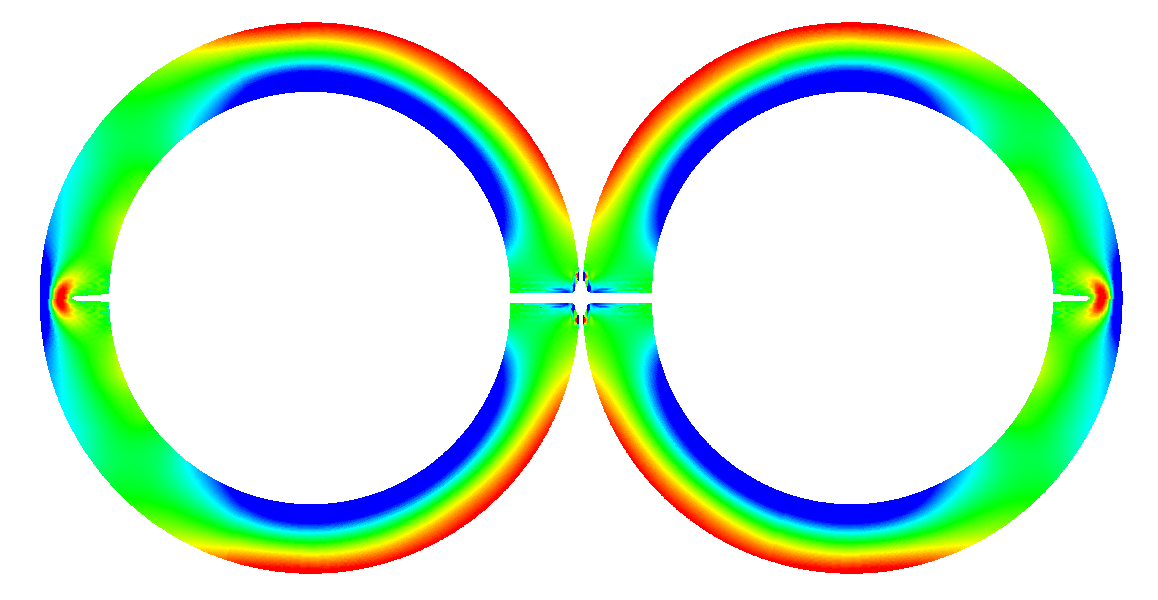}} \\
		\subfloat[\label{fig:CTR_HS_V10_200sec}]{
			\includegraphics[width=0.3\columnwidth]{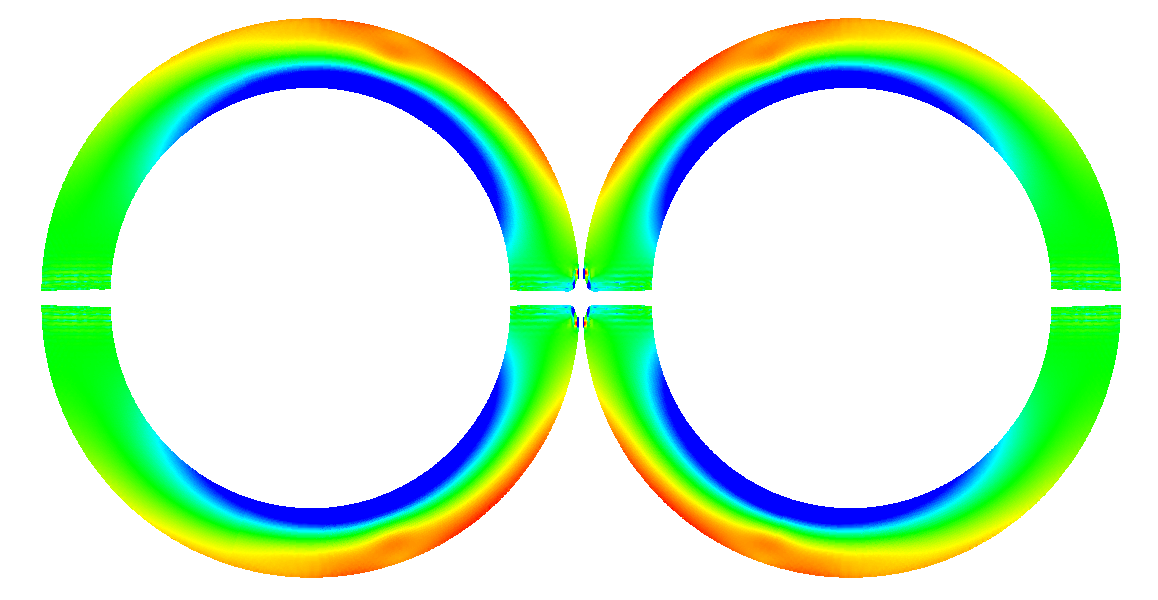}} &
		\subfloat[\label{fig:CTR_HS_V10_400sec}]{
			\includegraphics[width=0.3\columnwidth]{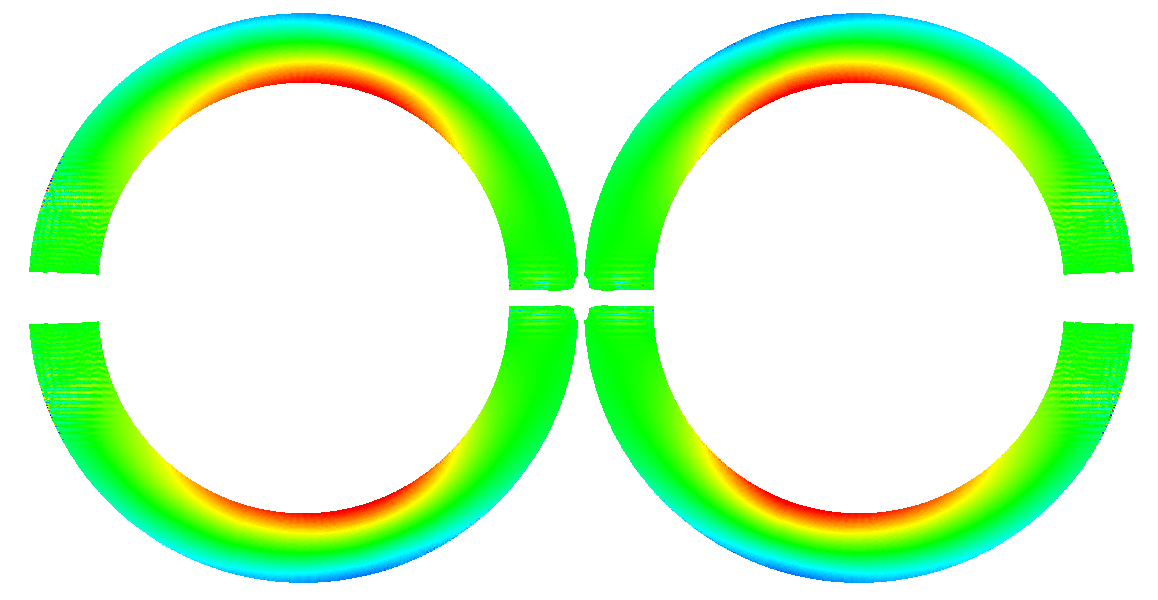}} &
		\subfloat[\label{fig:CTR_HS_V10_625sec}]{
			\includegraphics[width=0.3\columnwidth]{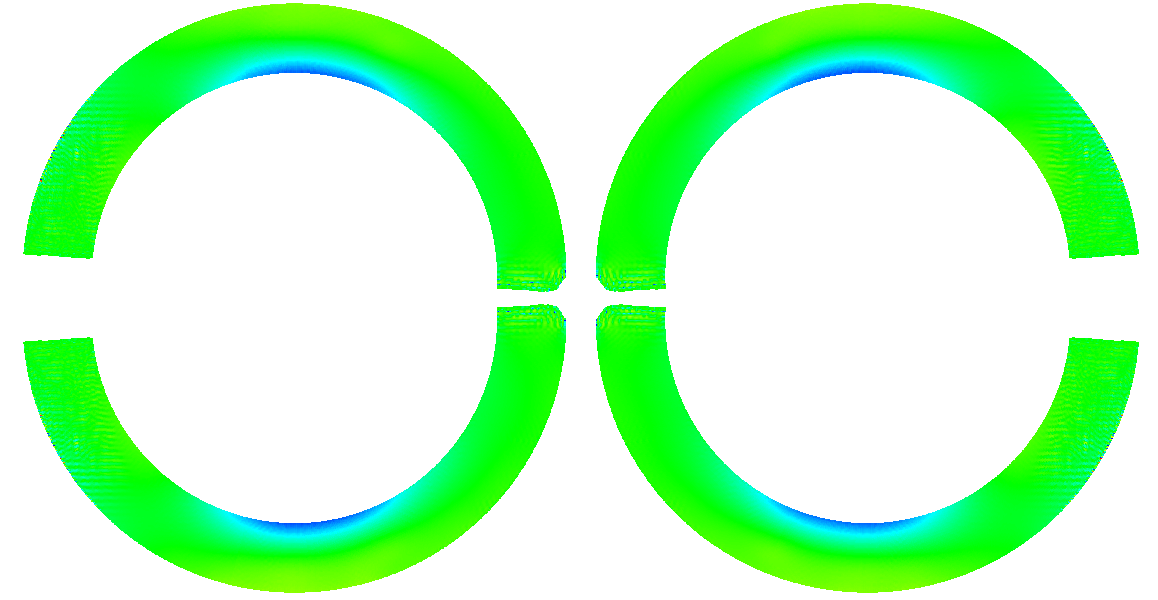} } \\
		\multicolumn{3}{c}{\subfloat{
				\includegraphics[width=1\columnwidth]{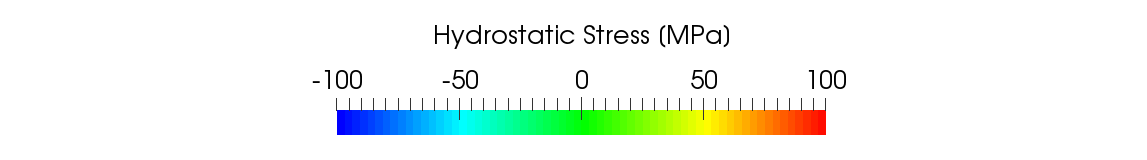}}} 
	\end{tabular}
	\caption[]{Collision of two rings - Case(i): Hydrostatic stresses for time steps \subref{fig:CTR_HS_V10_0sec} t=0 $\mu$s \subref{fig:CTR_HS_V10_75sec} t=75 $\mu$s \subref{fig:CTR_HS_V10_95sec} t=95 $\mu$s \subref{fig:CTR_HS_V10_200sec} t=200 $\mu$s \subref{fig:CTR_HS_V10_400sec} t=400 $\mu$s and \subref{fig:CTR_HS_V10_625sec} t=625 $\mu$s. Material points with $c_p<0.05$ have been removed.}
	\label{fig:CTR_HS_V10}
\end{figure}

In case (i), a crack initiates at the contact surface of the two rings due to their initial impact (see Fig. \ref{fig:CTR_PF_V10_75sec}) followed by a second crack that initiates and fully propagates on the opposite side of each ring (see Fig. \ref{fig:CTR_PF_V10_95sec} and Fig. \ref{fig:CTR_PF_V10_200sec}, respectively). Material degradation also occurs on the top and bottom surfaces of each ring. Crack opening gradually increases (see Fig. \ref{fig:CTR_PF_V10_400sec}) and eventually both rings are fully separated in two fragments (see Fig. \ref{fig:CTR_PF_V10_625sec}). As also shown in Fig. \ref{fig:CTR_FractEngs} the fracture process has been fully developed by point (4), hence the fracture energy remains constant along the path (4)-(6).

\subsubsection{Case (ii): $\dot{\mathbf{u}}_{{(0)}}=0.02$ {\normalfont mm/$\mu$s}}

\begin{figure}
	\centering
	\begin{tabular}{ccc}
		\subfloat[\label{fig:CTR_PF_V20_0sec}]{
			\includegraphics[width=0.28\columnwidth]{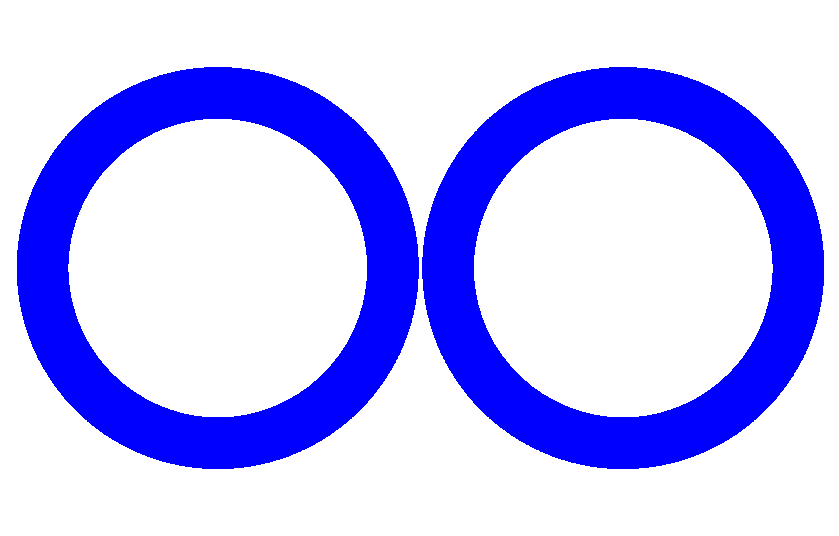}} &
		\subfloat[\label{fig:CTR_PF_V20_30sec}]{
			\includegraphics[width=0.28\columnwidth]{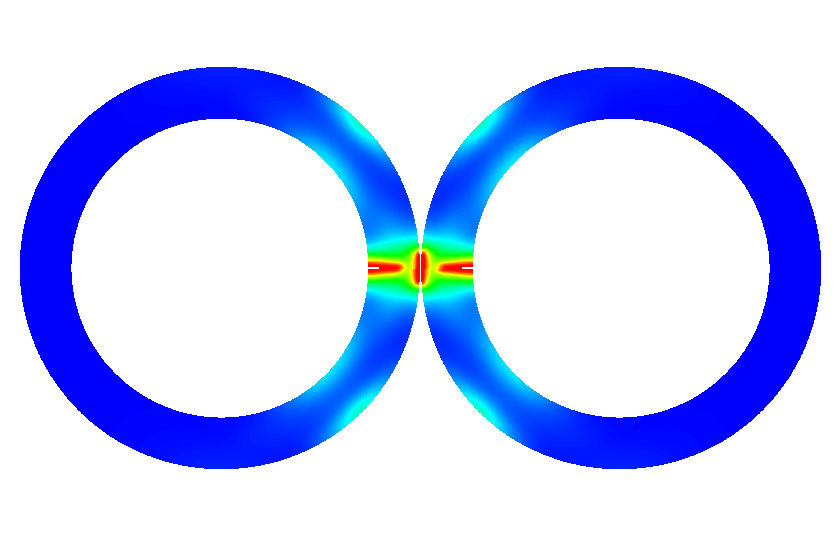}} &
		\subfloat[\label{fig:CTR_PF_V20_45sec}]{
			\includegraphics[width=0.28\columnwidth]{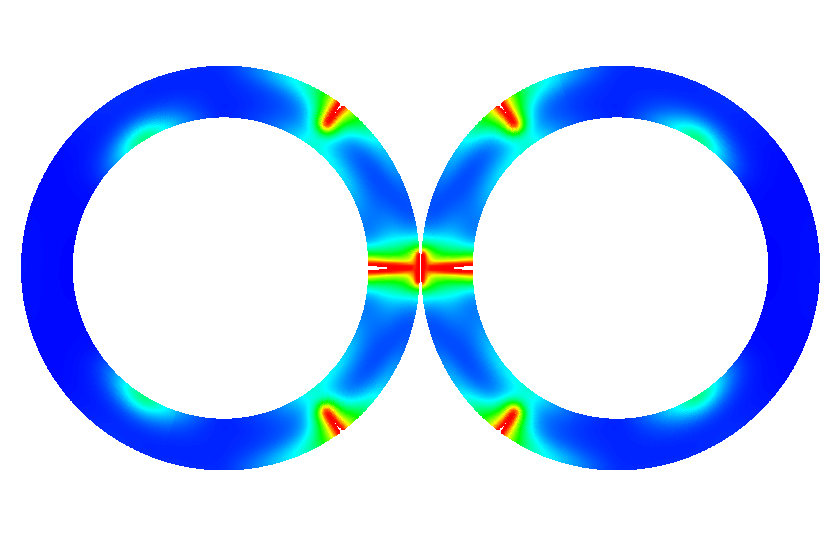}} \\
		\subfloat[\label{fig:CTR_PF_V20_60sec}]{
			\includegraphics[width=0.28\columnwidth]{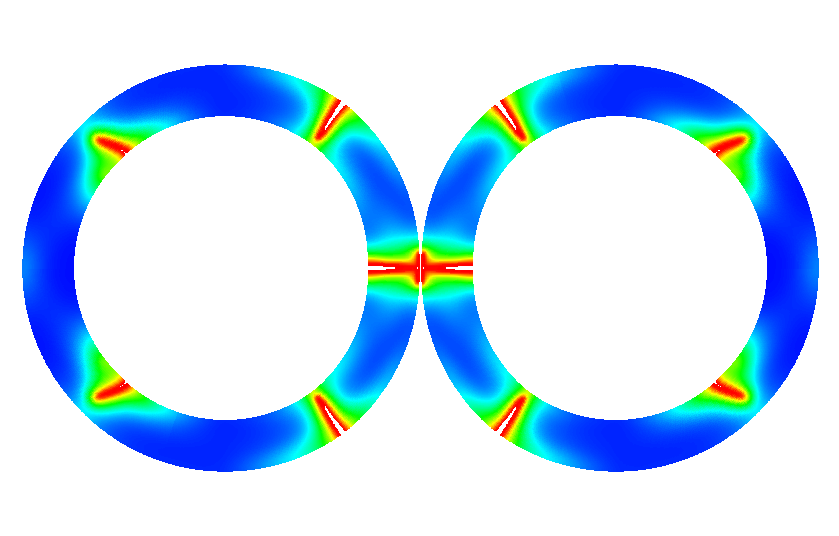}} &
		\subfloat[\label{fig:CTR_PF_V20_400sec}]{
			\includegraphics[width=0.28\columnwidth]{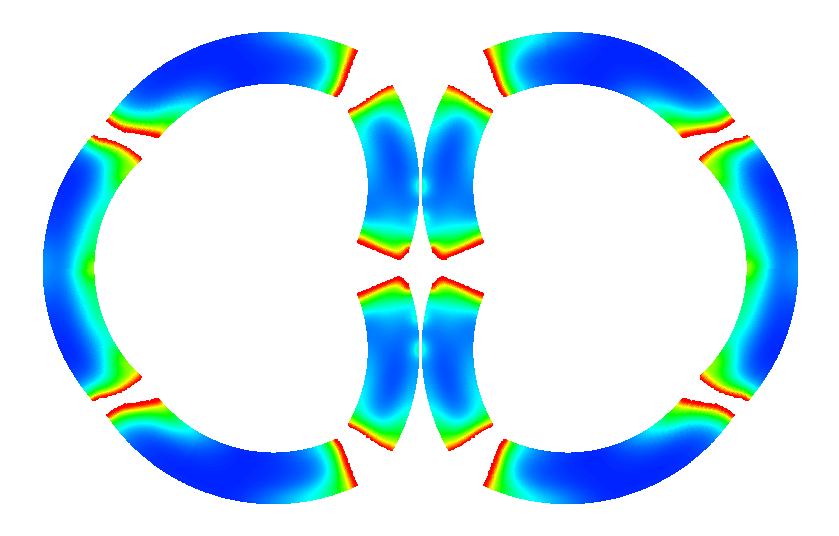}} &
		\subfloat[\label{fig:CTR_PF_V20_625sec}]{
			\includegraphics[width=0.28\columnwidth]{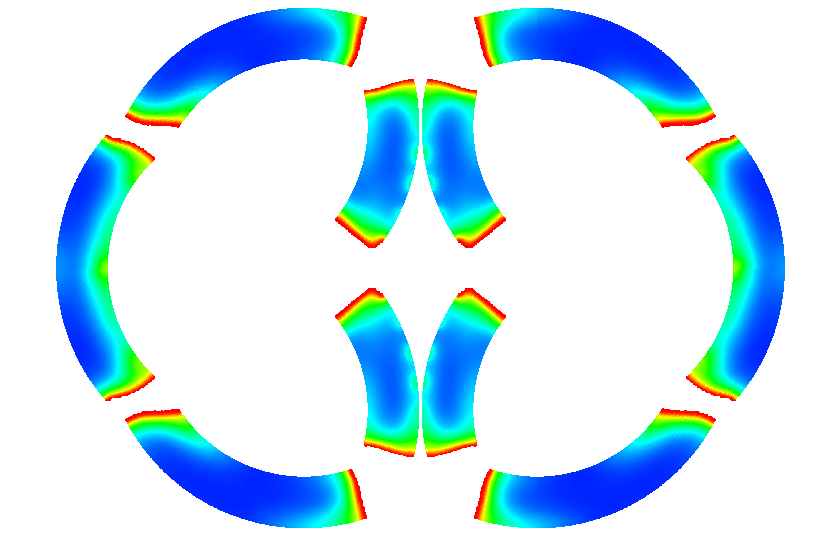} } \\
		\multicolumn{3}{c}{\subfloat{
				\includegraphics[width=1\columnwidth]{PF_ColorBar.png}}} 
	\end{tabular}
	\caption[]{Collision of two rings - Case (ii): Phase field for time steps \subref{fig:CTR_PF_V20_0sec} t=0 $\mu$s \subref{fig:CTR_PF_V20_30sec} t=30 $\mu$s \subref{fig:CTR_PF_V20_45sec} t=45 $\mu$s \subref{fig:CTR_PF_V20_60sec} t=60 $\mu$s \subref{fig:CTR_PF_V20_400sec} t=400 $\mu$s and \subref{fig:CTR_PF_V20_625sec} t=625 $\mu$s.}
	\label{fig:CTR_PF_V20}
\end{figure}

Similar to case (i), in case (ii) a crack initiates at the contact surface of the two rings due to their initial impact (see Fig. \ref{fig:CTR_PF_V20_30sec}). Next, and contrary to case (i), two additional cracks simultaneously propagate of the top right (left) and bottom right (left) of each ring (see Fig. \ref{fig:CTR_PF_V20_45sec}). This is due to the increased impact velocity compared to case (i) where the corresponding points underwent material degradation only. Two more cracks are observed on the top left (right) and bottom left (right) at each ring (see Fig. \ref{fig:CTR_PF_V20_60sec}). The complete crack paths are presented in Fig. \ref{fig:CTR_PF_V20_400sec}. After that point, the fracture energy remains constant; existing cracks do not propagate and new cracks are not initiated.

\begin{figure}
	\centering
	\begin{tabular}{ccc}
		\subfloat[\label{fig:CTR_HS_V20_0sec}]{
			\includegraphics[width=0.28\columnwidth]{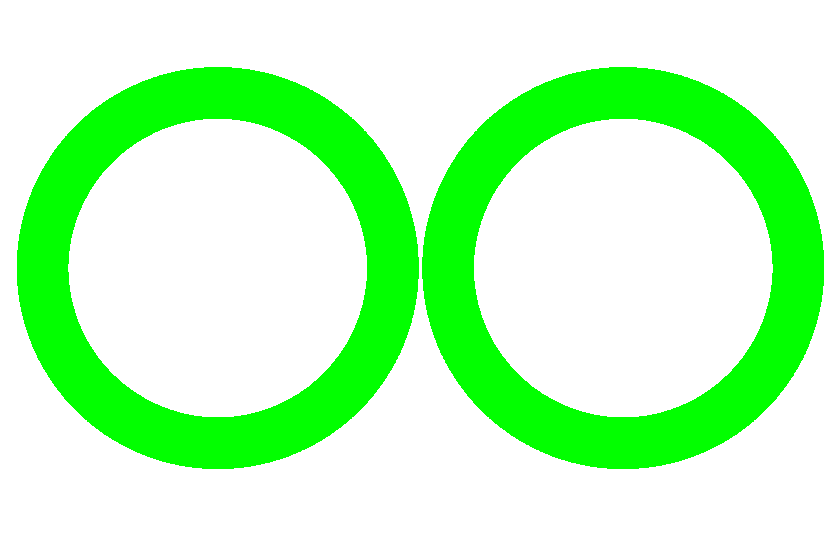}} &
		\subfloat[\label{fig:CTR_HS_V20_30sec}]{
			\includegraphics[width=0.28\columnwidth]{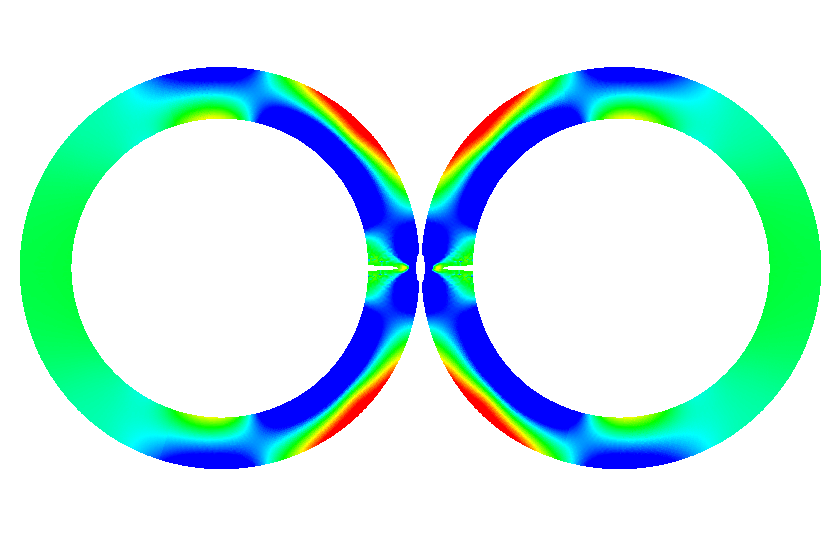}} &
		\subfloat[\label{fig:CTR_HS_V20_45sec}]{
			\includegraphics[width=0.28\columnwidth]{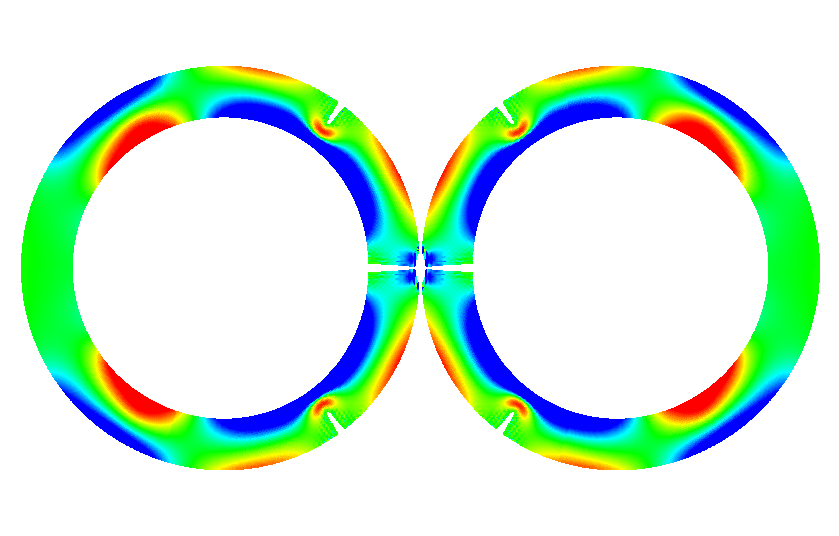}} \\
		\subfloat[\label{fig:CTR_HS_V20_60sec}]{
			\includegraphics[width=0.28\columnwidth]{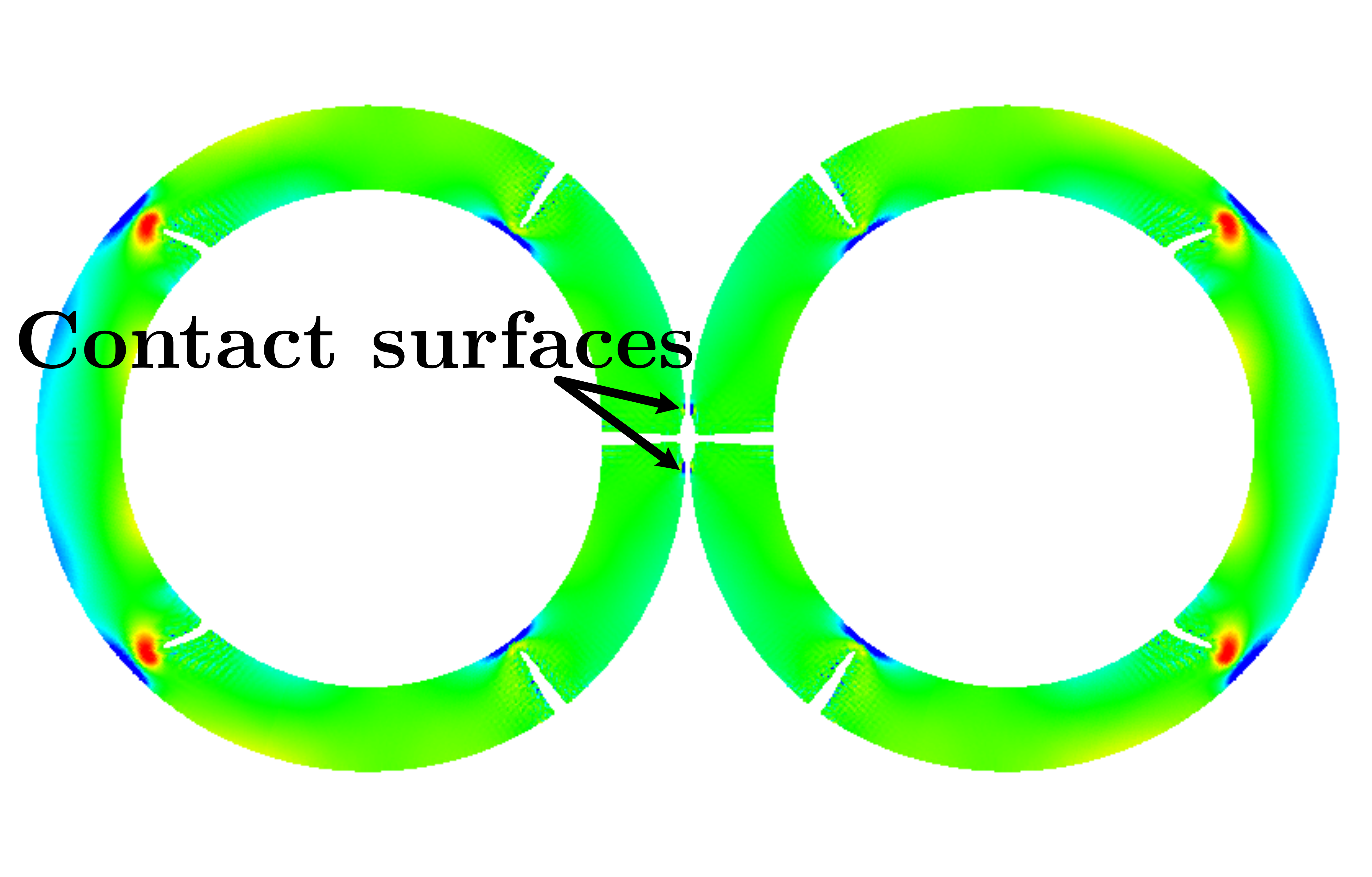}} &
		\subfloat[\label{fig:CTR_HS_V20_400sec}]{
			\includegraphics[width=0.28\columnwidth]{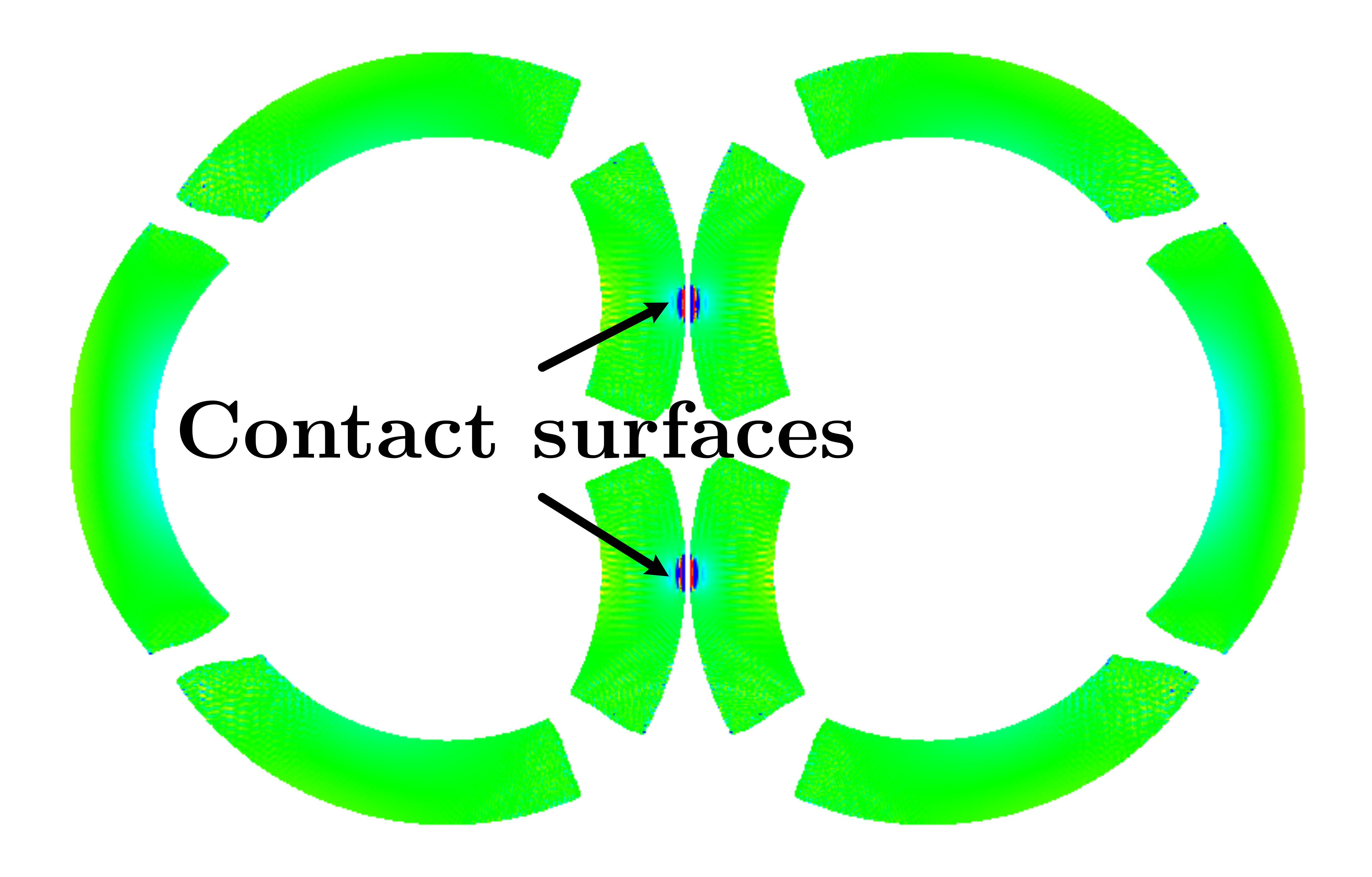}} &
		\subfloat[\label{fig:CTR_HS_V20_625sec}]{
			\includegraphics[width=0.28\columnwidth]{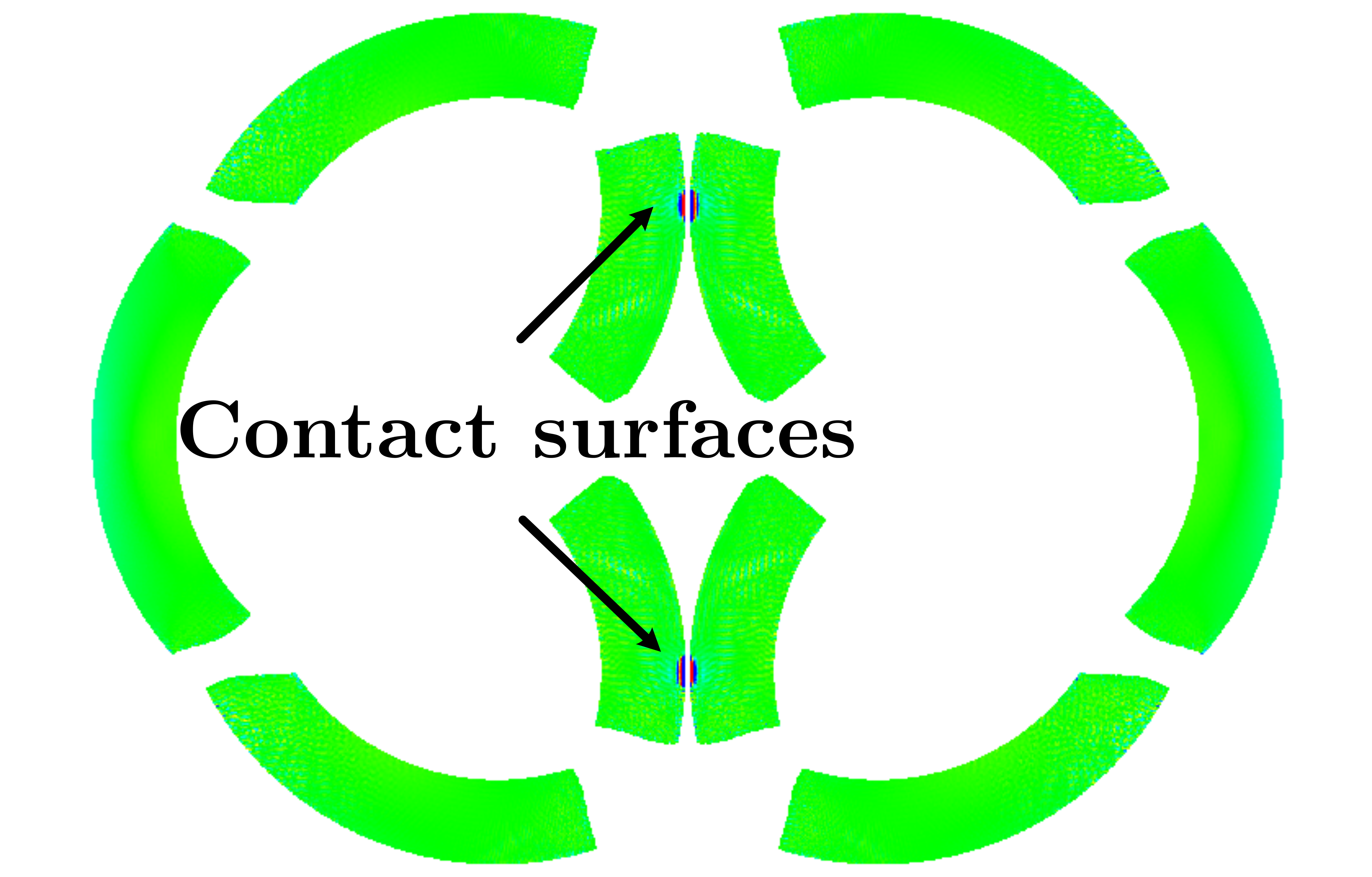} } \\
		\multicolumn{3}{c}{\subfloat{
				\includegraphics[width=1\columnwidth]{CTR_HS_ColorBar.png}}} 
	\end{tabular}
	\caption[]{Collision of two rings - Case (ii): Hydrostatic stresses for time steps \subref{fig:CTR_HS_V20_0sec} t=0 $\mu$s \subref{fig:CTR_HS_V20_30sec} t=30 $\mu$s \subref{fig:CTR_HS_V20_45sec} t=45 $\mu$s \subref{fig:CTR_HS_V20_60sec} t=60 $\mu$s \subref{fig:CTR_HS_V20_400sec} t=400 $\mu$s and \subref{fig:CTR_HS_V20_625sec} t=625 $\mu$s. Material points with $c_p<0.05$ have been removed.}
	\label{fig:CTR_HS_V20}
\end{figure}

The final deformed configuration of the problem is shown in Fig. \ref{fig:CTR_PF_V20_625sec}  where each ring is split into five fragments. The PF-MPM method naturally resolves the large displacement motion of the fragments, accounting also for the non-stationarity of the contact surfaces (see, also, Figs. \ref{fig:CTR_HS_V20_60sec}, \ref{fig:CTR_HS_V20_400sec} and \ref{fig:CTR_HS_V20_625sec}). Using a phase field driven fracture approximation allows both, the crack paths and the contact surfaces to not be tracked algorithmically during the simulation process. Furthermore, this is accomplished with no mesh distortion induced errors, contrary to a FEM based approach.

\subsection{Sphere-beam impact fracture problem} \label{lbl:SBI}

\begin{figure}
	\centering
	\includegraphics[width=0.40\columnwidth]{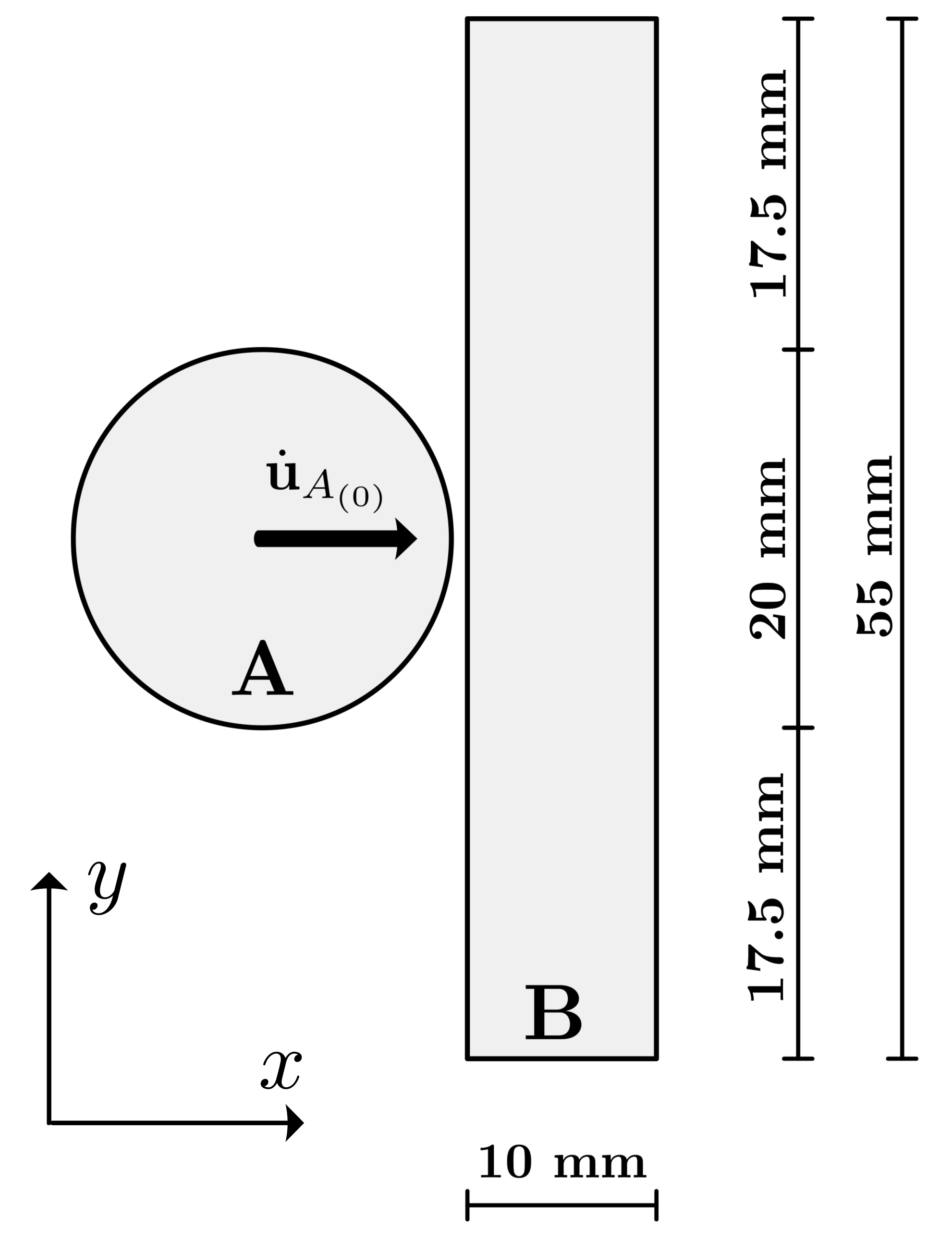}
	\caption[]{Geometry and initial conditions of the sphere-beam impact fracture problem. All boundaries are free.}
	\label{fig:SBI_GeoBCs}
\end{figure}

In this case, a sphere to beam impact fracture problem is examined. The geometry and boundary conditions of the problem are presented in Fig. \ref{fig:SBI_GeoBCs}. To examine the dependence of the resulting crack patterns and overall response of the beam on the level of assumed material anisotropy, three cases are examined, namely case (i), case (ii) and case (iii) with different surface energy densities. The material orientation is considered to be $\phi=+45^{o}$ with respect to $x$ axis (clockwise) for all cases.

The cell (patch) spacing is $h=0.125$ mm and plane stress conditions are assumed with thickness $2$ mm. The grid is formed by two knot vectors $\splitatcommas{\Xi=\{0,0,0,0.0025,0.0050,...,0.9950,0.9975,1,1,1\}}$ and $\splitatcommas{H = \{0,0,0,0.00167,0.00333,...,0.99667,0.99833,1,1,1\}}$, $242004$ control points and $400x600=240000$ cells. Two discrete fields are considered in this example, namely A for the sphere and B for the beam with the corresponding friction coefficient being $\mu_{f}=0.65$. The total number of material points is $536796$.

An initial velocity is applied to all material points in the sphere $\dot{\mathbf{u}}_{Ap_{(0)}}=0.02$ mm/$\mu$s while the beam is at rest at this stage. The initial distance between the sphere and the beam is considered to be $h=0.125$ mm. The elastic material parameters are chosen to be $E=190000$ N/mm\textsuperscript{2}, $\nu=0.30$ and $\rho=8000$ kg/m\textsuperscript{3} for both bodies. The solution procedure is implemented with a time step $\Delta t = 0.0125$ $\mu$s. The critical time step is $\tilde{\Delta t_{cr}}=0.018$ $\mu$s.

In case (i), the second order isotropic phase field model (i.e. $\gamma_{ijkl}=0$) is chosen with surface energy density $\mathcal{G}_{c} \left( \theta \right)=\bar{\mathcal{G}}_c=\mathcal{G}_{c_{max}}=\mathcal{G}_{c_{min}}=10.6066$ N/mm for the beam. In case (ii) , the second order isotropic phase field model is chosen again, but with a reduced surface energy density $\mathcal{G}_{c} \left( \theta \right)=\bar{\mathcal{G}}_c=\mathcal{G}_{c_{max}}=\mathcal{G}_{c_{min}}=9.75$ N/mm for the beam. Finally, in case (iii) the fourth order orthotropic model is utilized with anisotropic parameters $\bar{\mathcal{G}}_c=7.50$ N/mm, $\gamma_{1111}=80.00$, $\gamma_{2222}=1.00$, $\gamma_{1122}=0.00$ and $\gamma_{1212}=74.00$. These parameters result in maximum and minimum surface energy densities $\mathcal{G}_{c_{max}}=23.6892$ N/mm and $\mathcal{G}_{c_{min}}=10.6066$ N/mm, respectively for the beam. The surface energy density of the sphere is taken us sufficiently large (i.e. $\mathcal{G}_{c_{A}} \left( \theta \right) = 100 \mathcal{G}_{c_{B}} \left( \theta \right)$) so that the sphere remains undamaged in all cases. The length scale parameter is $l_0=0.25$ mm and $k_f=0.00$ in all cases. In all cases reported in this section, the projectile does not penetrate the beam, rather it bounces back and the beam undergoes free vibrations.

\subsubsection{Case (i): Isotropy - $\mathcal{G}_{c} \left( \theta \right)=10.6066$ {\normalfont N/mm}}

The time history of the total fracture energy is shown in Fig. \ref{fig:SBI_Fract_Iso1Cracks}. In Fig. \ref{fig:SBI_Fract_Iso1Cracks}, the path points (1-7) are labelled to facilitate discussion on the material response. Phase field and hydrostatic stress snapshots corresponding to points (1-6) are shown in Figs. \ref{fig:SBI_PF_C1SecondMPMIso1Cracks} and \ref{fig:SBI_HS_C1SecondMPMIso1Cracks}, respectively. 

The sphere initially comes into contact with the beam and fracture initiates at the contact surface (see. Fig. \ref{fig:SBI_PF_C1SecondMPMIso1Cracks_12sec} and point (2) in Fig. \ref{fig:SBI_Fract_Iso1Cracks}). Next, the right edge of the beam gradually degrades (see. \ref{fig:SBI_PF_C1SecondMPMIso1Cracks_40sec} and point (3) in Fig. \ref{fig:SBI_Fract_Iso1Cracks}) just before a crack initiates at the middle right-edge point. However, as the beam vibrates, the degradation continues at the left edge of the beam (see Fig. \ref{fig:SBI_PF_C1SecondMPMIso1Cracks_56sec} and point (4) in Fig. \ref{fig:SBI_Fract_Iso1Cracks}) a median crack develops and propagates just below the crack nucleation region (see Fig. \ref{fig:SBI_PF_C1SecondMPMIso1Cracks_70sec} and point (5) in Fig. \ref{fig:SBI_Fract_Iso1Cracks}). The complete crack path is shown in Fig. \ref{fig:SBI_PF_C1SecondMPMIso1Cracks_88sec}.

The results of Fig. \ref{fig:SBI_Fract_Iso1Cracks} can be further examined in view of the total fracture energy evolution. The evolution of the total fracture energy from point (1) to (2) corresponds to damage initiating between the sphere and the beam at their contact surface. Degradation at the right edge of the beam results in a further increase of the fracture energy corresponding to the path (2-3). Finally, the crack rapidly propagates from point (4) to (6). Hence, the total fracture energy corresponding to crack propagation is
\begin{equation*}
\prescript{(6)}{(4)}{  \Psi_{f} }^{}_{} = 311.82 - 94.29 = 217.53 \text{ mJ}.
\end{equation*}

This is in very good agreement with the analytical prediction as $A_{f} \cdot \mathcal{G}_{c} \left( \theta \right)=10 \cdot 2 \cdot 10.6066 = 212.13$ mJ, where $A_{f}$ stands for the fracture surface. The slight increase of the total fracture energy from point (6) to (7) corresponds to the marginal degradation of the beam material during the free vibration regime of its response.

\subsubsection{Case (ii): Isotropy - $\mathcal{G}_{c} \left( \theta \right)=9.75$ {\normalfont N/mm}}

Even though the variation in $\mathcal{G}_{c}$ is small compared to case (i), it results in a significantly different material response. The total fracture energy time-history for case (ii) is shown in Fig. \ref{fig:SBI_Fract_Iso2Cracks}. The evolution of the phase field and the hydrostatic stress for points ((1)-(6)) labelled in Fig. \ref{fig:SBI_Fract_Iso2Cracks} is shown in Figs. \ref{fig:SBI_PF_C1SecondMPMIso2Cracks} and \ref{fig:SBI_HS_C1SecondMPMIso2Cracks}, respectively.

Similar to case (i), the sphere initially comes into contact with the beam and causes damage at their contact surface. As a result, material degradation is observed at the right edge of the beam (see. Fig. \ref{fig:SBI_PF_C1SecondMPMIso2Cracks_30sec} and point (2) in Fig. \ref{fig:SBI_Fract_Iso2Cracks}) as in case (i). Contrary to case (i) however, a flexural crack initiates at the middle right-edge point of the beam due to maximum principal tensile stresses developing at the tensile fibre of the beam.   

As the beam oscillates the maximum tension region alternates between the two edges and the crack arrests (see Fig. \ref{fig:SBI_PF_C1SecondMPMIso2Cracks_40sec} and point (3) in Fig. \ref{fig:SBI_Fract_Iso2Cracks}). A second crack then initiates at the left edge (see Fig. \ref{fig:SBI_PF_C1SecondMPMIso2Cracks_68sec} and point (4) in Fig. \ref{fig:SBI_Fract_Iso2Cracks}) and propagates (see Fig. \ref{fig:SBI_PF_C1SecondMPMIso2Cracks_72sec} and point (5) in Fig. \ref{fig:SBI_Fract_Iso2Cracks}) until the two cracks finally merge as shown in Fig. \ref{fig:SBI_PF_C1SecondMPMIso2Cracks_78sec}.

As in case (i), the evolution of the fracture energy (shown shown in Fig. \ref{fig:SBI_Fract_Iso2Cracks}) is consistent with the observed response. The first crack (right crack) initiates at point (2) and stops at point (3). The second crack (left crack) propagates from point (4) to (6). Therefore, the total fracture energy is
\begin{equation*}
\Psi_{f} = \prescript{(3)}{(2)}{  \Psi_{f} }^{}_{} + \prescript{(6)}{(4)}{  \Psi_{f} }^{}_{} = (196.55-87.37) + (306.81-226.39) \text{ mJ}.
\end{equation*}

This is again in very good agreement with the analytical prediction as in this case $A_{f} \cdot \mathcal{G}_{c} \left( \theta \right)=10 \cdot 2 \cdot 9.75 = 195$ mJ.

\subsubsection{Case (iii): Orthotropy}

\begin{figure}
	\centering
	\begin{tabular}{lcr}
		{\subfloat[\label{fig:SBI_Fract_Iso1Cracks}]{
				\includegraphics[width=0.37\columnwidth]{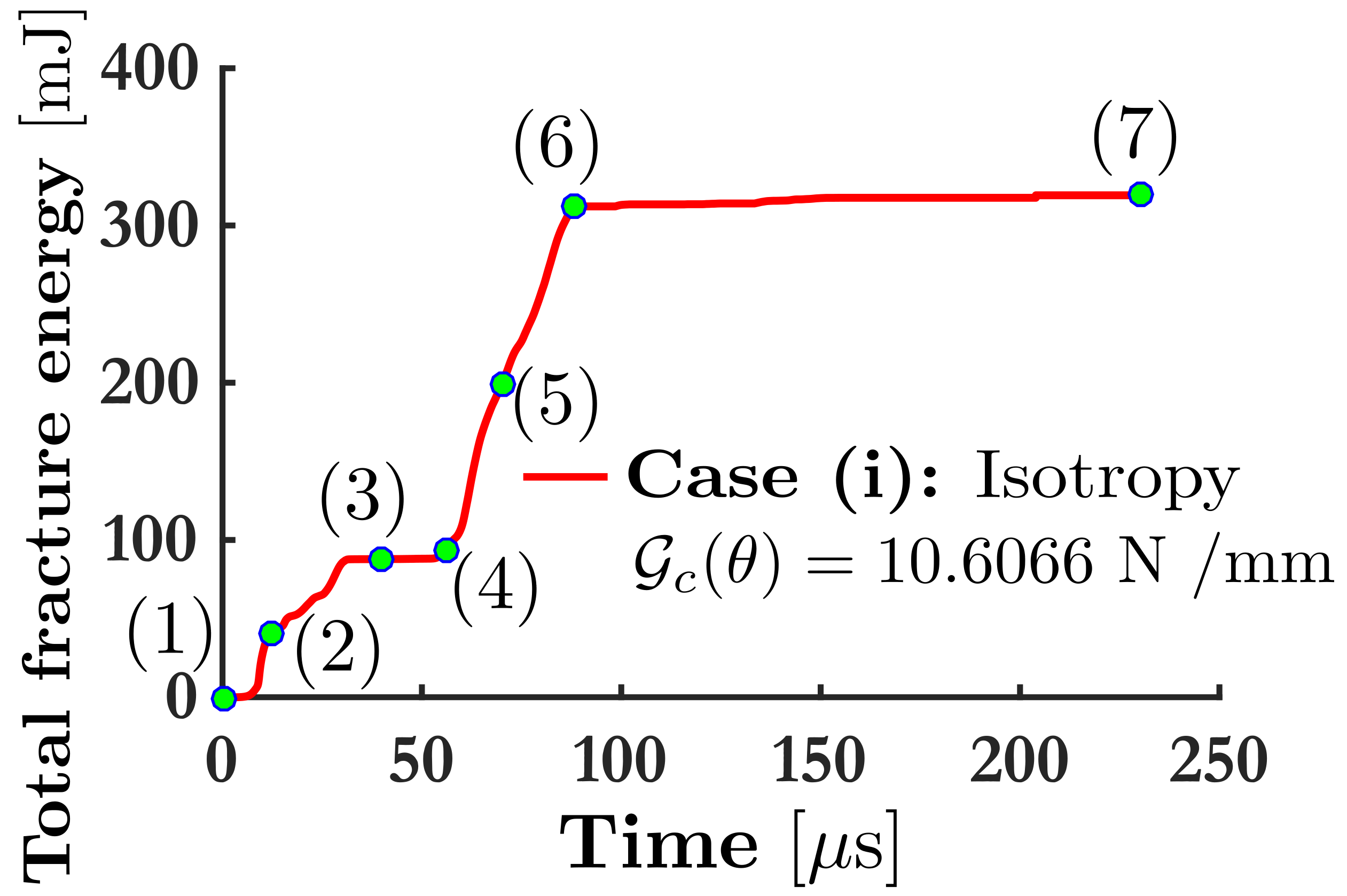}}} &
		\subfloat[\label{fig:SBI_Fract_Iso2Cracks}]{
			\includegraphics[width=0.37\columnwidth]{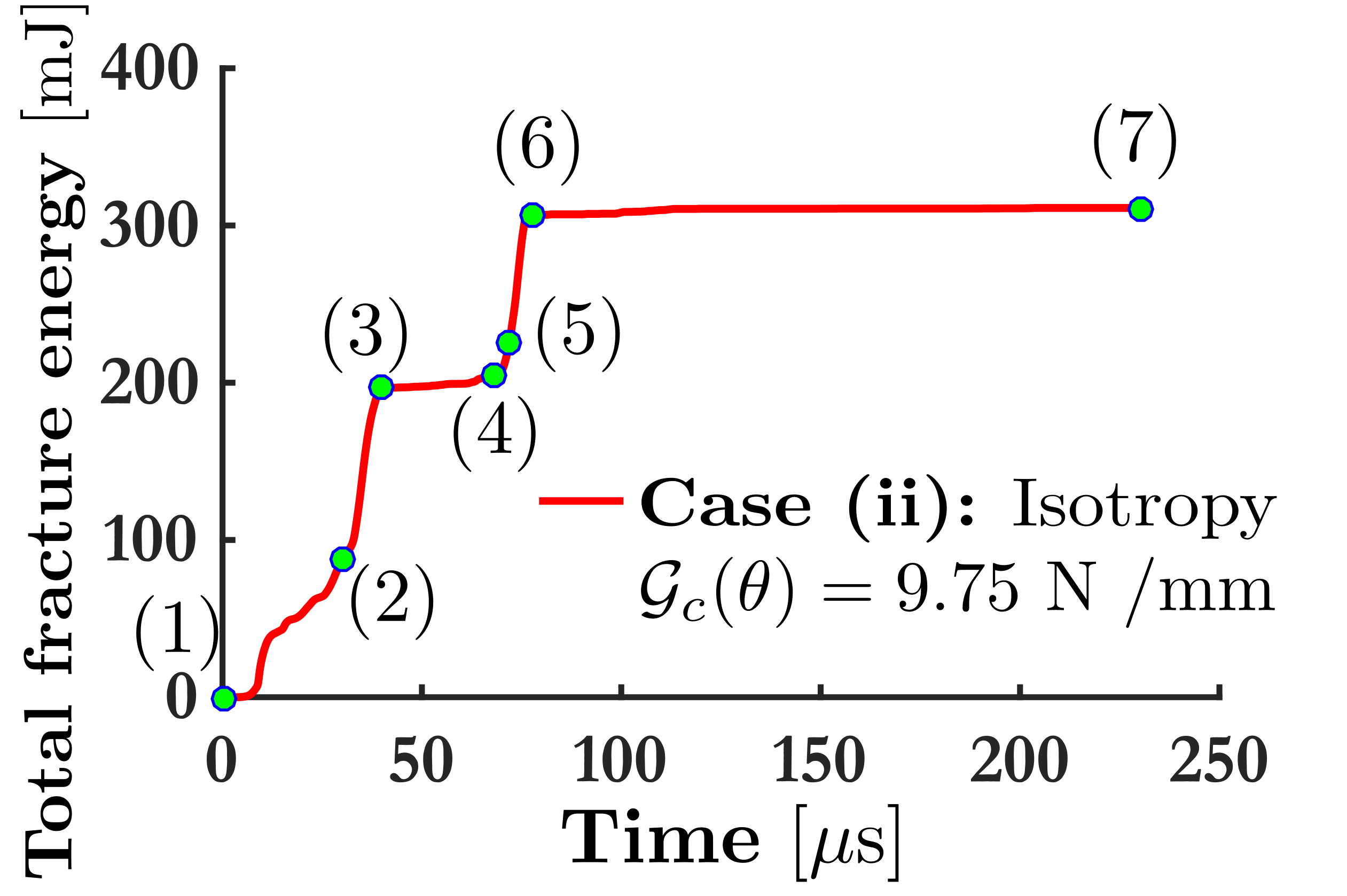}} &
		{\subfloat[\label{fig:SBI_Fract_Ortho45}]{
				\includegraphics[width=0.37\columnwidth]{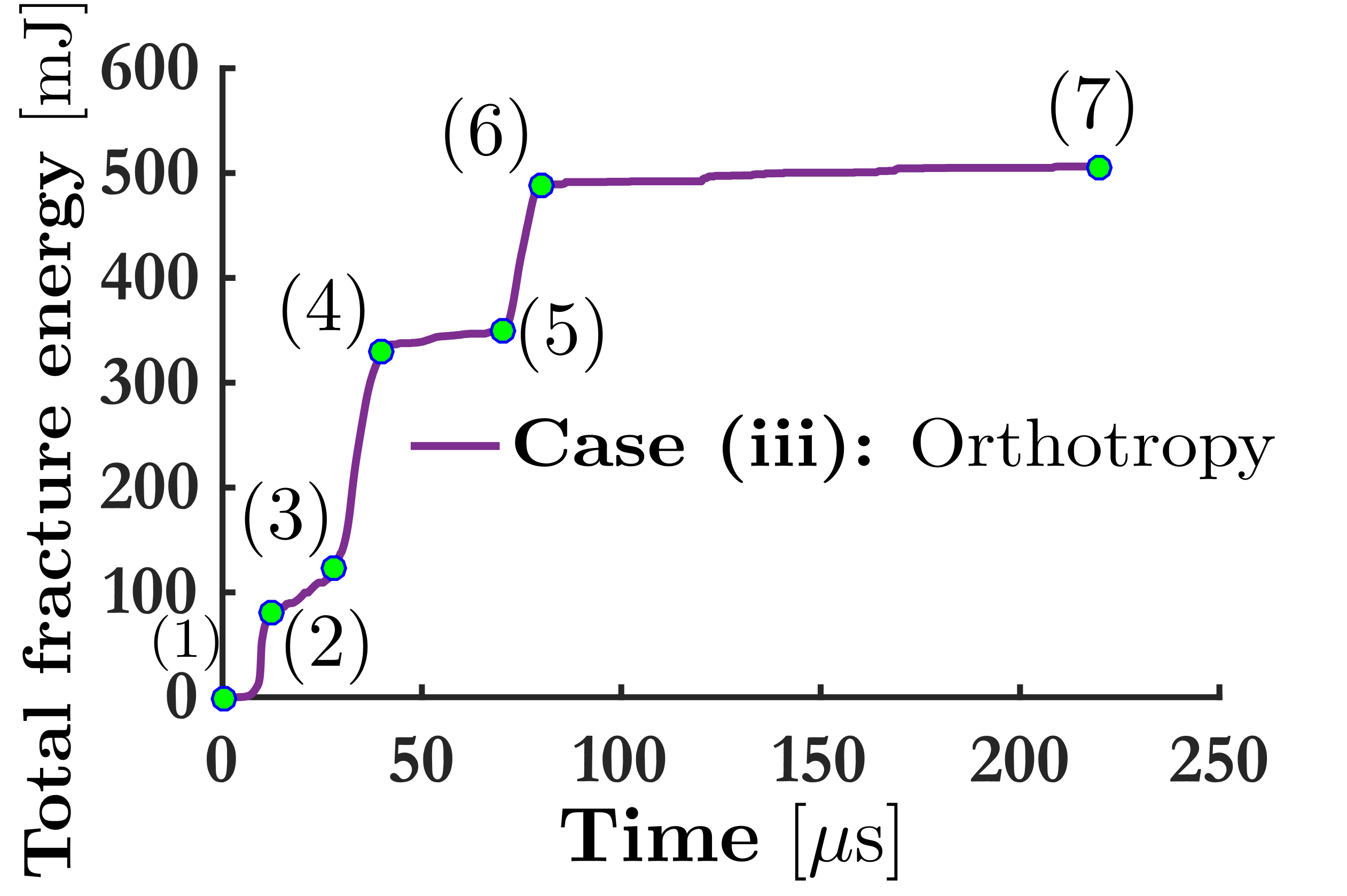}}}
	\end{tabular}
	\caption[]{Sphere-beam impact fracture problem: Total fracture energy time history for \subref{fig:SBI_Fract_Iso1Cracks} case (i): PF-MPM 2nd order isotropic model and $\mathcal{G}_{c} \left( \theta \right) = 10.6066$ N/mm \subref{fig:SBI_Fract_Iso2Cracks} case (ii): PF-MPM 2nd order isotropic model and $\mathcal{G}_{c} \left( \theta \right) = 9.75$ N/mm and \subref{fig:SBI_Fract_Ortho45} case (iii): PF-MPM 4th order orthotropic model for the beam.}
	\label{fig:SBI_Fract}
\end{figure}

Orthotropic anisotropy with a material orientation $\phi = +45^{o}$ results in two cracks at each beam edge (right and left) that do not coincide with the horizontal axis as in case (ii). The evolution of the phase field and the hydrostatic stress are represented for several time steps in Figs. \ref{fig:SBI_PF_C1FourthMPMOrtho} and \ref{fig:SBI_HS_C1FourthMPMOrtho}, respectively. The characteristic points ((1)-(7)) of that analysis are shown in Fig. \ref{fig:SBI_Fract_Ortho45}. 

Similarly to the previous cases, damage initiation is observed at the contact surface (see Fig. \ref{fig:SBI_PF_C1FourthMPMOrtho_12sec} and point (2) in Fig. \ref{fig:SBI_Fract_Ortho45}). Next, degradation occurs at the left edge of the beam (see Fig. \ref{fig:SBI_PF_C1FourthMPMOrtho_28sec} and point (3) in Fig. \ref{fig:SBI_Fract_Ortho45}). The first crack (right crack) initiates at middle right-edge point of the beam and propagates along the material’s week direction until it arrests in the vicinity of the beam’s neutral axis (see Fig. \ref{fig:SBI_PF_C1FourthMPMOrtho_40sec} and point (4) in Fig. \ref{fig:SBI_Fract_Ortho45}). After impact, further degradation occurs due to the beam’s free vibration resulting in degradation to its left edge (see Fig. \ref{fig:SBI_PF_C1FourthMPMOrtho_70sec} and point (5) in Fig. \ref{fig:SBI_Fract_Ortho45}). Finally, a second crack (left crack) initiates at the middle left-edge of the beam and propagates along the material’s week direction (see Fig. \ref{fig:SBI_PF_C1FourthMPMOrtho_80sec} and point (6) in Fig. \ref{fig:SBI_Fract_Ortho45}). Similar to the first crack, the second crack arrests in the vicinity of the beam’s neutral axis. The final crack paths are shown in Fig. \ref{fig:SBI_HS_C1FourthMPMOrtho_80sec} where the two cracks do not merge as in case (ii).

\begin{figure}
	\centering
	\begin{tabular}{ccc}
		\subfloat[\label{fig:SBI_PF_C1SecondMPMIso1Cracks_0sec}]{
			\includegraphics[width=0.25\columnwidth]{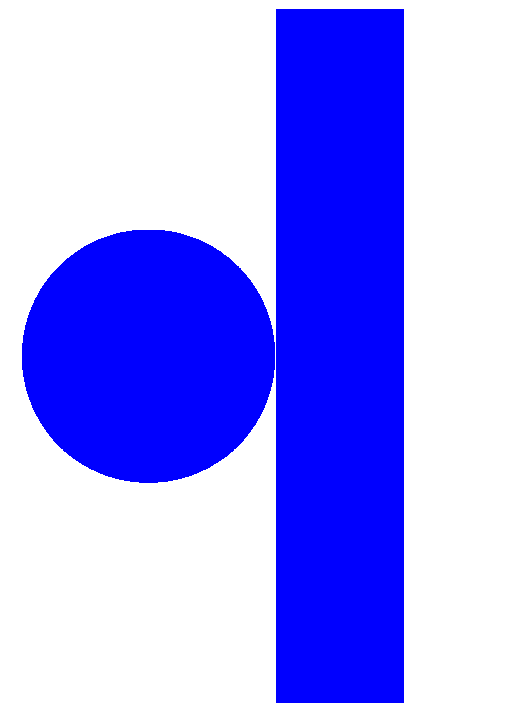}} &
		\subfloat[\label{fig:SBI_PF_C1SecondMPMIso1Cracks_12sec}]{
			\includegraphics[width=0.25\columnwidth]{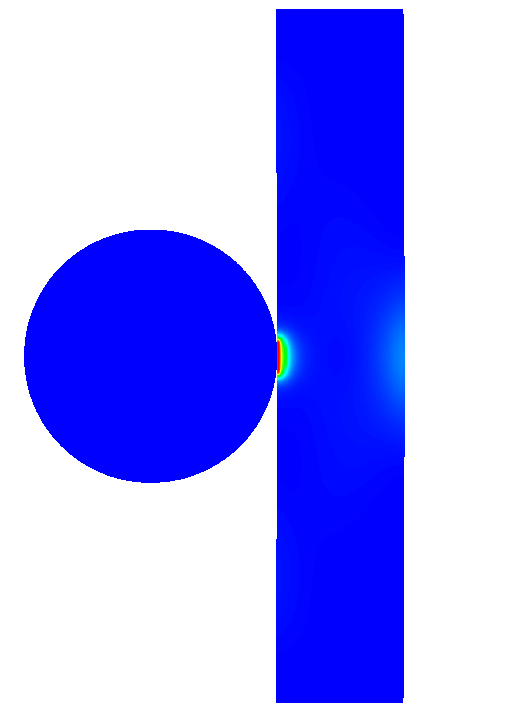}} &
		\subfloat[\label{fig:SBI_PF_C1SecondMPMIso1Cracks_40sec}]{
			\includegraphics[width=0.25\columnwidth]{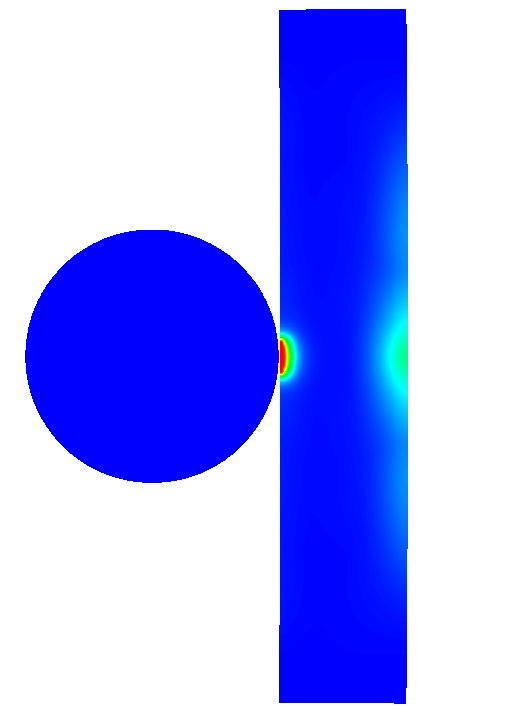}} \\
		\subfloat[\label{fig:SBI_PF_C1SecondMPMIso1Cracks_56sec}]{
			\includegraphics[width=0.25\columnwidth]{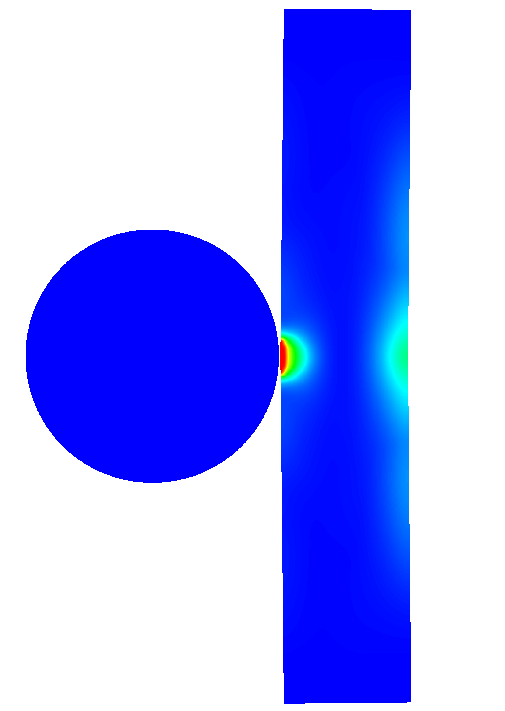}} &
		\subfloat[\label{fig:SBI_PF_C1SecondMPMIso1Cracks_70sec}]{
			\includegraphics[width=0.25\columnwidth]{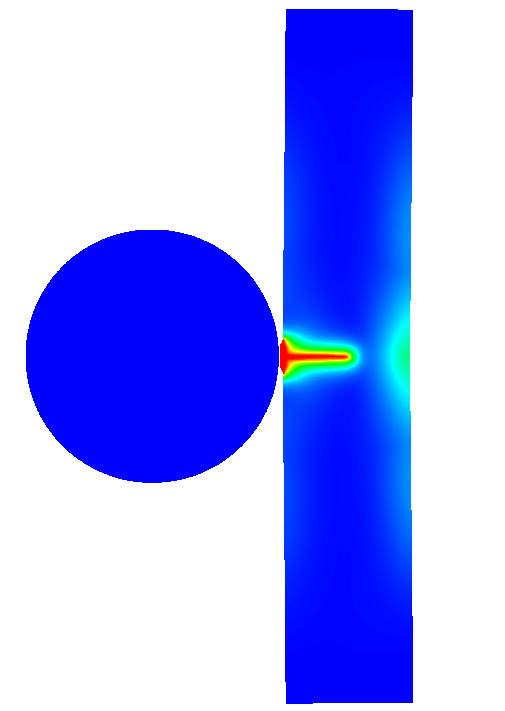}} &
		\subfloat[\label{fig:SBI_PF_C1SecondMPMIso1Cracks_88sec}]{
			\includegraphics[width=0.25\columnwidth]{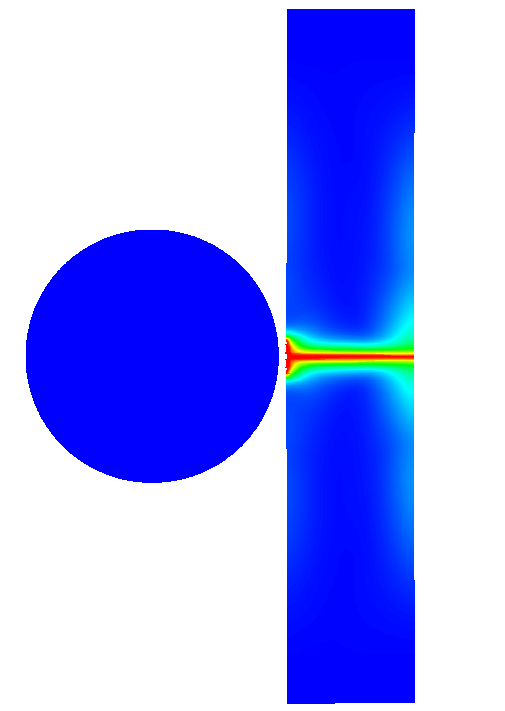} } \\
		\multicolumn{3}{c}{\subfloat{
				\includegraphics[width=0.75\columnwidth]{PF_ColorBar.png}}} 		
	\end{tabular}
	\caption[]{Sphere-beam impact fracture problem: Phase field for time steps \subref{fig:SBI_PF_C1SecondMPMIso1Cracks_0sec} t=0 $\mu$s \subref{fig:SBI_PF_C1SecondMPMIso1Cracks_12sec} t=12 $\mu$s \subref{fig:SBI_PF_C1SecondMPMIso1Cracks_40sec} t=40 $\mu$s \subref{fig:SBI_PF_C1SecondMPMIso1Cracks_56sec} t=56 $\mu$s \subref{fig:SBI_PF_C1SecondMPMIso1Cracks_70sec} t=70 $\mu$s and \subref{fig:SBI_PF_C1SecondMPMIso1Cracks_88sec} t=88 $\mu$s. Results for case (i): PF-MPM 2nd order isotropic model and $\mathcal{G}_{c} \left( \theta \right) = 10.6066$ N/mm for the beam.}
	\label{fig:SBI_PF_C1SecondMPMIso1Cracks}
\end{figure}

\begin{figure}
	\centering
	\begin{tabular}{ccc}
		\subfloat[\label{fig:SBI_HS_C1SecondMPMIso1Cracks_0sec}]{
			\includegraphics[width=0.25\columnwidth]{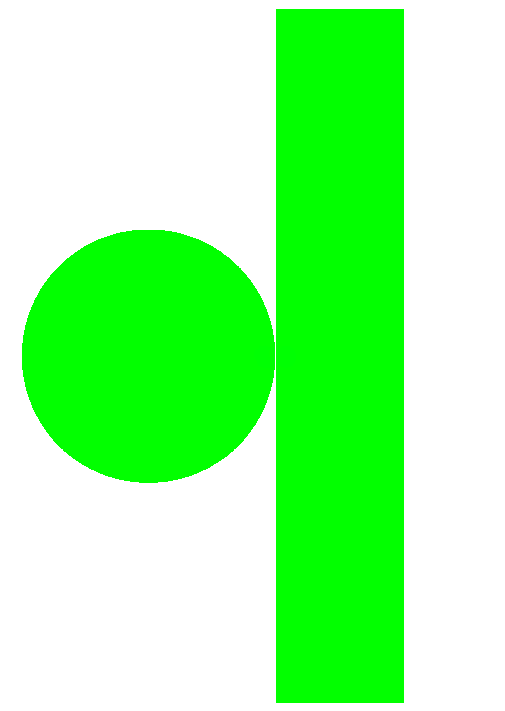}} &
		\subfloat[\label{fig:SBI_HS_C1SecondMPMIso1Cracks_12sec}]{
			\includegraphics[width=0.25\columnwidth]{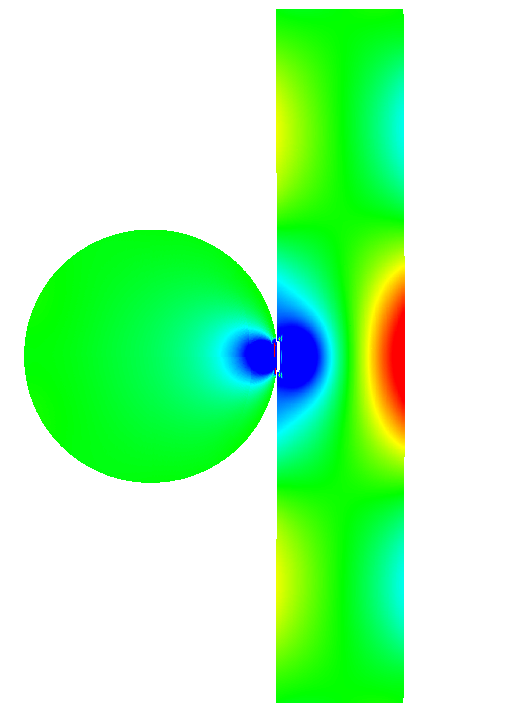}} &
		\subfloat[\label{fig:SBI_HS_C1SecondMPMIso1Cracks_40sec}]{
			\includegraphics[width=0.25\columnwidth]{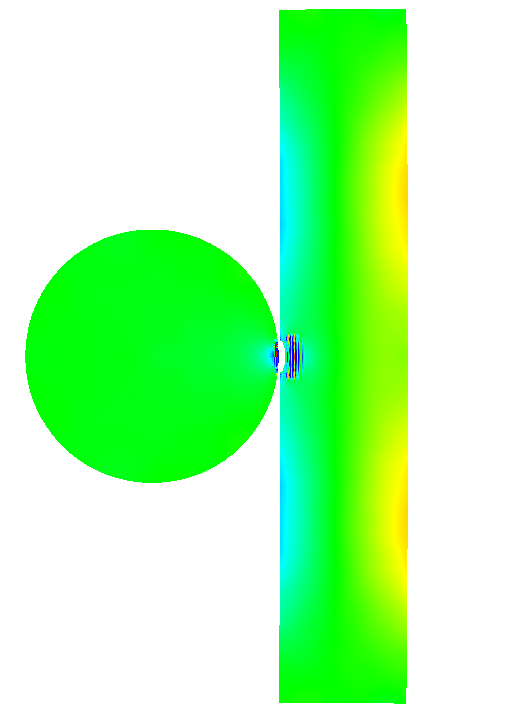}} \\
		\subfloat[\label{fig:SBI_HS_C1SecondMPMIso1Cracks_56sec}]{
			\includegraphics[width=0.25\columnwidth]{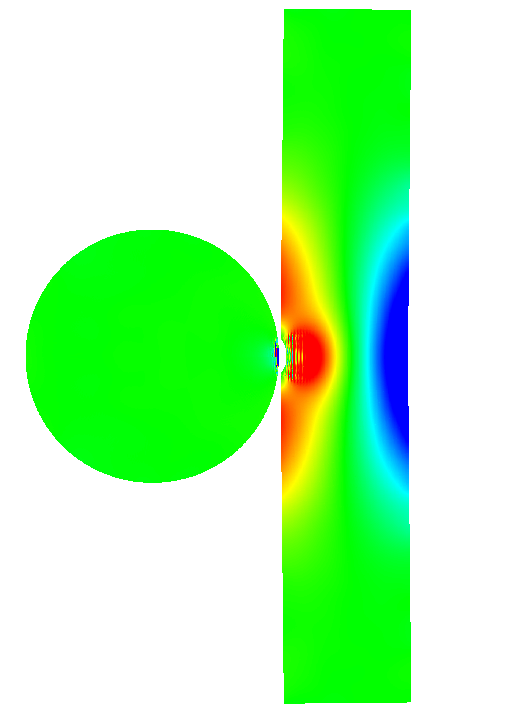}} &
		\subfloat[\label{fig:SBI_HS_C1SecondMPMIso1Cracks_70sec}]{
			\includegraphics[width=0.25\columnwidth]{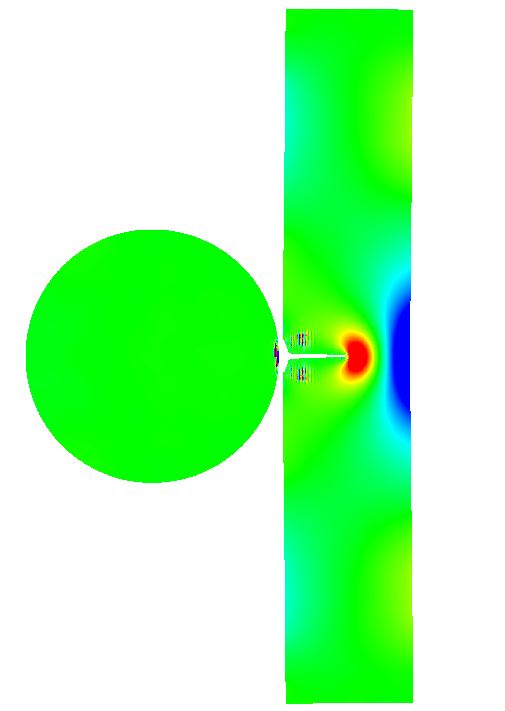}} &
		\subfloat[\label{fig:SBI_HS_C1SecondMPMIso1Cracks_88sec}]{
			\includegraphics[width=0.25\columnwidth]{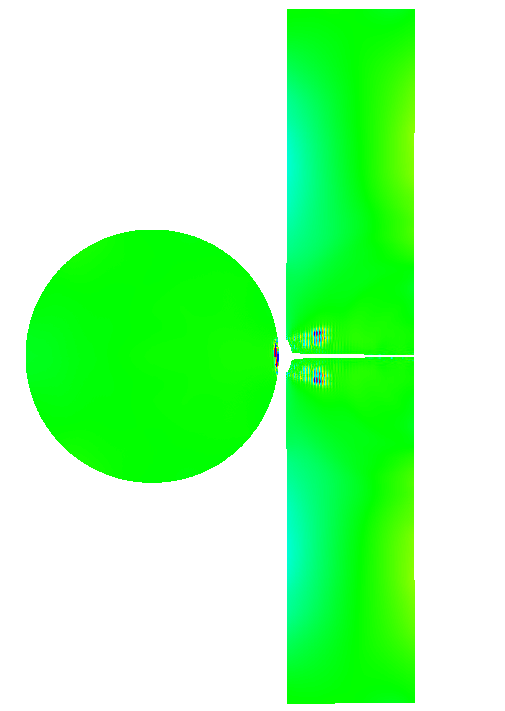} } \\
		\multicolumn{3}{c}{\subfloat{
				\includegraphics[width=0.75\columnwidth]{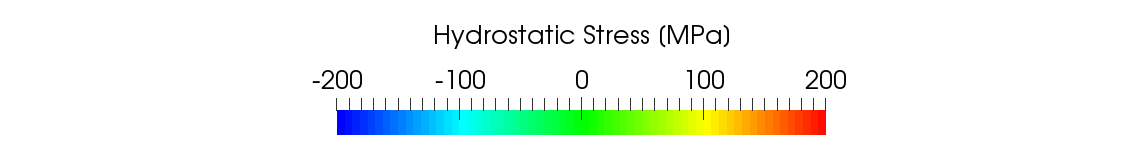}}} 		
	\end{tabular}
	\caption[]{Sphere-beam impact fracture problem: Hydrostatic stresses for time steps \subref{fig:SBI_HS_C1SecondMPMIso1Cracks_0sec} t=0 $\mu$s \subref{fig:SBI_HS_C1SecondMPMIso1Cracks_12sec} t=12 $\mu$s \subref{fig:SBI_HS_C1SecondMPMIso1Cracks_40sec} t=40 $\mu$s \subref{fig:SBI_HS_C1SecondMPMIso1Cracks_56sec} t=56 $\mu$s \subref{fig:SBI_HS_C1SecondMPMIso1Cracks_70sec} t=70 $\mu$s and \subref{fig:SBI_HS_C1SecondMPMIso1Cracks_88sec} t=88 $\mu$s. Results for case (i): PF-MPM 2nd order isotropic model and $\mathcal{G}_{c} \left( \theta \right) = 10.6066$ N/mm for the beam. Material points with $c_p<0.08$ have been removed.}
	\label{fig:SBI_HS_C1SecondMPMIso1Cracks}
\end{figure}

\begin{figure}
	\centering
	\begin{tabular}{ccc}
		\subfloat[\label{fig:SBI_PF_C1SecondMPMIso2Cracks_0sec}]{
			\includegraphics[width=0.25\columnwidth]{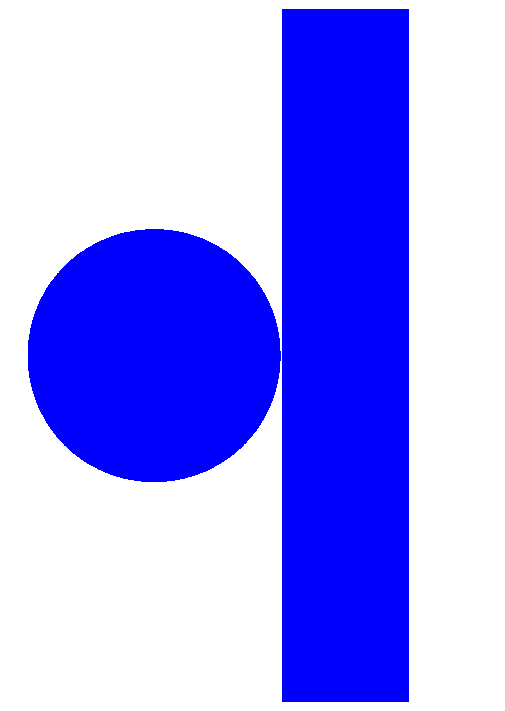}} &
		\subfloat[\label{fig:SBI_PF_C1SecondMPMIso2Cracks_30sec}]{
			\includegraphics[width=0.25\columnwidth]{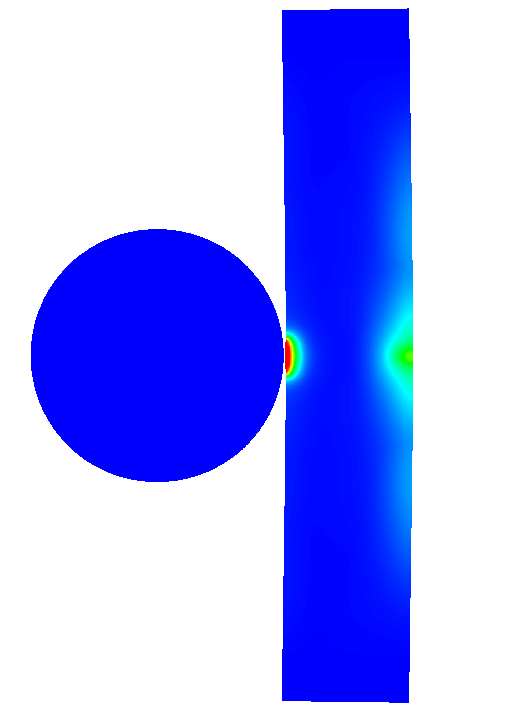}} &
		\subfloat[\label{fig:SBI_PF_C1SecondMPMIso2Cracks_40sec}]{
			\includegraphics[width=0.25\columnwidth]{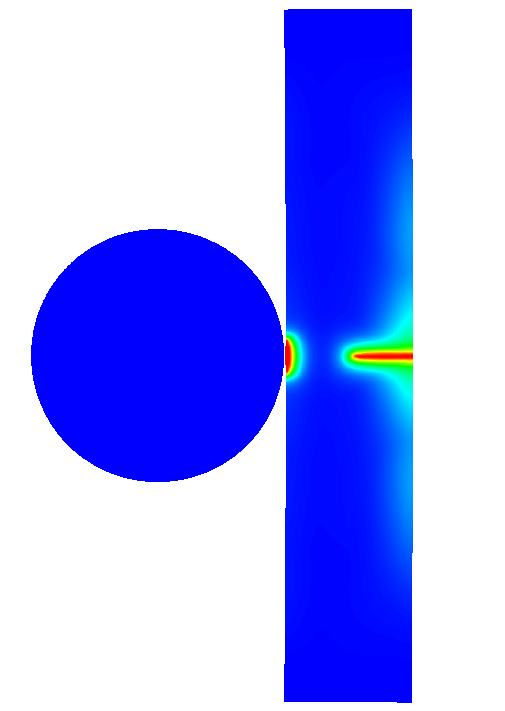}} \\
		\subfloat[\label{fig:SBI_PF_C1SecondMPMIso2Cracks_68sec}]{
			\includegraphics[width=0.25\columnwidth]{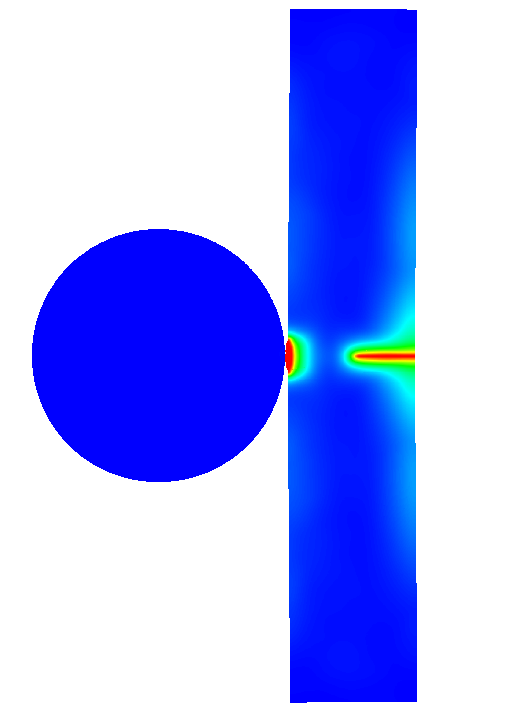}} &
		\subfloat[\label{fig:SBI_PF_C1SecondMPMIso2Cracks_72sec}]{
			\includegraphics[width=0.25\columnwidth]{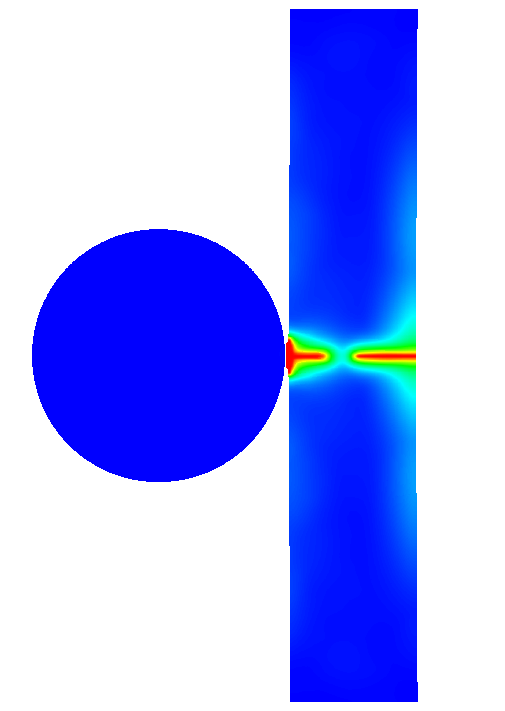}} &
		\subfloat[\label{fig:SBI_PF_C1SecondMPMIso2Cracks_78sec}]{
			\includegraphics[width=0.25\columnwidth]{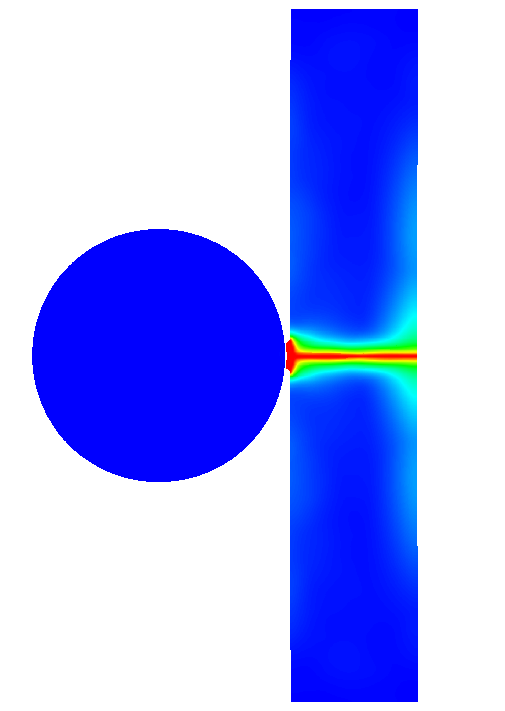} } \\
		\multicolumn{3}{c}{\subfloat{
				\includegraphics[width=0.75\columnwidth]{PF_ColorBar.png}}} 		
	\end{tabular}
	\caption[]{Sphere-beam impact fracture problem: Phase field for time steps \subref{fig:SBI_PF_C1SecondMPMIso2Cracks_0sec} t=0 $\mu$s \subref{fig:SBI_PF_C1SecondMPMIso2Cracks_30sec} t=30 $\mu$s \subref{fig:SBI_PF_C1SecondMPMIso2Cracks_40sec} t=40 $\mu$s \subref{fig:SBI_PF_C1SecondMPMIso2Cracks_68sec} t=68 $\mu$s \subref{fig:SBI_PF_C1SecondMPMIso2Cracks_72sec} t=72 $\mu$s and \subref{fig:SBI_PF_C1SecondMPMIso2Cracks_78sec} t=78 $\mu$s. Results for case (ii): PF-MPM 2nd order isotropic model and $\mathcal{G}_{c} \left( \theta \right) = 9.75$ N/mm for the beam.}
	\label{fig:SBI_PF_C1SecondMPMIso2Cracks}
\end{figure}

\begin{figure}
	\centering
	\begin{tabular}{ccc}
		\subfloat[\label{fig:SBI_HS_C1SecondMPMIso2Cracks_0sec}]{
			\includegraphics[width=0.25\columnwidth]{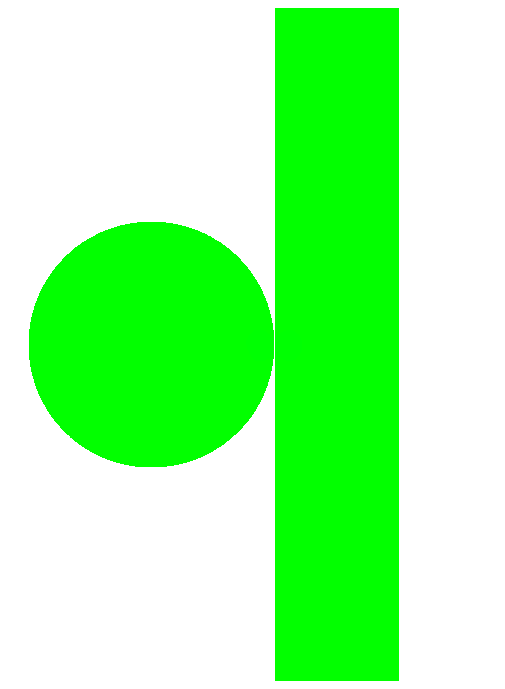}} &
		\subfloat[\label{fig:SBI_HS_C1SecondMPMIso2Cracks_30sec}]{
			\includegraphics[width=0.25\columnwidth]{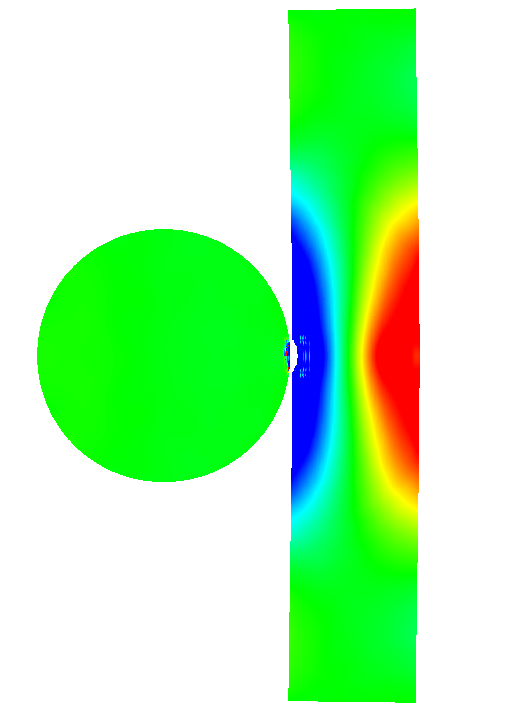}} &
		\subfloat[\label{fig:SBI_HS_C1SecondMPMIso2Cracks_40sec}]{
			\includegraphics[width=0.25\columnwidth]{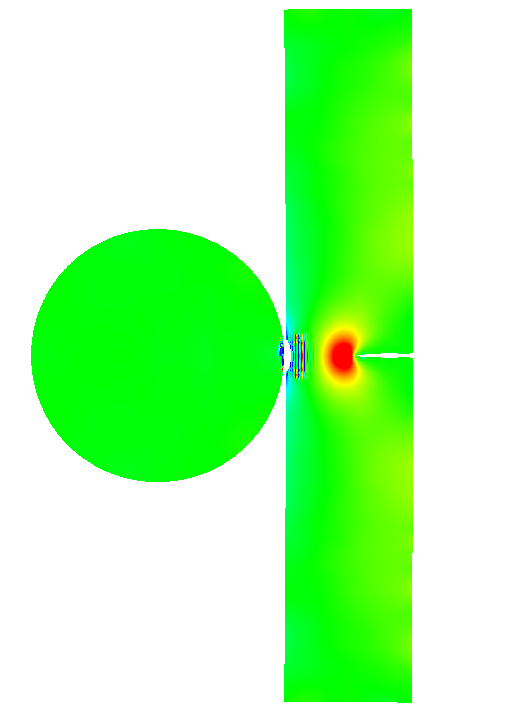}} \\
		\subfloat[\label{fig:SBI_HS_C1SecondMPMIso2Cracks_68sec}]{
			\includegraphics[width=0.25\columnwidth]{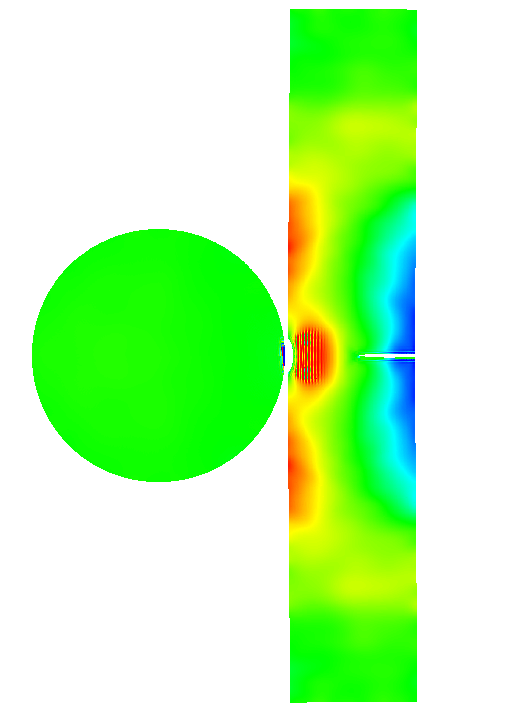}} &
		\subfloat[\label{fig:SBI_HS_C1SecondMPMIso2Cracks_72sec}]{
			\includegraphics[width=0.25\columnwidth]{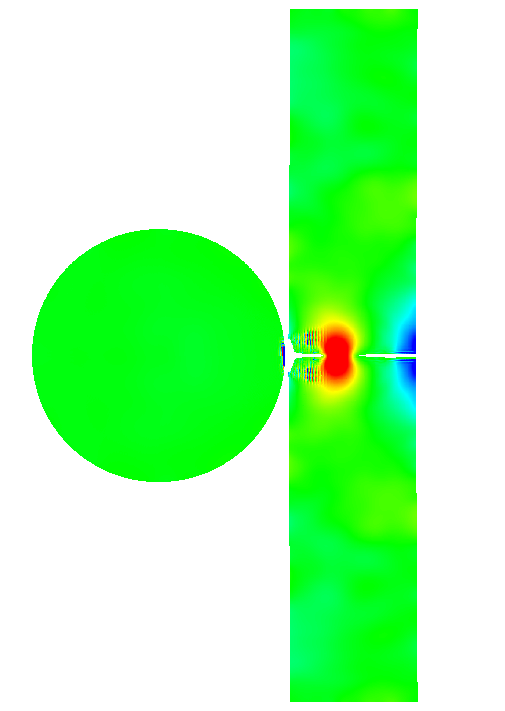}} &
		\subfloat[\label{fig:SBI_HS_C1SecondMPMIso2Cracks_78sec}]{
			\includegraphics[width=0.25\columnwidth]{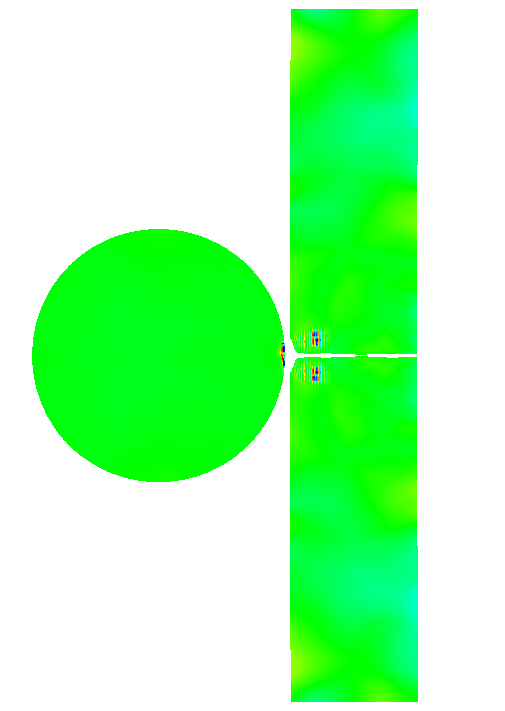} } \\
		\multicolumn{3}{c}{\subfloat{
				\includegraphics[width=0.75\columnwidth]{HS_ColorBar1.png}}} 		
	\end{tabular}
	\caption[]{Sphere-beam impact fracture problem: Hydrostatic stresses for time steps \subref{fig:SBI_HS_C1SecondMPMIso2Cracks_0sec} t=0 $\mu$s \subref{fig:SBI_HS_C1SecondMPMIso2Cracks_30sec} t=30 $\mu$s \subref{fig:SBI_HS_C1SecondMPMIso2Cracks_40sec} t=40 $\mu$s \subref{fig:SBI_HS_C1SecondMPMIso2Cracks_68sec} t=68 $\mu$s \subref{fig:SBI_HS_C1SecondMPMIso2Cracks_72sec} t=72 $\mu$s and \subref{fig:SBI_HS_C1SecondMPMIso2Cracks_78sec} t=78 $\mu$s. Results for case (ii): PF-MPM 2nd order isotropic model and $\mathcal{G}_{c} \left( \theta \right) = 9.75$ N/mm for the beam. Material points with $c_p<0.08$ have been removed.}
	\label{fig:SBI_HS_C1SecondMPMIso2Cracks}
\end{figure}

\begin{figure}
	\centering
	\begin{tabular}{lcr}
		\subfloat[\label{fig:SBI_PF_C1FourthMPMOrtho_0sec}]{
			\includegraphics[width=0.25\columnwidth]{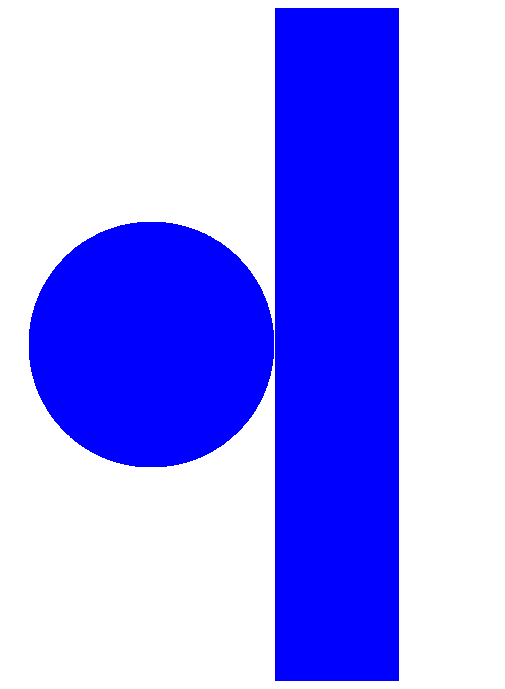}} &
		\subfloat[\label{fig:SBI_PF_C1FourthMPMOrtho_12sec}]{
			\includegraphics[width=0.25\columnwidth]{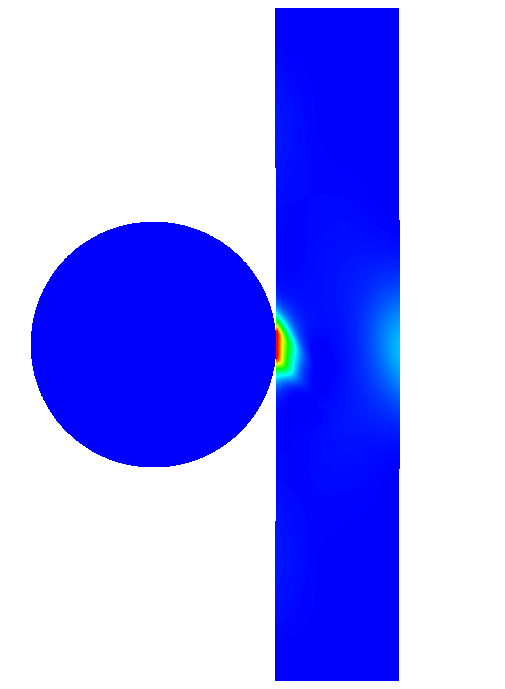}} &
		\subfloat[\label{fig:SBI_PF_C1FourthMPMOrtho_28sec}]{
			\includegraphics[width=0.25\columnwidth]{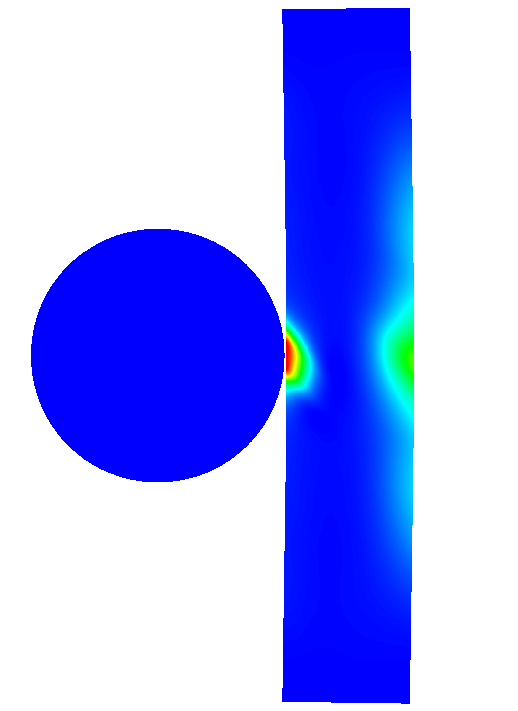}} \\
		\subfloat[\label{fig:SBI_PF_C1FourthMPMOrtho_40sec}]{
			\includegraphics[width=0.25\columnwidth]{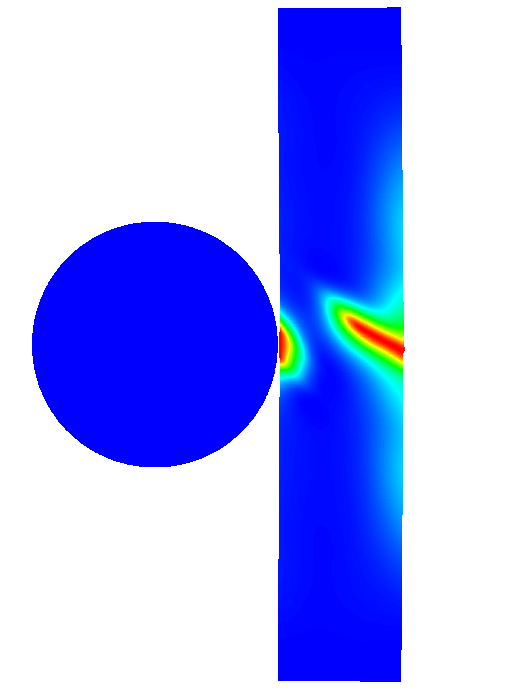}} &
		\subfloat[\label{fig:SBI_PF_C1FourthMPMOrtho_70sec}]{
			\includegraphics[width=0.25\columnwidth]{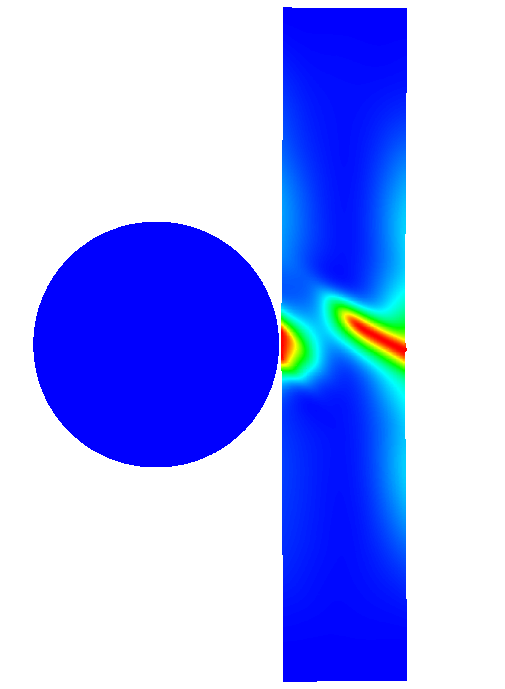}} &
		\subfloat[\label{fig:SBI_PF_C1FourthMPMOrtho_80sec}]{
			\includegraphics[width=0.25\columnwidth]{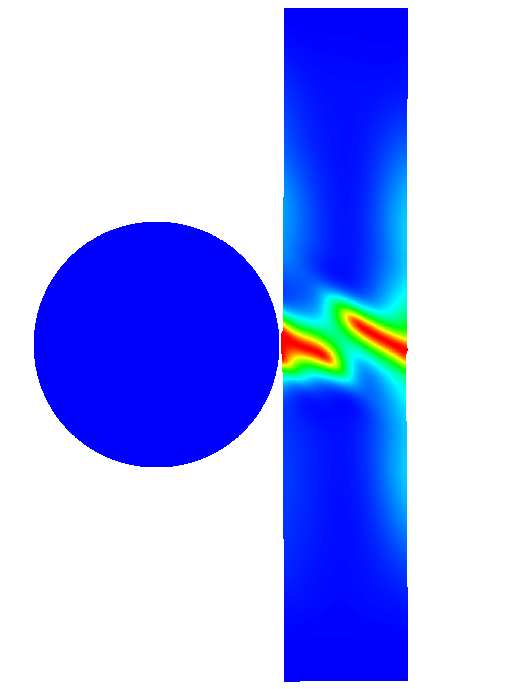} } \\
		\multicolumn{3}{c}{\subfloat{
				\includegraphics[width=0.75\columnwidth]{PF_ColorBar.png}}} 		
	\end{tabular}
	\caption[]{Sphere-beam impact fracture problem: Phase field for time steps \subref{fig:SBI_PF_C1FourthMPMOrtho_0sec} t=0 $\mu$s \subref{fig:SBI_PF_C1FourthMPMOrtho_12sec} t=12 $\mu$s \subref{fig:SBI_PF_C1FourthMPMOrtho_28sec} t=28 $\mu$s \subref{fig:SBI_PF_C1FourthMPMOrtho_40sec} t=40 $\mu$s \subref{fig:SBI_PF_C1FourthMPMOrtho_70sec} t=70 $\mu$s and \subref{fig:SBI_PF_C1FourthMPMOrtho_80sec} t=80 $\mu$s. Results for case (iii): PF-MPM 4th order orthotropic model.}
	\label{fig:SBI_PF_C1FourthMPMOrtho}
\end{figure}

\begin{figure}
	\centering
	\begin{tabular}{lcr}
		\subfloat[\label{fig:SBI_HS_C1FourthMPMOrtho_0sec}]{
			\includegraphics[width=0.25\columnwidth]{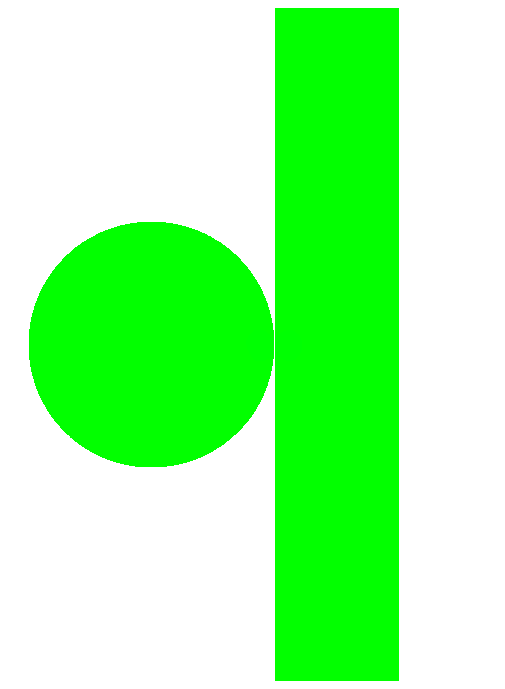}} &
		\subfloat[\label{fig:SBI_HS_C1FourthMPMOrtho_12sec}]{
			\includegraphics[width=0.25\columnwidth]{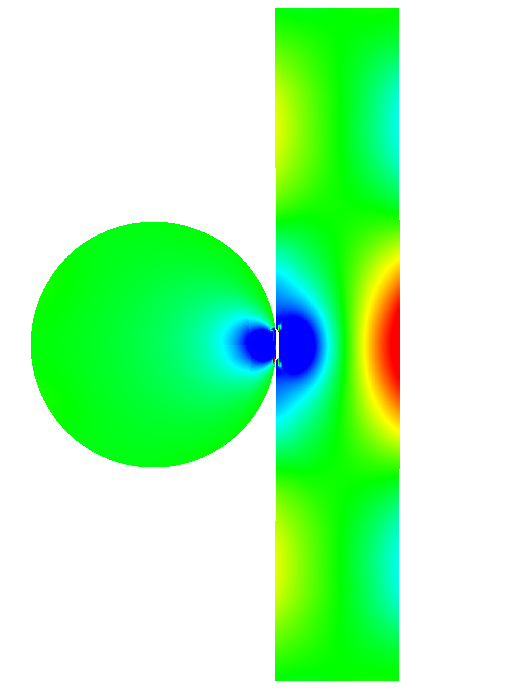}} &
		\subfloat[\label{fig:SBI_HS_C1FourthMPMOrtho_28sec}]{
			\includegraphics[width=0.25\columnwidth]{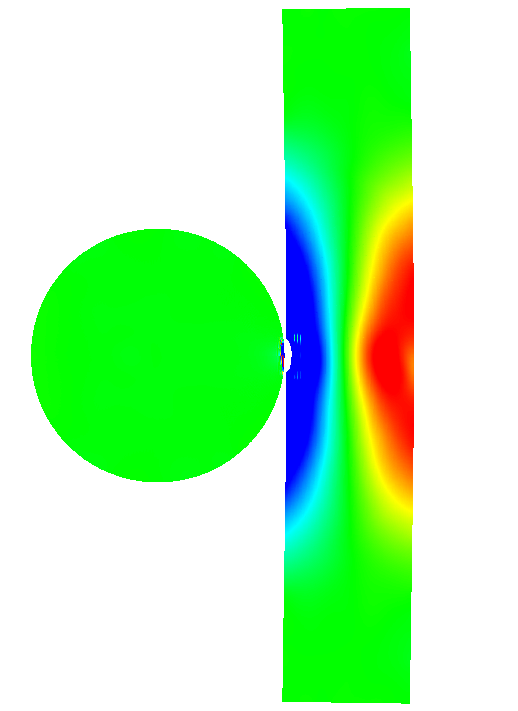}} \\
		\subfloat[\label{fig:SBI_HS_C1FourthMPMOrtho_40sec}]{
			\includegraphics[width=0.25\columnwidth]{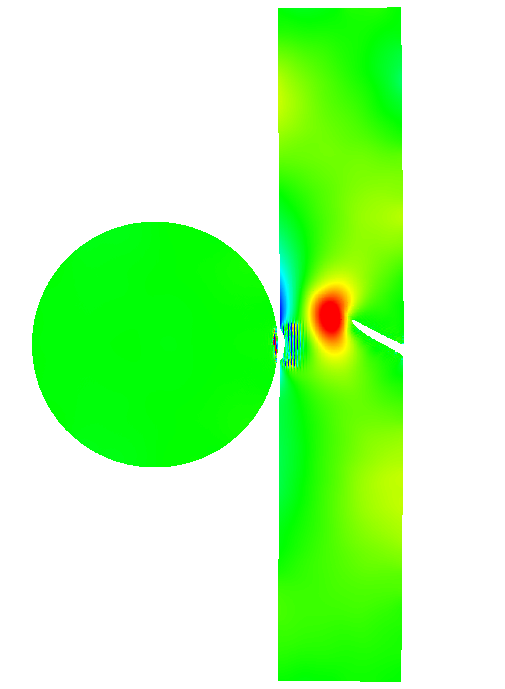}} &
		\subfloat[\label{fig:SBI_HS_C1FourthMPMOrtho_70sec}]{
			\includegraphics[width=0.25\columnwidth]{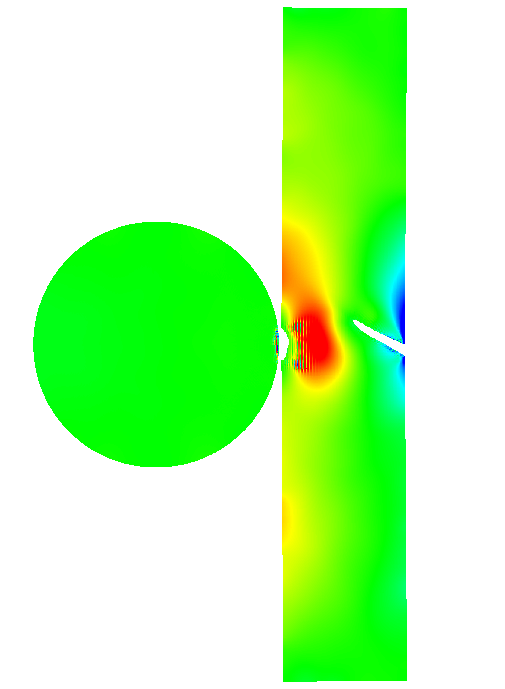}} &
		\subfloat[\label{fig:SBI_HS_C1FourthMPMOrtho_80sec}]{
			\includegraphics[width=0.25\columnwidth]{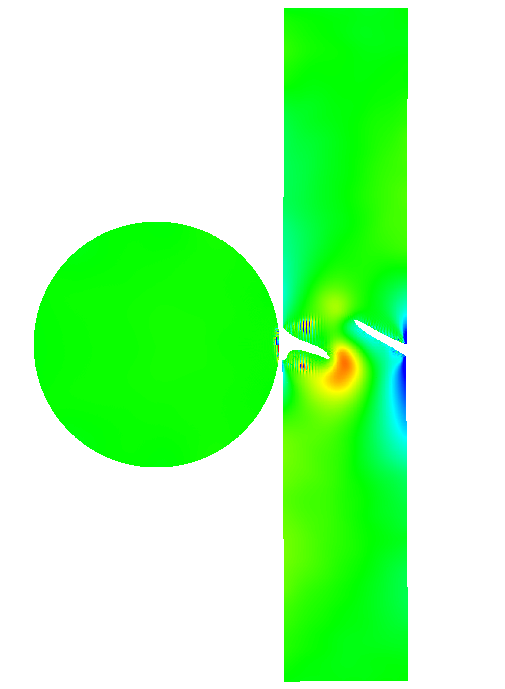} } \\
		\multicolumn{3}{c}{\subfloat{
				\includegraphics[width=0.75\columnwidth]{HS_ColorBar1.png}}} 		
	\end{tabular}
	\caption[]{Sphere-beam impact fracture problem: Hydrostatic stresses for time steps \subref{fig:SBI_HS_C1FourthMPMOrtho_0sec} t=0 $\mu$s \subref{fig:SBI_HS_C1FourthMPMOrtho_12sec} t=12 $\mu$s \subref{fig:SBI_HS_C1FourthMPMOrtho_28sec} t=28 $\mu$s \subref{fig:SBI_HS_C1FourthMPMOrtho_40sec} t=40 $\mu$s \subref{fig:SBI_HS_C1FourthMPMOrtho_70sec} t=70 $\mu$s and \subref{fig:SBI_HS_C1FourthMPMOrtho_80sec} t=80 $\mu$s. Results for case (iii): PF-MPM 4th order orthotropic model. Material points with $c_p<0.08$ have been removed.}
	\label{fig:SBI_HS_C1FourthMPMOrtho}
\end{figure}

\subsubsection{Discussion on observed fracture patterns}

Fracture patterns emerging from static indentation of a practically rigid sphere against a deformable solid as well as from low and high speed impact tests have been the focus of extensive experimental investigations, see, e.g., \cite{ball1994low, jelagin2008indentation}. The failure modes observed vary considerably with the velocity of the projectile, the flexibility of the impacted beam, and the interface properties \cite{jelagin2008indentation}.

We focus here in case (i) with isotropic fracture energy equal to $\mathcal{G}_{c} \left( \theta \right)=10.6066$ {\normalfont N/mm}. The fracture patterns shown in Fig. \ref{fig:SBI_PF_C1SecondMPMIso1Cracks} correspond to a median type of crack, with the crack at the left end nucleating due to impact and then propagating towards the right edge driven by the principal tensile stresses at the mid-span. To investigate the effect of the projectile velocity on the induced fracture pattern, a total of 33 analysis cases is performed keeping the geometry, the elastic and the fracture properties of the beam similar to those reported in section \ref{lbl:SBI}. In each case the projectile velocity is varied from 0.02  mm/$\mu$s to 0.18  mm/$\mu$s at a step size of 0.005 mm/$\mu$s. All analysis parameters, the background cell size and the cell density are similar to section \ref{lbl:SBI}.

\begin{figure}
	\centering
	\includegraphics[width=0.40\columnwidth]{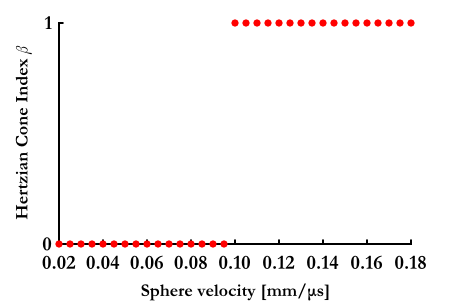}
	\caption[]{Hertzian cone index versus projectile velocity}
	\label{fig:SBI_Hertz}
\end{figure}

\begin{figure}
	\centering
	\begin{tabular}{lcr}
		\subfloat[\label{fig:Hertz1}]{
			\includegraphics[width=0.25\columnwidth]{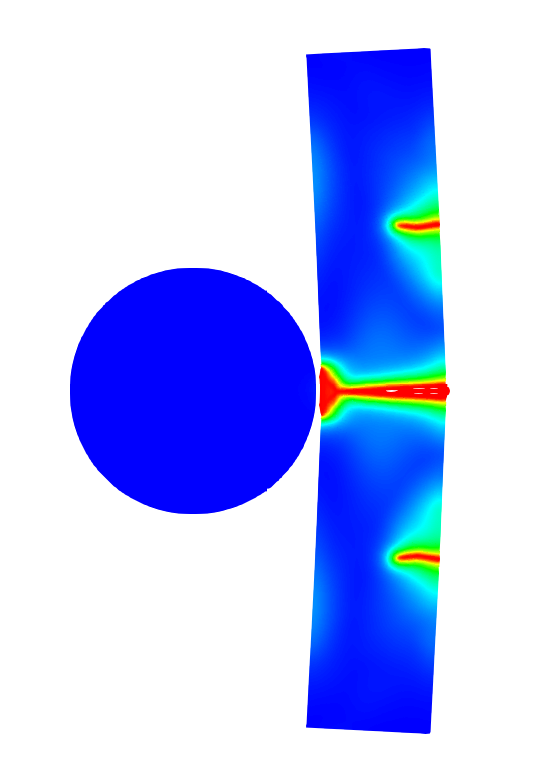}} &
		\subfloat[\label{fig:Hertz2}]{
			\includegraphics[width=0.25\columnwidth]{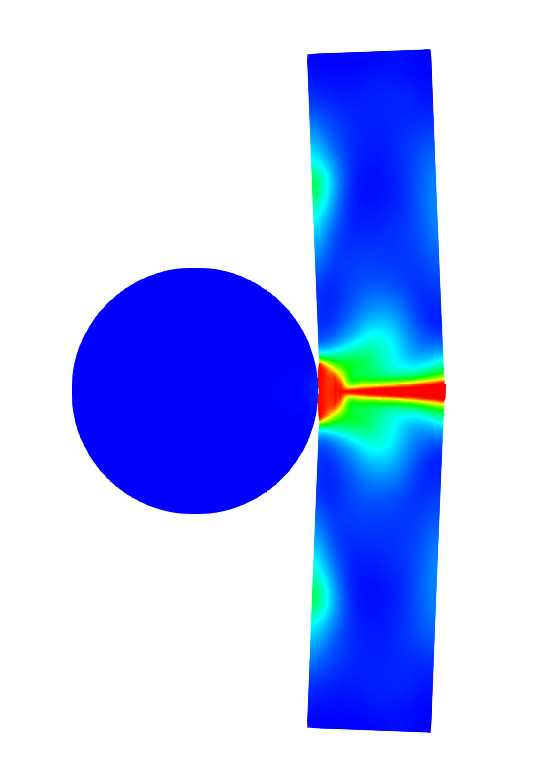}} &
		\subfloat[\label{fig:Hertz3}]{
			\includegraphics[width=0.25\columnwidth]{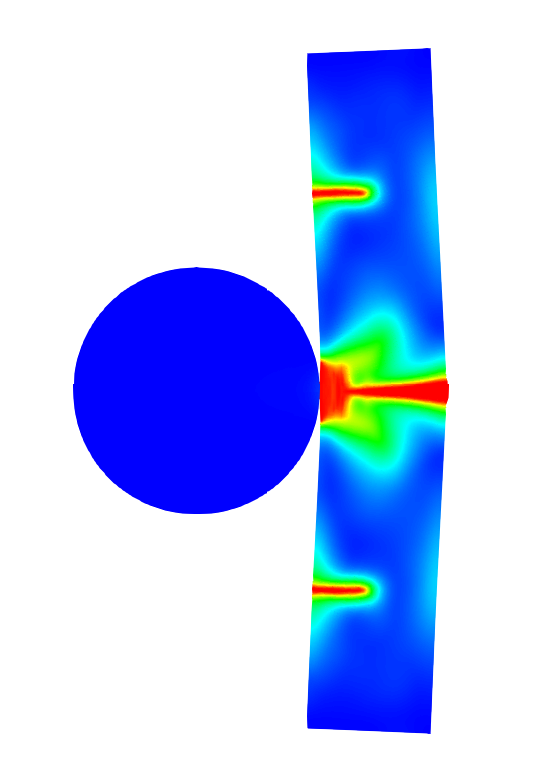}} \\
		\subfloat[\label{fig:Hertz4}]{
			\includegraphics[width=0.25\columnwidth]{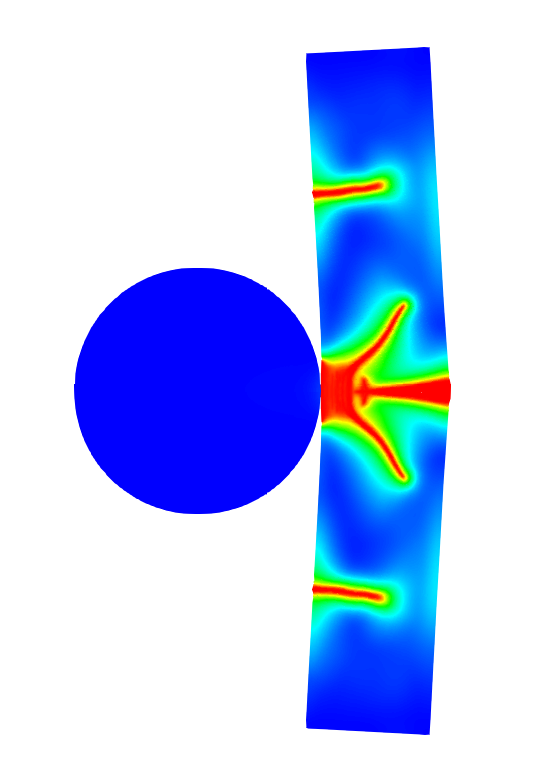}} &
		\subfloat[\label{fig:Hertz5}]{
			\includegraphics[width=0.25\columnwidth]{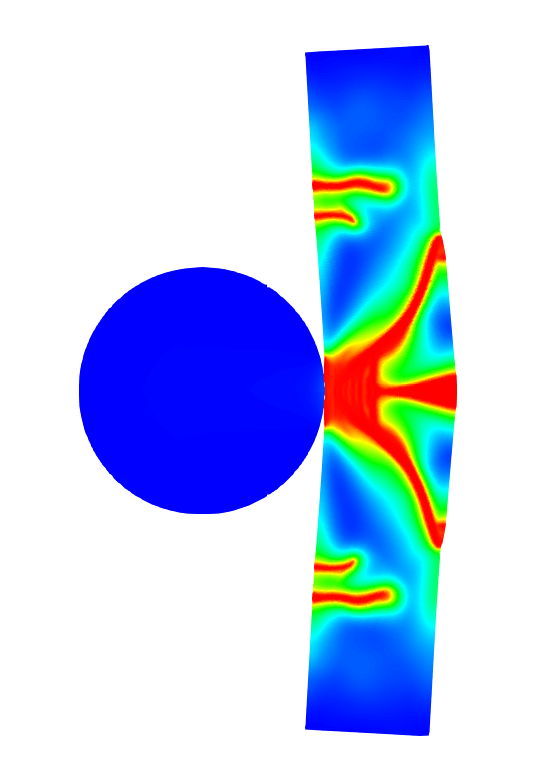}} &
		\subfloat[\label{fig:Hertz6}]{
			\includegraphics[width=0.25\columnwidth]{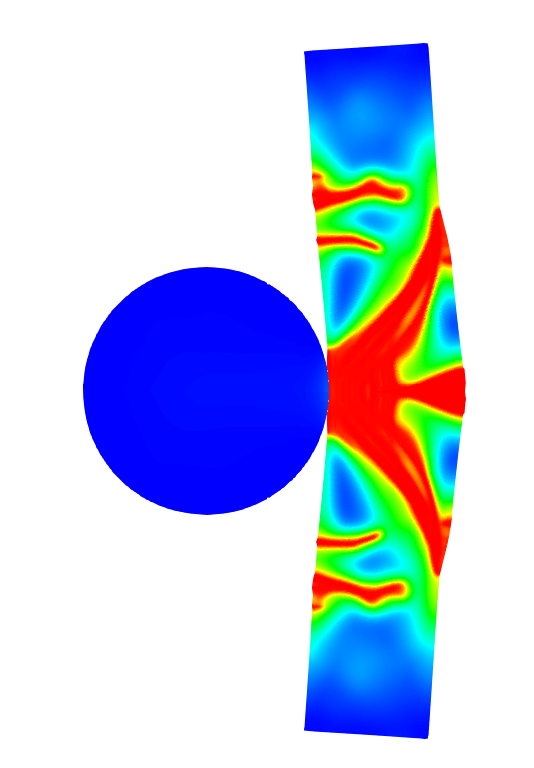} } \\
		\multicolumn{3}{c}{\subfloat{
				\includegraphics[width=0.75\columnwidth]{HS_ColorBar1.png}}} 		
	\end{tabular}
	\caption[]{Fracture patterns for varying projectile velocities \subref{fig:Hertz1}  $\dot{\mathbf{u}}_{Ap_{(0)}}=0.05$ mm/$\mu$s \subref{fig:Hertz2} $\dot{\mathbf{u}}_{Ap_{(0)}}=0.08$ mm/$\mu$s \subref{fig:Hertz3}  $\dot{\mathbf{u}}_{Ap_{(0)}}=0.095$ mm/$\mu$s   \subref{fig:Hertz4} $\dot{\mathbf{u}}_{Ap_{(0)}}=0.10$ mm/$\mu$s  \subref{fig:Hertz5} $\dot{\mathbf{u}}_{Ap_{(0)}}=0.14$ mm/$\mu$s \subref{fig:Hertz6} $\dot{\mathbf{u}}_{Ap_{(0)}}=0.18$ mm/$\mu$s}
	\label{fig:Hertz}
\end{figure}

We define the Hertzian cone index $\beta$ with a value $\beta=0$ corresponding to a cone not developing and $\beta=1$ when a cone develops. This is plotted versus the projectile velocity in Fig. \ref{fig:SBI_Hertz}; a cone fracture pattern occurs for velocities larger than 0.10 mm/$\mu$s. Conversely, for velocities smaller than 0.10 mm/$\mu$s the crack pattern is consistent with the flexure failure mode described in section \ref{lbl:SBI}. In Fig. \ref{fig:Hertz}, phase field snapshots are shown for six particular cases of projectile velocity. All snapshots correspond to time $t=16 \mu$ when the maximum value of fracture energy for all case has been attained. 

For the case of the lowest impact velocity considered in Fig. \ref{fig:Hertz1}, the crack pattern involves a median crack propagating from the left and towards the right edge of the beam. Secondary, flexural cracks appear at the right edge of the beam. Such a response is consistent with experimental observations on brittle materials at low impact loads where a plastic band initiates at the impact zone prior to crack formation, see, e.g., \cite{hagan1978origin}. In the framework presented herein, the material degradation prior to fracture assumes this role. 

Increasing impact velocities result in a Hertz cone formation at the vicinity of the impact zone. Secondary cracks also propagate from the left edge. Of interest is also the evidently smooth transition from a median to a Hertz cone fracture pattern from Fig. \ref{fig:Hertz3} to Fig. \ref{fig:Hertz4}.

\subsection{Anisotropic plate with centred crack} \label{lbl:AnisoPlate}

\begin{figure}
	\centering
	\includegraphics[width=0.45\columnwidth]{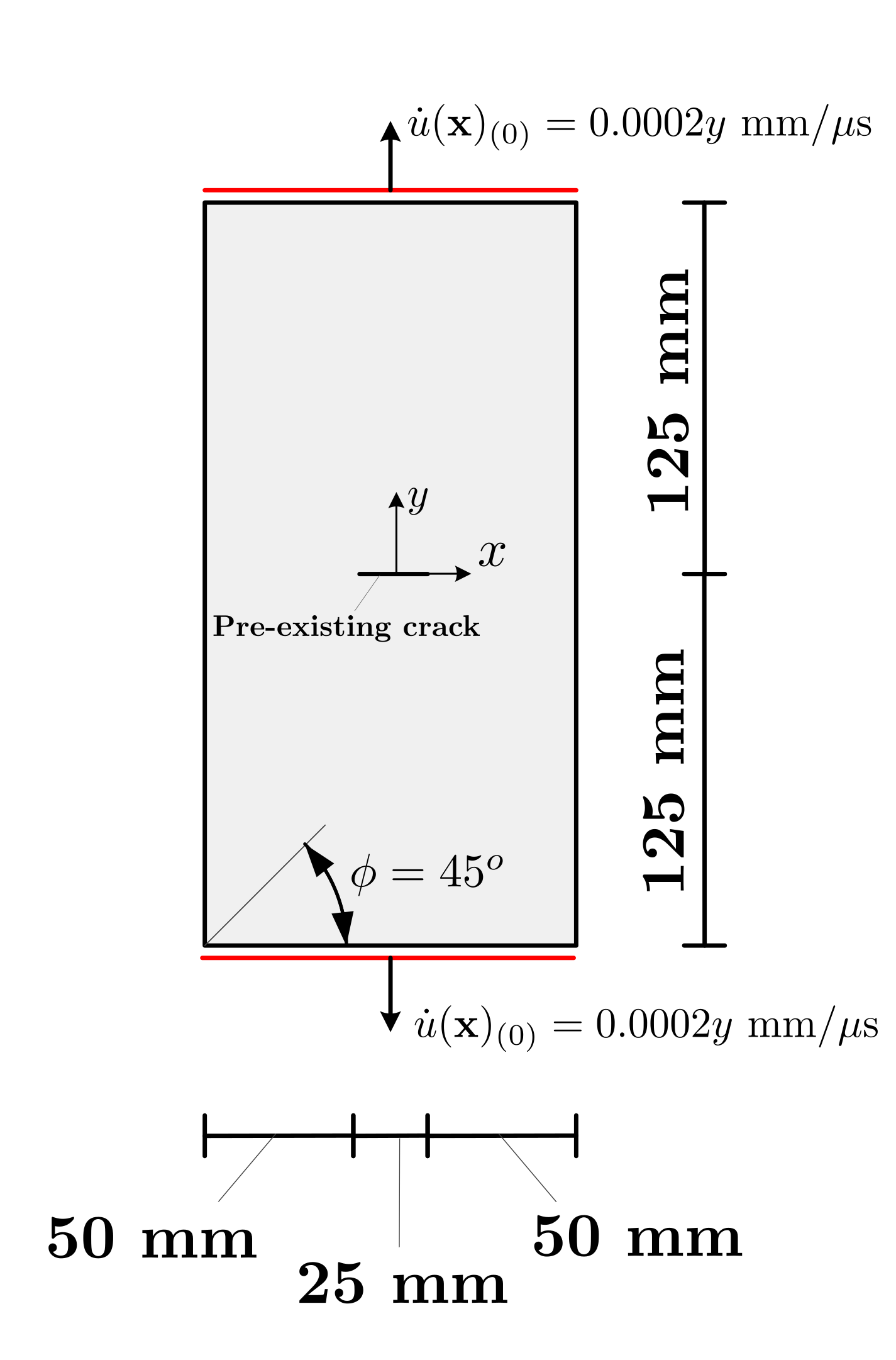}
	\caption[]{Anisotropic plate with centre crack: Geometry and boundary conditions.}
	\label{fig:AnisoPlate_GeoBCs}
\end{figure}

The case of the anisotropic rectangular plate shown in Fig. \ref{fig:AnisoPlate_GeoBCs} is examined, made from a unidirectional HTA/ 6376 composite laminate and subjected to an initial velocity field $\dot{u}\left(\mathbf{x}\right)_{(0)}=0.0002y$ mm/$\mu$s. The plate contains a pre-existing crack at its centre with length $25$ mm. The material properties of the composite are summarized in Table \ref{tab:1}. In addition, the following values hold, i.e., $\nu_{xy}$=0.3 and $G_{xy}=5500$ N/mm\textsuperscript{2}.

\begin{table}[]
	\begin{tabular}{lrrrr}
		& \multicolumn{1}{c}{HTA fibre} & \multicolumn{1}{c}{6376 epoxy} & \multicolumn{1}{c}{laminate ($\phi$ = 0$^o$)} & \multicolumn{1}{c}{laminate ($\phi$ = 90$^o$)} \\ \hline
		Young’s modulus {[}MPa{]}  & 235000                        & 3600                           & 136000                                        & 8750                                           \\ \hline
		Tensile strength {[}MPa{]} & 3920                          & 105                            & 1670                                          & 60                                             \\ \hline
		Density {[}kg/m$^3${]}     & 1770                          & 1310                           & 1586                                          & 1586                                           \\ \hline
	\end{tabular}
\label{tab:1}
\end{table}

In their experimental campaign, Cahil et al. \cite{cahill2014experimental} have shown that cracks grow parallel to the fibre direction hence indicating that the damage originates only through matrix failure. We consider herein the case of fibre orientation at $\phi=+45^o$ as shown in Fig. \ref{fig:AnisoPlate_GeoBCs}. This allows us to use the same Young's modulus and Poisson's ratio along $x$ and $y$. The elastic material properties considered are $E=14980$ N/mm\textsuperscript{2} and $\nu=0.36$.The  mass density is $\rho=1586$ kg/m\textsuperscript{3}. The length scale parameter is considered $l_0=1$ mm and the anisotropic parameters are taken to be $\gamma_{1111}=1.00$, $\gamma_{2222}=2900$, $\gamma_{1122}=0.00$ and $\gamma_{1222}=74.00$ and $\bar{\mathcal{G}}_c=4.175$ N/mm. These parameters correspond to an orthotropic surface energy with $\mathcal{G}_{c_{min}}=5.9067$ N/mm along the fibre orientation and $\mathcal{G}_{c_{max}}=30.9044$ N/mm normal to the fibre.

The grid is formed by two knot vectors $\splitatcommas{\Xi=\{0,0,0,0.00666,0.01333,...,0.98666,0.99333,1,1,1\}}$ and $\splitatcommas{H=\{0,0,0,0.0033,0.0066,0.01,...0.9899,0.9933,0.9966,1,1,1\}}$, $45904$ control points and $150x300=45000$ cells. The cell (patch) spacing is $h=1.00$ mm and plane stress conditions are assumed. The total number of material points is $281250$. The solution procedure is implemented with a time step $\Delta t = 0.0125$ $\mu$s for a total time of $25$ $\mu$s. The critical time step is $\tilde{\Delta t_{cr}}=0.201$ $\mu$s.

\begin{figure}
	\centering
	\begin{tabular}{cccc}
		\subfloat[\label{fig:AnisoPlate_45Deg_1}]{
			\includegraphics[width=0.20\columnwidth]{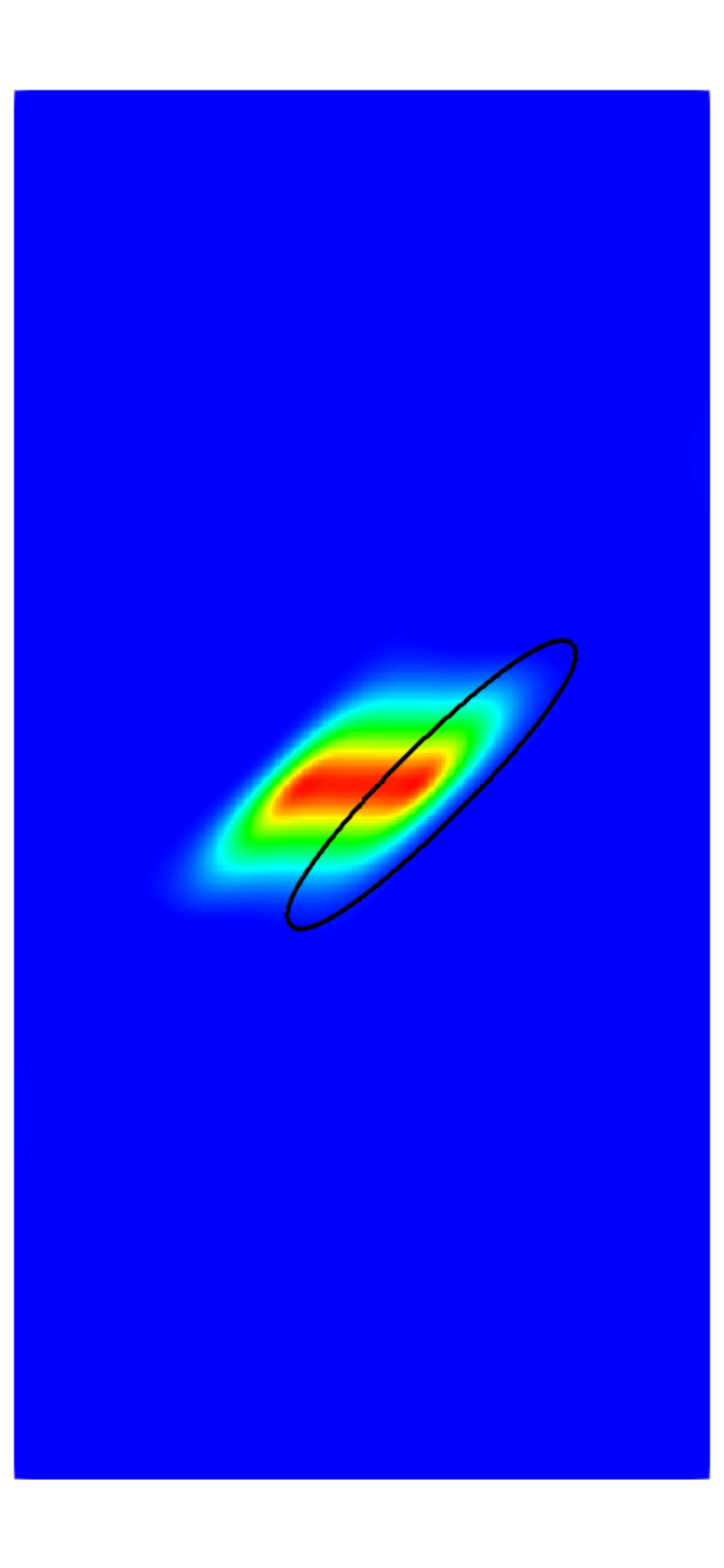}} &
		\subfloat[\label{fig:AnisoPlate_45Deg_2}]{
			\includegraphics[width=0.20\columnwidth]{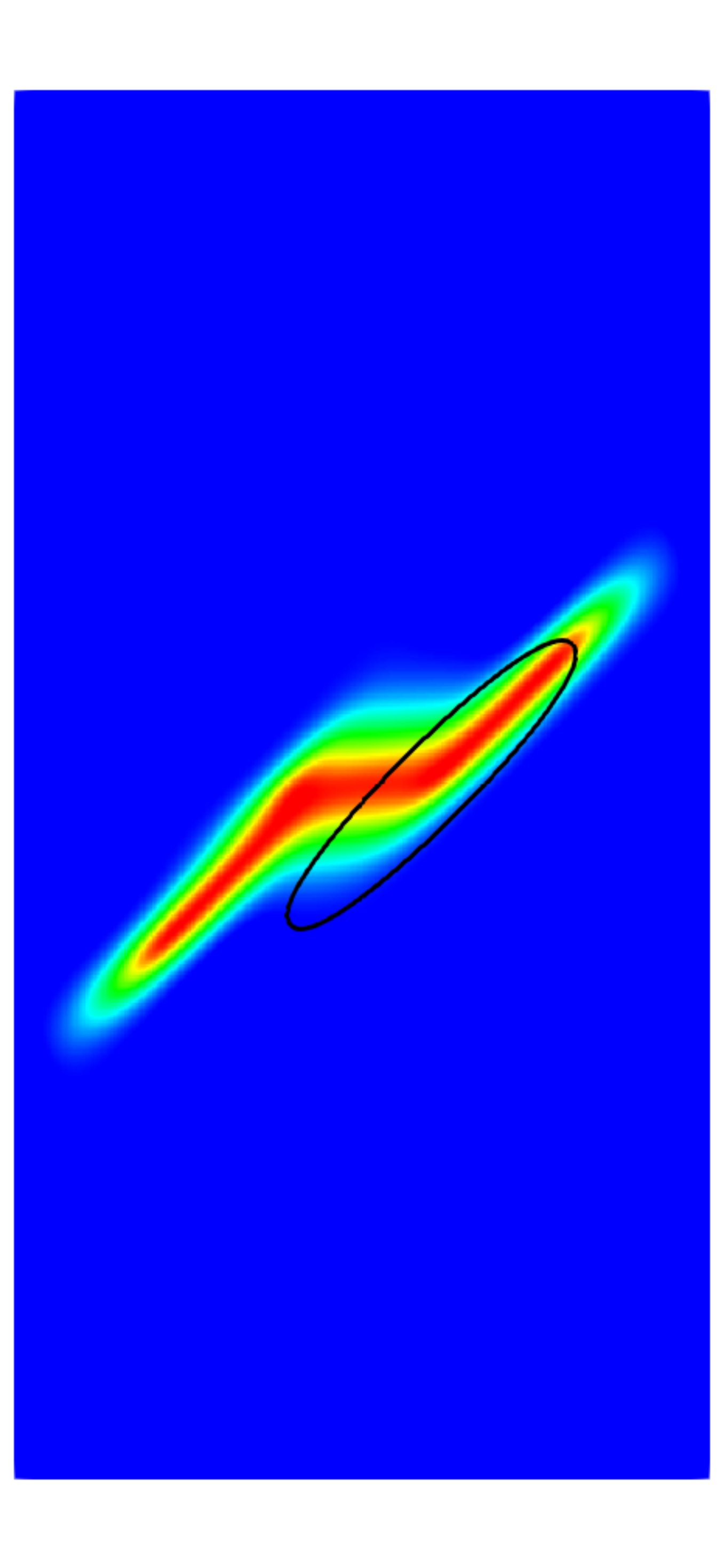}} &
		\subfloat[\label{fig:AnisoPlate_45Deg_3}]{
			\includegraphics[width=0.20\columnwidth]{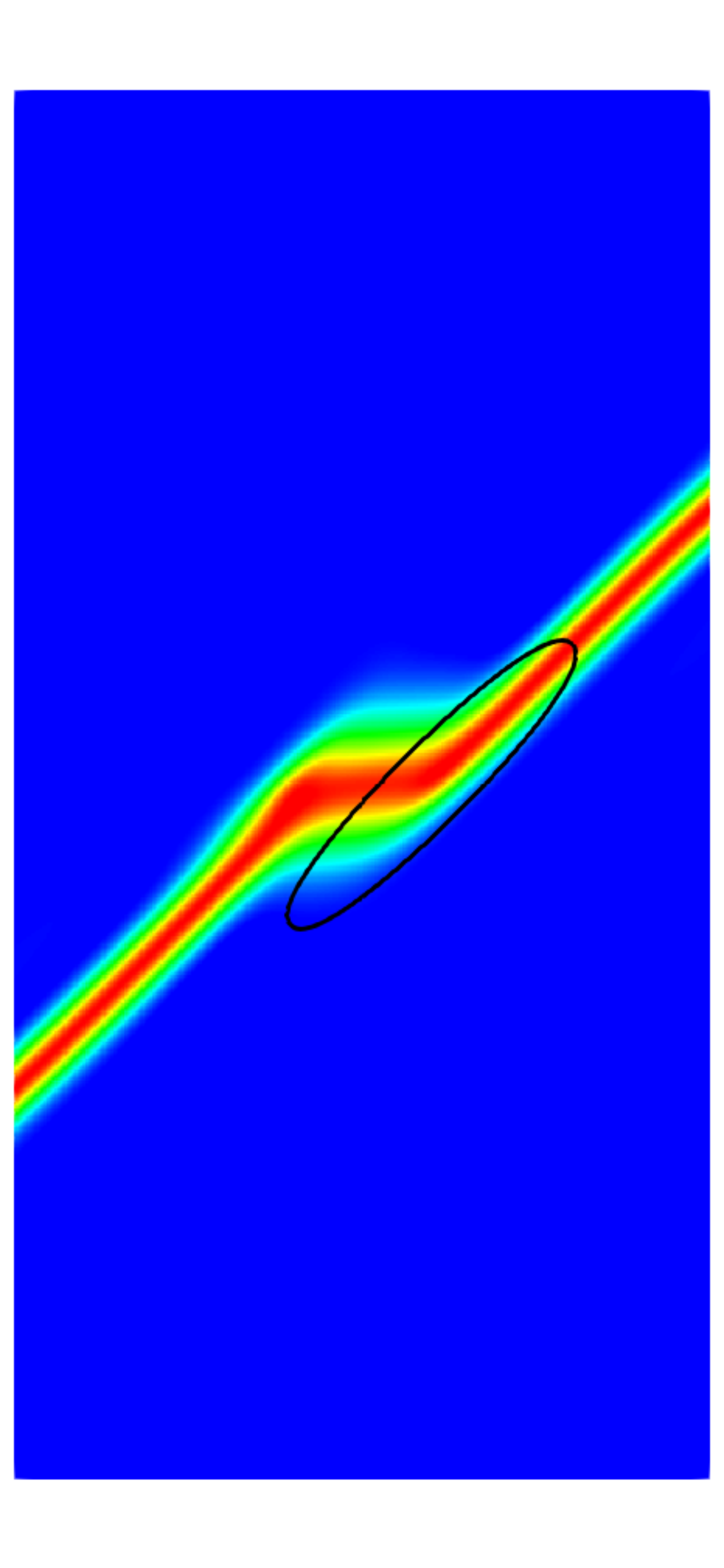}} &
		\subfloat[\label{fig:AnisoPlate_45Deg_Exp}]{
			\includegraphics[width=0.20\columnwidth]{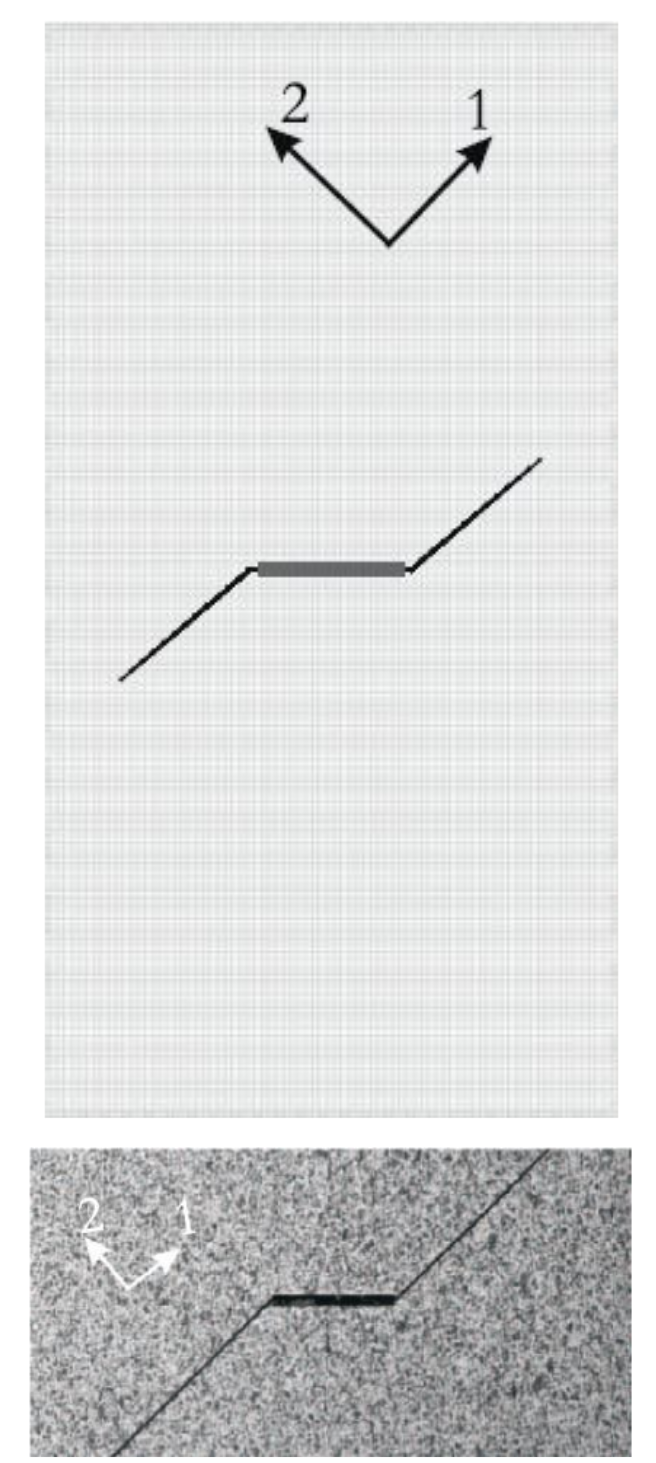}} \\
		\multicolumn{4}{c}{\subfloat{
				\includegraphics[width=0.75\columnwidth]{PF_ColorBar.png}}}		
	\end{tabular}
	\caption[]{Anisotropic plate with centre crack: Phase field for time steps \subref{fig:AnisoPlate_45Deg_1} t=0 $\mu$s \subref{fig:AnisoPlate_45Deg_2} t=22 $\mu$s and \subref{fig:AnisoPlate_45Deg_3} t=25 $\mu$s. The experimental observations are shown in \subref{fig:AnisoPlate_45Deg_Exp} (Cahil et al. \cite{cahill2014experimental}).}
	\label{fig:AnisoPlate_PF}
\end{figure}

The results for the numerical simulations together with the experimental observations are shown in Fig. \ref{fig:AnisoPlate_PF}. The reciprocal of the surface energy density (black eclipse) is also plotted on these snapshots. In Figs. \ref{fig:AnisoPlate_45Deg_2} and \ref{fig:AnisoPlate_45Deg_3}, the phase field evolves along the material orientation $\phi=+45^o$. The crack paths derived from our simulation agree well with the experimental crack paths presented in Fig. \ref{fig:AnisoPlate_45Deg_Exp}.

\vspace{-6pt}
\section{Conclusions} \label{sec:Conclusions}
\vspace{-2pt}

In this work, a novel numerical method is introduced for the treatment of dynamic brittle fracture in both isotropic and anisotropic media. The evolution of crack paths is represented by means of phase field models within a Material Point Method setting. Anisotropy is explicitly introduced in the fracture energy through a crack density functional. The method is further extended to account for frictional contact problems involving phase field fracture adopting a discrete field approach. A notable advantage of the proposed formulation is that both the equilibrium and phase field governing equations are solved independently for each discrete field rendering the method suitable for parallel implementation. 

The method is rather appealing for the case of phase field modelling where very fine meshes are commonly required due to the regularized crack topology. Rather than employing a uniform background mesh and material point density, multiple small scale problems can be solved separately for each discrete field at their corresponding background domain. In terms of contact driven fracture, contrary to standard FEM implementations that necessitate the algorithmic treatment of local contact features, these now naturally emerge from the interaction of material points within a fixed Eulerian mesh. Indeed, the fixed Eulerian grid is utilized to identify the contact surfaces using the material points’ projection on the grid.

A set of representative numerical examples is presented where the computational advantages of PF-MPM are demonstrated. The method is verified against the standard Phase Field Finite Element Method; the two methods are in good agreement. The influence of anisotropy is examined in terms of crack path, time history energy results and crack tip velocities. Benchmark problems with complex crack path i.e. crack branching and merging are considered and the robustness of the method is established. It is shown that different loading velocities and fracture material parameters strongly influence the dynamic failure response of the structure and the resulting crack paths. Fracture energies computed from the proposed method are compared and indeed verified against the corresponding analytical predictions. Finally, crack paths derived from the method are validated against experimental observations.

\section*{Acknowledgement} \label{sec:Acknowledgements}

The research described in this paper has been financed by the University of Nottingham through the Dean of Engineering Prize, a scheme for pump priming support for early career academic staff. The authors are grateful to the University of Nottingham for access to its high performance computing facility.

\bibliography{PFMPM}

\end{document}